\documentclass[11pt]{book}
\usepackage[a4paper,width=145mm,top=30mm,bottom=30mm]{geometry}
\usepackage[pagebackref=true]{hyperref}
\usepackage{amsmath}
\usepackage{graphicx}
\usepackage{mathrsfs}
\usepackage[latin1]{inputenc}
\usepackage{amssymb,wasysym,stmaryrd}
\usepackage{amsthm}
\usepackage{thmtools}
\usepackage{mathtools}
\usepackage{float}
\restylefloat{table}
\usepackage{multirow}
\usepackage{verbatim}
\usepackage{cmll}

\usepackage{tikz}
\usetikzlibrary{matrix,arrows,decorations.pathmorphing,backgrounds,decorations.markings}
\usetikzlibrary{arrows.meta}
\tikzset{2cell/.style={-implies,double,double equal sign distance,shorten >=0pt, shorten <=0pt}}
\tikzset{2cellr/.style={implies-,double,double equal sign distance,shorten >=0pt, shorten <=0pt}}
\tikzset{2cellnode/.style={-implies,double,double equal sign distance,shorten >=2pt, shorten <=3pt}}
\tikzset{2cellnoder/.style={implies-,double,double equal sign distance,shorten >=3pt, shorten <=2pt}}
\tikzset{3cell/.style={-implies,double,double distance=2.5pt,shorten >=2pt, shorten <=3pt}}

\tikzset{labelsize/.style={font=\scriptsize}}
\tikzset{myborder/.style={color=gray,thin}}

\usetikzlibrary{external} 
\tikzset{external/up to date check={md5}}
\tikzset{external/system call={pdflatex \tikzexternalcheckshellescape -halt-on-error -interaction=batchmode -jobname "\image" "\texsource"}}

\numberwithin{equation}{chapter}

\declaretheoremstyle[
spaceabove=6pt, spacebelow=6pt,
headfont=\normalfont\itshape,
notefont=\mdseries, notebraces={(}{)},
bodyfont=\normalfont,
postheadspace=0.5em,
qed=$\square$
]{myproofstyle}

\declaretheorem[style=plain,numberwithin=chapter,name=Theorem]{thm}
\declaretheorem[style=plain,sibling=thm,name=Lemma]{lem}
\declaretheorem[style=plain,sibling=thm,name=Proposition]{prop}

\declaretheorem[style=definition,qed=$\blacksquare$,sibling=thm,name=Definition]{defn}

\declaretheorem[style=myproofstyle,numbered=no,name={Proof}]{pf}

\mathchardef\mhyphen="2D

\newcommand{\defemph}[1]{{\bfseries #1}}

\newcommand{\const}[1]{\mathrm{const}_{#1}}

\newcommand{\TCat}{{\mathit{2}\mhyphen \Cat}}
\newcommand{\V}{\mathcal{V}}

\newcommand{\Res}[1]{\mathbf{Res}(#1)}

\newcommand{\Fib}[1]{\mathscr{F}\!ib(#1)}

\newcommand{\cat}[1]{\mathbb{#1}}
\newcommand{\tcat}[1]{\mathscr{#1}}
\newcommand{\st}[1]{\mathrm{st}(#1)}
\newcommand{\Mnd}[1]{\mathbf{Mnd}({#1})}
\newcommand{\CoMnd}[1]{\mathbf{Comnd}({#1})}
\newcommand{\Cat}{\mathscr{C}\!at}
\newcommand{\Epp}{\mathscr{E}^{++}}
\newcommand{\Epm}{\mathscr{E}^{+-}}
\newcommand{\Emp}{\mathscr{E}^{-+}}
\newcommand{\Emm}{\mathscr{E}^{--}}

\newcommand{\ppp}{p^{++}}
\newcommand{\ppm}{p^{+-}}
\newcommand{\pmp}{p^{-+}}
\newcommand{\pmm}{p^{--}}

\newcommand{\Inj}{\mathbf{Inj}}
\newcommand{\Set}{\mathbf{Set}}

\newcommand{\grst}{\gm{T}}
\newcommand{\vl}{{val}}
\newcommand{\wl}{{wal}}
\newcommand{\inst}{\im{T}}

\newcommand{\inr}{\mathrm{inr}}
\newcommand{\inl}{\mathrm{inl}}

\newcommand{\dprod}{{\prod\nolimits_{\cat{B}}}} 

\newcommand{\gm}[1]{\mathbf{#1}}
\newcommand{\im}[1]{\mathscr{#1}}
\newcommand{\il}[1]{\mathcal{#1}}
\newcommand{\Law}{\mathbf{Law}}
\newcommand{\mdl}[2]{\mathbf{Mod}(#1,#2)}

\newcommand{\imc}[1]{\im{T}_{{#1}}}  
\newcommand{\icc}[1]{\im{S}_{{#1}}}

\newcommand{\act}  {{\mathbin{\ast}}}
\newcommand{\sact} {{\mathbin{\circledcirc}}}

\newcommand{\emact}{{\mathbin{\circledast}}}
\newcommand{\klact}{{\mathbin{\odot}}}

\newcommand{\op}[1]{{\mathrm{op}(#1)}}
\newcommand{\opp}{\mathrm{op}}

\newcommand{\id}[1]{\mathrm{id}_{#1}}

\newcommand{\tensor}{{\mathbin{\otimes}}}
\newcommand{\e}{I}
\newcommand{\fork}{{\mathbin{\pitchfork}}}
\newcommand{\copow}{{\mathbin{\otimes}}}

\newcommand{\vc}{{\mathbin{.}}}
\newcommand{\hc}{{\mathbin{\ast}}}

\newcommand{\EM}[2]{{{#1}^{#2}}}
\newcommand{\EMg}{\EM{\cat{C}}{\gm{T}}}
\newcommand{\EMi}{\EM{\cat{C}}{\im{T}}}

\newcommand{\CoEMg}{\EM{\cat{C}}{\gm{S}}}
\newcommand{\CoEMi}{\EM{\cat{C}}{\im{S}}}

\newcommand{\Kl}[2]{{{#1}_{#2}}}
\newcommand{\Klg}{\Kl{\cat{C}}{\gm{T}}}
\newcommand{\Kli}{\Kl{\cat{C}}{\im{T}}}
\newcommand{\CoKlg}{\Kl{\cat{C}}{\gm{S}}}

\newcommand{\Gray}{\mathbf{Gray}}

\newcommand{\EKl}{\cat{E}}

\newcommand{\ecat}[1]{\mathcal{#1}}

\newcommand{\Sfunc}{\Sigma}  

\newcommand{\valg}[3]{
\begin{tikzpicture}[baseline=-\the\dimexpr\fontdimen22\textfont2\relax ]
      \node at (0,0.34)  [font=\footnotesize,inner sep=0.25pt]{$#1$};
      \node at (0,-0.34) [font=\footnotesize,inner sep=0.25pt]{$#3$};
      
      \draw [->] (0,0.15) to node [auto,font=\scriptsize,inner sep=0pt]{$\ #2$} (0,-0.15);
\end{tikzpicture}
} 
  
\newcommand{\valgp}[3]{
\left(
\valg{#1}{#2}{#3}
\right)
} 

\newcommand{\vcoalg}[3]{
\begin{tikzpicture}[baseline=-\the\dimexpr\fontdimen22\textfont2\relax ]
      \node at (0,0.34)  [font=\footnotesize,inner sep=0.25pt]{$#3$};
      \node at (0,-0.34) [font=\footnotesize,inner sep=0.25pt]{$#1$};
      
      \draw [->] (0,0.15) to node [auto,font=\scriptsize,inner sep=0pt]{$\ #2$} (0,-0.15);
\end{tikzpicture}
} 

\newcommand{\vcoalgp}[3]{
\left(
\vcoalg{#1}{#2}{#3}
\right)
} 
  
\newcommand{\ccol}[1]{\omit\hfill$#1$\hfill}

\title{A 2-Categorical Study of Graded and Indexed Monads}
\author{Soichiro Fujii}
\date{January 25, 2016}
\begin{document}
\frontmatter
\begin{titlepage}
    \begin{center}
        \vspace*{1cm}
 
        \Huge
        A 2-Categorical Study of 
        Graded and Indexed Monads
        
        \vspace{9cm}
        \LARGE{Soichiro Fujii}
 
        \vfill
 \LARGE{A Master Thesis\\ \vspace{1cm}
 Submitted to\\
 the Graduate School of the University of Tokyo\\
 on January 25, 2016\\
 in Partial Fulfillment of the Requirements\\
 for the Degree of Master of Information Science and
 Technology\\
 in Computer Science}
 
    \end{center}
\end{titlepage}
\chapter*{Abstract}
In the study of computational effects, it is important to 
consider the notion of 
\emph{computational effects with parameters}.
The need of such a notion arises when, for example, 
statically estimating the range of effects caused by a program, 
or studying the ways in which effects with local scopes 
are derived from effects with only the global scope.
Extending the classical observation that
computational effects can be modeled by monads,
these computational effects with parameters are modeled by 
various mathematical structures including 
\emph{graded monads} and \emph{indexed monads},
which are two different generalizations of 
ordinary monads.
The former has been employed in the semantics of effect systems,
whereas the latter in the study of the relationship 
between the local state monads and the global state monads,
each exemplifying the two situations mentioned above.
However, despite their importance, the mathematical theory of
graded and indexed monads is far less developed than
that of ordinary monads.

Here we develop the mathematical theory of graded 
and indexed monads from a 2-categorical viewpoint. 
We first introduce four 2-categories and observe that
in two of them graded
monads are in fact monads in the 2-categorical sense,
and similarly indexed monads are monads 
in the 2-categorical sense in the other two.
We then construct explicitly the Eilenberg--Moore and 
the Kleisli objects of graded monads, and 
the Eilenberg--Moore objects of indexed monads
in the sense of Street in appropriate 2-categories among
these four.
The corresponding results for \emph{graded} and 
\emph{indexed comonads} also follow.

We expect that the current work will provide a theoretical
foundation to a unified study of computational effects with
parameters, or dually (using the comonad variants), of
\emph{computational resources with parameters},
arising for example in Bounded Linear Logic.

\chapter*{Acknowledgements}
First of all I thank my supervisor Ichiro Hasuo for his patience and 
generous support.
Thanks to his assistance I was able to interact with 
many top-level researchers and up-and-coming students, 
and this meant a lot to me.

I am perhaps intellectually most indebted to Paul-Andr\'e
Melli\`es during my master's course, and the joint work with him 
and Shin-ya Katsumata forms the basis of this thesis.
Discussions with him in his office, and also in nice caf\'es in
Paris were always delightful.

I was fortunate enough to be able to learn portions of 
3-category theory from John Power himself.
It was indeed impressive to observe that he drew various
higher-dimensional diagrams in an experienced manner,
and I thank him for his kindness.

Let me finally thank all the people who have 
helped me broaden and deepen my knowledge 
through precious discussions:
Kazuyuki Asada, 
Wataru Hino,
Toshiki Kataoka,
Shin-ya Katsumata,
Kenji Maillard,
Tokio Sasaki, and
Takeshi Tsukada,
to name just a few.

\tableofcontents 

\mainmatter 

\chapter{Introduction}

\section{Background and related work}
Roughly speaking, there are two lines of research
on which the current work is based,
one \emph{computational} and one \emph{mathematical}.
From a computational point of view, this work can be 
considered as a contribution to a theoretical study of 
\emph{computational effects}, originating in
Moggi's work in the late 1980's; this perspective 
is explained in Section~\ref{subsec-intro-comp}.
The mathematical aspect of the thesis relies on
the \emph{2-categorical theory of monads}
developed by Street in the 1970's, as recalled in
Section~\ref{subsec-intro-math}.

\subsection{Monads}
Let us begin with a brief discussion on
the basic theory of \emph{monads}, which forms the common basis
for both of the two lines of research mentioned above.

The notion of monad originates in pure mathematics,
and has been one of the main subjects of research in category theory.
Monads are so fundamental that any reasonable introduction to
category theory contains some account of them; see for example 
\cite{maclane}.
Among the earliest important constructions in the theory of 
monads are those of 
Eilenberg--Moore~\cite{eilenberg-moore}
and Kleisli~\cite{kleisli} obtained in mid 1960's, 
each building a category~$\cat{D}$ out of a monad~$T$
on a category~$\cat{C}$, together with an
adjunction 
\[
    \begin{tikzpicture}
      \node (A)  at (0,0) {$\cat{C}$};
      \node (C) at (4,0) {$\cat{D}$};
      
      \draw [->] (A) to [bend left=23]node (dom)[auto,labelsize] {$F$}(C);
      \draw [->] (C) to [bend left=23]node (cod)[auto,labelsize] {$U$}(A);
      
      \node at (2,0) [rotate=90,font=\large]{$\vdash$};
    \end{tikzpicture}
\]
that \emph{generates} $T$ in the sense that $T=U\circ F$
(and similarly for the other components).

What we achieve in this thesis is, intuitively, 
the generalization of these constructions of 
Eilenberg--Moore and Kleisli from
ordinary monads to generalized notions of monad called 
\emph{graded} or \emph{indexed monads}.

\subsection{Computational effects and monads}
\label{subsec-intro-comp}
Examples of \emph{computational effects} include those behaviors of 
programs involving \emph{global state}, 
\emph{nondeterministic branch}, 
or \emph{input/output}.
Moggi~\cite{moggi:comp-lambda}
brought breakthrough to the theory of computational effects
by advocating that these
diverse kinds of computational effects can be 
uniformly modeled by \emph{monads}.
For example, corresponding to the notion of global state
there is a monad~$T$ (say, on $\Set$) called the
\emph{global state monad}, defined as 
$TX\coloneqq S\Rightarrow (S\times X)$
where $S$ is the set of \emph{states}.
Similarly, there are monads for nondeterministic branch,
input/output, etc.
What is especially interesting in Moggi's approach is that, 
not only mathematical concepts that model these computational effects 
possess the suitable monad structures, the \emph{theory of
monads} developed in pure mathematics can be fruitfully applied
to the study of computational effects.
Indeed, Moggi~\cite{moggi:comp-lambda,moggi:notions}
constructed categorical models of his calculus for computational
effects via the \emph{Kleisli construction}.

Nowadays, \emph{computational effects with parameters} are becoming
increasingly important.
For example, there are type systems called 
\emph{effect systems}~\cite{lucassen-gifford},
which statically estimate the range of effects
caused by a program.
Effect systems achieve this estimation by introducing,
instead of a single effect (e.g., global state),
a family of effects parametrized by {ranges} (e.g., global state
together with a set of registers specified, with the intuition being
that it contains those registers that are 
manipulated in an execution of the program).
On the other hand, there is an active line of research on the nature of
\emph{local state}~\cite{plotkin-power,power:indexed,power:lce,maillard-mellies}
(which, in addition to manipulating registers, 
is able to increase or decrease 
the number of registers in use),
especially on its relation to global state.
One attracting view~\cite{maillard-mellies} is 
that local state is obtained by 
``gluing'' a family of global states with different numbers of registers,
and there again the idea of parameters shows itself.

Corresponding to this \emph{parametrization} on the side of 
computational effects, various notions of 
\emph{monads with parameters} have been employed in the 
study of those computational effects with parameters.
The notion of \emph{graded monad} has been used to give
semantics of effect systems~\cite{katsumata},
whereas the notion of \emph{indexed monad} arose and has been applied
in the study of local state~\cite{power:indexed,power:lce,maillard-mellies}.

This thesis, by developing the fundamental mathematical theory of 
graded and indexed monads, aims to
lay the foundation of the 
\emph{unified study of computational effects
with parameters},
just as mathematical theory of monads has been the basis of
the unified study of computational effects initiated by Moggi.

\subsection{Street's 2-categorical study of monads}
\label{subsec-intro-math}
In 1972, Street~\cite{street:formal-theory-of-monads} 
presented a highly abstract \emph{2-categorical} framework 
that collects main results in the theory of monads 
obtained until then (including the 
Eilenberg--Moore and Kleisli constructions), and reconstructs them
inside an \emph{arbitrary 2-category} (subject to certain
completeness and cocompleteness requirements).
For the basic 2-categorical notions, see~\cite{kelly-street}.
In particular, he defined the
\emph{Eilenberg--Moore} and \emph{Kleisli objects} of 
a monad in a 2-category, indicating how the
Eilenberg--Moore and Kleisli constructions should be
generalized to 2-categories other than $\Cat$
(the 2-category of categories, functors, and natural transformations).

This framework of Street is particularly suited for
our study of graded and indexed monads, since it turns out that
the framework is general enough to include 
theories of graded and indexed monads 
as particular instances.
Indeed, the main theorems of this thesis
(Theorems~\ref{thm-univ-em-gm}, \ref{thm-univ-kl-gm} and
\ref{thm-univ-em-im}) amount to state that certain concrete
constructions defined in the thesis produce 
Eilenberg--Moore or Kleisli objects in the sense of Street,
providing a compelling justification for our constructions.

\section{Chapter overview}
We begin with the definition of graded and indexed monads in
Chapter~\ref{chap-param-monad}, together with 
examples illustrating them: the \emph{graded state monad}
(Section~\ref{subsec-gr-st-mnd})
and the \emph{indexed state monad}
(Section~\ref{subsec-in-st-mnd}), each extending the 
usual \emph{global state monad} in different directions.
The dual notions of graded and indexed comonad are also introduced.

Then, in Chapter~\ref{chap-four-2-cats},
we introduce four 2-categories by ``enlarging''
the familiar 2-category~$\Cat$ in various ways.
These 2-categories are now ``large'' enough to incorporate
the notions of graded and indexed monads as mere \emph{monads}
(in the 2-categorical sense) inside them.

The above observation has a key consequence that
it now makes perfect sense to speculate on 
\emph{Eilenberg--Moore} and \emph{Kleisli objects} of graded and 
indexed monads in the sense of Street,
by working inside these ``larger'' 2-categories instead of $\Cat$.
We show that indeed there do exist the following constructions
that produce the Eilenberg--Moore or
Kleisli objects:
\begin{itemize}
\item The Eilenberg--Moore construction for graded monads.
\item The Kleisli construction for graded monads.
\item The Eilenberg--Moore construction for indexed monads.
\end{itemize}
The definition, and establishment of the relevant 
universal property of them are the 
main tasks of Chapter~\ref{chap-main-constr},
which are also the main contributions of the current thesis.
We also present the corresponding constructions for
graded and indexed comonads.
The quest for the following construction is left as future work:
\begin{itemize}
\item The Kleisli construction for indexed monads.
\end{itemize}

We then present two applications of these constructions 
in Chapter~\ref{chap-appl}.
In Section~\ref{sec-decomp-lax-action} we apply 
the first two constructions,
the Eilenberg--Moore and Kleisli constructions for graded monads,
to the study of \emph{lax monoidal actions}.
We show how they give solutions to 
the problem of \emph{decomposing} lax monoidal actions
into \emph{strict monoidal actions} and \emph{adjunctions}, a
situation which generalizes the one that constitute the very
origin of the classical Eilenberg--Moore
and Kleisli constructions for ordinary 
monads~\cite{eilenberg-moore,kleisli}.
The second application, presented in 
Section~\ref{sec-constr-maillard-mellies},
employs the Eilenberg--Moore construction for indexed monads.
The story begins with the realization that the categories 
identical to our  
Eilenberg--Moore categories of indexed monads has already 
appeared in a recent paper~\cite{maillard-mellies} by
Maillard and Melli\`es.
They did not introduce these categories as Eilenberg--Moore
categories, but established a nice connection to the 
categories of \emph{models of indexed Lawvere theories}
defined by Power~\cite{power:lce,power:indexed}.
We show how the perspective provided by our construction enables one to
understand more conceptually this result of Maillard and Melli\`es. 

Finally in Chapter~\ref{chap-concl} we conclude the thesis and 
indicate possible directions for future work.

\section{Notation}
For the various compositions inside a 2-category, we adopt the following
notation: we denote composition of 1-cells by $\circ$,
vertical composition of 2-cells by $\vc$, 
and horizontal composition of 2-cells by $\hc$.
We often abbreviate \emph{whiskerings} 
$\id{g}\hc\alpha$ and $\beta\hc \id{f}$ by
$g\hc \alpha$ and $\beta\hc f$ respectively.
We also omit various composition symbols 
and write compositions simply by concatenation
when confusion is unlikely.

Given a 2-category~$\tcat{K}$, by $\tcat{K}^\op{1}$,
$\tcat{K}^\op{2}$, and $\tcat{K}^\op{1,2}$, we mean
the 2-categories obtained from $\tcat{K}$ by reversing, respectively,
1-cells but not 2-cells, 2-cells but not 1-cells,
and both 1-cells and 2-cells.
When $\tcat{K}$ is actually a category, we also write
$\tcat{K}^\op{1}$ as $\tcat{K}^\opp$. 

We denote an adjunction of type
\[
    \begin{tikzpicture}
      \node (A)  at (0,0) {$A$};
      \node (C) at (4,0) {$B$};
      
      \draw [->] (A) to [bend left=23]node (dom)[auto,labelsize] {$L$}(C);
      \draw [->] (C) to [bend left=23]node (cod)[auto,labelsize] {$R$}(A);
      
      \node at (2,0) [rotate=90,font=\large]{$\vdash$};

    \end{tikzpicture}
\]
by $L\dashv R\colon B\to A$.

\chapter{Adding parameters to monads}
\label{chap-param-monad}
In this chapter, we introduce two different ways to generalize
the notion of monad on a category 
so as to incorporate the intuitive idea of
\emph{parameters} into it.
These two ways give rise to the notions of \emph{graded monad}
and \emph{indexed monad} respectively.
The transition from monads to graded or
indexed monads is motivated by the various ways to
add parameters to computational effects.
We illustrate this point of view by two typical examples:
the \emph{graded state monad} and 
the \emph{indexed state monad}.

One can also argue that the notions of graded and indexed monad are
mathematically natural and are in some sense derived from a single
general principle;
see Appendix~\ref{apx-enriched-monad-tensor} for the discussion 
on this unified view.

\section{Graded monads}
A \emph{graded monad} takes its parameters from a 
\emph{monoidal category}.

\subsection{The definition}
\label{subsec-gm-defn}
\begin{defn}
Let $\cat{M}=(\cat{M},\tensor,\e)$ be a strict monoidal category
and $\cat{C}$ be a category.
An \defemph{$\cat{M}$-graded monad on $\cat{C}$}
is a (lax) monoidal functor 
\[
\gm{T}\quad\colon\quad (\cat{M},\tensor,\e)\quad\longrightarrow\quad
([\cat{C},\cat{C}],\circ,\id{\cat{C}}). \qedhere
\]
\end{defn}

To spell out the definition,
an $\cat{M}$-graded monad on $\cat{C}$
consists of the following data:
\begin{itemize}
\item A functor $\gm{T}_m\colon\cat{C}\to\cat{C}$ for each object $m$ 
of $\cat{M}$.
\item A natural transformation $\gm{T}_u\colon \gm{T}_m\Rightarrow \gm{T}_{m^\prime}$ for each
morphism~$u\colon m\to m^\prime$ of $\cat{M}$.
\item A natural transformation 
$\eta^{(\gm{T})}=\eta\colon \id{\cat{C}}\Rightarrow \gm{T}_\e$.
\item A natural transformation 
$\mu_{m,n}^{(\gm{T})}=\mu_{m,n}\colon \gm{T}_m\circ \gm{T}_n\Rightarrow \gm{T}_{m\tensor n}$ for each
pair of objects $m,n$ of $\cat{M}$.
\end{itemize}
These data are subject to the following axioms:
\begin{description}
\item[(GM1)] $\id{\gm{T}_m}=\gm{T}_{\id{m}}$ for each object~$m$ of $\cat{M}$.
\item[(GM2)] $\gm{T}_{u^\prime}\vc \gm{T}_u=\gm{T}_{u^\prime\circ u}$ for each composable
pair of morphisms~$u,u^\prime$ of $\cat{M}$. 
\item[(GM3)] $\begin{tikzpicture}[baseline=-\the\dimexpr\fontdimen22\textfont2\relax ]
      \node (TL)  at (0,0.75)  {$\gm{T}_m\circ \gm{T}_n$};
      \node (TR)  at (3,0.75)  {$\gm{T}_{m\tensor n}$};
      \node (BL)  at (0,-0.75) {$\gm{T}_{m^\prime}\circ \gm{T}_{n^\prime}$};
      \node (BR)  at (3,-0.75) {$\gm{T}_{m^\prime\tensor n^\prime}$};
      
      \draw [2cell] (TL) to node (T) [auto,labelsize]      {$\mu_{m,n}$} (TR);
      \draw [2cell] (TR) to node (R) [auto,labelsize]      {$\gm{T}_{u\tensor v}$} (BR);
      \draw [2cell] (TL) to node (L) [auto,swap,labelsize] {$\gm{T}_u\circ \gm{T}_v$}(BL);
      \draw [2cell] (BL) to node (B) [auto,swap,labelsize] {$\mu_{m^\prime,n^\prime}$}(BR);
\end{tikzpicture}$
commutes for each pair of morphisms $u\colon m\to m^\prime,
v\colon n\to n^\prime$ of $\cat{M}$.
\item[(GM4)] 
$\begin{tikzpicture}[baseline=-\the\dimexpr\fontdimen22\textfont2\relax ]
      \node (TL)  at (0,0.75)  {$\gm{T}_m$};
      \node (TR)  at (2,0.75)  {$\gm{T}_\e\circ \gm{T}_{m}$};
      \node (BR)  at (2,-0.75) {$\gm{T}_{m}$};
      
      \draw [2cell] (TL) to node (T) [auto,labelsize]      {$\eta\hc\gm{T}_m$} (TR);
      \draw [2cell] (TR) to node (R) [auto,labelsize]      {$\mu_{\e, m}$} (BR);
      \draw [2cell] (TL) to node (L) [auto,swap,labelsize] {$\id{\gm{T}_m}$}(BR);
\end{tikzpicture}$
commutes for each object~$m$ of $\cat{M}$.
\item[(GM5)]
$\begin{tikzpicture}[baseline=-\the\dimexpr\fontdimen22\textfont2\relax ]
      \node (TL)  at (0,0.75)  {$\gm{T}_m$};
      \node (TR)  at (2,0.75)  {$\gm{T}_m\circ \gm{T}_{\e}$};
      \node (BR)  at (2,-0.75) {$\gm{T}_{m}$};
      
      \draw [2cell] (TL) to node (T) [auto,labelsize]      {$\gm{T}_m\hc\eta$} (TR);
      \draw [2cell] (TR) to node (R) [auto,labelsize]      {$\mu_{m, \e}$} (BR);
      \draw [2cell] (TL) to node (L) [auto,swap,labelsize] {$\id{\gm{T}_m}$}(BR);
\end{tikzpicture}$
commutes for each object~$m$ of $\cat{M}$.
\item[(GM6)]
$\begin{tikzpicture}[baseline=-\the\dimexpr\fontdimen22\textfont2\relax ]
      \node (TL)  at (0,0.75)  {$\gm{T}_l\circ \gm{T}_m\circ \gm{T}_n$};
      \node (TR)  at (4,0.75)  {$\gm{T}_l\circ \gm{T}_{m\tensor n}$};
      \node (BL)  at (0,-0.75) {$\gm{T}_{l\tensor m}\circ \gm{T}_{n}$};
      \node (BR)  at (4,-0.75) {$\gm{T}_{l\tensor m\tensor n}$};
      
      \draw [2cell] (TL) to node (T) [auto,labelsize]      {$\gm{T}_l\hc\mu_{m,n}$} (TR);
      \draw [2cell] (TR) to node (R) [auto,labelsize]      {$\mu_{l,m\tensor n}$} (BR);
      \draw [2cell] (TL) to node (L) [auto,swap,labelsize] {$\mu_{l,m}\hc\gm{T}_n$}(BL);
      \draw [2cell] (BL) to node (B) [auto,swap,labelsize] {$\mu_{l\tensor m,n}$}(BR);
\end{tikzpicture}$
commutes for each triple of objects~$l$, $m$, $n$ of $\cat{M}$.
\end{description}

Graded monads generalize ordinary monads (monads in $\Cat$),
in the following sense.
\begin{prop}
Let $1$ be the terminal monoidal category and $\cat{C}$
a category.
A $1$-graded monad on $\cat{C}$ is nothing but a 
monad on $\cat{C}$.
\end{prop}

That monoidal functors from $1$ to $[\cat{C},\cat{C}]$ is the same as
monads on $\cat{C}$ has been known for a while; see for example \cite{benabou:bicat}.

\medbreak

Graded monads, which were previously also known as 
\emph{parametric monads}~\cite{mellies:tl-3,mellies:tl-8,katsumata},
arose in Melli\`es' study~\cite{mellies:tl-3} of 
\emph{Tensorial Logic}, a refinement of Linear Logic.
This notion was used by Katsumata~\cite{katsumata} 
in the study of semantics of 
\emph{effect systems}~\cite{lucassen-gifford,kammer-plotkin}, which 
are type systems to statically estimate the range of effects caused 
by a program (e.g., the set of registers that may be manipulated
during an execution of the program).
The intuition underlying Katsumata's work is that graded monads 
model \emph{computational effects with parameters},
extending Moggi's classical 
observation~\cite{moggi:comp-lambda,moggi:notions} that
monads model computational effects.

On the other hand,
the dual notion of \emph{graded comonad}, which will be introduced below,
have appeared in the study of \emph{computational resources with parameters},
as seen in e.g., the work of 
Petricek, Orchard and Mycroft~\cite{petricek-orchard-mycroft}.
Graded comonads also arise in the study of 
\emph{higher-order model checking}; see e.g.,
Grellois--Melli\`es~\cite{grellois-mellies}
and Tsukada--Ong~\cite{tsukada-ong}.

\subsection{An example: the graded state monad}
\label{subsec-gr-st-mnd}
One of the motivating examples of graded monads is the 
\emph{graded state monad}~\cite{fujii-katsumata-mellies},
which is an $\Inj$-graded monad on $\Set$.
Let us see its definition in detail.
We first define the monoidal category $\Inj=(\Inj,+,0)$.

\begin{defn}
Define the strict monoidal category~$(\Inj,+,0)$ as follows:
\begin{itemize}
\item An object of $\Inj$ is a natural number~$m=\{0,\dots,m-1\}$.
\item A morphism~$u$ of $\Inj$ from $m$ to $m^\prime$ is an injective 
function between sets
\[
u\quad\colon\quad\{0,\dots,m-1\}\quad\longrightarrow\quad
\{0,\dots,m^\prime-1\}.
\]
\item The monoidal product~$+$ on $\Inj$ is the restriction of
the binary coproduct on $\Set$.
\end{itemize}
Observe that the standard choice of 
$+$ indeed makes $(\Inj,+,0)$ a \emph{strict} monoidal category.
\end{defn}

\begin{defn}
Let $V$ be a set, interpreted as the \emph{set of values} that can be 
stored in one register.
The \defemph{graded state monad} with the 
set of values $V$ is 
the $\Inj$-graded monad~$\grst$ on $\Set$ defined as follows:
\begin{itemize}
\item For a natural number~$m$, define the functor~$
\grst_m\colon \Set\to\Set$ by 
\[
X\quad\longmapsto \quad V^m\Rightarrow (V^m\times X)\quad\cong\quad 
(V^m\Rightarrow V^m)\times(V^m\Rightarrow X),
\]
where $\Rightarrow$ denotes the exponential in $\Set$.
\item For an injective function~$u\colon m\to m^\prime$,
define the natural transformation $\grst_u\colon
\grst_m\Rightarrow \grst_{m^\prime}$ by setting 
its $X$-component, 
\[
\grst_{u,X}\quad\colon\quad
(V^m\Rightarrow V^m)\times(V^m\Rightarrow X)\quad
\longrightarrow\quad
(V^{m^\prime}\Rightarrow V^{m^\prime})\times(V^{m^\prime}\Rightarrow X),
\]
as $(\tau,\xi)\mapsto (u\bullet \tau, \xi\circ V^f)$, where
given a function $\tau\colon V^m\to V^m$, the function
$u\bullet \tau\colon V^{m^\prime}\to V^{m^\prime}$ is 
defined as 
\[
\big(\,(u\bullet \tau)(\vl_{k^\prime})_{k^\prime\in m^\prime}\big)_{l^\prime}
\ \coloneqq\ 
 \begin{cases} \big(\,\tau(\vl_{u(k)})_{k\in m}\big)_l 
 &\mbox{if } l^\prime = u(l)\mbox{ for some }l\in m;\\ 
 \,\vl_{l^\prime} & \mbox{otherwise.}\end{cases} 
\]
\item Define the natural transformation~$\eta^{(\grst)}\colon \id{\Set}
\Rightarrow \grst_0$ by setting its $X$-component,
\[
\eta^{(\grst)}_X\quad\colon\quad
X\quad\longrightarrow\quad(V^0\Rightarrow V^0)\times 
(V^0\Rightarrow X),
\]
to be the canonical bijection.
\item For a pair of natural numbers~$m,n$, define the natural 
transformation $\mu^{(\grst)}_{m,n}\colon \grst_m\circ \grst_n\Rightarrow\grst_{m+n}$,
by setting its $X$-component,
\begin{multline*}
\mu^{(\grst)}_{m,n,X}\ \colon\ 
(V^m\Rightarrow V^m)\times
(V^{m+n}\Rightarrow V^n)\times
(V^{m+n}\Rightarrow X)\\
\longrightarrow\quad
(V^{m+n}\Rightarrow V^{m+n})\times
(V^{m+n}\Rightarrow X),
\end{multline*}
as $(\tau, \sigma, \xi)\mapsto
(\tau\star \sigma, \xi)$, where
given functions~$\tau\colon V^m\to V^m$ and $\sigma\colon
V^{m+n}\to V^n$, 
the function~$\tau\star \sigma\colon V^{m+n}\to V^{m+n}$ is defined as
\begin{multline*}
\big(\,(\tau\star \sigma) (\vl_k, \wl_{k^\prime})_{k\in m, k^\prime\in n}\big)_{l}\\
 \coloneqq\quad
  \begin{cases} \big(\,\tau(\vl_{k})_{k\in m}\big)_l 
  &\mbox{if } l\in m;\\ 
  \big(\,\sigma(\vl_{k},\wl_{k^\prime})_{k\in m,k^\prime\in n}\big)_l  
  & \mbox{if } l\in n.\end{cases}
\tag*{\qedhere}
\end{multline*}
\end{itemize}
\end{defn}

\subsection{Graded comonads}
The dual notion of \emph{graded comonad} is of course
\begin{defn}
Let $\cat{M}=(\cat{M},\tensor,\e)$ be a strict monoidal category
and $\cat{C}$ be a category.
An \defemph{$\cat{M}$-graded comonad on $\cat{C}$}
is an oplax monoidal functor 
\[
\gm{S}\quad\colon \quad (\cat{M},\tensor,\e)\quad\longrightarrow\quad
([\cat{C},\cat{C}],\circ,\id{\cat{C}}). \qedhere
\]
\end{defn}

We write down the explicit description of an $\cat{M}$-graded comonad on $\cat{C}$
for the sake of concreteness; it consists of the following data:
\begin{itemize}
\item A functor $\gm{S}_m\colon\cat{C}\to\cat{C}$ for each object $m$ 
of $\cat{M}$.
\item A natural transformation $\gm{S}_u\colon \gm{S}_m\Rightarrow \gm{S}_{m^\prime}$ for each
morphism~$u\colon m\to m^\prime$ of $\cat{M}$.
\item A natural transformation 
$\varepsilon^{(\gm{S})}=\varepsilon\colon \gm{S}_\e\Rightarrow\id{\cat{C}}$.
\item A natural transformation 
$\delta_{m,n}^{(\gm{S})}=\delta_{m,n}\colon \gm{S}_{m\tensor n}\Rightarrow\gm{S}_m\circ \gm{S}_n$ for each
pair of objects $m,n$ of $\cat{M}$.
\end{itemize}
These data are subject to the following axioms:
\begin{description}
\item[(GC1)] $\id{\gm{S}_m}=\gm{S}_{\id{m}}$ for each object~$m$ of $\cat{M}$.
\item[(GC2)] $\gm{S}_{u^\prime}\vc \gm{S}_u=\gm{S}_{u^\prime\circ u}$ for each composable
pair of morphisms~$u,u^\prime$ of $\cat{M}$. 
\item[(GC3)] $\begin{tikzpicture}[baseline=-\the\dimexpr\fontdimen22\textfont2\relax ]
      \node (TL)  at (0,0.75)  {$\gm{S}_m\circ \gm{S}_n$};
      \node (TR)  at (-3,0.75)  {$\gm{S}_{m\tensor n}$};
      \node (BL)  at (0,-0.75) {$\gm{S}_{m^\prime}\circ \gm{S}_{n^\prime}$};
      \node (BR)  at (-3,-0.75) {$\gm{S}_{m^\prime\tensor n^\prime}$};
      
      \draw [2cellr] (TL) to node (T) [auto,swap,labelsize]      {$\delta_{m,n}$} (TR);
      \draw [2cell] (TR) to node (R) [auto,swap,labelsize]      {$\gm{S}_{u\tensor v}$} (BR);
      \draw [2cell] (TL) to node (L) [auto,labelsize] {$\gm{S}_u\circ \gm{S}_v$}(BL);
      \draw [2cellr] (BL) to node (B) [auto,labelsize] {$\delta_{m^\prime,n^\prime}$}(BR);
\end{tikzpicture}$
commutes for each pair of morphisms $u\colon m\to m^\prime,
v\colon n\to n^\prime$ of $\cat{M}$.
\item[(GC4)] 
$\begin{tikzpicture}[baseline=-\the\dimexpr\fontdimen22\textfont2\relax ]
      \node (TL)  at (0,0.75)  {$\gm{S}_m$};
      \node (TR)  at (2,0.75)  {$\gm{S}_\e\circ \gm{S}_{m}$};
      \node (BR)  at (2,-0.75) {$\gm{S}_{m}$};
      
      \draw [2cell] (TL) to node (T) [auto,labelsize]      {$\delta_{\e, m}$} (TR);
      \draw [2cell] (TR) to node (R) [auto,labelsize]      {$\varepsilon\hc\gm{S}_m$} (BR);
      \draw [2cell] (TL) to node (L) [auto,swap,labelsize] {$\id{\gm{S}_m}$}(BR);
\end{tikzpicture}$
commutes for each object~$m$ of $\cat{M}$.
\item[(GC5)]
$\begin{tikzpicture}[baseline=-\the\dimexpr\fontdimen22\textfont2\relax ]
      \node (TL)  at (0,0.75)  {$\gm{S}_m$};
      \node (TR)  at (2,0.75)  {$\gm{S}_m\circ \gm{S}_{\e}$};
      \node (BR)  at (2,-0.75) {$\gm{S}_{m}$};
      
      \draw [2cell] (TL) to node (T) [auto,labelsize]      {$\delta_{m, \e}$} (TR);
      \draw [2cell] (TR) to node (R) [auto,labelsize]      {$\gm{S}_m\hc\varepsilon$} (BR);
      \draw [2cell] (TL) to node (L) [auto,swap,labelsize] {$\id{\gm{S}_m}$}(BR);
\end{tikzpicture}$
commutes for each object~$m$ of $\cat{M}$.
\item[(GC6)]
$\begin{tikzpicture}[baseline=-\the\dimexpr\fontdimen22\textfont2\relax ]
      \node (TL)  at (0,-0.75)  {$\gm{S}_l\circ \gm{S}_m\circ \gm{S}_n$};
      \node (TR)  at (-4,-0.75)  {$\gm{S}_l\circ \gm{S}_{m\tensor n}$};
      \node (BL)  at (0,0.75) {$\gm{S}_{l\tensor m}\circ \gm{S}_{n}$};
      \node (BR)  at (-4,0.75) {$\gm{S}_{l\tensor m\tensor n}$};
      
      \draw [2cellr] (TL) to node (T) [auto,labelsize]      {$\gm{S}_l\hc\delta_{m,n}$} (TR);
      \draw [2cellr] (TR) to node (R) [auto,labelsize]      {$\delta_{l,m\tensor n}$} (BR);
      \draw [2cellr] (TL) to node (L) [auto,swap,labelsize] {$\delta_{l,m}\hc\gm{S}_n$}(BL);
      \draw [2cellr] (BL) to node (B) [auto,swap,labelsize] {$\delta_{l\tensor m,n}$}(BR);
\end{tikzpicture}$
commutes for each triple of objects~$l,m,n$ of $\cat{M}$.
\end{description}

It should be clear that $(\cat{M},\tensor,\e)$-graded comonads 
on $\cat{C}$ are
nothing but $(\cat{M}^\opp,\tensor,\e)$-graded monads on
$\cat{C}^\opp$.
Note that
we reversed the orientations of morphisms of $\cat{M}$ and $\cat{C}$,
but not the orientation of the monoidal product of $\cat{M}$.

\section{Indexed monads}
An \emph{indexed monad} takes its parameters from a 
\emph{category}.
\subsection{The definition}
We begin with a preliminary definition.
\begin{defn}
Let $\cat{C}$ be a category.
Define the category~$\Mnd{\cat{C}}$ of monads on $\cat{C}$ as follows:
\begin{itemize}
\item An object of $\Mnd{\cat{C}}$ is a monad~$T=(T,\eta,\mu)$ on 
$\cat{C}$.
\item A morphism of $\Mnd{\cat{C}}$ from $(T,\eta,\mu)$ to
$(T^\prime,\eta^\prime,\mu^\prime)$ is a 
natural transformation~$\tau\colon T^\prime\to T$
commuting with the monad structures, i.e., such that the diagrams
\[
\begin{tikzpicture}[baseline=-\the\dimexpr\fontdimen22\textfont2\relax ]
      \node (TL)  at (0,0.75)  {$\id{\cat{C}}$};
      \node (TR)  at (3,0.75)  {$\id{\cat{C}}$};
      \node (BL)  at (0,-0.75) {$T^\prime$};
      \node (BR)  at (3,-0.75) {$T$};
      
      \draw [2cell] (TL) to node (T) [auto,labelsize]      {$\id{\id{\cat{C}}}$} (TR);
      \draw [2cell] (TR) to node (R) [auto,labelsize]      {$\eta$} (BR);
      \draw [2cell] (TL) to node (L) [auto,swap,labelsize] {$\eta^\prime$}(BL);
      \draw [2cell] (BL) to node (B) [auto,swap,labelsize] {$\tau$}(BR);
\end{tikzpicture}
\qquad\quad
\begin{tikzpicture}[baseline=-\the\dimexpr\fontdimen22\textfont2\relax ]
      \node (TL)  at (0,0.75)  {$T^\prime\circ T^\prime$};
      \node (TR)  at (3,0.75)  {$T\circ T$};
      \node (BL)  at (0,-0.75) {$T^\prime$};
      \node (BR)  at (3,-0.75) {$T$};
      
      \draw [2cell] (TL) to node (T) [auto,labelsize]      {$\tau\hc\tau$} (TR);
      \draw [2cell] (TR) to node (R) [auto,labelsize]      {$\mu$} (BR);
      \draw [2cell] (TL) to node (L) [auto,swap,labelsize] {$\mu^\prime$}(BL);
      \draw [2cell] (BL) to node (B) [auto,swap,labelsize] {$\tau$}(BR);
\end{tikzpicture}
\]
commute.\qedhere
\end{itemize}
\end{defn}

The direction of morphisms of $\Mnd{\cat{C}}$ follows that of 
Street~\cite{street:formal-theory-of-monads}.
Note that this direction makes the ordinary Eilenberg--Moore construction
a covariant functor
\[
\EM{\cat{C}}{(-)}\quad\colon\quad\Mnd{\cat{C}}\quad\longrightarrow
\quad\Cat,
\]
with $\cat{C}^\tau\colon\cat{C}^T\to\cat{C}^{T^\prime}$
given by $\valgp{Tc}{\chi}{c}\mapsto
\valgp{T^\prime c}{\chi\circ \tau_c}{c}$.

\begin{defn}
Let $\cat{B}$ and $\cat{C}$ be categories.
A \defemph{$\cat{B}$-indexed monad on $\cat{C}$} is a functor
\[
\im{T}\quad\colon \quad\cat{B}^\opp\quad\longrightarrow\quad\Mnd{\cat{C}}.
\qedhere
\] 
\end{defn}

Explicitly, a $\cat{B}$-indexed monad on $\cat{C}$ consists of the 
following data:
\begin{itemize}
\item A functor~$\imc{b}\colon \cat{C}\to\cat{C}$ for each object~$b$
of $\cat{B}$.
\item A natural transformation~$\imc{u}\colon \imc{b}\Rightarrow\imc{b^\prime}$
for each morphism~$u\colon b\to b^\prime$ of $\cat{B}$.
\item A natural transformation~$\eta^{(\im{T})}_b
=\eta_b\colon \id{\cat{C}}
\Rightarrow\imc{b}$ for each object~$b$ of $\cat{B}$.
\item A natural transformation~$\mu^{(\im{T})}_b
=\mu_b\colon\imc{b}\circ\imc{b}\Rightarrow
\imc{b}$ for each object~$b$ of $\cat{B}$.
\end{itemize}
These data are subject to the following axioms:
\begin{description}
\item[(IM1)] $\id{\imc{b}}=\imc{\id{b}}$ for each object~$b$ of $\cat{B}$.
\item[(IM2)] $\imc{u^\prime}\vc\imc{u}=\imc{u^\prime\circ u}$ for each composable 
pair of morphisms~$u,u^\prime$ of $\cat{B}$.
\item[(IM3)]
$\begin{tikzpicture}[baseline=-\the\dimexpr\fontdimen22\textfont2\relax ]
      \node (TL)  at (0,0.75)  {$\id{\cat{C}}$};
      \node (TR)  at (3,0.75)  {$\id{\cat{C}}$};
      \node (BL)  at (0,-0.75) {$\imc{b}$};
      \node (BR)  at (3,-0.75) {$\imc{b^\prime}$};
      
      \draw [2cell] (TL) to node (T) [auto,labelsize]      {$\id{\id{\cat{C}}}$} (TR);
      \draw [2cell] (TR) to node (R) [auto,labelsize]      {$\eta_{b^\prime}$} (BR);
      \draw [2cell] (TL) to node (L) [auto,swap,labelsize] {$\eta_b$}(BL);
      \draw [2cell] (BL) to node (B) [auto,swap,labelsize] {$\imc{u}$}(BR);
\end{tikzpicture}$
commutes for each morphism~$u\colon b\to b^\prime$ of $\cat{B}$.
\item[(IM4)]
$\begin{tikzpicture}[baseline=-\the\dimexpr\fontdimen22\textfont2\relax ]
      \node (TL)  at (0,0.75)  {$\imc{b}\circ\imc{b}$};
      \node (TR)  at (3,0.75)  {$\imc{b^\prime}\circ\imc{b^\prime}$};
      \node (BL)  at (0,-0.75) {$\imc{b}$};
      \node (BR)  at (3,-0.75) {$\imc{b^\prime}$};
      
      \draw [2cell] (TL) to node (T) [auto,labelsize]      {$\imc{u}\hc\imc{u}$} (TR);
      \draw [2cell] (TR) to node (R) [auto,labelsize]      {$\mu_{b^\prime}$} (BR);
      \draw [2cell] (TL) to node (L) [auto,swap,labelsize] {$\mu_b$}(BL);
      \draw [2cell] (BL) to node (B) [auto,swap,labelsize] {$\imc{u}$}(BR);
\end{tikzpicture}$
\!\!commutes for each morphism~$u\colon b\to b^\prime$ of $\cat{B}$.
\item[(IM5)] $\begin{tikzpicture}[baseline=-\the\dimexpr\fontdimen22\textfont2\relax ]
      \node (TL)  at (0,0.75)  {$\imc{b}$};
      \node (TR)  at (2,0.75)  {$\imc{b}\circ\imc{b}$};
      \node (BR)  at (2,-0.75) {$\imc{b}$};
      
      \draw [2cell] (TL) to node (T) [auto,labelsize]      {$\eta_b\hc\imc{b}$} (TR);
      \draw [2cell] (TR) to node (R) [auto,labelsize]      {$\mu_b$} (BR);
      \draw [2cell] (TL) to node (L) [auto,swap,labelsize] {$\id{\imc{b}}$}(BR);
\end{tikzpicture}$
commutes for each object~$b$ of $\cat{B}$.
\item[(IM6)]
$\begin{tikzpicture}[baseline=-\the\dimexpr\fontdimen22\textfont2\relax ]
      \node (TL)  at (0,0.75)  {$\imc{b}$};
      \node (TR)  at (2,0.75)  {$\imc{b}\circ \imc{b}$};
      \node (BR)  at (2,-0.75) {$\imc{b}$};
      
      \draw [2cell] (TL) to node (T) [auto,labelsize]      {$\imc{b}\hc\eta_b$} (TR);
      \draw [2cell] (TR) to node (R) [auto,labelsize]      {$\mu_b$} (BR);
      \draw [2cell] (TL) to node (L) [auto,swap,labelsize] {$\id{\imc{b}}$}(BR);
\end{tikzpicture}$
commutes for each object~$b$ of $\cat{B}$.
\item[(IM7)] $\begin{tikzpicture}[baseline=-\the\dimexpr\fontdimen22\textfont2\relax ]
      \node (TL)  at (0,0.75)  {$\imc{b}\circ \imc{b}\circ \imc{b}$};
      \node (TR)  at (4,0.75)  {$\imc{b}\circ \imc{b}$};
      \node (BL)  at (0,-0.75) {$\imc{b}\circ \imc{b}$};
      \node (BR)  at (4,-0.75) {$\imc{b}$};
      
      \draw [2cell] (TL) to node (T) [auto,labelsize]      {$\imc{b}\hc\mu_{b}$} (TR);
      \draw [2cell] (TR) to node (R) [auto,labelsize]      {$\mu_{b}$} (BR);
      \draw [2cell] (TL) to node (L) [auto,swap,labelsize] {$\mu_{b}\hc\imc{b}$}(BL);
      \draw [2cell] (BL) to node (B) [auto,swap,labelsize] {$\mu_{b}$}(BR);
\end{tikzpicture}$
commutes for each object~$b$ of $\cat{B}$.
\end{description}

Indexed monads also provide a generalization of 
ordinary monads:
\begin{prop}
Let $1$ be the terminal category and $\cat{C}$
a category.
A $1$-indexed monad on $\cat{C}$ is nothing but a 
monad on $\cat{C}$.
\end{prop}

Indexed monads
have also arisen in the study of computational effects recently.
They appeared, in the form of \emph{indexed Lawvere theory},
in the work of Power~\cite{power:indexed,power:lce} on
the \emph{local state monad}.
Maillard and Melli\`es~\cite{maillard-mellies} shed new light on
Power's work using \emph{indexed monads} in their sense;
our notion of \emph{$\cat{B}$-indexed monad on $\cat{C}$}
is a restriction of their notion of \emph{$\cat{B}$-indexed monad} 
(without further qualification).
The primary example of indexed monads used in this series of 
study is the \emph{indexed state monad} which we recall shortly,
and can be intuitively thought of as a family of the state monads
indexed by (parametrized by) the number of registers;
so indexed monads also model the idea of \emph{computational
effects with parameters}.
We review in more detail Power's and Maillard and Melli\`es'
work in Section~\ref{sec-constr-maillard-mellies},
where we relate their results to ours.

\subsection{An example: the indexed state monad}
\label{subsec-in-st-mnd}
The indexed state monad~\cite{power:indexed,power:lce,maillard-mellies}
is an $\Inj$-indexed monad on $\Set$.
Note that this time we regard $\Inj$ as a mere category
rather than a strict monoidal category~$(\Inj,+,0)$.

\begin{defn}
Let $V$ be a set, interpreted again as the set of values that 
can be stored in one register.
The \defemph{indexed state monad} with the set of values~$V$
is the $\Inj$-indexed monad $\inst$ defined as follows:
\begin{itemize}
\item For a natural number~$m$, define the 
functor~$\inst_m\colon\Set\to\Set$ by
$\inst_m=\grst_m$.
\item For an injective function~$u\colon m\to m^\prime$, define the 
natural transformation
$\inst_u\colon\inst_m\Rightarrow\inst_{m^\prime}$ by
$\inst_u=\grst_u$.
\item For a natural number~$m$, define the natural transformation~$
\eta^{(\inst)}_m\colon \id{\Set}\Rightarrow\inst_m$ by
setting its $X$-component,
\[
\eta^{(\inst)}_{m,X}\quad\colon\quad
X\quad\longrightarrow\quad(V^m\Rightarrow V^m)\times (V^m\Rightarrow
X),
\] 
as $x\mapsto (\id{V^m},\const{x}=(-\mapsto x))$.
\item For a natural number~$m$, define the natural transformation~$
\mu^{(\inst)}_m\colon \inst_m\circ\inst_m\Rightarrow\inst_m$
by setting its $X$-component,
\begin{multline*}
\mu^{(\inst)}_{m,X}\ \colon\ 
(V^m\Rightarrow V^m)\times
(V^{m}\Rightarrow (V^m\Rightarrow V^m))\times
(V^{m}\Rightarrow (V^m\Rightarrow X))\\
\longrightarrow\quad
(V^{m}\Rightarrow V^{m})\times
(V^{m}\Rightarrow X),
\end{multline*}
by $(\tau,\sigma,\xi)\mapsto (\tau\rhd\sigma,\tau\blacktriangleright\xi)$,
where
\begin{itemize}
\item given functions~$\tau\colon V^m\to V^m$ and $\sigma\colon V^m\to
(V^m\Rightarrow V^m)$, the function~$\tau\rhd\sigma\colon V^m\to V^m$
is defined as 
\[
(\tau\rhd\sigma)(v)\quad\coloneqq\quad \sigma (v) (\tau v);
\] 
\item given functions~$\tau\colon V^m\to V^m$ and 
$\xi\colon V^m\to (V^m\Rightarrow X)$, the function~$
\tau\blacktriangleright\xi\colon V^m\to X$ is defined as 
\[
(\tau\blacktriangleright\xi)(v)\quad\coloneqq\quad
\xi(v)(\tau v).\qedhere
\]
\end{itemize}
\end{itemize}
\end{defn}

Power~\cite{power:lce,power:indexed} derives the indexed Lawvere theory
equivalent to the indexed state monad abstractly from the global state monad,
using the universality of $(\Inj,+,0)$ as the 
\emph{free monoidal category with the initial unit}, 
and the \emph{tensor product} operation naturally defined on Lawvere
theories.
Now the reader might noticed the obvious
similarity between the indexed state monad and the 
graded state monad in Section~\ref{subsec-gr-st-mnd}.
Indeed, one can also derive the graded state monad from the indexed
state monad abstractly; the requirement is again 
that $(\Inj,+,0)$ has the initial unit.

\begin{prop}
\label{prop-from-im-to-gm}
Let $\cat{M}=(\cat{M},\tensor,\e)$ be a strict monoidal category
such that the monoidal unit~$\e$ is the initial object of $\cat{M}$,
and $\cat{C}$ be a category.
Then every $\cat{M}$-indexed monad~$\im{T}=(\im{T},\eta^{(\im{T})},
\mu^{(\im{T})})$ on $\cat{C}$ naturally induces an $(\cat{M},\tensor,\e)$-graded
monad~$\gm{T}=(\gm{T},\eta^{(\gm{T})},\mu^{(\gm{T})})$
on $\cat{C}$.
\end{prop}
\begin{pf}
First define $\gm{T}_m\coloneqq \im{T}_m$ and $\gm{T}_u\coloneqq \im{T}_u$
for each object~$m$ and morphism~$u$ of $\cat{M}$.

Then define $\eta^{(\gm{T})}\colon \id{\cat{C}}\Rightarrow\gm{T}_\e$
by $\eta^{(\gm{T})}\coloneqq \eta^{(\im{T})}_\e$.

Finally, we have to define $\mu^{(\gm{T})}_{m,n}\colon \gm{T}_m\circ
\gm{T}_n\Rightarrow\gm{T}_{m\tensor n}$.
To this purpose, first observe that there are morphisms
\begin{align*}
\inl_{m,n}\quad&\colon\quad m\quad=\quad m\tensor \e \quad\xrightarrow{\id{m}\tensor !_n}\quad m\tensor n
\\
\inr_{m,n}\quad&\colon\quad n\quad=\quad \e\tensor n \quad\xrightarrow{!_m\tensor \id{n}}\quad m\tensor n
\end{align*}
of $\cat{M}$.
Now define
\[
\mu^{(\gm{T})}_{m,n}\quad\coloneqq\quad
\im{T}_m\circ\im{T}_n\xRightarrow{\im{T}_{\inl_{m,n}}\hc\im{T}_{\inl_{m,n}}}
\im{T}_{m\tensor n}\circ\im{T}_{m\tensor n}
\xRightarrow{\mu^{(\im{T})}_{m\tensor n}}\im{T}_{m\tensor n}.\qedhere
\]
\end{pf}

This gives an abstract explanation of how our leading examples 
of ``monads with parameters'', the graded and indexed 
state monads, are obtained from the global state monad.

\subsection{Indexed comonads}

\begin{defn}
Let $\cat{C}$ be a category.
Define the category~$\CoMnd{\cat{C}}$ of comonads on $\cat{C}$ as follows:
\begin{itemize}
\item An object of $\CoMnd{\cat{C}}$ is a comonad~$S=(S,\varepsilon,\delta)$ on 
$\cat{C}$.
\item A morphism of $\CoMnd{\cat{C}}$ from $(S,\varepsilon,\delta)$ to
$(S^\prime,\varepsilon^\prime,\delta^\prime)$ is a 
natural transformation~$\sigma\colon S\to S^\prime$
commuting with the comonad structures, i.e., such that the diagrams
\[
\begin{tikzpicture}[baseline=-\the\dimexpr\fontdimen22\textfont2\relax ]
      \node (TL)  at (0,-0.75)  {$\id{\cat{C}}$};
      \node (TR)  at (-3,-0.75)  {$\id{\cat{C}}$};
      \node (BL)  at (0,0.75) {$S^\prime$};
      \node (BR)  at (-3,0.75) {$S$};
      
      \draw [2cellr] (TL) to node (T) [auto,labelsize]      {$\id{\id{\cat{C}}}$} (TR);
      \draw [2cellr] (TR) to node (R) [auto,labelsize]      {$\varepsilon$} (BR);
      \draw [2cellr] (TL) to node (L) [auto,swap,labelsize] {$\varepsilon^\prime$}(BL);
      \draw [2cellr] (BL) to node (B) [auto,swap,labelsize] {$\sigma$}(BR);
\end{tikzpicture}
\qquad\quad
\begin{tikzpicture}[baseline=-\the\dimexpr\fontdimen22\textfont2\relax ]
      \node (TL)  at (0,-0.75)  {$S^\prime\circ S^\prime$};
      \node (TR)  at (-3,-0.75)  {$S\circ S$};
      \node (BL)  at (0,0.75) {$S^\prime$};
      \node (BR)  at (-3,0.75) {$S$};
      
      \draw [2cellr] (TL) to node (T) [auto,labelsize]      {$\sigma\hc\sigma$} (TR);
      \draw [2cellr] (TR) to node (R) [auto,labelsize]      {$\delta$} (BR);
      \draw [2cellr] (TL) to node (L) [auto,swap,labelsize] {$\delta^\prime$}(BL);
      \draw [2cellr] (BL) to node (B) [auto,swap,labelsize] {$\sigma$}(BR);
\end{tikzpicture}
\]
commute.\qedhere
\end{itemize}
\end{defn}

This definition makes the co-Eilenberg--Moore construction for comonads
a covariant functor
\[
\EM{\cat{C}}{(-)}\quad\colon\quad\CoMnd{\cat{C}}\quad\longrightarrow
\quad\Cat,
\]
with $\cat{C}^\sigma\colon\cat{C}^S\to\cat{C}^{S^\prime}$
given by $\vcoalgp{Sc}{\chi}{c}\mapsto
\vcoalgp{S^\prime c}{\sigma_c\circ\chi}{c}$.

\begin{defn}
Let $\cat{B}$ and $\cat{C}$ be categories.
A \defemph{$\cat{B}$-indexed comonad on $\cat{C}$} is a functor
\[
\im{S}\quad\colon \quad\cat{B}\quad\longrightarrow\quad\CoMnd{\cat{C}}.
\qedhere
\] 
\end{defn}

A $\cat{B}$-indexed comonad on $\cat{C}$ consists of the 
following data:
\begin{itemize}
\item A functor~$\icc{b}\colon \cat{C}\to\cat{C}$ for each object~$b$
of $\cat{B}$.
\item A natural transformation~$\icc{u}\colon \icc{b}\Rightarrow\icc{b^\prime}$
for each morphism~$u\colon b\to b^\prime$ of $\cat{B}$.
\item A natural transformation~$\varepsilon^{(\im{T})}_b
=\varepsilon_b\colon \icc{b}
\Rightarrow\id{\cat{C}}$ for each object~$b$ of $\cat{B}$.
\item A natural transformation~$\delta^{(\im{T})}_b
=\delta_b\colon\icc{b}\Rightarrow
\icc{b}\circ\icc{b}$ for each object~$b$ of $\cat{B}$.
\end{itemize}
These data are subject to the following axioms:
\begin{description}
\item[(IC1)] $\id{\icc{b}}=\icc{\id{b}}$ for each object~$b$ of $\cat{B}$.
\item[(IC2)] $\icc{u^\prime}\vc\icc{u}=\icc{u^\prime\circ u}$ for each composable 
pair of morphisms~$u,u^\prime$ of $\cat{B}$.
\item[(IC3)]
$\begin{tikzpicture}[baseline=-\the\dimexpr\fontdimen22\textfont2\relax ]
      \node (TL)  at (0,-0.75)  {$\id{\cat{C}}$};
      \node (TR)  at (3,-0.75)  {$\id{\cat{C}}$};
      \node (BL)  at (0,0.75) {$\icc{b}$};
      \node (BR)  at (3,0.75) {$\icc{b^\prime}$};
      
      \draw [2cell] (TL) to node (T) [auto,swap,labelsize]      {$\id{\id{\cat{C}}}$} (TR);
      \draw [2cellr] (TR) to node (R) [auto,swap,labelsize]      {$\varepsilon_{b^\prime}$} (BR);
      \draw [2cellr] (TL) to node (L) [auto,labelsize] {$\varepsilon_b$}(BL);
      \draw [2cell] (BL) to node (B) [auto,labelsize] {$\icc{u}$}(BR);
\end{tikzpicture}$
commutes for each morphism~$u\colon b\to b^\prime$ of $\cat{B}$.
\item[(IC4)]
$\begin{tikzpicture}[baseline=-\the\dimexpr\fontdimen22\textfont2\relax ]
      \node (TL)  at (0,-0.75)  {$\icc{b}\circ\icc{b}$};
      \node (TR)  at (3,-0.75)  {$\icc{b^\prime}\circ\icc{b^\prime}$};
      \node (BL)  at (0,0.75) {$\icc{b}$};
      \node (BR)  at (3,0.75) {$\icc{b^\prime}$};
      
      \draw [2cell] (TL) to node (T) [auto,swap,labelsize]      {$\icc{u}\hc\icc{u}$} (TR);
      \draw [2cellr] (TR) to node (R) [auto,swap,labelsize]      {$\delta_{b^\prime}$} (BR);
      \draw [2cellr] (TL) to node (L) [auto,labelsize] {$\delta_b$}(BL);
      \draw [2cell] (BL) to node (B) [auto,labelsize] {$\icc{u}$}(BR);
\end{tikzpicture}$
\!\!commutes for each morphism~$u\colon b\to b^\prime$ of $\cat{B}$.
\item[(IC5)] $\begin{tikzpicture}[baseline=-\the\dimexpr\fontdimen22\textfont2\relax ]
      \node (TL)  at (0,0.75)  {$\icc{b}$};
      \node (TR)  at (2,0.75)  {$\icc{b}\circ\icc{b}$};
      \node (BR)  at (2,-0.75) {$\icc{b}$};
      
      \draw [2cell] (TL) to node (T) [auto,labelsize]      {$\delta_b$} (TR);
      \draw [2cell] (TR) to node (R) [auto,labelsize]      {$\varepsilon_b\hc\icc{b}$} (BR);
      \draw [2cell] (TL) to node (L) [auto,swap,labelsize] {$\id{\icc{b}}$}(BR);
\end{tikzpicture}$
commutes for each object~$b$ of $\cat{B}$.
\item[(IC6)]
$\begin{tikzpicture}[baseline=-\the\dimexpr\fontdimen22\textfont2\relax ]
      \node (TL)  at (0,0.75)  {$\icc{b}$};
      \node (TR)  at (2,0.75)  {$\icc{b}\circ \icc{b}$};
      \node (BR)  at (2,-0.75) {$\icc{b}$};
      
      \draw [2cell] (TL) to node (T) [auto,labelsize]      {$\delta_b$} (TR);
      \draw [2cell] (TR) to node (R) [auto,labelsize]      {$\icc{b}\hc\varepsilon_b$} (BR);
      \draw [2cell] (TL) to node (L) [auto,swap,labelsize] {$\id{\icc{b}}$}(BR);
\end{tikzpicture}$
commutes for each object~$b$ of $\cat{B}$.
\item[(IC7)] $\begin{tikzpicture}[baseline=-\the\dimexpr\fontdimen22\textfont2\relax ]
      \node (TL)  at (0,-0.75)  {$\icc{b}\circ \icc{b}\circ \icc{b}$};
      \node (TR)  at (-4,-0.75)  {$\icc{b}\circ \icc{b}$};
      \node (BL)  at (0,0.75) {$\icc{b}\circ \icc{b}$};
      \node (BR)  at (-4,0.75) {$\icc{b}$};
      
      \draw [2cellr] (TL) to node (T) [auto,labelsize]      {$\icc{b}\hc\delta_{b}$} (TR);
      \draw [2cellr] (TR) to node (R) [auto,labelsize]      {$\delta_{b}$} (BR);
      \draw [2cellr] (TL) to node (L) [auto,swap,labelsize] {$\delta_{b}\hc\icc{b}$}(BL);
      \draw [2cellr] (BL) to node (B) [auto,swap,labelsize] {$\delta_{b}$}(BR);
\end{tikzpicture}$
commutes for each object~$b$ of $\cat{B}$.
\end{description}

Observe that $\cat{B}$-indexed comonads on $\cat{C}$ are equivalent to
$\cat{B}^\opp$-indexed monads on $\cat{C}^\opp$.

\section*{Notes}
For the reason why we have decided to change the terminology from
\emph{parametric monads} to \emph{graded monads},
see Section~2 of \cite{fujii-katsumata-mellies}.
The graded state monad, which appears in \cite{fujii-katsumata-mellies},
was introduced by Paul-Andr\'e Melli\`es.
I learned the abstract construction of the graded state monad 
from the indexed state monad (Proposition~\ref{prop-from-im-to-gm})
from the anonymous reviewer of the paper~\cite{fujii-katsumata-mellies}.

\chapter{The four 2-categories}
\label{chap-four-2-cats}
In this chapter, we introduce four 2-categories
$\Epp$, $\Epm$, $\Emp$ and $\Emm$.
These 2-categories enable us to regard
graded and indexed monads as mere monads in the 2-categorical 
sense, thus paving the way to a mathematical theory of 
graded and indexed monads inside the celebrated abstract framework 
of Street~\cite{street:formal-theory-of-monads}.
More precisely, the relationship of 
the notions of graded and indexed (co)monad
and the 2-categories $\Epp$, $\Epm$, $\Emp$ and $\Emm$ is 
summarized in the following table:
\begin{table}[H]
	\centering
	\begin{tabular}{c|c c c c}    
    & $\Epp$ & $\Epm$ & $\Emp$ & $\Emm$\\ \hline
    Graded monads    & \checkmark &    &    & \checkmark \\
    Indexed monads   &    & \checkmark & \checkmark & \\
    Graded comonads  &    & \checkmark & \checkmark & \\
    Indexed comonads & \checkmark & & & \checkmark
	\end{tabular}
	\vspace{-15pt}
\end{table}
\noindent
The \checkmark~mark indicates that the  
row notion can be seen 
as mere (co)monads in the column 2-category.




\section{The formal theory of monads}
\label{sec-ftm}
In his seminal and influential paper~\cite{street:formal-theory-of-monads}, Street developed an
abstract theory of monads relative to an arbitrary 2-category~$\tcat{K}$, so
that the usual theory of monads is regained by instantiating $\tcat{K}$ by
$\Cat$, the 2-category of categories, functors, and natural
transformations.
Since our principal justification of the constructions 
presented in Chapter~\ref{chap-main-constr} will be that 
they produce \emph{Eilenberg--Moore} or \emph{Kleisli objects},
which are among the key notions introduced in \cite{street:formal-theory-of-monads},
we begin this chapter with a quick review of Street's work.

\begin{defn}
\label{defn-monad}
Let $\tcat{K}$ be a 2-category and $K$ a 0-cell of $\tcat{K}$.
  A \defemph{monad~${T}$ in $\tcat{K}$ on $K$} is given by 
  a 1-cell $T\colon K\to K$, and 2-cells 
  $\eta\colon \id{K}\Rightarrow T$ and
  $\mu\colon T\circ T\Rightarrow T$ of $\tcat{K}$, satisfying
  the usual axioms
  $\mu\vc\eta T=\id{T}=\mu\vc T\eta$ and
  $\mu\vc T\mu=\mu\vc \mu T$.
\end{defn}
We fix a monad~$T$ in $\tcat{K}$ on $K$.
For the definition of Eilenberg--Moore and 
Kleisli objects, we adopt the following 
\emph{representable} one, which appears 
for instance in~\cite{street:two-constructions}:
\begin{defn}
\label{defn-EM-obj}
  The \defemph{Eilenberg--Moore object} of ${T}$
  is a 0-cell $K^{T}$ such that 
  there is a family of
  isomorphisms of categories
  \[
  \tcat{K}(X,K^{T})\quad\cong\quad \tcat{K}(X,K)^{\tcat{K}(X,{T})}
  \]
  2-natural in $X\in \tcat{K}$.
  Here the category on the right hand side is the 
  (usual) Eilenberg--Moore category of the 
  monad~$\tcat{K}(X,{T})$ on the category~$\tcat{K}(X,K)$.
\end{defn}
\begin{defn}
\label{defn-Kl-obj}
  The \defemph{Kleisli object} of ${T}$
  is a 0-cell $K_{T}$ such that
  there is a family of
  isomorphisms of categories
  \[
  \tcat{K}(K_{T},X)\quad\cong\quad \tcat{K}(K,X)^{\tcat{K}({T},X)}
  \]
  2-natural in $X\in \tcat{K}$.
  Here the category on the right hand side is the 
  (usual) Eilenberg--Moore category of the 
  monad~$\tcat{K}({T},X)$ on the category~$\tcat{K}(K,X)$.
\end{defn}

Noting that a monad in $\tcat{K}$ is the same thing as 
a monad in $\tcat{K}^\op{1}$, the Kleisli 
object of $T$ in $\tcat{K}$ can be equivalently defined as the 
Eilenberg--Moore object of $T$ in $\tcat{K}^\op{1}$.

Now a remarkable point is that from this simple
and abstract definition, one can reconstruct 
a fair amount of the well-known property of 
Eilenberg--Moore or Kleisli categories,
including the existence of adjunctions which
generate the monads, and the existence and uniqueness
of comparison 1-cells.
The interested reader should consult
\cite{street:formal-theory-of-monads}
for an ingenious 2-categorical manipulation
achieving this reconstruction.

\section{Enlarging $\Cat$}

The four 2-categories~$\Epp$, $\Epm$, $\Emp$ and $\Emm$ we now introduce
are obtained via suitable lax comma constructions applied to the 
3-category~$\TCat$ of 2-categories,
2-functors, 2-natural transformations and modifications.  
The intuitive idea is that these four 2-categories are obtained by
``enlarging'' the familiar 2-category~$\Cat$;
actually, they ``contain'' an arbitrary 2-category, in the sense 
that for any 2-category~$\tcat{K}$ (subject to a certain
size condition) there are canonical 
inclusion 2-functors $\tcat{K}\hookrightarrow\Epp$, 
$\tcat{K}\hookrightarrow\Epm$, $\tcat{K}\hookrightarrow\Emp$
and $\tcat{K}\hookrightarrow\Emm$.

In the following definitions, 
we denote the terminal 2-category (the 2-category consisting of 
a single 0-cell, a single 1-cell and a single 2-cell) by $1$.

\subsection{The 2-category $\Epp$}
\begin{defn}
  We define the 2-category~$\Epp$ by the following data.
  \begin{itemize}
  \item A 0-cell of $\Epp$ is a 2-functor
    $A\colon  1\to \tcat{A}$ where $\tcat{A}$ is a 2-category;
    equivalently, it is a pair $(\tcat{A},A)$ where 
    $A$ is a 0-cell of $\tcat{A}$.
  \item A 1-cell of $\Epp$ from $(\tcat{A},A)$ to 
    $(\tcat{B},B)$ is a diagram in $\TCat$
    \[
    \begin{tikzpicture}
      \node (C)  at (1,1.5) {$1$};
      \node (C2)  at (0,0) {$\tcat{A}$};
      \node (D)  at (2,0) {$\tcat{B}$};
      
      \draw [->] (C)  to [bend left=0]  node [auto,swap,labelsize] {$A$} (C2);
      \draw [->] (C) to [bend left=0]node [auto,labelsize] {$B$}(D);
      \draw [->] (C2) to [bend left=0]node [auto,swap,labelsize] {$F$}(D);
      
      \draw [2cell] (0.6,0.5) to node [auto,labelsize] {$f$} (1.4,0.5);
    \end{tikzpicture}
    \]
    where $F$ is a 2-functor and 
    $f$ a 2-natural transformation;
    equivalently, it is a pair $(F,f)$ where
    $f\colon  FA\to B$ is a 1-cell of
    $\tcat{B}$.
  \item A 2-cell of $\Epp$ from $(F,f)$ to 
    $(F^\prime,f^\prime)$ is a diagram in $\TCat$
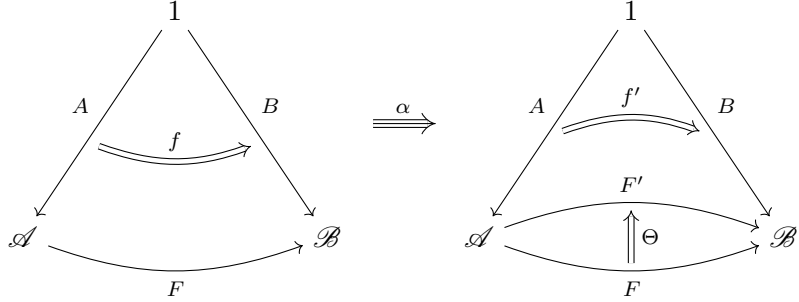
    \[
    \begin{tikzpicture}
      \node (C)  at (2,3) {$1$};
      \node (C2)  at (0,0) {$\tcat{A}$};
      \node (D)  at (4,0) {$\tcat{B}$};
      
      \draw [->] (C)  to [bend left=0]  node [auto,swap,labelsize] {$A$} (C2);
      \draw [->] (C) to [bend left=0]node [auto,labelsize] {$B$}(D);
      \draw [->] (C2) to [bend right=20]node [auto,swap,labelsize] {$F$}(D);
      
      \draw [2cell] (1,1.2) to [bend right=20] node [auto,labelsize] {$f$} (3,1.2);
      
      
      \draw [3cell] (4.5,1.5) to [bend left=0] node [auto,labelsize] {$\alpha$} (5.5,1.5);
            \draw (4.5,1.5)+(3pt,0pt) to  (5.4,1.5);
      
      \begin{scope}[shift={(6,0)}]
	  \node (C)  at (2,3) {$1$};
      \node (C2)  at (0,0) {$\tcat{A}$};
      \node (D)  at (4,0) {$\tcat{B}$};
      
      \draw [->] (C)  to [bend left=0]  node [auto,swap,labelsize] {$A$} (C2);
      \draw [->] (C) to [bend left=0]node [auto,labelsize] {$B$}(D);
      \draw [->] (C2) to [bend right=20]node [auto,swap,labelsize] {$F$}(D);
      \draw [->] (C2) to [bend left=20]node [auto,labelsize] {$F^\prime$}(D);
      
      \draw [2cell] (1.1,1.4) to [bend left=20] node [auto,labelsize] {$f^\prime$} (2.9,1.4);
      \draw [2cell] (2,-0.35) to [bend left=0] node [auto,swap,labelsize] {$\Theta$} (2,0.35);
      
      \end{scope}
    \end{tikzpicture}
    \]
    where $\Theta$ is a 2-natural transformation and
    $\alpha$ a modification;
    equivalently, it is a pair $(\Theta,\alpha)$
    where $\alpha$ is a 2-cell of $\tcat{B}$
    of the following type:
    
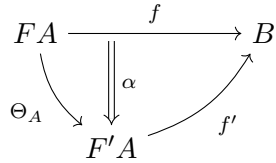
\[
    \begin{tikzpicture}
      \node (a)  at (3.5,2) {$FA$};
      \node (ap) at (6.5,2) {$B$};
      \node (b)  at (4.5,0.5) {$F^\prime A$};
      
      \draw [->] (a) to [bend right=0]node [auto,labelsize] {$f$}(ap);
      \draw [->] (a)  to [bend right=20]node [auto,swap,labelsize] {$\Theta_A$}(b);
      \draw [->] (b) to [bend right=20]node [auto,swap,labelsize] {$f^\prime$}(ap);
      \draw [2cell] (4.5,1.9) to node [auto,labelsize]{$\alpha$}(b);
    \end{tikzpicture}
    \qedhere
    \]
  \end{itemize}
\end{defn}

The compositions in $\Epp$ are defined abstractly by obvious
pasting diagrams in $\TCat$.
We also provide a concrete description of the compositions 
in Appendix~\ref{apx-composition}.

The first projection of the data defines a
2-functor~$\ppp\colon \Epp\to \TCat_2$,
where $\TCat_2$ is the 2-category of 2-categories,
2-functors, and 2-natural transformations.
We take a fibrational viewpoint~\cite{bakovic:2-Grothendieck,hermida:descent,hermida:fib} 
and say a notion~$X$ is \emph{above} $I$ if $\ppp (X)=I$.
Observe that the fiber over the object~$\tcat{K}$ of $\TCat_2$
is $\tcat{K}$ itself; 
this is why one can claim that $\Epp$ ``contains''
an arbitrary 2-category.

\subsection{The 2-category $\Epm$}
\begin{defn}
  We define the 2-category~$\Epm$ by the following data.
  \begin{itemize}
  \item A 0-cell of $\Epm$ is a 2-functor
    $A\colon  1\to \tcat{A}$ where $\tcat{A}$ is a 2-category;
    equivalently, it is a pair $(\tcat{A},A)$ where 
    $A$ is a 0-cell of $\tcat{A}$.
  \item A 1-cell of $\Epm$ from $(\tcat{A},A)$ to 
    $(\tcat{B},B)$ is a diagram in $\TCat$
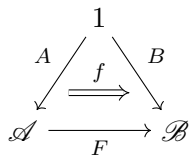
    \[
    \begin{tikzpicture}
      \node (C)  at (1,1.5) {$1$};
      \node (C2)  at (0,0) {$\tcat{A}$};
      \node (D)  at (2,0) {$\tcat{B}$};
      
      \draw [->] (C)  to [bend left=0]  node [auto,swap,labelsize] {$A$} (C2);
      \draw [->] (C) to [bend left=0]node [auto,labelsize] {$B$}(D);
      \draw [->] (C2) to [bend left=0]node [auto,swap,labelsize] {$F$}(D);
      
      \draw [2cell] (0.6,0.5) to node [auto,labelsize] {$f$} (1.4,0.5);
    \end{tikzpicture}
    \]
    where $F$ is a 2-functor and 
    $f$ a 2-natural transformation;
    equivalently, it is a pair $(F,f)$ where
    $f\colon  FA\to B$ is a 1-cell of
    $\tcat{B}$.
  \item A 2-cell of $\Epm$ from $(F,f)$ to 
    $(F^\prime,f^\prime)$ is a diagram in $\TCat$
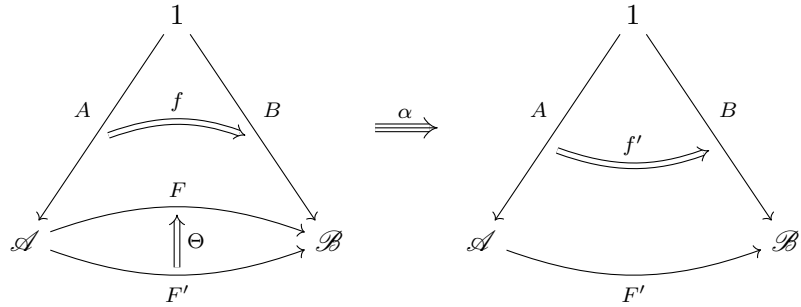
    \[
    \begin{tikzpicture}
      \node (C)  at (2,3) {$1$};
      \node (C2)  at (0,0) {$\tcat{A}$};
      \node (D)  at (4,0) {$\tcat{B}$};
      
      \draw [->] (C)  to [bend left=0]  node [auto,swap,labelsize] {$A$} (C2);
      \draw [->] (C) to [bend left=0]node [auto,labelsize] {$B$}(D);
      \draw [->] (C2) to [bend right=20]node [auto,swap,labelsize] {$F^\prime$}(D);
      
      \draw [2cell] (1,1.2) to [bend right=20] node [auto,labelsize] {$f^\prime$} (3,1.2);
      
      \begin{scope}[shift={(-6,0)}]
      \draw [3cell] (4.5,1.5) to [bend left=0] node [auto,labelsize] {$\alpha$} (5.5,1.5);
            \draw (4.5,1.5)+(3pt,0pt) to  (5.4,1.5);
      \end{scope}
      
      \begin{scope}[shift={(-6,0)}]
	  \node (C)  at (2,3) {$1$};
      \node (C2)  at (0,0) {$\tcat{A}$};
      \node (D)  at (4,0) {$\tcat{B}$};
      
      \draw [->] (C)  to [bend left=0]  node [auto,swap,labelsize] {$A$} (C2);
      \draw [->] (C) to [bend left=0]node [auto,labelsize] {$B$}(D);
      \draw [->] (C2) to [bend right=20]node [auto,swap,labelsize] {$F^\prime$}(D);
      \draw [->] (C2) to [bend left=20]node [auto,labelsize] {$F$}(D);
      
      \draw [2cell] (1.1,1.4) to [bend left=20] node [auto,labelsize] {$f$} (2.9,1.4);
      \draw [2cell] (2,-0.35) to [bend left=0] node [auto,swap,labelsize] {$\Theta$} (2,0.35);
      
      \end{scope}
    \end{tikzpicture}
    \]
    where $\Theta$ is a 2-natural transformation and
    $\alpha$ a modification;
    equivalently, it is a pair $(\Theta,\alpha)$
    where $\alpha$ is a 2-cell of $\tcat{B}$
    of the following type:
    
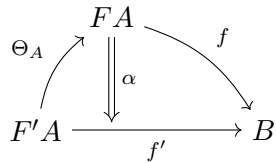
\[
    \begin{tikzpicture}
      \node (a)  at (3.5,2) {$F^\prime A$};
      \node (ap) at (6.5,2) {$B$};
      \node (b)  at (4.5,3.5) {$F A$};
      
      \draw [->] (a) to [bend right=0]node [auto,swap,labelsize] {$f^\prime$}(ap);
      \draw [->] (a)  to [bend left=20]node [auto,labelsize] {$\Theta_A$}(b);
      \draw [->] (b) to [bend left=20]node [auto,labelsize] {$f$}(ap);
      \draw [2cellr] (4.5,2.1) to node [auto,swap,labelsize]{$\alpha$}(b);
    \end{tikzpicture}
    \qedhere
    \]
  \end{itemize}
\end{defn}

For this 2-category, we may define by the first projection
the projection 2-functor~$\ppm\colon$ $\Epm\to\TCat_2^\op{2}$.

\subsection{The 2-category $\Emp$}
\begin{defn}
  We define the 2-category~$\Emp$ by the following data.
  \begin{itemize}
  \item A 0-cell of $\Emp$ is a 2-functor
    $A\colon  1\to \tcat{A}$ where $\tcat{A}$ is a 2-category;
    equivalently, it is a pair $(\tcat{A},A)$ where 
    $A$ is a 0-cell of $\tcat{A}$.
  \item A 1-cell of $\Emp$ from $(\tcat{A},A)$ to 
    $(\tcat{B},B)$ is a diagram in $\TCat$
    \[
    \begin{tikzpicture}
      \node (C)  at (1,1.5) {$1$};
      \node (C2)  at (0,0) {$\tcat{A}$};
      \node (D)  at (2,0) {$\tcat{B}$};
      
      \draw [->] (C)  to [bend left=0]  node [auto,swap,labelsize] {$A$} (C2);
      \draw [->] (C) to [bend left=0]node [auto,labelsize] {$B$}(D);
      \draw [<-] (C2) to [bend left=0]node [auto,swap,labelsize] {$F$}(D);
      
      \draw [2cell] (0.6,0.5) to node [auto,labelsize] {$f$} (1.4,0.5);
    \end{tikzpicture}
    \]
    where $F$ is a 2-functor and 
    $f$ a 2-natural transformation;
    equivalently, it is a pair $(F,f)$ where
    $f\colon  A\to FB$ is a 1-cell of
    $\tcat{A}$.
  \item A 2-cell of $\Emp$ from $(F,f)$ to 
    $(F^\prime,f^\prime)$ is a diagram in $\TCat$
	\[
    \begin{tikzpicture}
      \node (C)  at (2,3) {$1$};
      \node (C2)  at (0,0) {$\tcat{A}$};
      \node (D)  at (4,0) {$\tcat{B}$};
      
      \draw [->] (C)  to [bend left=0]  node [auto,swap,labelsize] {$A$} (C2);
      \draw [->] (C) to [bend left=0]node [auto,labelsize] {$B$}(D);
      \draw [<-] (C2) to [bend right=20]node [auto,swap,labelsize] {$F^\prime$}(D);
      
      \draw [2cell] (1,1.2) to [bend right=20] node [auto,labelsize] {$f^\prime$} (3,1.2);
      
      \begin{scope}[shift={(-6,0)}]
      \draw [3cell] (4.5,1.5) to [bend left=0] node [auto,labelsize] {$\alpha$} (5.5,1.5);
            \draw (4.5,1.5)+(3pt,0pt) to  (5.4,1.5);
      \end{scope}
      
      \begin{scope}[shift={(-6,0)}]
	  \node (C)  at (2,3) {$1$};
      \node (C2)  at (0,0) {$\tcat{A}$};
      \node (D)  at (4,0) {$\tcat{B}$};
      
      \draw [->] (C)  to [bend left=0]  node [auto,swap,labelsize] {$A$} (C2);
      \draw [->] (C) to [bend left=0]node [auto,labelsize] {$B$}(D);
      \draw [<-] (C2) to [bend right=20]node [auto,swap,labelsize] {$F^\prime$}(D);
      \draw [<-] (C2) to [bend left=20]node [auto,labelsize] {$F$}(D);
      
      \draw [2cell] (1.1,1.4) to [bend left=20] node [auto,labelsize] {$f$} (2.9,1.4);
      \draw [2cellr] (2,-0.35) to [bend left=0] node [auto,swap,labelsize] {$\Theta$} (2,0.35);
      
      \end{scope}
    \end{tikzpicture}
    \]
    where $\Theta$ is a 2-natural transformation and
    $\alpha$ a modification;
    equivalently, it is a pair $(\Theta,\alpha)$
    where $\alpha$ is a 2-cell of $\tcat{A}$
    of the following type:
    \[
    \begin{tikzpicture}
      \node (a)  at (3.5,2) {$A$};
      \node (ap) at (6.5,2) {$F^\prime B$};
      \node (b)  at (5.5,3.5) {$F B$};
      
      \draw [->] (a) to [bend right=0]node [auto,swap,labelsize] {$f^\prime$}(ap);
      \draw [->] (a)  to [bend left=20]node [auto,labelsize] {$f$}(b);
      \draw [->] (b) to [bend left=20]node [auto,labelsize] {$\Theta_{B}$}(ap);
      \draw [2cellr] (5.5,2.1) to node [auto,swap,labelsize]{$\alpha$}(b);
    \end{tikzpicture}
    \qedhere
    \]
  \end{itemize}
\end{defn}
The first projection defines a 2-functor~$\pmp\colon\Emp\to\TCat_2^\op{1}$.

\subsection{The 2-category $\Emm$}
\begin{defn}
  We define the 2-category~$\Emm$ by the following data.
  \begin{itemize}
  \item A 0-cell of $\Emm$ is a 2-functor
    $A\colon  1\to \tcat{A}$ where $\tcat{A}$ is a 2-category;
    equivalently, it is a pair $(\tcat{A},A)$ where 
    $A$ is a 0-cell of $\tcat{A}$.
  \item A 1-cell of $\Emm$ from $(\tcat{A},A)$ to 
    $(\tcat{B},B)$ is a diagram in $\TCat$
    \[
    \begin{tikzpicture}
      \node (C)  at (1,1.5) {$1$};
      \node (C2)  at (0,0) {$\tcat{A}$};
      \node (D)  at (2,0) {$\tcat{B}$};
      
      \draw [->] (C)  to [bend left=0]  node [auto,swap,labelsize] {$A$} (C2);
      \draw [->] (C) to [bend left=0]node [auto,labelsize] {$B$}(D);
      \draw [<-] (C2) to [bend left=0]node [auto,swap,labelsize] {$F$}(D);
      
      \draw [2cell] (0.6,0.5) to node [auto,labelsize] {$f$} (1.4,0.5);
    \end{tikzpicture}
    \]
    where $F$ is a 2-functor and 
    $f$ a 2-natural transformation;
    equivalently, it is a pair $(F,f)$ where
    $f\colon  A\to FB$ is a 1-cell of
    $\tcat{A}$.
  \item A 2-cell of $\Emm$ from $(F,f)$ to 
    $(F^\prime,f^\prime)$ is a diagram in $\TCat$
    \[
    \begin{tikzpicture}
      \node (C)  at (2,3) {$1$};
      \node (C2)  at (0,0) {$\tcat{A}$};
      \node (D)  at (4,0) {$\tcat{B}$};
      
      \draw [->] (C)  to [bend left=0]  node [auto,swap,labelsize] {$A$} (C2);
      \draw [->] (C) to [bend left=0]node [auto,labelsize] {$B$}(D);
      \draw [<-] (C2) to [bend right=20]node [auto,swap,labelsize] {$F$}(D);
      
      \draw [2cell] (1,1.2) to [bend right=20] node [auto,labelsize] {$f$} (3,1.2);
      
      
      \draw [3cell] (4.5,1.5) to [bend left=0] node [auto,labelsize] {$\alpha$} (5.5,1.5);
            \draw (4.5,1.5)+(3pt,0pt) to  (5.4,1.5);
      
      \begin{scope}[shift={(6,0)}]
	  \node (C)  at (2,3) {$1$};
      \node (C2)  at (0,0) {$\tcat{A}$};
      \node (D)  at (4,0) {$\tcat{B}$};
      
      \draw [->] (C)  to [bend left=0]  node [auto,swap,labelsize] {$A$} (C2);
      \draw [->] (C) to [bend left=0]node [auto,labelsize] {$B$}(D);
      \draw [<-] (C2) to [bend right=20]node [auto,swap,labelsize] {$F$}(D);
      \draw [<-] (C2) to [bend left=20]node [auto,labelsize] {$F^\prime$}(D);
      
      \draw [2cell] (1.1,1.4) to [bend left=20] node [auto,labelsize] {$f^\prime$} (2.9,1.4);
      \draw [2cellr] (2,-0.35) to [bend left=0] node [auto,swap,labelsize] {$\Theta$} (2,0.35);
      
      \end{scope}
    \end{tikzpicture}
    \]
    where $\Theta$ is a 2-natural transformation and
    $\alpha$ a modification;
    equivalently, it is a pair $(\Theta,\alpha)$
    where $\alpha$ is a 2-cell of $\tcat{A}$
    of the following type:
    \[
    \begin{tikzpicture}
      \node (a)  at (3.5,2) {$A$};
      \node (ap) at (6.5,2) {$F B$};
      \node (b)  at (5.5,0.5) {$F^\prime B$};
      
      \draw [->] (a) to [bend right=0]node [auto,labelsize] {$f$}(ap);
      \draw [->] (a)  to [bend right=20]node [auto,swap,labelsize] {$f^\prime$}(b);
      \draw [->] (b) to [bend right=20]node [auto,swap,labelsize] {$\Theta_{B}$}(ap);
      \draw [2cell] (5.5,1.9) to node [auto,labelsize]{$\alpha$}(b);
    \end{tikzpicture}
    \qedhere
    \]
  \end{itemize}
\end{defn}
The first projection defines a 2-functor~$\pmm\colon\Emm\to\TCat_2^\op{1,2}$.

\section{Reducing monads with parameters to mere monads}

Let us proceed to see how the 
notions of graded and indexed monad can be seen as mere \emph{monads}
in the 2-categorical sense of Definition~\ref{defn-monad}
inside the appropriate 2-categories among 
$\Epp$, $\Epm$, $\Emp$ and $\Emm$;
this also helps digress and motivate the definitions 
of these 2-categories, which might look 
a little bit complicated at first sight.

\subsection{Graded monads as monads in $\Epp$}
\label{subsec-gm-mnd-epp}
The first reduction for graded monads takes place in the 2-category~$\Epp$,
and will later be employed to perform the
\emph{Eilenberg--Moore construction} for graded monads.

Let us fix a strict monoidal category~$\cat{M}=(\cat{M},\tensor,\e)$.
First we need a preliminary definition.
\begin{defn}
\label{defn-M-times}
There is a 2-monad~$\cat{M}\times(-)=
(\cat{M}\times(-),H^\cat{M},M^\cat{M})$ on $\Cat$ defined as follows.
\begin{itemize}
\item The 2-functor~$\cat{M}\times(-)\colon\Cat\to\Cat$
is the product category construction.
\item The 2-natural transformation~$H^\cat{M}\colon \id{\Cat}\Rightarrow
\cat{M}\times(-)$ has as its $\cat{X}$-component the functor~$
H^\cat{M}_\cat{X}=H_\cat{X}\colon\cat{X}\to\cat{M}\times\cat{X}$ given by
$x\mapsto(\e,x)$.
\item The 2-natural transformation~$M^\cat{M}\colon \cat{M}\times\cat{M}
\times(-)\Rightarrow\cat{M}\times(-)$ has as its $\cat{X}$-component
the functor~$M^\cat{M}_\cat{X}=M_\cat{X}\colon\cat{M}\times\cat{M}\times\cat{X}
\to\cat{M}\times\cat{X}$ given by $(m,n,x)\mapsto
(m\tensor n,x)$.\qedhere
\end{itemize}
\end{defn}
See Appendix~\ref{apx-sec-enriched} for an abstract way to understand 
this 2-monad, as well as other 2-(co)monads introduced subsequently.

Now we have the following result:
\begin{prop}
Let $\cat{C}$ be a category.
Then an $\cat{M}$-graded monad on 
$\cat{C}$ is the same thing as a monad in $\Epp$
on $(\Cat,\cat{C})$ above the 2-monad~$\cat{M}\times(-)\colon
\Cat\to\Cat$.
\end{prop}
\begin{pf}
Let us spell out what \emph{a monad in $\Epp$ on $(\Cat,\cat{C})$ above
$\cat{M}\times(-)$}, actually is.
Such a thing consists of the following data:
\begin{itemize}
\item A 1-cell~$(\cat{M}\times(-),\gm{T})\colon
(\Cat,\cat{C})\to(\Cat,\cat{C})$ of $\Epp$, where $\gm{T}$ is a functor of type
\[
\gm{T}\quad\colon\quad\cat{M}\times\cat{C}
\quad\longrightarrow\quad\cat{C}.
\]
\item A 2-cell~$(H^\cat{M},\eta)\colon (\id{\Cat},\id{\cat{C}})\Rightarrow
(\cat{M}\times(-),\gm{T})$ of $\Epp$, where $\eta$ is a natural transformation of type
\[
    \begin{tikzpicture}
      \node (a)  at (3.5,2) {$\cat{C}$};
      \node (ap) at (6.5,2) {$\cat{C}$};
      \node (b)  at (4.5,0.5) {$\cat{M}\times\cat{C}$};
      
      \draw [->] (a) to [bend right=0]node [auto,labelsize] {$\id{\cat{C}}$}(ap);
      \draw [->] (a)  to [bend right=20]node [auto,swap,labelsize] {$H_\cat{C}$}(b);
      \draw [->] (b) to [bend right=20]node [auto,swap,labelsize] {$\gm{T}$}(ap);
      \draw [2cell] (4.5,1.9) to node [auto,labelsize]{$\eta$}(b);
    \end{tikzpicture}
\]
\item A 2-cell~$(M^\cat{M},\mu)\colon (\cat{M}\times(-),\gm{T})\circ
(\cat{M}\times(-),\gm{T})\Rightarrow(\cat{M}\times(-),\gm{T})$ of 
$\Epp$, where $\mu$ is a natural transformation of type
\[
\begin{tikzpicture}[baseline=-\the\dimexpr\fontdimen22\textfont2\relax ]
        \begin{scope} 
      \node (L)  at (0,0)  {$\cat{M}\times\cat{M}\times\cat{C}$};
      \node (M)  at (3,0)  {$\cat{M}\times\cat{C}$};
      \node (R)  at (6,0) {$\cat{C}$};
      \node (B)  at (2,-2) {$\cat{M}\times\cat{C}$};
      \draw [->] (L) to node  [auto,labelsize]      {$\cat{M}\times\gm{T}$} (M);
      \draw [->] (M) to node  [auto,labelsize]      {$\gm{T}$} (R);
      \draw [->] (L) to [bend right=20]node  [auto,swap,labelsize] {$M_\cat{C}$}(B);
      \draw [->] (B) to [bend right=20]node  [auto,swap,labelsize] {$\gm{T}$}(R);
      \draw [2cellnode] (2,0) to node[auto,labelsize] {$\mu$}(B);
        \end{scope}
\end{tikzpicture}
\]
\end{itemize}
These data satisfy the following axioms:
{\allowdisplaybreaks
\begin{align}
\label{eqn-epp-gm-1}
\begin{tikzpicture}[baseline=-\the\dimexpr\fontdimen22\textfont2\relax ]
        \begin{scope} 
      \node (A)  at (0,0) {$\cat{M}\times \cat{M}\times\cat{C}$};
      \node (B) at (3,0) {$\cat{M}\times\cat{C}$};
      \node (C) at (5,0) {$\cat{C}$};
      \node (D)  at (1.5,-1.5) {$\cat{M}\times\cat{C}$};
      \node (E) at (3.5,1.5) {$\cat{C}$};
      \node (F) at (0.5,1.5) {$\cat{M}\times\cat{C}$};
      \draw [->] (A) to [bend left=0]node [auto,labelsize] {$\cat{M}\times \gm{T}$}(B);
      \draw [->] (B) to [bend left=0]node [auto,labelsize] {$\gm{T}$}(C);
      \draw [->] (A) to [bend right=20]node [auto,swap,labelsize] {$M_\cat{C}$}(D);
      \draw [->] (D) to [bend right=20]node [auto,swap,labelsize] {$\gm{T}$}(C);
      \draw [->] (E) to [bend right=20]node [auto,swap,labelsize] {$H_\cat{C}$}(B);
      \draw [->] (F) to [bend right=0]node [auto,labelsize] {$\gm{T}$}(E);
      \draw [->] (E) to [bend left=20]node [auto,labelsize] {$\id{\cat{C}}$}(C);
      \draw [->] (F) to [bend right=20]node [auto,swap,labelsize] {$H_{\cat{M}\times\cat{C}}$}(A);
      \draw [2cell] (1.5,-0.1) to node [auto,labelsize]{$\mu$}(1.5,-1.1);
      \draw [2cell] (3.5,1.2) to node [auto,labelsize]{$\eta$}(3.5,0.3);
        \end{scope}
\end{tikzpicture}
\quad&=\quad
\begin{tikzpicture}[baseline=-\the\dimexpr\fontdimen22\textfont2\relax ]
        \begin{scope} 
       \node (A)  at (0,0) {$\cat{M}\times\cat{C}$};
       \node (C) at (2,0) {$\cat{C}$};
       \draw [->] (A) to [bend left=0]node [auto,labelsize] {$\gm{T}$}(C);
        \end{scope}
\end{tikzpicture}\\
\label{eqn-epp-gm-2}
\begin{tikzpicture}[baseline=-\the\dimexpr\fontdimen22\textfont2\relax ]
        \begin{scope} 
      \node (A)  at (0,0) {$\cat{M}\times \cat{M}\times\cat{C}$};
      \node (B) at (3,0) {$\cat{M}\times\cat{C}$};
      \node (C) at (5,0) {$\cat{C}$};
      \node (D)  at (1.5,-1.5) {$\cat{M}\times\cat{C}$};
      \node (F) at (1,1.5) {$\cat{M}\times\cat{C}$};
      \draw [->] (A) to [bend left=0]node [auto,labelsize] {$\cat{M}\times \gm{T}$}(B);
      \draw [->] (B) to [bend left=0]node [auto,labelsize] {$\gm{T}$}(C);
      \draw [->] (A) to [bend right=20]node [auto,swap,labelsize] {$M_\cat{C}$}(D);
      \draw [->] (D) to [bend right=20]node [auto,swap,labelsize] {$\gm{T}$}(C);
      \draw [->] (F) to [bend left=20]node [auto,labelsize] {$\cat{M}\times\id{\cat{C}}$}(B);
      \draw [->] (F) to [bend right=20]node [auto,swap,labelsize] {$\cat{M}\times H_{\cat{C}}$}(A);
      \draw [2cell] (1.5,-0.1) to node [auto,labelsize]{$\mu$}(1.5,-1.1);
      \draw [2cell] (1,1.2) to node [auto,labelsize]{$\cat{M}\times\eta$}(1,0.3);
        \end{scope}
\end{tikzpicture}
\quad&=\quad
\begin{tikzpicture}[baseline=-\the\dimexpr\fontdimen22\textfont2\relax ]
        \begin{scope} 
       \node (A)  at (0,0) {$\cat{M}\times\cat{C}$};
       \node (C) at (2,0) {$\cat{C}$};
       \draw [->] (A) to [bend left=0]node [auto,labelsize] {$\gm{T}$}(C);
        \end{scope}
\end{tikzpicture}
\end{align}
}%
\vspace{-20pt}
\begin{multline}
\label{eqn-epp-gm-3}
\begin{tikzpicture}[baseline=-\the\dimexpr\fontdimen22\textfont2\relax ]
        \begin{scope} 
      \node (A)  at (0,0) {$\cat{M}\times \cat{M}\times\cat{C}$};
      \node (B) at (4,0) {$\cat{M}\times\cat{C}$};
      \node (C) at (7,0) {$\cat{C}$};
      \node (D)  at (2,-1.5) {$\cat{M}\times\cat{C}$};
      \node (E) at (4,1.5) {$\cat{M}\times\cat{M}\times\cat{C}$};
      \node (F) at (0,1.5) {$\cat{M}\times\cat{M}\times \cat{M}\times\cat{C}$};
      \node (G) at (7,1.5) {$\cat{M}\times\cat{C}$};
      \draw [->] (A) to [bend left=0]node [auto,labelsize] {$\cat{M}\times \gm{T}$}(B);
      \draw [->] (B) to [bend left=0]node [auto,swap,labelsize] {$\gm{T}$}(C);
      \draw [->] (A) to [bend right=20]node [auto,swap,labelsize] {$M_\cat{C}$}(D);
      \draw [->] (D) to [bend right=15]node [auto,swap,labelsize] {$\gm{T}$}(C);
      \draw [->] (E) to [bend right=0]node [auto,swap,labelsize] {$M_\cat{C}$}(B);
      \draw [->] (G) to [bend right=0]node [auto,labelsize] {$\gm{T}$}(C);
      \draw [->] (E) to [bend right=0]node [auto,labelsize] {$\cat{M}\times \gm{T}$}(G);
      \draw [->] (F) to [bend right=0]node [auto,labelsize] {$\cat{M}\times\cat{M}\times \gm{T}$}(E);
      \draw [->] (F) to [bend right=0]node [auto,swap,labelsize] {$M_{\cat{M}\times\cat{C}}$}(A);
      \draw [2cell] (2,-0.1) to node [auto,labelsize]{$\mu$}(2,-1.1);
      \draw [2cell] (5.5,1.4) to node [auto,labelsize]{$\mu$}(5.5,0.1);
        \end{scope}
\end{tikzpicture}\\
=\quad
\begin{tikzpicture}[baseline=-\the\dimexpr\fontdimen22\textfont2\relax ]
        \begin{scope} 
      \node (A)  at (0,0) {$\cat{M}\times \cat{M}\times\cat{C}$};
      \node (B) at (4,0) {$\cat{M}\times\cat{C}$};
      \node (C) at (7,0) {$\cat{C}$};
      \node (D)  at (2,-1.5) {$\cat{M}\times\cat{C}$};
      \node (E) at (4,1.5) {$\cat{M}\times\cat{M}\times\cat{C}$};
      \node (F) at (0,1.5) {$\cat{M}\times\cat{M}\times\cat{M}\times\cat{C}$};
      \draw [->] (A) to [bend left=0]node [auto,labelsize] {$\cat{M}\times \gm{T}$}(B);
      \draw [->] (B) to [bend left=0]node [auto,labelsize] {$\gm{T}$}(C);
      \draw [->] (A) to [bend right=20]node [auto,swap,labelsize] {$M_\cat{C}$}(D);
      \draw [->] (D) to [bend right=15]node [auto,swap,labelsize] {$\gm{T}$}(C);
      \draw [->] (E) to [bend right=0]node [auto,labelsize] {$\cat{M}\times\gm{T}$}(B);
      \draw [->] (F) to [bend right=0]node [auto,labelsize] {$\cat{M}\times\cat{M}\times \gm{T}$}(E);
      \draw [->] (F) to [bend right=0]node [auto,swap,labelsize] {$\cat{M}\times M_\cat{C}$}(A);
      \draw [2cell] (2,-0.1) to node [auto,labelsize]{$\mu$}(2,-1.1);
      \draw [2cell] (2,1.4) to node [auto,labelsize]{$\cat{M}\times\mu$}(2,0.5);
        \end{scope}
\end{tikzpicture}
\end{multline}

To see the equivalence of this notion and that of $\cat{M}$-graded monad
on $\cat{C}$, first note that by adjointness
the functor~$\gm{T}\colon\cat{M}\times\cat{C}\to\cat{C}$ is equivalent to 
its transpose, also denoted by
$\gm{T}\colon\cat{M}\to[\cat{C},\cat{C}]$.
%
Now it is routine to see that in the definition of graded monads, 
giving the data $\gm{T}_m$ and $\gm{T}_u$ satisfying (GM1) and (GM2)
is equivalent to giving $\gm{T}$, and giving $\mu_{m,n}$ satisfying
(GM3) is equivalent to giving $\mu$.
Finally observe that (GM4), (GM5) and (GM6) are equivalent to
(\ref{eqn-epp-gm-1}), (\ref{eqn-epp-gm-2}) and (\ref{eqn-epp-gm-3})
respectively. 
\end{pf}

\subsection{Graded monads as monads in $\Emm$}
\label{subsec-gm-mnd-emm}
In this section we observe that there is another way to 
understand graded monads as monads in the 2-categorical sense,
by working inside the 2-category~$\Emm$;
this will be used for the \emph{Kleisli construction} for graded monads.

Again we fix a strict monoidal category~$\cat{M}=(\cat{M},\tensor,\e)$.
\begin{defn}
\label{defn-M-to}
There is a 2-comonad~$[\cat{M},-]=([\cat{M},-],H_\cat{M},M_\cat{M})$ 
on $\Cat$ defined as follows.
\begin{itemize}
\item The 2-functor~$[\cat{M},-]\colon\Cat\to\Cat$
is the functor category construction.
\item The 2-natural transformation~$H_\cat{M}\colon [\cat{M},-]\Rightarrow
\id{\Cat}$ has as its $\cat{X}$-component the functor~$
H_{\cat{M},\cat{X}}=H_\cat{X}\colon
[\cat{M},\cat{X}]\to\cat{X}$ given by
$X\mapsto X(\e)$.
\item The 2-natural transformation~$M_\cat{M}\colon [\cat{M},-]
\Rightarrow[\cat{M},[\cat{M},-]]$ has as its $\cat{X}$-component
the functor~$M_{\cat{M},\cat{X}}=M_\cat{X}\colon
[\cat{M},\cat{X}]\to[\cat{M},[\cat{M},\cat{X}]]$ given by 
$X\mapsto (m\mapsto X (-\tensor m) )$.\qedhere
\end{itemize}
\end{defn}

\begin{prop}
Let $\cat{C}$ be a category.
Then an $\cat{M}$-graded monad on $\cat{C}$ is the same thing as 
a monad in $\Emm$ on $(\Cat,\cat{C})$ above the 2-comonad~$
[\cat{M},-]\colon\Cat\to\Cat$.
\end{prop}
\begin{pf}
Indeed, the latter notion is given by the following data:
\begin{itemize}
\item A 1-cell~$([\cat{M},-],\gm{T})\colon(\Cat,\cat{C})\to
(\Cat,\cat{C})$ of $\Emm$, where $\gm{T}$ is a functor of type
\[
\gm{T}\quad\colon\quad
\cat{C}\quad\longrightarrow\quad[\cat{M},\cat{C}].
\]
\item A 2-cell~$(H_\cat{M},\eta)\colon(\id{\Cat},\id{\cat{C}})
\Rightarrow([\cat{M},-],\gm{T})$ of $\Emm$, where 
$\eta$ is a natural transformation of type
\[
    \begin{tikzpicture}
      \node (a)  at (3.5,2) {$\cat{C}$};
      \node (ap) at (6.5,2) {$\cat{C}$};
      \node (b)  at (5.5,0.5) {$[\cat{M},\cat{C}]$};
      
      \draw [->] (a) to [bend right=0]node [auto,labelsize] {$\id{\cat{C}}$}(ap);
      \draw [->] (a)  to [bend right=20]node [auto,swap,labelsize] {$\gm{T}$}(b);
      \draw [->] (b) to [bend right=20]node [auto,swap,labelsize] {$H_{\cat{C}}$}(ap);
      \draw [2cell] (5.5,1.9) to node [auto,labelsize]{$\eta$}(b);
    \end{tikzpicture}
\]
\item A 2-cell~$(M_\cat{M},\mu)\colon([\cat{M},-],\gm{T})\circ([\cat{M},-],\gm{T})
\Rightarrow([\cat{M},-],\gm{T})$ of $\Emm$, where $\mu$ is a natural 
transformation of type
\[
\begin{tikzpicture}[baseline=-\the\dimexpr\fontdimen22\textfont2\relax ]
        \begin{scope} [shift={(0,0.75)}]
      \node (A)  at (0,0) {$[\cat{M},[\cat{M},\cat{C}]]$};
      \node (B) at (-3,0) {$[\cat{M},\cat{C}]$};
      \node (C) at (-6,0) {$\cat{C}$};
      \node (D)  at (-2,-2) {$[\cat{M},\cat{C}]$};
      \draw [->] (B) to [bend left=0]node [auto,labelsize] {$[\cat{M},\gm{T}]$}(A);
      \draw [->] (C) to [bend left=0]node [auto,labelsize] {$\gm{T}$}(B);
      \draw [->] (D) to [bend right=20]node [auto,swap,labelsize] {$M_\cat{C}$}(A);
      \draw [->] (C) to [bend right=20]node [auto,swap,labelsize] {$\gm{T}$}(D);
      \draw [2cell] (-2,-0.1) to node [auto,labelsize]{$\mu$}(-2,-1.6);
        \end{scope}
\end{tikzpicture}
\]
\end{itemize}
These data satisfy the following axioms:
{\allowdisplaybreaks
\begin{align*}
\begin{tikzpicture}[baseline=-\the\dimexpr\fontdimen22\textfont2\relax ]
        \begin{scope} 
      \node (A)  at (0,0) {$[\cat{M},[\cat{M},\cat{C}]]$};
      \node (B) at (-3,0) {$[\cat{M},\cat{C}]$};
      \node (C) at (-5,0) {$\cat{C}$};
      \node (D)  at (-1.5,-1.5) {$[\cat{M},\cat{C}]$};
      \node (E) at (-3.5,1.5) {$\cat{C}$};
      \node (F) at (-0.5,1.5) {$[\cat{M},\cat{C}]$};
      \draw [->] (B) to [bend left=0]node [auto,labelsize] {$[\cat{M},\gm{T}]$}(A);
      \draw [->] (C) to [bend left=0]node [auto,labelsize] {$\gm{T}$}(B);
      \draw [->] (D) to [bend right=20]node [auto,swap,labelsize] {$M_\cat{C}$}(A);
      \draw [->] (C) to [bend right=20]node [auto,swap,labelsize] {$\gm{T}$}(D);
      \draw [->] (B) to [bend right=20]node [auto,swap,labelsize] {$H_\cat{C}$}(E);
      \draw [->] (E) to [bend right=0]node [auto,labelsize] {$\gm{T}$}(F);
      \draw [->] (C) to [bend left=20]node [auto,labelsize] {$\id{\cat{C}}$}(E);
      \draw [->] (A) to [bend right=20]node [auto,swap,labelsize] {$H_{[\cat{M},\cat{C}]}$}(F);
      \draw [2cell] (-1.5,-0.1) to node [auto,labelsize]{$\mu$}(-1.5,-1.1);
      \draw [2cell] (-3.5,1.2) to node [auto,labelsize]{$\eta$}(-3.5,0.3);
        \end{scope}
\end{tikzpicture}
\quad&=\quad
\begin{tikzpicture}[baseline=-\the\dimexpr\fontdimen22\textfont2\relax ]
        \begin{scope} 
       \node (A)  at (0,0) {$[\cat{M},\cat{C}]$};
       \node (C) at (-2,0) {$\cat{C}$};
       \draw [->] (C) to [bend left=0]node [auto,labelsize] {$\gm{T}$}(A);
        \end{scope}
\end{tikzpicture}\\
\begin{tikzpicture}[baseline=-\the\dimexpr\fontdimen22\textfont2\relax ]
        \begin{scope} 
      \node (A)  at (0,0) {$[\cat{M},[\cat{M},\cat{C}]]$};
      \node (B) at (-3,0) {$[\cat{M},\cat{C}]$};
      \node (C) at (-5,0) {$\cat{C}$};
      \node (D)  at (-1.5,-1.5) {$[\cat{M},\cat{C}]$};
      \node (F) at (-1,1.5) {$[\cat{M},\cat{C}]$};
      \draw [->] (B) to [bend left=0]node [auto,labelsize] {$[\cat{M},\gm{T}]$}(A);
      \draw [->] (C) to [bend left=0]node [auto,labelsize] {$\gm{T}$}(B);
      \draw [->] (D) to [bend right=20]node [auto,swap,labelsize] {$M_\cat{C}$}(A);
      \draw [->] (C) to [bend right=20]node [auto,swap,labelsize] {$\gm{T}$}(D);
      \draw [->] (B) to [bend left=20]node [auto,labelsize] {$[\cat{M},\id{\cat{C}}]$}(F);
      \draw [->] (A) to [bend right=20]node [auto,swap,labelsize] {$[\cat{M},H_{\cat{C}}]$}(F);
      \draw [2cell] (-1.5,-0.1) to node [auto,labelsize]{$\mu$}(-1.5,-1.1);
      \draw [2cell] (-1,1.2) to node [auto,swap,labelsize]{$[\cat{M},\eta]$}(-1,0.3);
        \end{scope}
\end{tikzpicture}
\quad&=\quad
\begin{tikzpicture}[baseline=-\the\dimexpr\fontdimen22\textfont2\relax ]
        \begin{scope} 
       \node (A)  at (0,0) {$[\cat{M},\cat{C}]$};
       \node (C) at (-2,0) {$\cat{C}$};
       \draw [->] (C) to [bend left=0]node [auto,labelsize] {$\gm{T}$}(A);
        \end{scope}
\end{tikzpicture}
\end{align*}
}
\vspace{-20pt}
\begin{multline*}
\begin{tikzpicture}[baseline=-\the\dimexpr\fontdimen22\textfont2\relax ]
        \begin{scope} 
      \node (A)  at (0,0) {$[\cat{M},[\cat{M},\cat{C}]]$};
      \node (B) at (-4,0) {$[\cat{M},\cat{C}]$};
      \node (C) at (-7,0) {$\cat{C}$};
      \node (D)  at (-2,-1.5) {$[\cat{M},\cat{C}]$};
      \node (E) at (-4,1.5) {$[\cat{M},[\cat{M},\cat{C}]]$};
      \node (F) at (0,1.5) {$[\cat{M},[\cat{M},[\cat{M},\cat{C}]]]$};
      \node (G) at (-7,1.5) {$[\cat{M},\cat{C}]$};
      \draw [->] (B) to [bend left=0]node [auto,labelsize] {$[\cat{M},\gm{T}]$}(A);
      \draw [->] (C) to [bend left=0]node [auto,swap,labelsize] {$\gm{T}$}(B);
      \draw [->] (D) to [bend right=20]node [auto,swap,labelsize] {$M_\cat{C}$}(A);
      \draw [->] (C) to [bend right=15]node [auto,swap,labelsize] {$\gm{T}$}(D);
      \draw [->] (B) to [bend right=0]node [auto,swap,labelsize] {$M_\cat{C}$}(E);
      \draw [->] (C) to [bend right=0]node [auto,labelsize] {$\gm{T}$}(G);
      \draw [->] (G) to [bend right=0]node [auto,labelsize] {$[\cat{M},\gm{T}]$}(E);
      \draw [->] (E) to [bend right=0]node [auto,labelsize] {$[\cat{M},[\cat{M},\gm{T}]]$}(F);
      \draw [->] (A) to [bend right=0]node [auto,swap,labelsize] {$M_{[\cat{M},\cat{C}]}$}(F);
      \draw [2cell] (-2,-0.1) to node [auto,labelsize]{$\mu$}(-2,-1.1);
      \draw [2cell] (-5.5,1.4) to node [auto,labelsize]{$\mu$}(-5.5,0.1);
        \end{scope}
\end{tikzpicture}\\
=\quad
\begin{tikzpicture}[baseline=-\the\dimexpr\fontdimen22\textfont2\relax ]
        \begin{scope} 
      \node (A)  at (0,0) {$[\cat{M},[\cat{M},\cat{C}]]$};
      \node (B) at (-4,0) {$[\cat{M},\cat{C}]$};
      \node (C) at (-7,0) {$\cat{C}$};
      \node (D)  at (-2,-1.5) {$[\cat{M},\cat{C}]$};
      \node (E) at (-4,1.5) {$[\cat{M},[\cat{M},\cat{C}]]$};
      \node (F) at (0,1.5) {$[\cat{M},[\cat{M}, [\cat{M},\cat{C}]]]$};
      \draw [->] (B) to [bend left=0]node [auto,labelsize] {$[\cat{M}, \gm{T}]$}(A);
      \draw [->] (C) to [bend left=0]node [auto,labelsize] {$\gm{T}$}(B);
      \draw [->] (D) to [bend right=20]node [auto,swap,labelsize] {$M_\cat{C}$}(A);
      \draw [->] (C) to [bend right=15]node [auto,swap,labelsize] {$\gm{T}$}(D);
      \draw [->] (B) to [bend right=0]node [auto,labelsize] {$[\cat{M},\gm{T}]$}(E);
      \draw [->] (E) to [bend right=0]node [auto,labelsize] {$[\cat{M},[\cat{M},\gm{T}]]$}(F);
      \draw [->] (A) to [bend right=0]node [auto,swap,labelsize] {$[\cat{M},M_\cat{C}]$}(F);
      \draw [2cell] (-2,-0.1) to node [auto,labelsize]{$\mu$}(-2,-1.1);
      \draw [2cell] (-2,1.4) to node [auto,labelsize]{$[\cat{M},\mu]$}(-2,0.5);
        \end{scope}
\end{tikzpicture}
\end{multline*}
The axioms correspond respectively to (GM4), (GM5) and (GM6).
\end{pf}

\subsection{Indexed monads as monads in $\Epm$}
\label{subsec-im-mnd-epm}
Indexed monads also admit similar reductions to mere monads 
in the 2-categorical sense.
The first such reduction is achieved in the 2-category~$\Epm$, and 
will provide a basis for the 
\emph{Eilenberg--Moore construction} for indexed monads.

Let us fix a category~$\cat{B}$.

\begin{defn}
\label{defn-B-times}
There is a 2-comonad~$\cat{B}\times(-)=
(\cat{B}\times(-),H^\cat{B},M^\cat{B})$ on $\Cat$ defined as follows.
\begin{itemize}
\item The 2-functor~$\cat{B}\times(-)$ is the product category 
construction.
\item The 2-natural transformation~$H^\cat{B}\colon
\cat{B}\times(-)\Rightarrow\id{\Cat}$ has as its $\cat{X}$-component
the functor~$H^\cat{B}_\cat{X}
=H_\cat{X}\colon\cat{B}\times\cat{X}\to\cat{X}$
given by $(b,x)\mapsto x$.
\item The 2-natural transformation~$M^\cat{B}\colon
\cat{B}\times(-)\Rightarrow\cat{B}\times\cat{B}\times(-)$ has as its
$\cat{X}$-component the functor~$M^\cat{B}_\cat{X}=M_\cat{X}\colon
\cat{B}\times\cat{X}\to\cat{B}\times\cat{B}\times\cat{X}$ given by
$(b,x)\mapsto (b,b,x)$.\qedhere
\end{itemize}
\end{defn}

\begin{prop}
Let $\cat{C}$ be a category.
Then a $\cat{B}$-indexed monad on $\cat{C}$ is the same thing as a monad 
in $\Epm$ on $(\Cat,\cat{C})$ above the 2-comonad~$\cat{B}\times(-)
\colon\Cat\to\Cat$.
\end{prop}
\begin{pf}
The latter notion is given by the following data:
\begin{itemize}
\item A 1-cell~$(\cat{B}\times(-),\im{T})\colon
(\Cat,\cat{C})\to(\Cat,\cat{C})$ of $\Epm$, where $\im{T}$ is a functor of type
\[
\im{T}\quad\colon\quad\cat{B}\times\cat{C}
\quad\longrightarrow\quad\cat{C}.
\]
\item A 2-cell~$(H^\cat{B},\eta)\colon (\id{\Cat},\id{\cat{C}})\Rightarrow
(\cat{B}\times(-),\im{T})$ of $\Epm$, where $\eta$ is a natural transformation of type
\[
    \begin{tikzpicture}
      \node (a)  at (3.5,2) {$\cat{B}\times\cat{C}$};
      \node (ap) at (6.5,2) {$\cat{C}$};
      \node (b)  at (4.5,3.5) {$\cat{C}$};
      
      \draw [->] (a) to [bend right=0]node [auto,swap,labelsize] {$\im{T}$}(ap);
      \draw [->] (a)  to [bend left=20]node [auto,labelsize] {$H_\cat{C}$}(b);
      \draw [->] (b) to [bend left=20]node [auto,labelsize] {$\id{\cat{C}}$}(ap);
      \draw [2cellr] (4.5,2.1) to node [auto,swap,labelsize]{$\eta$}(b);
    \end{tikzpicture}
\]
\item A 2-cell~$(M^\cat{B},\mu)\colon (\cat{B}\times(-),\im{T})\circ
(\cat{B}\times(-),\im{T})\Rightarrow(\cat{B}\times(-),\im{T})$ of 
$\Epm$, where $\mu$ is a natural transformation of type
\[
\begin{tikzpicture}[baseline=-\the\dimexpr\fontdimen22\textfont2\relax ]
        \begin{scope} 
      \node (L)  at (0,0)  {$\cat{B}\times\cat{C}$};
      \node (R)  at (6,0) {$\cat{C}$};
      \node (B)  at (2,2) {$\cat{B}\times\cat{B}\times\cat{C}$};
      \node (N) at (4.7,1.3) {$\cat{B}\times\cat{C}$};
      %
      \draw [->] (L) to node  [auto,swap,labelsize]      {$\im{T}$} (R);
      \draw [->] (L) to [bend left=20]node  [auto,labelsize] {$M_\cat{C}$}(B);
      \draw [->] (B) to [bend left=7]node  [auto,labelsize] {$\cat{B}\times\im{T}$}(N);
      \draw [->] (N) to [bend left=5]node  [auto,labelsize] {$\im{T}$}(R);
      \draw [2cellnode] (B) to node[auto,labelsize] {$\mu$}(2,0);
        \end{scope}
\end{tikzpicture}
\]
\end{itemize}
These data satisfy the following axioms:
{\allowdisplaybreaks
\begin{align*}
\begin{tikzpicture}[baseline=-\the\dimexpr\fontdimen22\textfont2\relax ]
        \begin{scope} 
      \node (A)  at (0,0) {$\cat{B}\times\cat{B}\times\cat{C}$};
      \node (B) at (3,0) {$\cat{B}\times\cat{C}$};
      \node (C) at (5,0) {$\cat{C}$};
      \node (D)  at (1.5,-1.5) {$\cat{B}\times\cat{C}$};
      \node (E) at (3.5,1.5) {$\cat{C}$};
      \node (F) at (0.5,1.5) {$\cat{B}\times\cat{C}$};
      \draw [->] (A) to [bend left=0]node [auto,labelsize] {$\cat{B}\times \im{T}$}(B);
      \draw [->] (B) to [bend left=0]node [auto,labelsize] {$\im{T}$}(C);
      \draw [<-] (A) to [bend right=20]node [auto,swap,labelsize] {$M_\cat{C}$}(D);
      \draw [->] (D) to [bend right=20]node [auto,swap,labelsize] {$\im{T}$}(C);
      \draw [<-] (E) to [bend right=20]node [auto,swap,labelsize] {$H_\cat{C}$}(B);
      \draw [->] (F) to [bend right=0]node [auto,labelsize] {$\im{T}$}(E);
      \draw [->] (E) to [bend left=20]node [auto,labelsize] {$\id{\cat{C}}$}(C);
      \draw [<-] (F) to [bend right=20]node [auto,swap,labelsize] {$H_{\cat{B}\times\cat{C}}$}(A);
      \draw [2cell] (1.5,-0.1) to node [auto,labelsize]{$\mu$}(1.5,-1.1);
      \draw [2cell] (3.5,1.2) to node [auto,labelsize]{$\eta$}(3.5,0.3);
        \end{scope}
\end{tikzpicture}
\quad&=\quad
\begin{tikzpicture}[baseline=-\the\dimexpr\fontdimen22\textfont2\relax ]
        \begin{scope} 
       \node (A)  at (0,0) {$\cat{B}\times\cat{C}$};
       \node (C) at (2,0) {$\cat{C}$};
       \draw [->] (A) to [bend left=0]node [auto,labelsize] {$\im{T}$}(C);
        \end{scope}
\end{tikzpicture}\\
\begin{tikzpicture}[baseline=-\the\dimexpr\fontdimen22\textfont2\relax ]
        \begin{scope} 
      \node (A)  at (0,0) {$\cat{B}\times\cat{B}\times\cat{C}$};
      \node (B) at (3,0) {$\cat{B}\times\cat{C}$};
      \node (C) at (5,0) {$\cat{C}$};
      \node (D)  at (1.5,-1.5) {$\cat{B}\times\cat{C}$};
      \node (F) at (1,1.5) {$\cat{B}\times\cat{C}$};
      \draw [->] (A) to [bend left=0]node [auto,labelsize] {$\cat{B}\times \im{T}$}(B);
      \draw [->] (B) to [bend left=0]node [auto,labelsize] {$\im{T}$}(C);
      \draw [<-] (A) to [bend right=20]node [auto,swap,labelsize] {$M_\cat{C}$}(D);
      \draw [->] (D) to [bend right=20]node [auto,swap,labelsize] {$\im{T}$}(C);
      \draw [->] (F) to [bend left=20]node [auto,labelsize] {$\cat{B}\times\id{\cat{C}}$}(B);
      \draw [<-] (F) to [bend right=20]node [auto,swap,labelsize] {$\cat{B}\times H_{\cat{C}}$}(A);
      \draw [2cell] (1.5,-0.1) to node [auto,labelsize]{$\mu$}(1.5,-1.1);
      \draw [2cell] (1,1.2) to node [auto,labelsize]{$\cat{B}\times\eta$}(1,0.3);
        \end{scope}
\end{tikzpicture}
\quad&=\quad
\begin{tikzpicture}[baseline=-\the\dimexpr\fontdimen22\textfont2\relax ]
        \begin{scope} 
       \node (A)  at (0,0) {$\cat{B}\times\cat{C}$};
       \node (C) at (2,0) {$\cat{C}$};
       \draw [->] (A) to [bend left=0]node [auto,labelsize] {$\im{T}$}(C);
        \end{scope}
\end{tikzpicture}
\end{align*}
} %
\vspace{-15pt}
\begin{multline*}
\begin{tikzpicture}[baseline=-\the\dimexpr\fontdimen22\textfont2\relax ]
        \begin{scope} 
      \node (A)  at (0,0) {$\cat{B}\times \cat{B}\times\cat{C}$};
      \node (B) at (4,0) {$\cat{B}\times\cat{C}$};
      \node (C) at (7,0) {$\cat{C}$};
      \node (D)  at (2,-1.5) {$\cat{B}\times\cat{C}$};
      \node (E) at (4,1.5) {$\cat{B}\times\cat{B}\times\cat{C}$};
      \node (F) at (0,1.5) {$\cat{B}\times\cat{B}\times \cat{B}\times\cat{C}$};
      \node (G) at (7,1.5) {$\cat{B}\times\cat{C}$};
      \draw [->] (A) to [bend left=0]node [auto,labelsize] {$\cat{B}\times \im{T}$}(B);
      \draw [->] (B) to [bend left=0]node [auto,swap,labelsize] {$\im{T}$}(C);
      \draw [<-] (A) to [bend right=20]node [auto,swap,labelsize] {$M_\cat{C}$}(D);
      \draw [->] (D) to [bend right=15]node [auto,swap,labelsize] {$\im{T}$}(C);
      \draw [<-] (E) to [bend right=0]node [auto,swap,labelsize] {$M_\cat{C}$}(B);
      \draw [->] (G) to [bend right=0]node [auto,labelsize] {$\im{T}$}(C);
      \draw [->] (E) to [bend right=0]node [auto,labelsize] {$\cat{B}\times \im{T}$}(G);
      \draw [->] (F) to [bend right=0]node [auto,labelsize] {$\cat{B}\times\cat{B}\times \im{T}$}(E);
      \draw [<-] (F) to [bend right=0]node [auto,swap,labelsize] {$M_{\cat{B}\times\cat{C}}$}(A);
      \draw [2cell] (2,-0.1) to node [auto,labelsize]{$\mu$}(2,-1.1);
      \draw [2cell] (5.5,1.4) to node [auto,labelsize]{$\mu$}(5.5,0.1);
        \end{scope}
\end{tikzpicture}\\
=\quad
\begin{tikzpicture}[baseline=-\the\dimexpr\fontdimen22\textfont2\relax ]
        \begin{scope} 
      \node (A)  at (0,0) {$\cat{B}\times \cat{B}\times\cat{C}$};
      \node (B) at (4,0) {$\cat{B}\times\cat{C}$};
      \node (C) at (7,0) {$\cat{C}$};
      \node (D)  at (2,-1.5) {$\cat{B}\times\cat{C}$};
      \node (E) at (4,1.5) {$\cat{B}\times\cat{B}\times\cat{C}$};
      \node (F) at (0,1.5) {$\cat{B}\times\cat{B}\times\cat{B}\times\cat{C}$};
      \draw [->] (A) to [bend left=0]node [auto,labelsize] {$\cat{B}\times \im{T}$}(B);
      \draw [->] (B) to [bend left=0]node [auto,labelsize] {$\im{T}$}(C);
      \draw [<-] (A) to [bend right=20]node [auto,swap,labelsize] {$M_\cat{C}$}(D);
      \draw [->] (D) to [bend right=15]node [auto,swap,labelsize] {$\im{T}$}(C);
      \draw [->] (E) to [bend right=0]node [auto,labelsize] {$\cat{B}\times\im{T}$}(B);
      \draw [->] (F) to [bend right=0]node [auto,labelsize] {$\cat{B}\times\cat{B}\times \im{T}$}(E);
      \draw [<-] (F) to [bend right=0]node [auto,swap,labelsize] {$\cat{B}\times M_\cat{C}$}(A);
      \draw [2cell] (2,-0.1) to node [auto,labelsize]{$\mu$}(2,-1.1);
      \draw [2cell] (2,1.4) to node [auto,labelsize]{$\cat{B}\times\mu$}(2,0.5);
        \end{scope}
\end{tikzpicture}
\end{multline*}
The axioms correspond respectively to (IM5), (IM6) and (IM7).
\end{pf}

\subsection{Indexed monads as monads in $\Emp$}
\label{subsec-im-mnd-emp}
We fix a category~$\cat{B}$.

\begin{defn}
There is a 2-monad~$[\cat{B},-]=([\cat{B},-],H_\cat{B},M_\cat{B})$
on $\Cat$ defined as follows.
\begin{itemize}
\item The 2-functor~$[\cat{B},-]\colon\Cat\to\Cat$ is the functor category
construction.
\item The 2-natural transformation~$H_\cat{B}\colon
\id{\Cat}\Rightarrow[\cat{B},-]$ has as its $\cat{X}$-component
the functor~$H_{\cat{B},\cat{X}}=H_\cat{X}\colon
\cat{X}\to[\cat{B},\cat{X}]$ given by $x\mapsto (b\mapsto x)$.
\item The 2-natural transformation~$M_\cat{B}\colon
[\cat{B},[\cat{B},-]]\Rightarrow[\cat{B},-]$ has as its $\cat{X}$-component
the functor~$M_{\cat{B},\cat{X}}=M_\cat{X}\colon
[\cat{B},[\cat{B},\cat{X}]]\to [\cat{B},\cat{X}]$ given by
$\Xi\mapsto (b\mapsto (\Xi b)b)$.\qedhere
\end{itemize}
\end{defn}

\begin{prop}
Let $\cat{C}$ be a category. Then a $\cat{B}$-indexed monad on $\cat{C}$
is the same thing as a monad in $\Emp$ on $(\Cat,\cat{C})$ above the
2-monad~$[\cat{B},-]\colon\Cat\to\Cat$.
\end{prop}
\begin{pf}
The latter notion is given by the following data:
\begin{itemize}
\item A 1-cell~$([\cat{B},-],\im{T})\colon(\Cat,\cat{C})\to
(\Cat,\cat{C})$ of $\Emp$, where $\im{T}$ is a functor of type
\[
\im{T}\quad\colon\quad
\cat{C}\quad\longrightarrow\quad[\cat{B},\cat{C}].
\]
\item A 2-cell~$(H_\cat{B},\eta)\colon(\id{\Cat},\id{\cat{C}})
\Rightarrow([\cat{B},-],\im{T})$ of $\Emp$, where 
$\eta$ is a natural transformation of type
\[
    \begin{tikzpicture}
      \node (a)  at (3.5,2) {$\cat{C}$};
      \node (ap) at (6.5,2) {$[\cat{B},\cat{C}]$};
      \node (b)  at (5.5,3.5) {$\cat{C}$};
      
      \draw [->] (a) to [bend right=0]node [auto,swap,labelsize] {$\im{T}$}(ap);
      \draw [->] (a)  to [bend left=20]node [auto,labelsize] {$\id{\cat{C}}$}(b);
      \draw [->] (b) to [bend left=20]node [auto,labelsize] {$H_{\cat{C}}$}(ap);
      \draw [2cellr] (5.5,2.1) to node [auto,swap,labelsize]{$\eta$}(b);
    \end{tikzpicture}
\]
\item A 2-cell~$(M_\cat{B},\mu)\colon([\cat{B},-],\im{T})\circ([\cat{B},-],\im{T})
\Rightarrow([\cat{B},-],\im{T})$ of $\Emp$, where $\mu$ is a natural 
transformation of type
\[
\begin{tikzpicture}[baseline=-\the\dimexpr\fontdimen22\textfont2\relax ]
        \begin{scope} [shift={(0,0.75)}]
      \node (M)  at (0,0)  {$[\cat{B},\cat{C}]$};
      \node (N)  at (-4.7,1.3) {$[\cat{B},\cat{C}]$};
      \node (RR)  at (-6,0) {$\cat{C}$};
      \node (B)  at (-2,2) {$[\cat{B},[\cat{B},\cat{C}]]$};
      \draw [->] (RR) to node  [auto,swap,labelsize]      {$\im{T}$} (M);
      \draw [->] (B) to [bend left=20]node  [auto,labelsize] {$M_\cat{C}$}(M);
      \draw [->] (N) to [bend left=8]node  [auto,labelsize] {$[\cat{B},\im{T}]$}(B);
      \draw [->] (RR) to [bend left=4]node  [auto,labelsize] {$\im{T}$}(N);
      \draw [2cellnode] (B) to node[auto,labelsize] {$\mu$}(-2,0);
        \end{scope}
\end{tikzpicture}
\]
\end{itemize}
These data satisfy the following axioms:
{\allowdisplaybreaks
\begin{align*}
\begin{tikzpicture}[baseline=-\the\dimexpr\fontdimen22\textfont2\relax ]
        \begin{scope} 
      \node (A)  at (0,0) {$[\cat{B},[\cat{B},\cat{C}]]$};
      \node (B) at (-3,0) {$[\cat{B},\cat{C}]$};
      \node (C) at (-5,0) {$\cat{C}$};
      \node (D)  at (-1.5,-1.5) {$[\cat{B},\cat{C}]$};
      \node (E) at (-3.5,1.5) {$\cat{C}$};
      \node (F) at (-0.5,1.5) {$[\cat{B},\cat{C}]$};
      \draw [->] (B) to [bend left=0]node [auto,labelsize] {$[\cat{B},\im{T}]$}(A);
      \draw [->] (C) to [bend left=0]node [auto,labelsize] {$\im{T}$}(B);
      \draw [<-] (D) to [bend right=20]node [auto,swap,labelsize] {$M_\cat{C}$}(A);
      \draw [->] (C) to [bend right=20]node [auto,swap,labelsize] {$\im{T}$}(D);
      \draw [<-] (B) to [bend right=20]node [auto,swap,labelsize] {$H_\cat{C}$}(E);
      \draw [->] (E) to [bend right=0]node [auto,labelsize] {$\im{T}$}(F);
      \draw [->] (C) to [bend left=20]node [auto,labelsize] {$\id{\cat{C}}$}(E);
      \draw [<-] (A) to [bend right=20]node [auto,swap,labelsize] {$H_{[\cat{B},\cat{C}]}$}(F);
      \draw [2cell] (-1.5,-0.1) to node [auto,labelsize]{$\mu$}(-1.5,-1.1);
      \draw [2cell] (-3.5,1.2) to node [auto,labelsize]{$\eta$}(-3.5,0.3);
        \end{scope}
\end{tikzpicture}
\quad&=\quad
\begin{tikzpicture}[baseline=-\the\dimexpr\fontdimen22\textfont2\relax ]
        \begin{scope} 
       \node (A)  at (0,0) {$[\cat{B},\cat{C}]$};
       \node (C) at (-2,0) {$\cat{C}$};
       \draw [->] (C) to [bend left=0]node [auto,labelsize] {$\im{T}$}(A);
        \end{scope}
\end{tikzpicture}\\
\begin{tikzpicture}[baseline=-\the\dimexpr\fontdimen22\textfont2\relax ]
        \begin{scope} 
      \node (A)  at (0,0) {$[\cat{B},[\cat{B},\cat{C}]]$};
      \node (B) at (-3,0) {$[\cat{B},\cat{C}]$};
      \node (C) at (-5,0) {$\cat{C}$};
      \node (D)  at (-1.5,-1.5) {$[\cat{B},\cat{C}]$};
      \node (F) at (-1,1.5) {$[\cat{B},\cat{C}]$};
      \draw [->] (B) to [bend left=0]node [auto,labelsize] {$[\cat{B},\im{T}]$}(A);
      \draw [->] (C) to [bend left=0]node [auto,labelsize] {$\im{T}$}(B);
      \draw [<-] (D) to [bend right=20]node [auto,swap,labelsize] {$M_\cat{C}$}(A);
      \draw [->] (C) to [bend right=20]node [auto,swap,labelsize] {$\im{T}$}(D);
      \draw [->] (B) to [bend left=20]node [auto,labelsize] {$[\cat{B},\id{\cat{C}}]$}(F);
      \draw [<-] (A) to [bend right=20]node [auto,swap,labelsize] {$[\cat{B},H_{\cat{C}}]$}(F);
      \draw [2cell] (-1.5,-0.1) to node [auto,labelsize]{$\mu$}(-1.5,-1.1);
      \draw [2cell] (-1,1.2) to node [auto,swap,labelsize]{$[\cat{B},\eta]$}(-1,0.3);
        \end{scope}
\end{tikzpicture}
\quad&=\quad
\begin{tikzpicture}[baseline=-\the\dimexpr\fontdimen22\textfont2\relax ]
        \begin{scope} 
       \node (A)  at (0,0) {$[\cat{B},\cat{C}]$};
       \node (C) at (-2,0) {$\cat{C}$};
       \draw [->] (C) to [bend left=0]node [auto,labelsize] {$\im{T}$}(A);
        \end{scope}
\end{tikzpicture}
\end{align*}
}
\vspace{-20pt}
\begin{multline*}
\begin{tikzpicture}[baseline=-\the\dimexpr\fontdimen22\textfont2\relax ]
        \begin{scope} 
      \node (A)  at (0,0) {$[\cat{B},[\cat{B},\cat{C}]]$};
      \node (B) at (-4,0) {$[\cat{B},\cat{C}]$};
      \node (C) at (-7,0) {$\cat{C}$};
      \node (D)  at (-2,-1.5) {$[\cat{B},\cat{C}]$};
      \node (E) at (-4,1.5) {$[\cat{B},[\cat{B},\cat{C}]]$};
      \node (F) at (0,1.5) {$[\cat{B},[\cat{B},[\cat{B},\cat{C}]]]$};
      \node (G) at (-7,1.5) {$[\cat{B},\cat{C}]$};
      \draw [->] (B) to [bend left=0]node [auto,labelsize] {$[\cat{B},\im{T}]$}(A);
      \draw [->] (C) to [bend left=0]node [auto,swap,labelsize] {$\im{T}$}(B);
      \draw [<-] (D) to [bend right=20]node [auto,swap,labelsize] {$M_\cat{C}$}(A);
      \draw [->] (C) to [bend right=15]node [auto,swap,labelsize] {$\im{T}$}(D);
      \draw [<-] (B) to [bend right=0]node [auto,swap,labelsize] {$M_\cat{C}$}(E);
      \draw [->] (C) to [bend right=0]node [auto,labelsize] {$\im{T}$}(G);
      \draw [->] (G) to [bend right=0]node [auto,labelsize] {$[\cat{B},\im{T}]$}(E);
      \draw [->] (E) to [bend right=0]node [auto,labelsize] {$[\cat{B},[\cat{B},\im{T}]]$}(F);
      \draw [<-] (A) to [bend right=0]node [auto,swap,labelsize] {$M_{[\cat{B},\cat{C}]}$}(F);
      \draw [2cell] (-2,-0.1) to node [auto,labelsize]{$\mu$}(-2,-1.1);
      \draw [2cell] (-5.5,1.4) to node [auto,labelsize]{$\mu$}(-5.5,0.1);
        \end{scope}
\end{tikzpicture}\\
=\quad
\begin{tikzpicture}[baseline=-\the\dimexpr\fontdimen22\textfont2\relax ]
        \begin{scope} 
      \node (A)  at (0,0) {$[\cat{B},[\cat{B},\cat{C}]]$};
      \node (B) at (-4,0) {$[\cat{B},\cat{C}]$};
      \node (C) at (-7,0) {$\cat{C}$};
      \node (D)  at (-2,-1.5) {$[\cat{B},\cat{C}]$};
      \node (E) at (-4,1.5) {$[\cat{B},[\cat{B},\cat{C}]]$};
      \node (F) at (0,1.5) {$[\cat{B},[\cat{B}, [\cat{B},\cat{C}]]]$};
      \draw [->] (B) to [bend left=0]node [auto,labelsize] {$[\cat{B}, \im{T}]$}(A);
      \draw [->] (C) to [bend left=0]node [auto,labelsize] {$\im{T}$}(B);
      \draw [<-] (D) to [bend right=20]node [auto,swap,labelsize] {$M_\cat{C}$}(A);
      \draw [->] (C) to [bend right=15]node [auto,swap,labelsize] {$\im{T}$}(D);
      \draw [->] (B) to [bend right=0]node [auto,labelsize] {$[\cat{B},\im{T}]$}(E);
      \draw [->] (E) to [bend right=0]node [auto,labelsize] {$[\cat{B},[\cat{B},\im{T}]]$}(F);
      \draw [<-] (A) to [bend right=0]node [auto,swap,labelsize] {$[\cat{B},M_\cat{C}]$}(F);
      \draw [2cell] (-2,-0.1) to node [auto,labelsize]{$\mu$}(-2,-1.1);
      \draw [2cell] (-2,1.4) to node [auto,labelsize]{$[\cat{B},\mu]$}(-2,0.5);
        \end{scope}
\end{tikzpicture}
\end{multline*}
The axioms correspond respectively to (IM5), (IM6) and (IM7).
\end{pf}

\subsection{Graded and indexed comonads as comonads}
\label{subsec-gicomnd-comnd}
Here we list the analogous results for graded and indexed
\emph{comonads}.
Let $\cat{M}=(\cat{M},\tensor,\e)$ be a strict monoidal category and 
$\cat{B}$ a category.

\begin{prop}
Let $\cat{C}$ be a category. Then an $\cat{M}$-graded comonad on $\cat{C}$
is the same thing as a comonad in $\Epm$ on $(\Cat,\cat{C})$ 
above the 2-monad~$\cat{M}\times(-)\colon\Cat\to\Cat$.
\end{prop}
\begin{pf}
The latter notion is given by the following data:
\begin{itemize}
\item A 1-cell~$(\cat{M}\times(-),\gm{S})\colon(\Cat,\cat{C})\to
(\Cat,\cat{C})$ of $\Epm$, where $\gm{S}$ is a functor of type
\[
\gm{S}\quad\colon\quad\cat{M}\times\cat{C}\quad\longrightarrow\quad\cat{C}.
\]
\item A 2-cell~$(H^\cat{M},\varepsilon)\colon(\cat{M}\times(-),\gm{S})
\Rightarrow(\id{\Cat},\id{\cat{C}})$ of $\Epm$, where $\varepsilon$ is 
a natural transformation of type
\[
    \begin{tikzpicture}
      \node (a)  at (3.5,2) {$\cat{C}$};
      \node (ap) at (6.5,2) {$\cat{C}$};
      \node (b)  at (4.5,3.5) {$\cat{M}\times\cat{C}$};
      
      \draw [->] (a) to [bend right=0]node [auto,swap,labelsize] {$\id{\cat{C}}$}(ap);
      \draw [->] (a)  to [bend left=20]node [auto,labelsize] {$H_\cat{C}$}(b);
      \draw [->] (b) to [bend left=20]node [auto,labelsize] {$\gm{S}$}(ap);
      \draw [2cellr] (4.5,2.1) to node [auto,swap,labelsize]{$\varepsilon$}(b);
    \end{tikzpicture}
\]
\item A 2-cell~$(M^\cat{M},\delta)\colon(\cat{M}\times(-),\gm{S})
\Rightarrow(\cat{M}\times(-),\gm{S})\circ(\cat{M}\times(-),\gm{S})$
of $\Epm$, where $\delta$ is a natural transformation of type
    \[
    \begin{tikzpicture}
      \node (L) at (0,0) {$\cat{M}\times\cat{M}\times\cat{C}$};
      \node (M) at (3,0) {$\cat{M}\times\cat{C}$};
      \node (R) at (6,0) {$\cat{C}$};
      \node (T) at (2,2) {$\cat{M}\times\cat{C}$};
      
      \draw [->] (L) to node [auto,swap,labelsize] {$\cat{M}\times\gm{S}$}(M);
      \draw [->] (M) to node [auto,swap,labelsize] {$\gm{S}$}(R);
      \draw [->] (L) to [bend left=20]node [auto,labelsize] {$M_\cat{C}$}(T);
      \draw [->] (T) to [bend left=20]node [auto,labelsize] {$\gm{S}$}(R);
            
      \draw [2cell] (T) to node [auto,labelsize]{$\delta$}(2,0.1);
    \end{tikzpicture}
    \]
\end{itemize}
They satisfy axioms corresponding to (GC4), (GC5) and (GC6).
\end{pf}

\begin{prop}
Let $\cat{C}$ be a category. Then an $\cat{M}$-graded comonad on $\cat{C}$
is the same thing as a comonad in $\Emp$ on $(\Cat,\cat{C})$ 
above the 2-comonad~$[\cat{M},-]\colon\Cat\to\Cat$.
\end{prop}
\begin{pf}
The latter notion is given by the following data:
\begin{itemize}
\item A 1-cell~$([\cat{M},-],\gm{S})\colon(\Cat,\cat{C})\to
(\Cat,\cat{C})$ of $\Emp$, where $\gm{S}$ is a functor of type
\[
\gm{S}\quad\colon\quad\cat{C}\quad\longrightarrow\quad[\cat{M},\cat{C}].
\]
\item A 2-cell~$(H_\cat{M},\varepsilon)\colon([\cat{M},-],\gm{S})
\Rightarrow(\id{\Cat},\id{\cat{C}})$ of $\Emp$, where $\varepsilon$ is 
a natural transformation of type
\[
    \begin{tikzpicture}
      \node (a)  at (3.5,2) {$\cat{C}$};
      \node (ap) at (6.5,2) {$\cat{C}$};
      \node (b)  at (5.5,3.5) {$[\cat{M},\cat{C}]$};
      
      \draw [->] (a) to [bend right=0]node [auto,swap,labelsize] {$\id{\cat{C}}$}(ap);
      \draw [->] (a)  to [bend left=20]node [auto,labelsize] {$\gm{S}$}(b);
      \draw [->] (b) to [bend left=20]node [auto,labelsize] {$H_\cat{C}$}(ap);
      \draw [2cellr] (5.5,2.1) to node [auto,swap,labelsize]{$\varepsilon$}(b);
    \end{tikzpicture}
\]
\item A 2-cell~$(M_\cat{M},\delta)\colon([\cat{M},-],\gm{S})
\Rightarrow([\cat{M},-],\gm{S})\circ([\cat{M},-],\gm{S})$
of $\Emp$, where $\delta$ is a natural transformation of type
    \[
    \begin{tikzpicture}
      \node (L) at (0,0) {$\cat{C}$};
      \node (M) at (3,0) {$[\cat{M},\cat{C}]$};
      \node (R) at (6,0) {$[\cat{M},[\cat{M},\cat{C}]]$};
      \node (T) at (4,2) {$[\cat{M},\cat{C}]$};
      
      \draw [->] (L) to node [auto,swap,labelsize] {$\gm{S}$}(M);
      \draw [->] (M) to node [auto,swap,labelsize] {$[\cat{M},\gm{S}]$}(R);
      \draw [->] (L) to [bend left=20]node [auto,labelsize] {$\gm{S}$}(T);
      \draw [->] (T) to [bend left=20]node [auto,labelsize] {$M_\cat{C}$}(R);
            
      \draw [2cell] (T) to node [auto,labelsize]{$\delta$}(4,0.1);
    \end{tikzpicture}
    \]
\end{itemize}
They satisfy axioms corresponding to (GC4), (GC5) and (GC6).
\end{pf}

\begin{prop}
Let $\cat{C}$ be a category. Then a $\cat{B}$-indexed comonad on $\cat{C}$
is the same thing as a comonad in $\Epp$ on $(\Cat,\cat{C})$ 
above the 2-comonad~$\cat{B}\times(-)\colon\Cat\to\Cat$.
\end{prop}
\begin{pf}
The latter notion is given by the following data:
\begin{itemize}
\item A 1-cell~$(\cat{B}\times(-),\im{S})\colon(\Cat,\cat{C})\to
(\Cat,\cat{C})$ of $\Epp$, where $\im{S}$ is a functor of type
\[
\im{S}\quad\colon\quad\cat{B}\times\cat{C}\quad\longrightarrow\quad\cat{C}.
\]
\item A 2-cell~$(H^\cat{B},\varepsilon)\colon(\cat{B}\times(-),\im{S})
\Rightarrow(\id{\Cat},\id{\cat{C}})$ of $\Epp$, where $\varepsilon$ is 
a natural transformation of type
\[
    \begin{tikzpicture}
      \node (a)  at (3.5,2) {$\cat{B}\times\cat{C}$};
      \node (ap) at (6.5,2) {$\cat{C}$};
      \node (b)  at (4.5,0.5) {$\cat{C}$};
      
      \draw [->] (a) to [bend right=0]node [auto,labelsize] {$\im{S}$}(ap);
      \draw [->] (a)  to [bend right=20]node [auto,swap,labelsize] {$H_\cat{C}$}(b);
      \draw [->] (b) to [bend right=20]node [auto,swap,labelsize] {$\id{\cat{C}}$}(ap);
      \draw [2cell] (4.5,1.9) to node [auto,labelsize]{$\varepsilon$}(b);
    \end{tikzpicture}
\]
\item A 2-cell~$(M^\cat{B},\delta)\colon(\cat{B}\times(-),\im{S})
\Rightarrow(\cat{B}\times(-),\im{S})\circ(\cat{B}\times(-),\im{S})$
of $\Epp$, where $\delta$ is a natural transformation of type
    \[
\begin{tikzpicture}[baseline=-\the\dimexpr\fontdimen22\textfont2\relax ]
        \begin{scope} 
      \node (L)  at (0,0)  {$\cat{B}\times\cat{C}$};
      \node (R)  at (6,0) {$\cat{C}$};
      \node (B)  at (2,-2) {$\cat{B}\times\cat{B}\times\cat{C}$};
      \node (N) at (4.7,-1.3) {$\cat{B}\times\cat{C}$};
      %
      \draw [->] (L) to node  [auto,labelsize]      {$\im{S}$} (R);
      \draw [->] (L) to [bend left=-20]node  [auto,swap,labelsize] {$M_\cat{C}$}(B);
      \draw [->] (B) to [bend left=-7]node  [auto,swap,labelsize] {$\cat{B}\times\im{S}$}(N);
      \draw [->] (N) to [bend left=-5]node  [auto,swap,labelsize] {$\im{S}$}(R);
      \draw [2cellnode] (2,0) to node[auto,labelsize] {$\delta$}(B);
        \end{scope}
\end{tikzpicture}
    \]
\end{itemize}
They satisfy axioms corresponding to (IC5), (IC6) and (IC7).
\end{pf}

\begin{prop}
Let $\cat{C}$ be a category. Then an $\cat{B}$-graded comonad on $\cat{C}$
is the same thing as a comonad in $\Emm$ on $(\Cat,\cat{C})$ 
above the 2-monad~$[\cat{B},-]\colon\Cat\to\Cat$.
\end{prop}
\begin{pf}
The latter notion is given by the following data:
\begin{itemize}
\item A 1-cell~$([\cat{B},-],\im{S})\colon(\Cat,\cat{C})\to
(\Cat,\cat{C})$ of $\Emm$, where $\im{S}$ is a functor of type
\[
\im{S}\quad\colon\quad\cat{C}\quad\longrightarrow\quad[\cat{B},\cat{C}].
\]
\item A 2-cell~$(H_\cat{B},\varepsilon)\colon([\cat{B},-],\im{S})
\Rightarrow(\id{\Cat},\id{\cat{C}})$ of $\Emm$, where $\varepsilon$ is 
a natural transformation of type
\[
    \begin{tikzpicture}
      \node (a)  at (3.5,-2) {$\cat{C}$};
      \node (ap) at (6.5,-2) {$[\cat{B},\cat{C}]$};
      \node (b)  at (5.5,-3.5) {$\cat{C}$};
      
      \draw [->] (a) to [bend right=0]node [auto,labelsize] {$\im{S}$}(ap);
      \draw [->] (a)  to [bend left=-20]node [auto,swap,labelsize] {$\id{\cat{C}}$}(b);
      \draw [->] (b) to [bend left=-20]node [auto,swap,labelsize] {$H_{\cat{C}}$}(ap);
      \draw [2cell] (5.5,-2.1) to node [auto,labelsize]{$\varepsilon$}(b);
    \end{tikzpicture}
\]
\item A 2-cell~$(M_\cat{B},\delta)\colon([\cat{B},-],\im{S})
\Rightarrow([\cat{B},-],\im{S})\circ([\cat{B},-],\im{S})$
of $\Emm$, where $\delta$ is a natural transformation of type
\[
\begin{tikzpicture}[baseline=-\the\dimexpr\fontdimen22\textfont2\relax ]
        \begin{scope} [shift={(0,0.75)}]
      \node (M)  at (0,0)  {$[\cat{B},\cat{C}]$};
      \node (N)  at (-4.7,-1.3) {$[\cat{B},\cat{C}]$};
      \node (RR)  at (-6,0) {$\cat{C}$};
      \node (B)  at (-2,-2) {$[\cat{B},[\cat{B},\cat{C}]]$};
      \draw [->] (RR) to node  [auto,labelsize]      {$\im{S}$} (M);
      \draw [->] (B) to [bend left=-20]node  [auto,swap,labelsize] {$M_\cat{C}$}(M);
      \draw [->] (N) to [bend left=-8]node  [auto,swap,labelsize] {$[\cat{B},\im{S}]$}(B);
      \draw [->] (RR) to [bend left=-4]node  [auto,swap,labelsize] {$\im{S}$}(N);
      \draw [2cellnoder] (B) to node[auto,swap,labelsize] {$\delta$}(-2,0);
        \end{scope}
\end{tikzpicture}
\]
\end{itemize}
They satisfy axioms corresponding to (IC5), (IC6) and (IC7).
\end{pf}

\section*{Notes}
The definitions of the 2-categories~$\Epp$, $\Epm$, $\Emp$ and $\Emm$
have occurred to me after I learned from Paul-Andr\'e Melli\`es
his key observation that, by ``enlarging'' $\Cat$ in a certain way,
one can regard graded monads as mere monads in the 2-categorical sense;
indeed, his 2-category was a full sub 2-category of $\Epp$. 
When actually writing down the definition of 
the 2-category~$\Epp$ for the first time, Kenji Maillard helped me
by telling me the (perhaps folklore) 
view of the Grothendieck construction as a certain comma construction.

The reduction of graded monads to mere monads are presented in 
\cite{fujii-katsumata-mellies}.

The observation that the notion of indexed monad can also be 
reduced to monads using these 2-categories seems to be new here.

\chapter{The main constructions}
\label{chap-main-constr}
In this chapter, we describe the Eilenberg--Moore and the Kleisli
constructions for graded and indexed monads
(except for the Kleisli construction for indexed monads).
These constructions are natural yet nontrivial generalization of
the classical Eilenberg--Moore and Kleisli constructions.
Moreover, our constructions satisfy the relevant 2-dimensional
universal properties in naturally arising 2-categories, i.e., 
they produce 
\emph{Eilenberg--Moore} and \emph{Kleisli objects} respectively in 
appropriate 2-categories introduced in 
Chapter~\ref{chap-four-2-cats}; we regard this fact as 
the major justification of the definitions given below. 
More precisely, the relationship of the 2-categories~$\Epp$,
$\Epm$, $\Emp$ and $\Emm$,   
graded and indexed (co)monads, and the generalized
(co)Eilenberg--Moore and (co)Kleisli constructions for
them is summarized in the following table, which is a refinement of 
the table appearing at the beginning of the previous chapter:
\vspace{0pt}
\begin{table}[H]
	\centering
	\begin{tabular}{c|c c c c}    
    & $\Epp$ & $\Epm$ & $\Emp$ & $\Emm$\\ \hline
    Graded monads    & EM &    &    & Kl\\
    Indexed monads   &    & EM & Kl$^\ast$ & \\
    Graded comonads  &    & coEM & coKl & \\
    Indexed comonads & coEM & & &coKl$^\ast$
	\end{tabular}
	\vspace{-15pt}
\end{table}
\noindent
So far, we have not been able to identify 
the (co)Kleisli constructions for indexed
(co)monads; the $^\ast$ mark indicates   
the conjectural status of the construction.

As we have already mentioned in Section~\ref{sec-ftm}, the 2-dimensional 
universality of Eilenberg--Moore and Kleisli objects is 
powerful enough to reconstruct some of the main development of 
the classical Eilenberg--Moore and Kleisli construction
abstractly.
However, in the current chapter we have chosen to start by
following more
closely the style of the classical theory, and 
construct adjunctions that generate monads
and comparison maps explicitly.
The discussion on 2-categorical properties of our
constructions, which seems to be largely of interest only to
the experts, is placed after that and the reader can 
harmlessly skip these parts.

We conclude this chapter by briefly indicating 
the suitable co-Eilenberg--Moore and co-Kleisli
constructions for graded and indexed comonads,
again except for the conjectural 
co-Kleisli construction for indexed comonads
which is left for future work.


\section{The Eilenberg--Moore construction for graded monads}
\label{sec-em-gm}
Let $\cat{M}=(\cat{M},\tensor,\e)$ be a strict monoidal category,
$\cat{C}$ a category, and 
$\gm{T}$ an $\cat{M}$-graded monad on $\cat{C}$.
Recall from Section~\ref{subsec-gm-mnd-epp}
that $\gm{T}$ may be seen as a monad in $\Epp$;
the Eilenberg--Moore adjunction for $\gm{T}$ 
lives in $\Epp$, and lies above
the Eilenberg--Moore adjunction for the 2-monad~$\cat{M}\times(-)$ on 
$\Cat$.
See the picture below for an illustration.
    
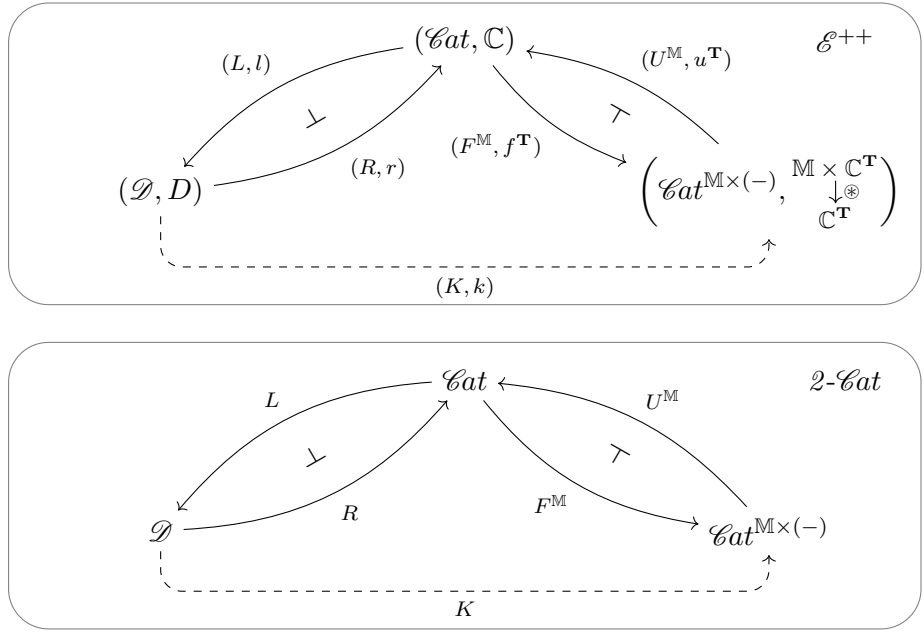
\[
    \begin{tikzpicture}
      \node (A)  at (8,0) {$\left(\EM{\Cat}{\cat{M}\times(-)},
            \valg{\cat{M}\times\EMg}{\emact}{\EMg}\right)$};
      \node (B) at (0,0) {$(\tcat{D},D)$};
      \node (C) at (4,2) {$(\Cat,\cat{C})$};
      
      \node at (2,1) [rotate=120,font=\large]{$\vdash$};
      \node at (6,1) [rotate=240,font=\large]{$\vdash$};
      
      \draw [->] (A) to [bend right=17]node [auto,labelsize,swap] {$(U^\cat{M},u^\gm{T})$}(C);
      \draw [->] (C) to [bend right=15]node [auto,labelsize,swap] {$(F^\cat{M},f^\gm{T})$}(A);
      
      \draw [->] (C) to [bend right=20]node [auto,swap,labelsize] {$(L,l)$}(B);
      \draw [->] (B) to [bend right=20]node [auto,swap,labelsize] {$(R,r)$}(C);
      
      \draw [->,dashed,rounded corners=8pt] (B)--(0,-1) to [bend right=0]node [auto,swap,labelsize]
       {$(K,k)$} (8,-1)--(A);
      
      \draw [rounded corners=15pt,myborder] (-2,-1.5) rectangle (10,2.5);
      \node at (9,2) {$\Epp$};

      \begin{scope}[shift={(0,-4.5)}]
       \node (A)  at (8,0) {$\EM{\Cat}{\cat{M}\times(-)}$};
       \node (B) at (0,0) {$\tcat{D}$};
       \node (C) at (4,2) {$\Cat$};
      
       \node at (2,1) [rotate=120,font=\large]{$\vdash$};
       \node at (6,1) [rotate=240,font=\large]{$\vdash$};
      
       \draw [->] (A) to [bend right=21]node [auto,labelsize,swap] {$U^\cat{M}$}(C);
       \draw [->] (C) to [bend right=21]node [auto,labelsize,swap] {$F^\cat{M}$}(A);
      
       \draw [->] (C) to [bend right=23]node [auto,swap,labelsize] (Fpp) {$L$}(B);
       \draw [->] (B) to [bend right=23]node [auto,swap,labelsize] (Fpp) {$R$}(C);
       
       \draw [->,dashed,rounded corners=8pt] (B)--(0,-0.8) to [bend right=0]node [auto,swap,labelsize]
              {$K$} (8,-0.8)--(A);
       \draw [rounded corners=15pt,myborder] (-2,-1.3) rectangle (10,2.5);
             \node at (9,2) {$\TCat$};
   
      \end{scope}
    \end{tikzpicture}
    \]
We will define the data appearing in the picture, one by one.
We write the functor part of $\gm{T}$
as $\act\colon \cat{M}\times\cat{C}\to\cat{C}$ 
as well and use the infix notation for it.

\subsection{The Eilenberg--Moore category}
Extending the classical construction of the Eilenberg--Moore category of 
an ordinary monad as the category of algebras,
the Eilenberg--Moore category of a graded monad is given as
the category of \emph{graded algebras}.

\begin{defn}
Define the category~$\EMg$ as follows:
\begin{itemize}
\item An object of $\EMg$ is a
\defemph{graded $\gm{T}$-algebra}, i.e., a pair
$(A,h)$ where
$A\colon \cat{M}\to\cat{C}$ is a functor and 
$h$ is a natural transformation of type
\[
\begin{tikzpicture}
      \node (TL)  at (0,0.75)  {$\cat{M}\times\cat{M}$};
      \node (TR)  at (3,0.75)  {$\cat{M}\times\cat{C}$};
      \node (BL)  at (0,-0.75) {$\cat{M}$};
      \node (BR)  at (3,-0.75) {$\cat{C}$};
      
      \draw [->] (TL) to node (T) [auto,labelsize]      {$\cat{M}\times A$} (TR);
      \draw [->] (TR) to node (R) [auto,labelsize]      {$\gm{T}$} (BR);
      \draw [->] (TL) to node (L) [auto,swap,labelsize] {$\tensor$}(BL);
      \draw [->] (BL) to node (B) [auto,swap,labelsize] {$A$}(BR);
      
      \draw [2cellnode] (R) to  [bend right=30] node [auto,swap,labelsize] {$h$}(B);
\end{tikzpicture}
\]
So the component of $h$ at $(m,n)\in\cat{M}\times\cat{M}$ is of type
\[
h_{m,n}\quad\colon \quad m\act A_n\quad \longrightarrow\quad A_{m\tensor n}.
\] 
These data are subject to the following axioms:
\begin{itemize}
\item $\begin{tikzpicture}[baseline=-\the\dimexpr\fontdimen22\textfont2\relax ]
      \node (TL)  at (0,0.75)  {$A_n$};
      \node (TR)  at (2,0.75)  {$\e\act A_n$};
      \node (BR)  at (2,-0.75) {$A_n$};
      
      \draw [->] (TL) to node (T) [auto,labelsize]      {$\eta_{A_n}$} (TR);
      \draw [->] (TR) to node (R) [auto,labelsize]      {$h_{\e, n}$} (BR);
      \draw [->] (TL) to node (L) [auto,swap,labelsize] {$\id{A_n}$}(BR);
\end{tikzpicture}$
commutes for each object~$n$ of $\cat{M}$.
\item 
$\begin{tikzpicture}[baseline=-\the\dimexpr\fontdimen22\textfont2\relax ]
      \node (TL)  at (0,0.75)  {$l\act m\act A_n$};
      \node (TR)  at (4,0.75)  {$(l\tensor m)\act A_n$};
      \node (BL)  at (0,-0.75) {$l\act A_{m\tensor n}$};
      \node (BR)  at (4,-0.75) {$A_{l\tensor m\tensor n}$};
      
      \draw [->] (TL) to node (T) [auto,labelsize]      {$\mu_{l,m,A_n}$} (TR);
      \draw [->] (TR) to node (R) [auto,labelsize]      {$h_{l\tensor m,n}$} (BR);
      \draw [->] (TL) to node (L) [auto,swap,labelsize] {$l\act h_{m,n}$}(BL);
      \draw [->] (BL) to node (B) [auto,swap,labelsize] {$h_{l,m\tensor n}$}(BR);
\end{tikzpicture}$
commutes for each triple of objects~$l$, $m$, $n$ of $\cat{M}$.
\end{itemize}
\item
A morphism of $\EMg$ from $(A,h)$ to $(A^\prime,h^\prime)$ is a 
\defemph{homomorphism} of graded $\gm{T}$-algebras between them,
i.e., a natural transformation
$\varphi\colon A\Rightarrow A^\prime$ making the diagram
\[
\begin{tikzpicture}
      \node (TL)  at (0,0.75)  {$m\act A_n$};
      \node (TR)  at (3,0.75)  {$m\act A^\prime_n$};
      \node (BL)  at (0,-0.75) {$A_{m\tensor n}$};
      \node (BR)  at (3,-0.75) {$A^\prime_{m\tensor n}$};
      
      \draw [->] (TL) to node (T) [auto,labelsize]      {$m\act \varphi_{n}$} (TR);
      \draw [->] (TR) to node (R) [auto,labelsize]      {$h^\prime_{m,n}$} (BR);
      \draw [->] (TL) to node (L) [auto,swap,labelsize] {$h_{m,n}$}(BL);
      \draw [->] (BL) to node (B) [auto,swap,labelsize] {$\varphi_{m\tensor n}$}(BR);
\end{tikzpicture}
\]
commute for each pair of objects~$m,n$ of $\cat{M}$.\qedhere
\end{itemize}
\end{defn}

Let us introduce a convenient notation for graded $\gm{T}$-algebras.
We write a graded $\gm{T}$-algebra~$(A,h)$ as 
$((A_n)_{n\in\cat{M}},\,(h_{m,n})_{m,n\in\cat{M}})$,
and use this notation to indicate definitions in what follows.
In principle we need to check the relevant functoriality or naturality
to validate such definitions, but these are all
completely routine and left to the interested reader.
Similarly, we denote a homomorphism~$\varphi\colon(A,h)\to(A^\prime,h^\prime)$
by $(\varphi_n)_{n\in\cat{M}}$.

The category~$\EMg$ becomes an object of the 
2-category~$\EM{\Cat}{\cat{M}\times(-)}$ by the following 
functor~$\emact$:

\begin{defn}
Define the functor
\[
\emact\quad\colon\quad\cat{M}\times\EMg\quad\longrightarrow\quad\EMg
\]
as follows:
\begin{itemize}
\item Given objects $p$ and $(A,h)$ of 
$\cat{M}$ and $\EMg$ respectively,
we define the graded $\gm{T}$-algebra~$p\emact (A,h)$ by
the precomposition of $(-)\tensor p\colon\cat{M}\to\cat{M}$:
\[
p\emact (A,h)\quad\coloneqq\quad\big(\,(A_{n\tensor p})_{n\in\cat{M}},\ 
(h_{m,n\tensor p})_{m,n\in\cat{M}}\,\big).
\]
\item Given morphisms $u\colon p\to p^\prime$ and $\varphi\colon(A,h)\to(A^\prime,h^\prime)$
of $\cat{M}$ and $\EMg$ respectively, 
we define the homomorphism~$u\emact \varphi\colon
p\emact(A,h)\to p^\prime\emact(A^\prime,h^\prime)$
by setting the component~$(u\emact \varphi)_{n}\colon
A_{n\tensor p}\to A^\prime_{n\tensor p^\prime}$ at $n\in\cat{M}$
to be either of the following two equivalent composites:
\[
\begin{tikzpicture}
      \node (TL)  at (0,0.75)  {$A_{n\tensor p}$};
      \node (TR)  at (3,0.75)  {$A^\prime_{n\tensor p}$};
      \node (BL)  at (0,-0.75) {$A_{n\tensor p^\prime}$};
      \node (BR)  at (3,-0.75) {$A^\prime_{n\tensor p^\prime}$};
      
      \draw [->] (TL) to node (T) [auto,labelsize]      {$\varphi_{n\tensor p}$} (TR);
      \draw [->] (TR) to node (R) [auto,labelsize]      {$A^\prime_{n\tensor u}$} (BR);
      \draw [->] (TL) to node (L) [auto,swap,labelsize] {$A_{n\tensor u}$}(BL);
      \draw [->] (BL) to node (B) [auto,swap,labelsize] {$\varphi_{n\tensor p^\prime}$}(BR);
\end{tikzpicture}
\qedhere
\]
\end{itemize}
\end{defn}

That this functor~$\emact\colon\cat{M}\times\EMg\to\EMg$ gives a 
$\EM{\Cat}{\cat{M}\times(-)}$-algebra structure, also 
known as a \emph{strict left action} of $(\cat{M},\tensor,\e)$, 
is easily verified;
a key step is the following:
\begin{align*}
p\emact (q\emact (A,h))\ &=\ 
p\emact ((A_{n\tensor q})_{n\in\cat{M}},\ (h_{m,n\tensor q})_{m,n\in\cat{M}})\\
&=\ ((A_{n\tensor p\tensor q})_{n\in\cat{M}},\ (h_{m,n\tensor p\tensor q})_{m,n\in\cat{M}})\\
&=\ (p\tensor q)\emact (A,h).
\end{align*}

\subsection{The Eilenberg--Moore adjunction}

\subsubsection{The left adjoint}
We define the 1-cell
\[
(F^\cat{M},f^\gm{T})\quad\colon\quad 
(\Cat,\cat{C})
\quad\longrightarrow\quad
\left(\EM{\Cat}{\cat{M}\times(-)},
\valg{\cat{M}\times\EMg}{\emact}{\EMg}
\right)
\] 
of $\Epp$ as follows:

\begin{defn}
The 2-functor~$F^\cat{M}\colon \Cat\to\EM{\Cat}{\cat{M\times(-)}}$
is the free 2-functor~$
\cat{X}\mapsto\valgp{\cat{M}\times\cat{M}\times\cat{X}}{M^\cat{M}_\cat{X}}{\cat{M}\times\cat{X}}$
where $M^\cat{M}_\cat{X}=M_\cat{X}$
is the one defined in Definition~\ref{defn-M-times}.
\end{defn}

\begin{defn}
The 1-cell~$
f^\gm{T}\colon 
\valgp{\cat{M}\times\cat{M}\times\cat{C}}{M_\cat{C}}{\cat{M}\times\cat{C}}
\to
\valgp{\cat{M}\times\EMg}{\emact}{\EMg}$
of $\EM{\Cat}{\cat{M}\times(-)}$
is the functor~$f^\gm{T}\colon \cat{M}\times\cat{C}\to\EMg$ defined as
\begin{multline*}
f^\gm{T}(p,c)\quad\coloneqq\quad
\big(\,((n\tensor p)\act c)_{n\in\cat{M}},\\
(\mu_{m,n\tensor p,c}\colon 
m\act((n\tensor p)\act c)\to(m\tensor n\tensor p)\act c
)_{m,n\in\cat{M}}\,\big)
\end{multline*}
on an object~$(p,c)$ and 
\[
f^\gm{T}(u,f)\quad\coloneqq\quad
\big(\,(n\tensor u)\act f\colon (n\tensor p)\act c\to(n\tensor p^\prime)\act c^\prime\,\big)_{n\in\cat{M}}
\]
on a morphism~$(u,f)\colon(p,c)\to(p^\prime,c^\prime)$.
\end{defn}

\subsubsection{The right adjoint}
We define the 1-cell
\[
(U^\cat{M},u^\gm{T})\quad\colon\quad 
\left(\EM{\Cat}{\cat{M}\times(-)},
\valg{\cat{M}\times\EMg}{\emact}{\EMg}
\right)
\quad\longrightarrow\quad
(\Cat,\cat{C})
\] 
of $\Epp$ as follows:
\begin{defn}
The 2-functor~$U^\cat{M}\colon \EM{\Cat}{\cat{M}\times(-)}\to\Cat$ is 
the forgetful 2-functor $
\valgp{\cat{M}\times\cat{A}}{\alpha}{\cat{A}}\mapsto\cat{A}$.
\end{defn}

\begin{defn}
The functor~$u^\gm{T}\colon \EMg\to\cat{C}$
is given by the evaluation at the monoidal unit $\e\in\cat{M}$:
$u^\gm{T}(A,h)\coloneqq A_\e$ for an object~$(A,h)$
and $u^\gm{T}(\varphi)\coloneqq \varphi_\e$ for a morphism~$\varphi$.
\end{defn}

\subsubsection{The unit}
We define the 2-cell
    \[
    \begin{tikzpicture}
      \node (A)  at (0,0) {$(\Cat,\cat{C})$};
      \node (B) at (8,0) {$(\Cat,\cat{C})$};
      \node (C) at (4,-2) {$\left(\EM{\Cat}{\cat{M}\times(-)},
      \valg{\cat{M}\times\EMg}{\emact}{\EMg}\right)$};
      
      \draw [->] (A) to [bend right=20]node [auto,swap,labelsize] {$(F^\cat{M},f^\gm{T})$}(C);
      \draw [->] (C) to [bend right=20]node [auto,swap,labelsize] {$(U^\cat{M},u^\gm{T})$}(B);

      \draw [->] (A) to [bend left=20]node [auto,labelsize] (Fpp) {$(\id{\Cat},\id{\cat{C}})$}(B);
      \draw [2cellnoder] (C) to node [auto,labelsize,swap]{$(H^\cat{M},\eta^\gm{T})$}(Fpp);

    \end{tikzpicture}
    \]
of $\Epp$ as follows:
\begin{defn}
The 2-natural transformation~$
H^\cat{M}\colon \id{\Cat}\Rightarrow U^\cat{M}\circ F^\cat{M}$
is the one defined in Definition~\ref{defn-M-times}.
\end{defn}

\begin{defn}
The natural transformation
    \[
    \begin{tikzpicture}
      \node (L) at (0,0) {$\cat{C}$};
      \node (R) at (6,0) {$\cat{C}$};
      \node (T) at (2,-2) {$\cat{M}\times\cat{C}$};
      \node (N) at (4.5,-1.35) {$\EMg$};
      
      \draw [->] (L) to node [auto,labelsize] {$\id{\cat{C}}$}(R);
      \draw [->] (L) to [bend right=20]node [auto,swap,labelsize] {$H_\cat{C}$}(T);
      \draw [->] (T) to [bend right=10]node [auto,swap,labelsize] {$f^\gm{T}$}(N);
      \draw [->] (N) to [bend right=10]node [auto,swap,labelsize] {$u^\gm{T}$}(R);
            
      \draw [2cellr] (T) to node [auto,labelsize,swap]{$\eta^\gm{T}$}(2,-0.1);
    \end{tikzpicture}
    \]
has components~$
\eta^\gm{T}_c\colon c\to\e\act c$ given by the data of the graded monad~$
\gm{T}$.
\end{defn}

\subsubsection{The counit}

We define the 2-cell
    \[
    \begin{tikzpicture}
      \node (A)  at (0,0) {$\left(\EM{\Cat}{\cat{M}\times(-)},
            \valg{\cat{M}\times\EMg}{\emact}{\EMg}\right)$};
      \node (B) at (8,0) {$\left(\EM{\Cat}{\cat{M}\times(-)},
            \valg{\cat{M}\times\EMg}{\emact}{\EMg}\right)$};
      \node (C) at (4,2) {$(\Cat,\cat{C})$};
      
      \draw [->] (A) to [bend left=20]node [auto,labelsize] {$(U^\cat{M},u^\gm{T})$}(C);
      \draw [->] (C) to [bend left=20]node [auto,labelsize] {$(F^\cat{M},f^\gm{T})$}(B);

      \draw [->] (A) to [bend right=20]node [auto,swap,labelsize] (Fpp) {$(\id{\EM{\Cat}{\cat{M}\times(-)}},\id{\emact})$}(B);
      \draw [2cellnoder] (Fpp) to node [auto,swap,labelsize]{$(E^\cat{M},\varepsilon^\gm{T})$}(C);

    \end{tikzpicture}
    \]
of $\Epp$ as follows:

\begin{defn}
The 2-natural transformation~$E^\cat{M}\colon F^\cat{M}\circ U^\cat{M}\Rightarrow \id{\EM{\Cat}{\cat{M}\times(-)}}$
has components~$
E^\cat{M}_{\alpha}\colon 
\valgp{\cat{M}\times\cat{M}\times\cat{A}}{M_\cat{A}}{\cat{M}\times\cat{A}}
\to
\valgp{\cat{M}\times\cat{A}}{\alpha}{\cat{A}}$
given by $E^\cat{M}_\alpha=E_\alpha\coloneqq\alpha\colon \cat{M}\times\cat{A}\to\cat{A}$.
\end{defn}

\begin{defn}
The 2-cell
    \[
    \begin{tikzpicture}
      \node (L) at (0,0) {$\valgp{\cat{M}\times\cat{M}\times\EMg}{M_\EMg}{\cat{M}\times\EMg}$};
      \node (R) at (9,0) {$\valgp{\cat{M}\times\EMg}{\emact}{\EMg}$};
      \node (T) at (3,-3) {$\valgp{\cat{M}\times\EMg}{\emact}{\EMg}$};
      
      \node (M) at (4.5,0) {$\valgp{\cat{M}\times\cat{M}\times\cat{C}}{M_\cat{C}}{\cat{M}\times\cat{C}}$};
      
      \draw [->] (L) to [bend right=0]node [auto,labelsize] {$\cat{M}\times u^\gm{T}$}(M);
      \draw [->] (M) to [bend right=0]node [auto,labelsize] {$f^\gm{T}$}(R);

      \draw [->] (L) to [bend right=20]node [auto,swap,labelsize] {$E_\emact$}(T);
      \draw [->] (T) to [bend right=20]node [auto,labelsize,swap] {$\id{\emact}$}(R);
            
      \draw [2cell] (3,-0.1) to node [auto,labelsize]{$\varepsilon^\gm{T}$}(T);
    \end{tikzpicture}
    \]
of $\EM{\Cat}{\cat{M}\times(-)}$ is the natural transformation
with its component at $(p,(A,h))\in\cat{M}\times\EMg$ of type
\begin{multline*}
\varepsilon^\gm{T}_{p,(A,h)}\quad\colon\quad 
\big(\,((n\tensor p)\act A_\e)_{n\in\cat{M}},\ 
(\mu_{m,n\tensor p,A_\e})_{m,n\in\cat{M}}\,\big)\\
\longrightarrow\quad 
\big(\,(A_{n\tensor p})_{n\in\cat{M}},\ 
(h_{m,n\tensor p})_{m,n\in\cat{M}}\,\big),
\end{multline*}
itself with the component at $n\in\cat{M}$
given by
$\varepsilon^\gm{T}_{p,(A,h),n}\coloneqq h_{n\tensor p,\e}
\colon 
(n\tensor p)\act A_\e\to (A_{n\tensor p})$.
\end{defn}
Observe that $\varepsilon^\gm{T}_{p,(A,h)}$ is
indeed a homomorphism of
graded $\gm{T}$-algebras, i.e., the diagram
\[
\begin{tikzpicture}
      \node (TL)  at (0,0.75)  {$m\act ((n\tensor p)\act A_\e)$};
      \node (TR)  at (4,0.75)  {$m\act A_{n\tensor p}$};
      \node (BL)  at (0,-0.75) {$(m\tensor n\tensor p)\act A_\e$};
      \node (BR)  at (4,-0.75) {$A_{m\tensor n\tensor p}$};
      
      \draw [->] (TL) to node (T) [auto,labelsize]      {$m\act h_{n\tensor p,\e}$} (TR);
      \draw [->] (TR) to node (R) [auto,labelsize]      {$h_{m,n\tensor p}$} (BR);
      \draw [->] (TL) to node (L) [auto,swap,labelsize] {$\mu_{m,n\tensor p,A_\e}$}(BL);
      \draw [->] (BL) to node (B) [auto,swap,labelsize] {$h_{m\tensor n\tensor p,\e}$}(BR);
\end{tikzpicture}
\]
commutes, thanks to one of the axioms of graded $\gm{T}$-algebras.

\subsection{Comparison maps}
Suppose we have an adjunction~$
(L,l)\dashv(R,r)\colon(\tcat{D},D)\to(\Cat,\cat{C})$ in $\Epp$
with unit~$(H,\eta)\colon$ $(\id{\Cat},\id{\cat{C}})\Rightarrow(R,r)\circ(L,l)$ 
and counit~$(E,\varepsilon)\colon(L,l)\circ(R,r)\Rightarrow(\id{\tcat{D}},\id{D})$,
which gives a resolution of the monad $(\cat{M}\times(-),\gm{T})$,
i.e., such that the following equations hold:
{\allowdisplaybreaks
\begin{align}
\label{eqn-comp-gem-1}
\Cat\xrightarrow{\cat{M}\times(-)}\Cat\quad
&=
\quad\Cat\xrightarrow{L}\tcat{D}\xrightarrow{R}\Cat\\
\label{eqn-comp-gem-2}
\cat{M}\times\cat{C}\xrightarrow{\gm{T}}\cat{C}\quad&=\quad
RL\cat{C}\xrightarrow{Rl}RD\xrightarrow{r}\cat{C}\\
\label{eqn-comp-gem-3}
H^\cat{M}\quad&=\quad H\\
\label{eqn-comp-gem-4}
\eta^\gm{T}\quad&=\quad\eta\\
\label{eqn-comp-gem-5}
\begin{tikzpicture}[baseline=-\the\dimexpr\fontdimen22\textfont2\relax ]
        \begin{scope} 
      \node (TL)  at (0,0)  {$\Cat$};
      \node (TT)  at (1.2,0.7) {$\Cat$};
      \node (RR)  at (2.4,0) {$\Cat$};
      \draw [->] (TL) to [bend left=20]node (L) [auto,labelsize] {$\cat{M}\times(-)$}(TT);
      \draw [->] (TT) to [bend left=20]node  [auto,labelsize] {$\cat{M}\times(-)$}(RR);
      \draw [->] (TL) to [bend right=20]node (B) [auto,swap,labelsize] {$\cat{M}\times(-)$} (RR);
      \draw [2cellnode] (TT) to node[auto,labelsize] {$M^\cat{M}$}(B);
        \end{scope}
\end{tikzpicture}
\quad&=\quad
\begin{tikzpicture}[baseline=-\the\dimexpr\fontdimen22\textfont2\relax ]
        \begin{scope} 
      \node (TL)  at (0,0)  {$\Cat$};
      \node (TR)  at (1.2,0)  {$\tcat{D}$};
      \node (TT)  at (2.4,0.7) {$\Cat$};
      \node (BR)  at (3.6,0) {$\tcat{D}$};
      \node (RR)  at (4.8,0) {$\Cat$};
      \draw [->] (TL) to node (T) [auto,labelsize]      {$L$} (TR);
      \draw [->] (TR) to [bend left=20]node (R) [auto,labelsize]      {$R$} (TT);
      \draw [->] (TT) to [bend left=20]node (L) [auto,labelsize] {$L$}(BR);
      \draw [->] (BR) to node  [auto,labelsize] {$R$}(RR);
      \draw [->] (TR) to [bend right=20]node (B) [auto,swap,labelsize] {$\id{\tcat{D}}$} (BR);
      \draw [2cellnode] (TT) to node[auto,labelsize] {$E$}(B);
        \end{scope}
\end{tikzpicture}
\end{align}
}%
\vspace{-20pt}
\begin{multline}
\label{eqn-comp-gem-6}
\begin{tikzpicture}[baseline=-\the\dimexpr\fontdimen22\textfont2\relax ]
        \begin{scope} 
      \node (L)  at (0,0)  {$\cat{M}\times\cat{M}\times\cat{C}$};
      \node (M)  at (3,0)  {$\cat{M}\times\cat{C}$};
      \node (R)  at (6,0) {$\cat{C}$};
      \node (B)  at (2,-2) {$\cat{M}\times\cat{C}$};
      \draw [->] (L) to node  [auto,labelsize]      {$\cat{M}\times\gm{T}$} (M);
      \draw [->] (M) to node  [auto,labelsize]      {$\gm{T}$} (R);
      \draw [->] (L) to [bend right=20]node  [auto,swap,labelsize] {$M_\cat{C}$}(B);
      \draw [->] (B) to [bend right=20]node  [auto,swap,labelsize] {$\gm{T}$}(R);
      \draw [2cellnode] (2,0) to node[auto,labelsize] {$\mu^\gm{T}$}(B);
        \end{scope}
\end{tikzpicture}\\
=\quad
\begin{tikzpicture}[baseline=-\the\dimexpr\fontdimen22\textfont2\relax ]
        \begin{scope} [shift={(0,0.75)}]
      \node (L)  at (0,0)  {$RLRL\cat{C}$};
      \node (M)  at (2.75,0)  {$RLRD$};
      \node (R)  at (5.2,0) {$RL\cat{C}$};
      \node (RR)  at (7.15,0) {$RD$};
      \node (RRR)  at (8.6,0) {$\cat{C}$};
      \node (B)  at (4.2,-1.5) {$RD$};
      \draw [->] (L) to node  [auto,labelsize]      {$RLRl$} (M);
      \draw [->] (M) to node  [auto,labelsize]      {$RLr$} (R);
      \draw [->] (R) to node  [auto,labelsize]      {$Rl$} (RR);
      \draw [->] (RR) to node  [auto,labelsize]      {$r$} (RRR);
      \draw [->] (M) to [bend right=20]node  [auto,swap,labelsize] {$RE_D$}(B);
      \draw [->] (B) to [bend right=20]node  [auto,swap,labelsize] {$\id{RD}$}(RR);
      \draw [2cellnode] (4.2,0) to node[auto,labelsize] {$R\varepsilon$}(B);
        \end{scope}
\end{tikzpicture}
\end{multline}

First note that 
equations~(\ref{eqn-comp-gem-1}), (\ref{eqn-comp-gem-3}) and
(\ref{eqn-comp-gem-5}) imply the existence of the comparison
2-functor~$K$
by a classical result of 2-monad theory (or rather enriched monad theory).

\begin{defn}
The 2-functor~$K\colon \tcat{D}\to\EM{\Cat}{\cat{M}\times(-)}$ 
is the comparison 2-functor 
$D^\prime\mapsto\valgp{\cat{M}\times RD^\prime}{RE_{D^\prime}}{RD^\prime}$.
\end{defn}

This 2-functor~$K$ becomes the 2-functor part of the 
comparison map~$(K,k)$
under construction.
So now it remains to construct an appropriate 1-cell~$
k$
of $\EM{\Cat}{\cat{M}\times(-)}$.

\begin{defn}
The 1-cell~$k\colon\valgp{\cat{M}\times RD}{RE_{D}}{RD}
\to\valgp{\cat{M}\times\EMg}{\emact}{\EMg}$
of $\EM{\Cat}{\cat{M}\times(-)}$
is the functor~$
k\colon RD\to\EMg$ defined as
\begin{equation*}
k(d)\quad\coloneqq\quad 
\big(\,(rRE_D(n,d))_{n\in\cat{M}},\ 
(rR\varepsilon_{m,RE_{D}(n,d)})_{m,n\in\cat{M}}\,\big)
\end{equation*}
on an object~$d$.
The types of the structure maps indeed match:
\begin{alignat*}{3}
&m\ast rRE_D(n,d)&
&\ccol{\ =\ }&
&rRl(m,rRE_D(n,d))\tag*{by (\ref{eqn-comp-gem-2})}\\
&&
&\ccol{\ =\ }&
&rRlRLr(m,RE_D(n,d))\\
&&
&\ccol{\ \xrightarrow{rR\varepsilon_{m,RE_{D}(n,d)}}\ }&
&rRE_D(m,RE_D(n,d))\\
&&
&\ccol{\ =\ }&
&rRE_DRLRE_{D}(m, n,d)\\
&&
&\ccol{\ =\ }&
&rRE_DRE_{LRD}(m, n,d)\\
&&
&\ccol{\ =\ }&
&rRE_D(m\tensor n,d)\tag*{by (\ref{eqn-comp-gem-5}).}
\end{alignat*}
$k$ is defined as
\[
k(w)\quad\coloneqq \quad 
\big(\,
rRE_D(n,w)\colon rRE_D(n,d)\to rRE_D(n,d^\prime)
\,\big)_{n\in\cat{M}}
\]
on a morphism~$w\colon d\to d^\prime$.
\end{defn}

\begin{prop}
\label{prop-comparison-em-gm}
The 1-cell~$(K,k)$ of $\Epp$ satisfies the equations~$
(K,k)\circ(L,l)=(F^\cat{M},f^\gm{T})$ and
$(R,r)=(U^\cat{M},u^\gm{T})\circ(K,k)$.
Moreover, it is the unique such.
\end{prop}

We omit a proof since it follows from the 2-dimensional universality
discussed below.

\subsection{The 2-dimensional universality}
\subsubsection{Statement of the theorem}
We will show that there is a family of
isomorphisms of categories
\begin{multline*}
  \Epp\left((\tcat{X},X),\ 
  \left(\EM{\Cat}{\cat{M}\times(-)},
  \valg{\cat{M}\times\EMg}{\emact}{\EMg}\right)\right)\\
  \cong\quad 
  \Epp\big((\tcat{X},X),\ (\Cat,\cat{C})\big)^{\Epp((\tcat{X},X),\ (\cat{M}\times(-),\gm{T}))}
\end{multline*}
2-natural in $(\tcat{X},X)\in \Epp$; cf.~Definition~\ref{defn-EM-obj}.
More precisely, we claim that the data 
    \[
    \begin{tikzpicture}
      \node (A)  at (0,0) {$\left(\EM{\Cat}{\cat{M}\times(-)},
            \valg{\cat{M}\times\EMg}{\emact}{\EMg}\right)$};
      \node (B) at (4,-2) {$\left(\EM{\Cat}{\cat{M}\times(-)},
            \valg{\cat{M}\times\EMg}{\emact}{\EMg}\right)$};
      \node (C) at (4,0) {$(\Cat,\cat{C})$};
      \node (D) at (4,-4) {$(\Cat,\cat{C})$};
      
      \draw [->] (A) to [bend left=0]node [auto,labelsize] {$(U^\cat{M},u^\gm{T})$}(C);
      \draw [->] (C) to [bend left=0]node [auto,labelsize] {$(F^\cat{M},f^\gm{T})$}(B);
      \draw [->] (B) to [bend left=0]node [auto,labelsize] {$(U^\cat{M},u^\gm{T})$}(D);

      \draw [->] (A) to [bend right=20]node [auto,swap,labelsize] (Fpp) {$(\id{\EM{\Cat}{\cat{M}\times(-)}},\id{\emact})$}(B);
      \draw [2cell] (2,-0.1) to node [auto,labelsize]{$(E^\cat{M},\varepsilon^\gm{T})$}(2,-1.6);

    \end{tikzpicture}
    \]
provides the universal left $\gm{T}$-module, in the sense that 
every left $\gm{T}$-module
    \[
    \begin{tikzpicture}
      \node (A)  at (0,0) {$(\tcat{X},X)$};
      \node (B) at (3,-1.5) {$(\Cat,\cat{C})$};
      \node (C) at (3,0) {$(\Cat,\cat{C})$};
      
      \draw [->] (A) to [bend left=0]node [auto,labelsize] {$(G,g)$}(C);
      \draw [->] (C) to [bend left=0]node [auto,labelsize] {$(\cat{M}\times(-),\gm{T})$}(B);

      \draw [->] (A) to [bend right=20]node [auto,swap,labelsize] (Fpp) {$(G,g)$}(B);
      \draw [2cell] (1.5,-0.1) to node [auto,labelsize]{$(\Gamma,\gamma)$}(1.5,-1.1);

    \end{tikzpicture}
    \]
i.e., an object of the category~$
\Epp((\tcat{X},X),\ (\Cat,\cat{C}))^{\Epp((\tcat{X},X),\ (\cat{M}\times(-),\gm{T}))}$,
factors uniquely as 
    \[
    \begin{tikzpicture}
      \node (E) at (-4,0) {$(\tcat{X},X)$};
      \node (A)  at (0,0) {$\left(\EM{\Cat}{\cat{M}\times(-)},
            \valg{\cat{M}\times\EMg}{\emact}{\EMg}\right)$};
      \node (B) at (4,-2) {$\left(\EM{\Cat}{\cat{M}\times(-)},
            \valg{\cat{M}\times\EMg}{\emact}{\EMg}\right)$};
      \node (C) at (4,0) {$(\Cat,\cat{C})$};
      \node (D) at (4,-4) {$(\Cat,\cat{C})$};

      \draw [->] (E) to [bend left=0]node [auto,labelsize] {$(\tilde{G},\tilde{g})$}(A);
      \draw [->] (A) to [bend left=0]node [auto,labelsize] {$(U^\cat{M},u^\gm{T})$}(C);
      \draw [->] (C) to [bend left=0]node [auto,labelsize] {$(F^\cat{M},f^\gm{T})$}(B);
      \draw [->] (B) to [bend left=0]node [auto,labelsize] {$(U^\cat{M},u^\gm{T})$}(D);

      \draw [->] (A) to [bend right=20]node [auto,swap,labelsize] (Fpp) {$(\id{\EM{\Cat}{\cat{M}\times(-)}},\id{\emact})$}(B);
      \draw [2cell] (2,-0.1) to node [auto,labelsize]{$(E^\cat{M},\varepsilon^\gm{T})$}(2,-1.6);

    \end{tikzpicture}
    \]
and similarly every morphism of left $\gm{T}$-modules
    \[
    \begin{tikzpicture}
      \node (A)  at (0,0) {$(\tcat{X},X)$};
      \node (C) at (4,0) {$(\Cat,\cat{C})$};
      
      \draw [->] (A) to [bend left=20]node (dom)[auto,labelsize] {$(G,g)$}(C);
      \draw [->] (A) to [bend right=20]node (cod)[auto,swap,labelsize] {$(G^\prime,g^\prime)$}(C);

      \draw [2cellnode] (dom) to node [auto,labelsize]{$(\Omega,\omega)$}(cod);

    \end{tikzpicture}
    \]
i.e., a morphism of the category~$
\Epp((\tcat{X},X),\ (\Cat,\cat{C}))^{\Epp((\tcat{X},X),\ (\cat{M}\times(-),\gm{T}))}$,
factors uniquely as 
    \[
    \begin{tikzpicture}
      \node (E) at (-5,0) {$(\tcat{X},X)$};
      \node (A)  at (0,0) {$\left(\EM{\Cat}{\cat{M}\times(-)},
            \valg{\cat{M}\times\EMg}{\emact}{\EMg}\right)$};
      \node (C) at (4,0) {$(\Cat,\cat{C})$};

      \draw [->] ([yshift=3mm]E.east) to [bend left=20]node(dom) [auto,labelsize] {$(\tilde{G},\tilde{g})$} ([yshift=3mm]A.west);
      \draw [->] ([yshift=-3mm]E.east) to [bend right=20]node(cod) [auto,swap,labelsize] {$(\tilde{G^\prime},\tilde{g^\prime})$}([yshift=-3mm]A.west);
      \draw [->] (A) to [bend left=0]node [auto,labelsize] {$(U^\cat{M},u^\gm{T})$}(C);

      \draw [2cellnode] (dom) to node [auto,labelsize]{$(\tilde{\Omega},\tilde{\omega})$}(cod);

    \end{tikzpicture}
    \]
    
\subsubsection{The 1-dimensional aspect}
Let us first verify the unique factorization for left $\gm{T}$-modules.
Suppose we have a left $\gm{T}$-modules, i.e., a piece of data
    \[
    \begin{tikzpicture}
      \node (A)  at (0,0) {$(\tcat{X},X)$};
      \node (B) at (3,-1.5) {$(\Cat,\cat{C})$};
      \node (C) at (3,0) {$(\Cat,\cat{C})$};
      
      \draw [->] (A) to [bend left=0]node [auto,labelsize] {$(G,g)$}(C);
      \draw [->] (C) to [bend left=0]node [auto,labelsize] {$(\cat{M}\times(-),\gm{T})$}(B);

      \draw [->] (A) to [bend right=20]node [auto,swap,labelsize] (Fpp) {$(G,g)$}(B);
      \draw [2cell] (1.5,-0.1) to node [auto,labelsize]{$(\Gamma,\gamma)$}(1.5,-1.1);

    \end{tikzpicture}
    \]
in $\Epp$ satisfying
%
%
%
{\allowdisplaybreaks
\begin{align}
\label{eqn-univ-gem-1}
\begin{tikzpicture}[baseline=-\the\dimexpr\fontdimen22\textfont2\relax ]
        \begin{scope} [shift={(0,0.75)}]
      \node (A)  at (0,0) {$\tcat{X}$};
      \node (B) at (3,-1.5) {$\Cat$};
      \node (C) at (3,0) {$\Cat$};
      \draw [->] (A) to [bend left=0]node [auto,labelsize] {$G$}(C);
      \draw [->] (C) to [bend left=0]node (t)[auto,swap,labelsize] {$\cat{M}\times(-)$}(B);
      \draw [->] (C) to [bend left=90]node(id) [auto,labelsize] {$\id{\Cat}$}(B);
      %
      \draw [->] (A) to [bend right=20]node [auto,swap,labelsize] (Fpp) {$G$}(B);
      \draw [2cell] (1.5,-0.1) to node [auto,swap,labelsize]{$\Gamma$}(1.5,-1.1);
      \draw [2cellnode] (id) to node [auto,swap,labelsize]{$H^\cat{M}$}(t);
        \end{scope}
\end{tikzpicture}
\quad&=\quad
\begin{tikzpicture}[baseline=-\the\dimexpr\fontdimen22\textfont2\relax ]
        \begin{scope} 
       \node (A)  at (0,0) {$\tcat{X}$};
       \node (C) at (2,0) {$\Cat$};
       \draw [->] (A) to [bend left=0]node [auto,labelsize] {$G$}(C);
        \end{scope}
\end{tikzpicture}\\
\label{eqn-univ-gem-2}
\begin{tikzpicture}[baseline=-\the\dimexpr\fontdimen22\textfont2\relax ]
        \begin{scope} 
      \node (A)  at (0,0) {$\cat{M}\times GX$};
      \node (B) at (3,0) {$\cat{M}\times\cat{C}$};
      \node (C) at (5,0) {$\cat{C}$};
      \node (D)  at (1.5,-1.5) {$GX$};
      \node (E) at (3.5,1.5) {$\cat{C}$};
      \node (F) at (0.5,1.5) {$GX$};
      \draw [->] (A) to [bend left=0]node [auto,labelsize] {$\cat{M}\times g$}(B);
      \draw [->] (B) to [bend left=0]node [auto,labelsize] {$\gm{T}$}(C);
      \draw [->] (A) to [bend right=20]node [auto,swap,labelsize] {$\Gamma_X$}(D);
      \draw [->] (D) to [bend right=20]node [auto,swap,labelsize] {$g$}(C);
      \draw [->] (E) to [bend right=20]node [auto,swap,labelsize] {$H_\cat{C}$}(B);
      \draw [->] (F) to [bend right=0]node [auto,labelsize] {$g$}(E);
      \draw [->] (E) to [bend left=20]node [auto,labelsize] {$\id{\cat{C}}$}(C);
      \draw [->] (F) to [bend right=20]node [auto,swap,labelsize] {$H_{GX}$}(A);
      \draw [2cell] (1.5,-0.1) to node [auto,labelsize]{$\gamma$}(1.5,-1.1);
      \draw [2cell] (3.5,1.2) to node [auto,labelsize]{$\eta$}(3.5,0.3);
        \end{scope}
\end{tikzpicture}
\quad&=\quad
\begin{tikzpicture}[baseline=-\the\dimexpr\fontdimen22\textfont2\relax ]
        \begin{scope} 
       \node (A)  at (0,0) {$G{X}$};
       \node (C) at (2,0) {$\cat{C}$};
       \draw [->] (A) to [bend left=0]node [auto,labelsize] {$g$}(C);
        \end{scope}
\end{tikzpicture}\\
\label{eqn-univ-gem-3}
\begin{tikzpicture}[baseline=-\the\dimexpr\fontdimen22\textfont2\relax ]
        \begin{scope} 
      \node (A)  at (0,0) {$\tcat{X}$};
      \node (B) at (3,-1.5) {$\Cat$};
      \node (C) at (4.5,0) {$\Cat$};
      \node (D) at (3,1.5) {$\Cat$};
      \draw [->] (D) to [bend left=0]node [auto,swap,labelsize] {$\cat{M}\times(-)$}(B);
      \draw [->] (D) to [bend left=20]node (t)[auto,labelsize] {$\cat{M}\times(-)$}(C);
      \draw [->] (C) to [bend left=20]node (t)[auto,labelsize] {$\cat{M}\times(-)$}(B);
      \draw [->] (A) to [bend left=20]node [auto,labelsize] (Fpp) {$G$}(D);
      \draw [->] (A) to [bend right=20]node [auto,swap,labelsize] (Fpp) {$G$}(B);
      \draw [2cell] (1.5,1) to node [auto,swap,labelsize]{$\Gamma$}(1.5,-1);
      \draw [2cell] (4.1,0) to node [auto,swap,labelsize]{$M^\cat{M}$}(3.1,0);
        \end{scope}
\end{tikzpicture}
\quad&=\quad 
\begin{tikzpicture}[baseline=-\the\dimexpr\fontdimen22\textfont2\relax ]
        \begin{scope} 
      \node (A)  at (0,0) {$\tcat{X}$};
      \node (B) at (3,-1.5) {$\Cat$};
      \node (C) at (4,0) {$\Cat$};
      \node (D) at (3,1.5) {$\Cat$};
      \draw [->] (A) to [bend left=0]node [auto,labelsize] {$G$}(C);
      \draw [->] (D) to [bend left=20]node (t)[auto,labelsize] {$\cat{M}\times(-)$}(C);
      \draw [->] (C) to [bend left=20]node (t)[auto,labelsize] {$\cat{M}\times(-)$}(B);
      \draw [->] (A) to [bend left=20]node [auto,labelsize] (Fpp) {$G$}(D);
      \draw [->] (A) to [bend right=20]node [auto,swap,labelsize] (Fpp) {$G$}(B);
      \draw [2cell] (3,1.2) to node [auto,labelsize]{$\Gamma$}(3,0.1);
      \draw [2cell] (3,-0.1) to node [auto,labelsize]{$\Gamma$}(3,-1.2);
        \end{scope}
\end{tikzpicture}
\end{align}
}%
\vspace{-20pt}
\begin{multline}
\label{eqn-univ-gem-4}
\begin{tikzpicture}[baseline=-\the\dimexpr\fontdimen22\textfont2\relax ]
        \begin{scope} 
      \node (A)  at (0,0) {$\cat{M}\times GX$};
      \node (B) at (4,0) {$\cat{M}\times\cat{C}$};
      \node (C) at (7,0) {$\cat{C}$};
      \node (D)  at (2,-1.5) {$GX$};
      \node (E) at (4,1.5) {$\cat{M}\times\cat{M}\times\cat{C}$};
      \node (F) at (0,1.5) {$\cat{M}\times\cat{M}\times GX$};
      \node (G) at (7,1.5) {$\cat{M}\times\cat{C}$};
      \draw [->] (A) to [bend left=0]node [auto,labelsize] {$\cat{M}\times g$}(B);
      \draw [->] (B) to [bend left=0]node [auto,swap,labelsize] {$\gm{T}$}(C);
      \draw [->] (A) to [bend right=20]node [auto,swap,labelsize] {$\Gamma_X$}(D);
      \draw [->] (D) to [bend right=15]node [auto,swap,labelsize] {$g$}(C);
      \draw [->] (E) to [bend right=0]node [auto,swap,labelsize] {$M_\cat{C}$}(B);
      \draw [->] (G) to [bend right=0]node [auto,labelsize] {$\gm{T}$}(C);
      \draw [->] (E) to [bend right=0]node [auto,labelsize] {$\cat{M}\times \gm{T}$}(G);
      \draw [->] (F) to [bend right=0]node [auto,labelsize] {$\cat{M}\times\cat{M}\times g$}(E);
      \draw [->] (F) to [bend right=0]node [auto,swap,labelsize] {$M_{GX}$}(A);
      \draw [2cell] (2,-0.1) to node [auto,labelsize]{$\gamma$}(2,-1.1);
      \draw [2cell] (5.5,1.4) to node [auto,labelsize]{$\mu$}(5.5,0.1);
        \end{scope}
\end{tikzpicture}\\
=\quad
\begin{tikzpicture}[baseline=-\the\dimexpr\fontdimen22\textfont2\relax ]
        \begin{scope} 
      \node (A)  at (0,0) {$\cat{M}\times GX$};
      \node (B) at (4,0) {$\cat{M}\times\cat{C}$};
      \node (C) at (7,0) {$\cat{C}$};
      \node (D)  at (2,-1.5) {$GX$};
      \node (E) at (4,1.5) {$\cat{M}\times\cat{M}\times\cat{C}$};
      \node (F) at (0,1.5) {$\cat{M}\times\cat{M}\times GX$};
      \draw [->] (A) to [bend left=0]node [auto,labelsize] {$\cat{M}\times g$}(B);
      \draw [->] (B) to [bend left=0]node [auto,labelsize] {$\gm{T}$}(C);
      \draw [->] (A) to [bend right=20]node [auto,swap,labelsize] {$\Gamma_X$}(D);
      \draw [->] (D) to [bend right=15]node [auto,swap,labelsize] {$g$}(C);
      \draw [->] (E) to [bend right=0]node [auto,labelsize] {$\cat{M}\times\gm{T}$}(B);
      \draw [->] (F) to [bend right=0]node [auto,labelsize] {$\cat{M}\times\cat{M}\times g$}(E);
      \draw [->] (F) to [bend right=0]node [auto,swap,labelsize] {$\cat{M}\times\Gamma_{X}$}(A);
      \draw [2cell] (2,-0.1) to node [auto,labelsize]{$\gamma$}(2,-1.1);
      \draw [2cell] (2,1.4) to node [auto,labelsize]{$\cat{M}\times\gamma$}(2,0.5);
        \end{scope}
\end{tikzpicture}
\end{multline}

First note that equations~(\ref{eqn-univ-gem-1}) and (\ref{eqn-univ-gem-3})
imply the unique factorization
    \[
    \begin{tikzpicture}[baseline=-\the\dimexpr\fontdimen22\textfont2\relax ]
            \begin{scope} [shift={(0,0.75)}]
      \node (A)  at (0,0) {$\tcat{X}$};
      \node (B) at (3,-1.5) {$\Cat$};
      \node (C) at (3,0) {$\Cat$};
          
      \draw [->] (A) to [bend left=0]node [auto,labelsize] {$G$}(C);
      \draw [->] (C) to [bend left=0]node [auto,labelsize] {$\cat{M}\times(-)$}(B);

      \draw [->] (A) to [bend right=20]node [auto,swap,labelsize] (Fpp) {$G$}(B);
      \draw [2cell] (1.5,-0.1) to node [auto,labelsize]{$\Gamma$}(1.5,-1.1);
            \end{scope}
    \end{tikzpicture}
    \quad=\quad
    \begin{tikzpicture}[baseline=-\the\dimexpr\fontdimen22\textfont2\relax ]
            \begin{scope} [shift={(0,1.5)}]
      \node (E) at (-3,0) {$\tcat{X}$};
      \node (A)  at (0,0) {$\EM{\Cat}{\cat{M}\times(-)}$};
      \node (B) at (3,-1.5) {$\EM{\Cat}{\cat{M}\times(-)}$};
      \node (C) at (3,0)  {$\Cat$};
      \node (D) at (3,-3) {$\Cat$};

      \draw [->] (E) to [bend left=0]node [auto,labelsize] {$\tilde{G}$}(A);
      \draw [->] (A) to [bend left=0]node [auto,labelsize] {$U^\cat{M}$}(C);
      \draw [->] (C) to [bend left=0]node [auto,labelsize] {$F^\cat{M}$}(B);
      \draw [->] (B) to [bend left=0]node [auto,labelsize] {$U^\cat{M}$}(D);

      \draw [->] (A) to [bend right=20]node [auto,swap,labelsize] (Fpp) {$\id{\EM{\Cat}{\cat{M}\times(-)}}$}(B);
      \draw [2cell] (1.5,-0.1) to node [auto,labelsize]{$E^\cat{M}$}(1.5,-1.1);
            \end{scope}
    \end{tikzpicture}
    \]
since the Eilenberg--Moore 2-category~$\EM{\Cat}{\cat{M}\times(-)}$
is the Eilenberg--Moore object in $\TCat$.
Concretely, the 2-functor~$\tilde{G}$ is defined as follows:
\begin{defn}
The 2-functor~$\tilde{G}\colon\tcat{X}\to\EM{\Cat}{\cat{M}\times(-)}$
is the mediating 2-functor 
$X^\prime\mapsto\valgp{\cat{M}\times GX^\prime}{\Gamma_{X^\prime}}{GX^\prime}$.
\end{defn}
This 2-functor is the only possible choice for the 2-functor
part of the desired factorization;
thus it remains to construct a 1-cell~$\tilde{g}\colon
\valgp{\cat{M}\times GX}{\Gamma_{X}}{GX}\to\valgp{\cat{M}\times\EMg}{\emact}{\EMg}$
of $\EM{\Cat}{\cat{M}\times(-)}$ which satisfies
\begin{align}
\label{eqn-univ-gem-5}
\begin{tikzpicture}[baseline=-\the\dimexpr\fontdimen22\textfont2\relax ]
        \begin{scope} 
      \node (A)  at (0,0) {$GX$};
      \node (B) at (2,0) {$\cat{C}$};
      \draw [->] (A) to [bend left=0]node [auto,labelsize] {${g}$}(B);
        \end{scope}
\end{tikzpicture}\quad
=\quad\begin{tikzpicture}[baseline=-\the\dimexpr\fontdimen22\textfont2\relax ]
        \begin{scope} 
      \node (A)  at (0,0) {$GX$};
      \node (B) at (2,0) {$\EMg$};
      \node (C) at (4,0) {$\cat{C}$};
      \draw [->] (A) to [bend left=0]node [auto,labelsize] {$\tilde{g}$}(B);
      \draw [->] (B) to [bend left=0]node [auto,labelsize] {$u^\gm{T}$}(C);
        \end{scope}
\end{tikzpicture}
\end{align}
\begin{multline}
\label{eqn-univ-gem-6}
\begin{tikzpicture}[baseline=-\the\dimexpr\fontdimen22\textfont2\relax ]
        \begin{scope} [shift={(0,0.75)}]
      \node (A)  at (0,0) {$\cat{M}\times GX$};
      \node (B) at (3,0) {$\cat{M}\times\cat{C}$};
      \node (C) at (5,0) {$\cat{C}$};
      \node (D)  at (1.5,-1.5) {$GX$};
      \draw [->] (A) to [bend left=0]node [auto,labelsize] {$\cat{M}\times g$}(B);
      \draw [->] (B) to [bend left=0]node [auto,labelsize] {$\gm{T}$}(C);
      \draw [->] (A) to [bend right=20]node [auto,swap,labelsize] {$\Gamma_X$}(D);
      \draw [->] (D) to [bend right=20]node [auto,swap,labelsize] {$g$}(C);
      \draw [2cell] (1.5,-0.1) to node [auto,labelsize]{$\gamma$}(1.5,-1.1);
        \end{scope}
\end{tikzpicture}\\
=\quad
\begin{tikzpicture}[baseline=-\the\dimexpr\fontdimen22\textfont2\relax ]
        \begin{scope} [shift={(0,0.75)}]  
      \node (E) at (-3,0) {$\cat{M}\times GX$};
      \node (F) at (6.5,0) {$\cat{C}$};
      \node (A)  at (0,0) {$\cat{M}\times \EMg$};
      \node (B) at (3,0) {$\cat{M}\times\cat{C}$};
      \node (C) at (5,0) {$\EMg$};
      \node (D)  at (1.5,-1.5) {$\EMg$};
      \draw [->] (A) to [bend left=0]node [auto,labelsize] {$\cat{M}\times u^\gm{T}$}(B);
      \draw [->] (E) to [bend left=0]node [auto,labelsize] {$\cat{M}\times\tilde{g}$}(A);
      \draw [->] (B) to [bend left=0]node [auto,labelsize] {$f^\gm{T}$}(C);
      \draw [->] (C) to [bend left=0]node [auto,labelsize] {$u^\gm{T}$}(F);
      \draw [->] (A) to [bend right=20]node [auto,swap,labelsize] {$E_\emact$}(D);
      \draw [->] (D) to [bend right=20]node [auto,swap,labelsize] {$\id{\EMg}$}(C);
      \draw [2cell] (1.5,-0.1) to node [auto,labelsize]{$\varepsilon^\gm{T}$}(1.5,-1.1);
        \end{scope}
\end{tikzpicture}
\end{multline}
and show its uniqueness.

\begin{defn}
The 1-cell~$\tilde{g}\colon\valgp{\cat{M}\times GX}{\Gamma_{X}}{GX}
\to\valgp{\cat{M}\times\EMg}{\emact}{\EMg}$
of $\EM{\Cat}{\cat{M}\times(-)}$ is a functor of type~$
\tilde{g}\colon GX\to\EMg$.
As a functor, it is defined as
\[
\tilde{g}(x)\quad\coloneqq\quad\big(\,(g\Gamma_X(n,x))_{n\in\cat{M}},\ 
(\gamma_{m,\Gamma_X(n,x)})_{m,n\in\cat{M}}\,\big)
\]
on an object~$x$.
Let us observe that the structure maps are well-typed:
\begin{alignat*}{3}
&m\act g\Gamma_X(n,x)&
&\ccol{\ =\ }&
&\gm{T}(\cat{M}\times g)(m,\Gamma_X(n,x))\\
&&
&\ccol{\ \xrightarrow{\gamma_{m,\Gamma_X(n,x)}}\ }&
&g\Gamma_X(m,\Gamma_X(n,x))\\
&&
&\ccol{\ =\ }&
&g\Gamma_XM_{GX}(m,n,x)\tag*{by (\ref{eqn-univ-gem-3})}\\
&&
&\ccol{\ =\ }&
&g\Gamma_X(m\tensor n,x).
\end{alignat*}
$\tilde{g}$ is defined as 
\[
\tilde{g}(z)\quad\coloneqq\quad\big(\,
g\Gamma_X(n,z)\colon g\Gamma_X(n,x)\to
g\Gamma_X(n,x^\prime)\,\big)_{n\in\cat{M}}
\]
on a morphism~$z\colon x\to x^\prime$.
\end{defn}

\begin{prop}
The functor~$\tilde{g}$ defined above is indeed a 1-cell 
of $\EM{\Cat}{\cat{M}\times(-)}$ which satisfies
(\ref{eqn-univ-gem-5}) and (\ref{eqn-univ-gem-6}).
Moreover, it is the unique such.
\end{prop}
\begin{pf}
First, $\tilde{g}$ is a 1-cell of $\EM{\Cat}{\cat{M}\times(-)}$, i.e.,
the diagram
\begin{equation}
\label{diag-1-cell-EM-M-times}
\begin{tikzpicture}[baseline=-\the\dimexpr\fontdimen22\textfont2\relax ]
      \node (TL)  at (0,0.75)  {$\cat{M}\times GX$};
      \node (TR)  at (3,0.75)  {$\cat{M}\times \EMg$};
      \node (BL)  at (0,-0.75) {$GX$};
      \node (BR)  at (3,-0.75) {$\EMg$};
      
      \draw [->] (TL) to node (T) [auto,labelsize]      {$\cat{M}\times \tilde{g}$} (TR);
      \draw [->] (TR) to node (R) [auto,labelsize]      {$\emact$} (BR);
      \draw [->] (TL) to node (L) [auto,swap,labelsize] {$\Gamma_X$}(BL);
      \draw [->] (BL) to node (B) [auto,swap,labelsize] {$\tilde{g}$}(BR);
\end{tikzpicture}
\end{equation}
commutes, because, 
\begin{align*}
\tilde{g}\Gamma_X (p,x)
\ &=\ 
((g\Gamma_X(n,\Gamma_X (p,x)))_{n\in\cat{M}}
,\,
(\gamma_{m,\Gamma_X(n,\Gamma_X(p,x))})_{m,n\in\cat{M}})\\
&=\ 
((g\Gamma_X(\cat{M}\times\Gamma_X)(n,p,x))_{n\in\cat{M}}
,\,
(\gamma_{m,\Gamma_X(\cat{M}\times\Gamma_X)(n,p,x)})_{m,n\in\cat{M}})\\
&=\ 
((g\Gamma_X M_{GX}(n,p,x))_{n\in\cat{M}}
,\,
(\gamma_{m,\Gamma_X M_{GX}(n,p,x)})_{m,n\in\cat{M}})
\tag*{by (\ref{eqn-univ-gem-3})}\\
&=\ 
((g\Gamma_X (n\tensor p,x))_{n\in\cat{M}}
,\,
(\gamma_{m,\Gamma_X (n\tensor p,x)})_{m,n\in\cat{M}})\\
&=\ 
p\emact 
((g\Gamma_X (n,x))_{n\in\cat{M}}
,\,
(\gamma_{m,\Gamma_X (n,x)})_{m,n\in\cat{M}})\\
&=\ 
p\emact \tilde{g}x
\end{align*}
on objects and 
\begin{align*}
\tilde{g}\Gamma_X (u,z)
\ &=\ 
(g\Gamma_X(n,\Gamma_X (u,z)))_{n\in\cat{M}}\\
&=\ 
(g\Gamma_X(\cat{M}\times\Gamma_X)(n,u,z))_{n\in\cat{M}}\\
&=\ 
(g\Gamma_X M_{GX}(n,u,z))_{n\in\cat{M}}
\tag*{by (\ref{eqn-univ-gem-3})}\\
&=\ 
(g\Gamma_X (n\tensor u,z))_{n\in\cat{M}}\\
&=\ 
u\emact 
(g\Gamma_X (n,z))_{n\in\cat{M}}\\
&=\ 
u\emact \tilde{g}z
\end{align*}
on morphisms.

$\tilde{g}$ satisfies (\ref{eqn-univ-gem-5}) since,
\begin{align*}
gx
\ &=\ 
g\Gamma_X H_{GX}x
\tag*{by (\ref{eqn-univ-gem-1})}\\
&=\ 
g\Gamma_X (\e,x)\\
&=\ 
u^\gm{T}
((g\Gamma_X(n,x))_{n\in\cat{M}}
,\,
(\gamma_{m,\Gamma_X(n,x)})_{m,n\in\cat{M}})\\
&=\ 
u^\gm{T}\tilde{g}x
\end{align*}
on objects and 
\begin{align*}
gz
\ &=\ 
g\Gamma_X H_{GX}z
\tag*{by (\ref{eqn-univ-gem-1})}\\
&=\ 
g\Gamma_X (\e,z)\\
&=\ 
u^\gm{T}
(g\Gamma_X(n,z))_{n\in\cat{M}}\\
&=\ 
u^\gm{T}\tilde{g}z
\end{align*}
on morphisms.

Finally, $\tilde{g}$ satisfies (\ref{eqn-univ-gem-6}) since,
\begin{align*}
\gamma_{m,x}
\ &=\ 
\gamma_{m,\Gamma_X(\e,x)}\tag*{by (\ref{eqn-univ-gem-1})}\\
&=\ 
\varepsilon^\gm{T}_{m,\tilde{g}x,\e}\\
&=\ 
u^\gm{T}\varepsilon^\gm{T}_{m,\tilde{g}x}
\end{align*}
on objects.

For the uniqueness, suppose that a functor~$\hat{g}\colon GX\to\EMg$ with
$\hat{g}x=((A^x_n)_{n\in\cat{M}},$ $(h^x_{m,n})_{m,n\in\cat{M}})$ on 
objects
and $\hat{g}z=(\varphi^z_n)_{n\in\cat{M}}$ on morphisms
also satisfy the conditions.
First, equation~(\ref{eqn-univ-gem-5}) forces 
$A^x_\e=gx$ and $\varphi^z_\e=gz$,
whereas equation~(\ref{eqn-univ-gem-6}) says that 
$h^x_{m,\e}=\gamma_{m,x}$.
Now the requirement that $\hat{g}$ is a 1-cell of 
$\EM{\Cat}{\cat{M}\times(-)}$ determines everything else.
By chasing diagram~(\ref{diag-1-cell-EM-M-times})
(with $\tilde{g}$ replaced by $\hat{g}$) starting from the object
$(n,x)$ and evaluating at $\e$ we may conclude
$A^x_n=\Gamma_X(n,x)$ and $h^x_{m,n}=\gamma_{m,\Gamma_X(n,x)}$;
whereas chasing it 
starting from the morphism $(\id{n},z)$ and evaluating at $\e$
enables us to conclude
$\varphi^z_n=g\Gamma_X(n,z)$.
\end{pf}

\subsubsection{The 2-dimensional aspect}
Let us proceed to the unique factorization for 
morphisms of left $\gm{T}$-modules.
Suppose we have a morphism of left $\gm{T}$-modules, i.e., 
a 2-cell
    \[
    \begin{tikzpicture}
      \node (A)  at (0,0) {$(\tcat{X},X)$};
      \node (C) at (4,0) {$(\Cat,\cat{C})$};
      
      \draw [->] (A) to [bend left=20]node (dom)[auto,labelsize] {$(G,g)$}(C);
      \draw [->] (A) to [bend right=20]node (cod)[auto,swap,labelsize] {$(G^\prime,g^\prime)$}(C);

      \draw [2cellnode] (dom) to node [auto,labelsize]{$(\Omega,\omega)$}(cod);

    \end{tikzpicture}
    \]
of $\Epp$ satisfying
\begin{align}
\label{eqn-univ-gem-7}
    \begin{tikzpicture}[baseline=-\the\dimexpr\fontdimen22\textfont2\relax ]
            \begin{scope} [shift={(0,0.75)}]
      \node (A)  at (0,0) {$\tcat{X}$};
      \node (B) at (3,-1.5) {$\Cat$};
      \node (C) at (3,0) {$\Cat$};
      \draw [->] (A) to [bend left=20]node [auto,labelsize] {$G$}(C);
      \draw [->] (A) to [bend right=20]node [auto,very near start,labelsize] {$G^\prime$}(C);
      \draw [->] (C) to [bend left=0]node [auto,labelsize] {$\cat{M}\times(-)$}(B);
      %
      \draw [->] (A) to [bend right=20]node [auto,swap,labelsize] (Fpp) {$G^\prime$}(B);
      \draw [2cell] (1.5,0.25) to node [auto,labelsize]{$\Omega$}(1.5,-0.25);
      \draw [2cell] (1.5,-0.5) to node [auto,labelsize]{$\Gamma^\prime$}(1.5,-1.1);
            \end{scope}
    \end{tikzpicture}
    \quad=\quad
    \begin{tikzpicture}[baseline=-\the\dimexpr\fontdimen22\textfont2\relax ]
            \begin{scope} [shift={(0,0.75)}]
      \node (A)  at (0,0) {$\tcat{X}$};
      \node (B) at (3,-1.5) {$\Cat$};
      \node (C) at (3,0) {$\Cat$};
      \draw [->] (A) to [bend left=20]node [auto,labelsize] {$G$}(C);
      \draw [->] (A) to [bend left=20]node [auto,near end,labelsize] {$G$}(B);
      \draw [->] (C) to [bend left=0]node [auto,labelsize] {$\cat{M}\times(-)$}(B);
      %
      \draw [->] (A) to [bend right=20]node [auto,swap,labelsize] (Fpp) {$G^\prime$}(B);
      \draw [2cell] (1.5,0.25) to node [auto,labelsize]{$\Gamma$}(1.5,-0.25);
      \draw [2cell] (1.5,-0.5) to node [auto,labelsize]{$\Omega$}(1.5,-1.1);
            \end{scope}
    \end{tikzpicture}
\end{align}
\begin{multline}
\label{eqn-univ-gem-8}
\begin{tikzpicture}[baseline=-\the\dimexpr\fontdimen22\textfont2\relax ]
        \begin{scope} 
      \node (A)  at (0,0) {$\cat{M}\times G^\prime X$};
      \node (B) at (3,0) {$\cat{M}\times\cat{C}$};
      \node (C) at (5,0) {$\cat{C}$};
      \node (D)  at (1.5,-1.5) {$G^\prime X$};
      \node (F) at (1,1.5) {$\cat{M}\times GX$};
      \draw [->] (A) to [bend left=0]node [auto,labelsize] {$\cat{M}\times g^\prime$}(B);
      \draw [->] (B) to [bend left=0]node [auto,labelsize] {$\gm{T}$}(C);
      \draw [->] (A) to [bend right=20]node [auto,swap,labelsize] {$\Gamma^\prime _X$}(D);
      \draw [->] (D) to [bend right=20]node [auto,swap,labelsize] {$g^\prime$}(C);
      \draw [->] (F) to [bend right=20]node [auto,swap,labelsize] {$\cat{M}\times\Omega_X$}(A);
      \draw [->] (F) to [bend left=20]node [auto,labelsize] {$\cat{M}\times g$}(B);
      \draw [2cell] (1.5,-0.1) to node [auto,labelsize]{$\gamma^\prime$}(1.5,-1.1);
      \draw [2cell] (1,1.2) to node [auto,labelsize]{$\cat{M}\times\omega$}(1,0.4);
        \end{scope}
\end{tikzpicture}\\
=\quad
    \begin{tikzpicture}[baseline=-\the\dimexpr\fontdimen22\textfont2\relax ]
            \begin{scope} 
          \node (M)  at (0,0)  {$GX$};
          \node (N)  at (3.8,1.05) {$\cat{M}\times\cat{C}$};
          \node (RR)  at (5,0) {$\cat{C}$};
          \node (B)  at (1.5,1.5) {$\cat{M}\times GX$};
          \node (E) at (1.5,-1.5) {$G^\prime X$};
          \draw [->] (M) to node  [auto,labelsize]      {$g$} (RR);
          \draw [<-] (M) to [bend left=20]node  [auto,labelsize] {$\Gamma_X$}(B);
          \draw [->] (B) to [bend left=8]node  [auto,labelsize] {$\cat{M}\times g$}(N);
          \draw [->] (N) to [bend left=5]node  [auto,labelsize] {$\gm{T}$}(RR);
          \draw [<-] (E) to [bend left=20]node  [auto,labelsize] {$\Omega_X$}(M);
          \draw [->] (E) to [bend right=20]node  [auto,swap,labelsize] {$g^\prime$}(RR);
          \draw [2cellnode] (B) to node[auto,labelsize] {$\gamma$}(1.5,0);
          \draw [2cellnode] (1.5,0) to node[auto,labelsize] {$\omega$}(E);
            \end{scope}
    \end{tikzpicture}
\end{multline}

Equation~(\ref{eqn-univ-gem-7}) implies the unique factorization
    \[
    \begin{tikzpicture}[baseline=-\the\dimexpr\fontdimen22\textfont2\relax ]
      \node (E) at (-3,0) {$\tcat{X}$};
      \node (A)  at (0,0) {$\Cat$};
      \draw [->] ([yshift=2mm]E.east) to [bend left=20]node(dom) [auto,labelsize] {${G}$} ([yshift=2mm]A.west);
      \draw [->] ([yshift=-2mm]E.east) to [bend right=20]node(cod) [auto,swap,labelsize] {${G^\prime}$}([yshift=-2mm]A.west);
      \draw [2cellnode] (dom) to node [auto,labelsize]{${\Omega}$}(cod);

    \end{tikzpicture}
    \quad=\quad
    \begin{tikzpicture}[baseline=-\the\dimexpr\fontdimen22\textfont2\relax ]
      \node (E) at (-3,0) {$\tcat{X}$};
      \node (A)  at (0,0) {$\EM{\Cat}{\cat{M}\times(-)}$};
      \node (C) at (3,0) {$\Cat$};
      \draw [->] ([yshift=2mm]E.east) to [bend left=20]node(dom) [auto,labelsize] {$\tilde{G}$} ([yshift=2mm]A.west);
      \draw [->] ([yshift=-2mm]E.east) to [bend right=20]node(cod) [auto,swap,labelsize] {$\tilde{G^\prime}$}([yshift=-2mm]A.west);
      \draw [->] (A) to [bend left=0]node [auto,labelsize] {$U^\cat{M}$}(C);
      \draw [2cellnode] (dom) to node [auto,labelsize]{$\tilde{\Omega}$}(cod);

    \end{tikzpicture}
    \]
where the 2-natural transformation~$\tilde{\Omega}$ is defined as follows:

\begin{defn}
The 2-natural transformation~$\tilde{\Omega}\colon\tilde{G}\Rightarrow\tilde{G^\prime}$
has, as components, 1-cells~$
\tilde{\Omega}_{X^\prime}\colon
\valgp{\cat{M}\times G{X^\prime}}{\Gamma_{X^\prime}}{G{X^\prime}}\to
\valgp{\cat{M}\times G^\prime {X^\prime}}{\Gamma^\prime_{X^\prime}}{G^\prime {X^\prime}}$
of $\EM{\Cat}{\cat{M}\times(-)}$
given by $\tilde{\Omega}_{X^\prime}\coloneqq\Omega_{X^\prime}$
as functors.
\end{defn}

Next we construct a 2-cell
    \[
    \begin{tikzpicture}
      \node (a)  at (0,0) {$\valgp{\cat{M}\times G{X}}{\Gamma_{X}}{G{X}}$};
      \node (ap) at (6,0) {$\valgp{\cat{M}\times\EMg}{\emact}{\EMg}$};
      \node (b)  at (2,-3) {$\valgp{\cat{M}\times G^\prime{X}}{\Gamma^\prime_{X}}{G^\prime{X}}$};
      
      \draw [->] (a) to [bend right=0]node [auto,labelsize] {$\tilde{g}$}(ap);
      \draw [->] (a)  to [bend right=20]node [auto,swap,labelsize] {$\tilde{\Omega}_X$}(b);
      \draw [->] (b) to [bend right=20]node [auto,swap,labelsize] {$\tilde{g^\prime}$}(ap);
      \draw [2cell] (2,-0.1) to node [auto,labelsize]{$\tilde{\omega}$}(b);
    \end{tikzpicture}
    \]
of $\EM{\Cat}{\cat{M}\times(-)}$ satisfying
\begin{align}
\label{eqn-univ-gem-9}
    \begin{tikzpicture}[baseline=-\the\dimexpr\fontdimen22\textfont2\relax ]
     \begin{scope} [shift={(0,0.75)}]
      \node (a)  at (0,0) {$GX$};
      \node (ap) at (3,0) {$\cat{C}$};
      \node (b)  at (1,-1.5) {$G^\prime X$};
      \draw [->] (a) to [bend right=0]node [auto,labelsize] {$g$}(ap);
      \draw [->] (a)  to [bend right=20]node [auto,swap,labelsize] {$\Omega_X$}(b);
      \draw [->] (b) to [bend right=20]node [auto,swap,labelsize] {$g^\prime$}(ap);
      \draw [2cell] (1,-0.1) to node [auto,labelsize]{$\omega$}(b);
     \end{scope}
    \end{tikzpicture}
   \quad=\quad
    \begin{tikzpicture}[baseline=-\the\dimexpr\fontdimen22\textfont2\relax ]
     \begin{scope} [shift={(0,0.75)}]
      \node (a)  at (0,0) {$GX$};
      \node (ap) at (3,0) {$\EMg$};
      \node (b)  at (1,-1.5) {$G^\prime X$};
      \node (c)  at (5,0) {$\cat{C}$};
      \draw [->] (a) to [bend right=0]node [auto,labelsize] {$\tilde{g}$}(ap);
      \draw [->] (a)  to [bend right=20]node [auto,swap,labelsize] {$\tilde{\Omega}_X$}(b);
      \draw [->] (b) to [bend right=20]node [auto,swap,labelsize] {$\tilde{g^\prime}$}(ap);
      \draw [->] (ap) to [bend right=0]node [auto,labelsize] {$u^\gm{T}$}(c);
      \draw [2cell] (1,-0.1) to node [auto,labelsize]{$\tilde{\omega}$}(b);
     \end{scope}
    \end{tikzpicture}
\end{align}

\begin{defn}
The 2-cell~$\tilde{\omega}$ of $\EM{\Cat}{\cat{M}\times(-)}$ 
is the natural transformation consisting of
components at $x\in GX$
\begin{multline*}
\tilde{\omega}_x\quad\colon\quad
\big(\,(g\Gamma_X(n,x))_{n\in\cat{M}},\ 
(\gamma_{m,\Gamma_X(n,x)})_{m,n\in\cat{M}}\,\big)\\
\longrightarrow\quad
\big(\,(g^\prime\Gamma^\prime_X(n,\Omega_X x))_{n\in\cat{M}},\ 
(\gamma^\prime_{m,\Gamma^\prime_X(n,\Omega_X x)})_{m,n\in\cat{M}}\,\big)
\end{multline*}
defined as
\[
\tilde{\omega}_{x,n}\coloneqq
\omega_{\Gamma_X(n,x)}\quad\colon\quad
g\Gamma_X(n,x)\quad\longrightarrow\quad
g^\prime\Omega_X\Gamma_X(n,x)\overset{\text{(\ref{eqn-univ-gem-7})}}= 
g^\prime\Gamma^\prime_X(n,\Omega_X x).
\qedhere
\]
\end{defn}
\begin{prop}
The natural transformation~$\tilde{\omega}$ defined above is 
indeed a 2-cell 
of $\EM{\Cat}{\cat{M}\times(-)}$ which satisfies
(\ref{eqn-univ-gem-9}).
Moreover, it is the unique such.
\end{prop}
\begin{pf}
$\tilde{\omega}$ is a 2-cell of $\EM{\Cat}{\cat{M}\times(-)}$, i.e.,
\begin{equation}
\label{diag-2-cell-EM-M-times}
    \begin{tikzpicture}[baseline=-\the\dimexpr\fontdimen22\textfont2\relax ]
     \begin{scope} [shift={(0,-1)}]
      \node (c)  at (0,2) {$\cat{M}\times GX$};
      \node (a)  at (0,0) {$GX$};
      \node (ap) at (4,0) {$\EMg$};
      \node (b)  at (2,-0.5) {$G^\prime X$};
      \draw [->] (c) to [bend right=0]node [auto,swap,labelsize] {$\Gamma_X$}(a);
      \draw [->] (a) to [bend left=25]node [auto,labelsize] {$\tilde{g}$}(ap);
      \draw [->] (a)  to [bend right=10]node [auto,swap,labelsize] {$\tilde{\Omega}_X$}(b);
      \draw [->] (b) to [bend right=10]node [auto,swap,labelsize] {$\tilde{g^\prime}$}(ap);
      \draw [2cell] (2,0.5) to node [auto,labelsize]{$\tilde{\omega}$}(b);
     \end{scope}
    \end{tikzpicture}
\ =\ 
\begin{tikzpicture}[baseline=-\the\dimexpr\fontdimen22\textfont2\relax ]
     \begin{scope} [shift={(0,1)}]
      \node (a)  at (0,0) {$\cat{M}\times GX$};
      \node (ap) at (4,0) {$\cat{M}\times\EMg$};
      \node (b)  at (2,-0.5) {$\cat{M}\times G^\prime X$};
      \node (c)  at (4,-2) {$\EMg$};
      \draw [->] (ap) to [bend right=0]node [auto,labelsize] {$\emact$}(c);
      \draw [->] (a) to [bend left=25]node [auto,labelsize] {$\cat{M}\times\tilde{g}$}(ap);
      \draw [->] (a)  to [bend right=10]node [auto,swap,labelsize] {$\cat{M}\times\tilde{\Omega}_X$}(b);
      \draw [->] (b) to [bend right=10]node [auto,swap,very near start,labelsize] {$\cat{M}\times\tilde{g^\prime}$}(ap);
      \draw [2cell] (2,0.5) to node [auto,labelsize]{$\cat{M}\times\tilde{\omega}$}(b);
     \end{scope}
\end{tikzpicture}
\end{equation}
holds, since 
\begin{align*}
\tilde{\omega}_{\Gamma_X(n,x),m}
\ &=\ 
\omega_{\Gamma_X(m,\Gamma_X(n,x))}\\
&=\ 
\omega_{\Gamma_X(\cat{M}\times\Gamma_X)(m,n,x)}\\
&=\ 
\omega_{\Gamma_XM_{GX}(m,n,x)}\tag*{by (\ref{eqn-univ-gem-3})}\\
&=\ 
\omega_{\Gamma_X(m\tensor n,x)}\\
&=\ 
\tilde{\omega}_{x,m\tensor n}\\
&=\ 
(n\emact \tilde{\omega}_x)_m
\end{align*}
on objects.

$\tilde{\omega}$ satisfies (\ref{eqn-univ-gem-9})
since,
\begin{align*}
\omega_x
\ &=\ 
\omega_{\Gamma_X H_{GX}x}
\tag*{by (\ref{eqn-univ-gem-1})}\\
&=\ 
\omega_{\Gamma_X(\e,x)}\\
&=\ 
u^\gm{T}
(\omega_{\Gamma_X(n,x)})_{n\in\cat{M}}\\
&=\ 
u^\gm{T}\tilde{\omega}_x
\end{align*}
on objects.

For the uniqueness, suppose a natural transformation~$\hat{\omega}
\colon \tilde{g^\prime}\circ\tilde{\Omega}_X\Rightarrow\tilde{g}$
with components~$\hat{\omega}_x=(\varphi^x_n)_{n\in\cat{M}}$
also satisfies the conditions.
Equation~(\ref{eqn-univ-gem-9}) implies that
$\varphi^x_\e=\omega_x$, and putting the object~$(n,x)$
into equation~(\ref{diag-2-cell-EM-M-times}) and evaluating the 
resulting morphism at $\e$ determines that 
$\varphi^x_n=\omega_{\Gamma_X(n,x)}$.
\end{pf}

Finally we may conclude:
\begin{thm}
\label{thm-univ-em-gm}
The object~$\left(\EM{\Cat}{\cat{M}\times(-)},
\valg{\cat{M}\times\EMg}{\emact}{\EMg}\right)$ of 
$\Epp$ is the Eilenberg--Moore object of the
graded monad~$\gm{T}$, considered as a monad
$(\cat{M}\times(-),\gm{T})$ in $\Epp$ on $(\Cat,\cat{C})$.
\end{thm}


\section{The Kleisli construction for graded monads}
\label{sec-kl-gm}
As in the previous section, suppose we have a strict monoidal 
category~$\cat{M}=(\cat{M},\tensor,\e)$, a category~$\cat{C}$, 
and an $\cat{M}$-graded monad~$\gm{T}$ on $\cat{C}$.
We continue to write the functor part of $\gm{T}$
as $\act\colon \cat{M}\times\cat{C}\to\cat{C}$,
by identifying the adjoint transposes.
Following the observation in Section~\ref{subsec-gm-mnd-emm}
this time,
the Kleisli adjunction for $\gm{T}$
lives in $\Emm$ and lies above
the \emph{co-Eilenberg--Moore adjunction} for the \emph{2-comonad}~$[\cat{M},-]$ on 
$\Cat$, as in the picture below:
    \[
    \begin{tikzpicture}
      \node (A)  at (0,0) {$\left(\EM{\Cat}{[\cat{M},-]},
            \vcoalg{[\cat{M},\Klg]}{\klact}{\Klg}\right)$};
      \node (B) at (8,0) {$(\tcat{D},D)$};
      \node (C) at (4,2) {$(\Cat,\cat{C})$};
      
      \node at (2,1) [rotate=120,font=\large]{$\vdash$};
      \node at (6,1) [rotate=240,font=\large]{$\vdash$};
      
      \draw [->] (A) to [bend right=18]node [auto,labelsize,swap] {$(U_\cat{M},u_\gm{T})$}(C);
      \draw [->] (C) to [bend right=16]node [auto,labelsize,swap] {$(F_\cat{M},f_\gm{T})$}(A);
      
      \draw [->] (C) to [bend right=20]node [auto,swap,labelsize] {$(L,l)$}(B);
      \draw [->] (B) to [bend right=20]node [auto,swap,labelsize] {$(R,r)$}(C);
      
      \draw [->,dashed,rounded corners=8pt] (A)--(0,-1) to [bend right=0]node [auto,swap,labelsize]
       {$(K,k)$} (8,-1)--(B);
      
      \draw [rounded corners=15pt,myborder] (-2,-1.5) rectangle (10,2.5);
      \node at (9,2) {$\Emm$};

      \begin{scope}[shift={(0,-4.5)}]
       \node (A)  at (0,0) {$\EM{\Cat}{[\cat{M},-]}$};
       \node (B) at (8,0) {$\tcat{D}$};
       \node (C) at (4,2) {$\Cat$};
      
       \node at (2,1) [rotate=120,font=\large]{$\vdash$};
       \node at (6,1) [rotate=240,font=\large]{$\vdash$};
      
       \draw [<-] (A) to [bend right=21]node [auto,labelsize,swap] {$U_\cat{M}$}(C);
       \draw [<-] (C) to [bend right=21]node [auto,labelsize,swap] {$F_\cat{M}$}(A);
      
       \draw [<-] (C) to [bend right=23]node [auto,swap,labelsize] (Fpp) {$L$}(B);
       \draw [<-] (B) to [bend right=23]node [auto,swap,labelsize] (Fpp) {$R$}(C);
       
       \draw [<-,dashed,rounded corners=8pt] (A)--(0,-0.8) to [bend right=0]node [auto,swap,labelsize]
              {$K$} (8,-0.8)--(B);
       \draw [rounded corners=15pt,myborder] (-2,-1.3) rectangle (10,2.5);
             \node at (9,2) {$\TCat$};
   
      \end{scope}
    \end{tikzpicture}
    \]

\subsection{The Kleisli category}
\begin{defn}
The category~$\Klg$ is defined as follows:
\begin{itemize}
\item An object of $\Klg$ is a pair~$(m,c)$ where $m$ and $c$
are objects of $\cat{M}$ and $\cat{C}$ respectively.
\item The set of morphisms from $(m,c)$ to $(m^\prime,c^\prime)$
is defined by the coend formula
\[
\Klg((m,c),(m^\prime,c^\prime))\quad\coloneqq\quad\int^{n\in \cat{M}}
\cat{M}(m\tensor n,m^\prime)\times\cat{C}(c,n\act c^\prime).
\]
Explicitly, a morphism~$(m,c)\to(m^\prime,c^\prime)$ is an 
equivalence class~$
[n,\,m\tensor n\xrightarrow{v}m^\prime,\,c\xrightarrow{f}n\act c^\prime]$
of tuples consisting of an object~$n\in\cat{M}$ and morphisms~$v,f$, 
where the equivalence relation is generated by 
\begin{multline*}
(n,\ m\tensor n\xrightarrow{m\tensor w}m\tensor n^\prime\xrightarrow{v}m^\prime,\ 
c\xrightarrow{f}n\act c^\prime)\\
\sim\quad(n^\prime,\ m\tensor n^\prime\xrightarrow{v}m^\prime,\ 
c\xrightarrow{f}n\act c^\prime\xrightarrow{w\act c^\prime}n^\prime\act c^\prime)
\end{multline*}
for each morphism~$w\colon n\to n^\prime$ of $\cat{M}$.
\item The identity morphism on $(m,c)$ is given by 
$[\e,\,m\tensor \e\xrightarrow{\id{m}}m,\,
c\xrightarrow{\eta_c}\e\act c]$.
\item For two composable morphisms
\begin{align*}
[n,\ m\tensor n\xrightarrow{v}m^\prime,\ c\xrightarrow{f}n\act c^\prime]
\quad&\colon\quad
(m,c)\quad\longrightarrow\quad(m^\prime,c^\prime),\\
[n^\prime,\ m^\prime\tensor n^\prime\xrightarrow{v^\prime}m^{\prime\prime},\ 
c^\prime\xrightarrow{f^\prime}n^\prime\act c^{\prime\prime}]
\quad&\colon\quad 
(m^\prime,c^\prime)\quad\longrightarrow\quad(m^{\prime\prime},c^{\prime\prime}),
\end{align*}
their composite is given by
\begin{multline*}
[n\tensor n^\prime,\ m\tensor n\tensor n^\prime\xrightarrow{v\tensor n^\prime}m^\prime\tensor n^\prime
\xrightarrow{v^\prime}m^{\prime\prime},\\
c\xrightarrow{f}n\act c^\prime\xrightarrow{n\act f^\prime}
n\act n^\prime\act c^{\prime\prime}\xrightarrow{\mu_{n,n^\prime,c^{\prime\prime}}}
(n\tensor n^\prime)\act c^{\prime\prime}]\\
\quad\colon\quad
(m,c)\quad\longrightarrow\quad(m^{\prime\prime},c^{\prime\prime}).
\tag*{\qedhere}
\end{multline*}
\end{itemize}
\end{defn}

We need to check well-definedness to validate such 
definitions as that of compositions in $\Klg$.
All these are established through straightforward calculation.

\begin{defn}
Define the functor
\[
\klact\quad\colon\quad\Klg\quad\longrightarrow\quad[\cat{M},\Klg]
\]
as follows;
we will use the infix notation~$l\klact (m,c)$
to denote the value of the functor $\klact (m,c)\colon \cat{M}\to\Klg$
applied to $l\in\cat{M}$ and similarly for morphisms.
\begin{itemize}
\item Given objects~$l$ and $(m,c)$ of $\cat{M}$ and $\Klg$ respectively,
we define $l\klact(m,c)\coloneqq(l\tensor m,c)$.
\item Given morphisms~$u\colon l\to l^\prime$ and 
$[n,\,m\tensor n\xrightarrow{v}m^\prime,\,c\xrightarrow{f} n\act c^\prime]
\colon (m,c)\to (m^\prime,c^\prime)$
of $\cat{M}$ and $\Klg$ respectively,
we define
\begin{multline*}
u\klact 
[n,\ m\tensor n\xrightarrow{v}m^\prime,\ c\xrightarrow{f} n\act c^\prime]\\
\coloneqq\quad
[n,\ l\tensor m\tensor n\xrightarrow{u\tensor v}l^\prime\tensor m^\prime,\ 
c\xrightarrow{f}n\act c^\prime]\\
\colon\quad
(l\tensor m,c)\quad\longrightarrow\quad(l^\prime \tensor m^\prime,c^\prime).\tag*{\qedhere}
\end{multline*}
\end{itemize}
\end{defn}

\subsection{The Kleisli adjunction}

\subsubsection{The left adjoint}
We define the 1-cell
\[
(F_\cat{M},f_\gm{T})\quad\colon\quad 
(\Cat,\cat{C})
\quad\longrightarrow\quad
\left(\EM{\Cat}{[\cat{M},-]},
\vcoalg{[\cat{M},\Klg]}{\klact}{\Klg}
\right)
\] 
of $\Emm$ as follows:

\begin{defn}
The 2-functor~$F_\cat{M}\colon \EM{\Cat}{[\cat{M},-]}\to\Cat$
is the forgetful 2-functor 
$\vcoalgp{[\cat{M},\cat{A}]}{\alpha}{\cat{A}}\mapsto \cat{A}$.
\end{defn}

\begin{defn}
The functor~$
f_\gm{T}\colon 
\cat{C}
\to
\Klg$
is defined as
$
f_\gm{T}(c)\coloneqq(\e,c)
$
on an object~$c$ and 
\[
f_\gm{T}(f)\quad\coloneqq\quad
[\e,\ \e\tensor\e\xrightarrow{\id{\e}}\e,\ 
c\xrightarrow{f}c^\prime\xrightarrow{\eta_{c^\prime}}\e\act c^\prime]
\quad\colon\quad (\e,c)\quad\longrightarrow\quad (\e,c^\prime)
\]
on a morphism~$f\colon c\to c^\prime$.
\end{defn}

\subsubsection{The right adjoint}

We define the 1-cell
\[
(U_\cat{M},u_\gm{T})\quad\colon\quad 
\left(\EM{\Cat}{[\cat{M},-]},
\vcoalg{[\cat{M},\Klg]}{\klact}{\Klg}\right)
\quad\longrightarrow\quad
(\Cat,\cat{C})
\]
of $\Emm$ as follows:
\begin{defn}
The 2-functor~$U_\cat{M}\colon\Cat\to\EM{\Cat}{[\cat{M},-]}$ is the cofree 
2-functor~$\cat{X}\mapsto
\vcoalgp{[\cat{M},[\cat{M},\cat{X}]]}{M_{\cat{M},\cat{X}}}{[\cat{M},\cat{X}]}$,
where $M_{\cat{M},\cat{X}}=M_\cat{X}$ is the one defined in
Definition~\ref{defn-M-to}.
\end{defn}

\begin{defn}
The 1-cell~$
u_\gm{T}\colon 
\vcoalgp{[\cat{M},\Klg]}{\klact}{\Klg}
\to
\vcoalgp{[\cat{M},[\cat{M},\cat{C}]]}{M_\cat{C}}{[\cat{M},\cat{C}]}
$
of $\EM{\Cat}{[\cat{M},-]}$ 
is the functor~$
u_\gm{T}\colon \Klg\to[\cat{M},\cat{C}]$
defined as
\[
u_\gm{T}(m,c)\quad\coloneqq\quad ((-)\tensor m)\act c
\quad\colon\quad
\cat{M}\quad\longrightarrow\quad \cat{C}
\]
on an object~$(m,c)$ and 
\begin{multline*}
\big(\,u_\gm{T}[n,\ m\tensor n\xrightarrow{v}m^\prime,\ c\xrightarrow{f} n\act c^\prime]\,\big)_l
\quad\coloneqq\\
(l\tensor m)\act c\xrightarrow{(l\tensor m)\act f}(l\tensor m)\act n\act c^\prime
\xrightarrow{\mu_{l\tensor m,n,c^\prime}}
(l\tensor m\tensor n)\act c^\prime\\
\xrightarrow{(l\tensor v)\act c^\prime}
(l\tensor m^\prime)\act c^\prime
\end{multline*}
on a morphism~$[n,\ m\tensor n\xrightarrow{v}m^\prime,\ c\xrightarrow{f} n\act c^\prime]
\colon (m,c)\to (m^\prime,c^\prime)$.
\end{defn}

\subsubsection{The unit}
We define the 2-cell
    \[
    \begin{tikzpicture}
      \node (A)  at (0,0) {$(\Cat,\cat{C})$};
      \node (B) at (8,0) {$(\Cat,\cat{C})$};
      \node (C) at (4,-2) {$\left(\EM{\Cat}{[\cat{M},-]},
      \vcoalg{[\cat{M},\Klg]}{\klact}{\Klg}\right)$};
      
      \draw [->] (A) to [bend right=20]node [auto,swap,labelsize] {$(F_\cat{M},f_\gm{T})$}(C);
      \draw [->] (C) to [bend right=20]node [auto,swap,labelsize] {$(U_\cat{M},u_\gm{T})$}(B);

      \draw [->] (A) to [bend left=20]node [auto,labelsize] (Fpp) {$(\id{\Cat},\id{\cat{C}})$}(B);
      \draw [2cellnoder] (C) to node [auto,labelsize,swap]{$(H_\cat{M},\eta_\gm{T})$}(Fpp);

    \end{tikzpicture}
    \]
of $\Emm$ as follows:

\begin{defn}
The 2-natural transformation~$
H_\cat{M}\colon F_\cat{M}\circ U_\cat{M}\Rightarrow\id{\Cat}$
is the one defined in Definition~\ref{defn-M-to}.
\end{defn}

\begin{defn}
The natural transformation
    \[
    \begin{tikzpicture}
      \node (L) at (0,0) {$\cat{C}$};
      \node (R) at (6,0) {$\cat{C}$};
      \node (N) at (1.5,-1.35) {$\Klg$};
      \node (T) at (4,-2) {$[\cat{M},\cat{C}]$};
      
      \draw [->] (L) to node [auto,labelsize] {$\id{\cat{C}}$}(R);
      \draw [->] (L) to [bend right=10]node [auto,swap,labelsize] {$f_\gm{T}$}(N);
      \draw [->] (N) to [bend right=10]node [auto,swap,labelsize] {$u_\gm{T}$}(T);
      \draw [->] (T) to [bend right=20]node [auto,swap,labelsize] {$H_{\cat{C}}$}(R);
            
      \draw [2cell] (4,-0.1) to node [auto,labelsize]{$\eta_\gm{T}$}(T);
    \end{tikzpicture}
    \]
has components~$\eta_{\gm{T},c}\colon c\to(\e\act c)$ 
given by the data of the graded monad.
\end{defn}

\subsubsection{The counit}
We define the 2-cell
    \[
    \begin{tikzpicture}
      \node (A)  at (0,0) {$\left(\EM{\Cat}{[\cat{M},-]},
            \vcoalg{[\cat{M},\Klg]}{\klact}{\Klg}\right)$};
      \node (B) at (8,0) {$\left(\EM{\Cat}{[\cat{M},-]},
            \vcoalg{[\cat{M},\Klg]}{\klact}{\Klg}\right)$};
      \node (C) at (4,2) {$(\Cat,\cat{C})$};
      
      \draw [->] (A) to [bend left=20]node [auto,labelsize] {$(U_\cat{M},u_\gm{T})$}(C);
      \draw [->] (C) to [bend left=20]node [auto,labelsize] {$(F_\cat{M},f_\gm{T})$}(B);

      \draw [->] (A) to [bend right=20]node [auto,swap,labelsize] (Fpp) {$(\id{\EM{\Cat}{[\cat{M},-]}},\id{\klact})$}(B);
      \draw [2cellnoder] (Fpp) to node [auto,swap,labelsize]{$(E_\cat{M},\varepsilon_\gm{T})$}(C);

    \end{tikzpicture}
    \]
of $\Emm$ as follows:

\begin{defn}
The 2-natural transformation~$E_\cat{M}\colon \id{\EM{\Cat}{[\cat{M},-]}}
\Rightarrow U_\cat{M}\circ F_\cat{M}$
has components~$
E_{\cat{M},\alpha}\colon 
\vcoalgp{[\cat{M},\cat{A}]}{\alpha}{\cat{A}}\to
\vcoalgp{[\cat{M},[\cat{M},\cat{A}]]}{M_\cat{A}}{[\cat{M},\cat{A}]}$
given by $E_{\cat{M},\alpha}=E_\alpha\coloneqq\alpha\colon [\cat{B},\cat{A}]\to\cat{A}$.
\end{defn}

\begin{defn}
The 2-cell
    \[
    \begin{tikzpicture}
      \node (L) at (0,0) {$\vcoalgp{[\cat{M},\Klg]}{\klact}{\Klg}$};
      \node (R) at (9,0) {$\vcoalgp{[\cat{M},[\cat{M},\Klg]]}{M_\Klg}{[\cat{M},\Klg]}$};
      \node (T) at (6,-3) {$\vcoalgp{[\cat{M},\Klg]}{\klact}{\Klg}$};
      
      \node (M) at (4.5,0) {$\vcoalgp{[\cat{M},[\cat{M},\cat{C}]]}{M_\cat{C}}{[\cat{M},\cat{C}]}$};
      
      \draw [->] (L) to [bend right=0]node [auto,labelsize] {$ u_\gm{T}$}(M);
      \draw [->] (M) to [bend right=0]node [auto,labelsize] {$[\cat{M},f_\gm{T}]$}(R);

      \draw [->] (L) to [bend right=20]node [auto,swap,labelsize] {$\id{\klact}$}(T);
      \draw [->] (T) to [bend right=20]node [auto,labelsize,swap] {$E_{\klact}$}(R);
            
      \draw [2cell] (6,-0.1) to node [auto,labelsize]{$\varepsilon_\gm{T}$}(T);
    \end{tikzpicture}
    \]
of $\EM{\Cat}{[\cat{M},-]}$ is the natural transformation
with the component at $(m,c)\in\Klg$ being
\[
\varepsilon_{\gm{T},(m,c)}\quad\colon\quad
(\e,((-)\tensor m)\act c)
\quad \Longrightarrow\quad
((-)\tensor m,c),
\]
itself with the component at $l\in\cat{M}$
given by
\begin{multline*}
\varepsilon_{\gm{T},(m,c),l}
\quad\coloneqq\quad
[l\tensor m,\ l\tensor m\xrightarrow{\id{}}l\tensor m,\ 
(l\tensor m)\act c\xrightarrow{\id{}}(l\tensor m)\act c]\\
\colon\quad(\e,(l\tensor m)\act c)\quad\longrightarrow\quad(l\tensor m,c).
\tag*{\qedhere}
\end{multline*}
\end{defn}

\subsection{Comparison maps}
Suppose we have an adjunction~$
(L,l)\dashv(R,r)\colon(\tcat{D},D)\to(\Cat,\cat{C})$ in $\Emm$
with unit~$(H,\eta)\colon$ $(\id{\Cat},\id{\cat{C}})\Rightarrow(R,r)\circ(L,l)$ 
and counit~$(E,\varepsilon)\colon(L,l)\circ(R,r)\Rightarrow(\id{\tcat{D}},\id{D})$,
which gives a resolution of the monad $([\cat{M},-],\gm{T})$,
i.e., such that the following equations hold:
{\allowdisplaybreaks
\begin{align}
\label{eqn-comp-gkl-1}
\Cat\xrightarrow{[\cat{M},-]}\Cat\quad
&=
\quad\Cat\xrightarrow{R}\tcat{D}\xrightarrow{L}\Cat\\
\label{eqn-comp-gkl-2}
\cat{C}\xrightarrow{\gm{T}}[\cat{M},\cat{C}]\quad&=\quad
\cat{C}\xrightarrow{l}LD\xrightarrow{Lr}LR\cat{C}\\
\label{eqn-comp-gkl-3}
H_\cat{M}\quad&=\quad H\\
\label{eqn-comp-gkl-4}
\eta_\gm{T}\quad&=\quad\eta\\
\label{eqn-comp-gkl-5}
\begin{tikzpicture}[baseline=-\the\dimexpr\fontdimen22\textfont2\relax ]
        \begin{scope} 
      \node (TL)  at (0,0)  {$\Cat$};
      \node (TT)  at (1.2,0.7) {$\Cat$};
      \node (RR)  at (2.4,0) {$\Cat$};
      \draw [->] (TL) to [bend left=20]node (L) [auto,labelsize] {$[\cat{M},-]$}(TT);
      \draw [->] (TT) to [bend left=20]node  [auto,labelsize] {$[\cat{M},-]$}(RR);
      \draw [->] (TL) to [bend right=20]node (B) [auto,swap,labelsize] {$[\cat{M},-]$} (RR);
      \draw [2cellnoder] (TT) to node[auto,labelsize] {$M_\cat{M}$}(B);
        \end{scope}
\end{tikzpicture}
\quad&=\quad
\begin{tikzpicture}[baseline=-\the\dimexpr\fontdimen22\textfont2\relax ]
        \begin{scope} 
      \node (TL)  at (0,0)  {$\Cat$};
      \node (TR)  at (1.2,0)  {$\tcat{D}$};
      \node (TT)  at (2.4,0.7) {$\Cat$};
      \node (BR)  at (3.6,0) {$\tcat{D}$};
      \node (RR)  at (4.8,0) {$\Cat$};
      \draw [->] (TL) to node (T) [auto,labelsize]      {$R$} (TR);
      \draw [->] (TR) to [bend left=20]node (R) [auto,labelsize]      {$L$} (TT);
      \draw [->] (TT) to [bend left=20]node (L) [auto,labelsize] {$R$}(BR);
      \draw [->] (BR) to node  [auto,labelsize] {$L$}(RR);
      \draw [->] (TR) to [bend right=20]node (B) [auto,swap,labelsize] {$\id{\tcat{D}}$} (BR);
      \draw [2cellnoder] (TT) to node[auto,labelsize] {$E$}(B);
        \end{scope}
\end{tikzpicture}
\end{align}
}%
\vspace{-20pt}
\begin{multline}
\label{eqn-comp-gkl-6}
\begin{tikzpicture}[baseline=-\the\dimexpr\fontdimen22\textfont2\relax ]
        \begin{scope} 
      \node (L)  at (0,0)  {$[\cat{M},[\cat{M},\cat{C}]]$};
      \node (M)  at (-3,0)  {$[\cat{M},\cat{C}]$};
      \node (R)  at (-6,0) {$\cat{C}$};
      \node (B)  at (-2,-2) {$[\cat{M},\cat{C}]$};
      \draw [<-] (L) to node  [auto,swap,labelsize]      {$[\cat{M},\gm{T}]$} (M);
      \draw [<-] (M) to node  [auto,swap,labelsize]      {$\gm{T}$} (R);
      \draw [<-] (L) to [bend left=20]node  [auto,labelsize] {$M_\cat{C}$}(B);
      \draw [<-] (B) to [bend left=20]node  [auto,labelsize] {$\gm{T}$}(R);
      \draw [2cellnode] (-2,0) to node[auto,labelsize] {$\mu_\gm{T}$}(B);
        \end{scope}
\end{tikzpicture}\\
=\quad
\begin{tikzpicture}[baseline=-\the\dimexpr\fontdimen22\textfont2\relax ]
        \begin{scope} [shift={(0,0.75)}]
      \node (L)  at (-0,0)  {$LRLR\cat{C}$};
      \node (M)  at (-2.75,0)  {$LRLD$};
      \node (R)  at (-5.2,0) {$LR\cat{C}$};
      \node (RR)  at (-7.15,0) {$LD$};
      \node (RRR)  at (-8.6,0) {$\cat{C}$};
      \node (B)  at (-4.2,-1.5) {$LD$};
      \draw [<-] (L) to node  [auto,swap,labelsize]      {$LRLr$} (M);
      \draw [<-] (M) to node  [auto,swap,labelsize]      {$LRl$} (R);
      \draw [<-] (R) to node  [auto,swap,labelsize]      {$Lr$} (RR);
      \draw [<-] (RR) to node  [auto,swap,labelsize]      {$l$} (RRR);
      \draw [<-] (M) to [bend left=20]node  [auto,labelsize] {$LE_D$}(B);
      \draw [<-] (B) to [bend left=20]node  [auto,labelsize] {$\id{LD}$}(RR);
      \draw [2cellnode] (-4.2,0) to node[auto,labelsize] {$L\varepsilon$}(B);
        \end{scope}
\end{tikzpicture}
\end{multline}

Equations~(\ref{eqn-comp-gkl-1}), (\ref{eqn-comp-gkl-3}) and
(\ref{eqn-comp-gkl-5}) imply the existence of the comparison
2-functor~$K$:

\begin{defn}
The 2-functor~$K\colon \tcat{D}\to\EM{\Cat}{[\cat{M},-]}$ 
is the comparison 2-functor
$D^\prime\mapsto\vcoalgp{[\cat{M},LD^\prime]}{LE_{D^\prime}}{LD^\prime}$.
\end{defn}

This provides the 2-functor part of the comparison map
under construction.
The remaining piece of data is given as follows:

\begin{defn}
The 1-cell~$k\colon\vcoalgp{[\cat{M},\Klg]}{\klact}{\Klg}
\to\vcoalgp{[\cat{M},LD]}{LE_D}{LD}$
of $\EM{\Cat}{[\cat{M},-]}$
is the functor~$
k\colon \Klg\to LD$ defined as
\[
k(m,c)\quad\coloneqq\quad (LE_Dlc)m
\]
on an object~$(m,c)$ and 
\[
k[n,\ m\tensor n\xrightarrow{v}m^\prime,\ c\xrightarrow{f} n\act c^\prime]
\quad\coloneqq\quad
((LE_Dlc^\prime)v)
\circ
(LE_DL\varepsilon_{lc^\prime,n})_m
\circ
(LE_Dlf)_m
\]
on a morphism~$
[n,\ m\tensor n\xrightarrow{v}m^\prime,\ c\xrightarrow{f} n\act c^\prime]
\colon (m,c)\to(m^\prime,c^\prime)$.
Let us check that the type of $k[n,\ m\tensor n\xrightarrow{v}m^\prime,\ 
c\xrightarrow{f} n\act c^\prime]$ is indeed the right one:
\begin{alignat*}{3}
&(LE_Dlc)m&
&\ccol{\ \xrightarrow{(LE_Dlf)_m}\ }&
&(LE_Dl(n\act c^\prime))m\\
&&
&\ccol{\ =\ }&
&(LE_D((LRlLrlc^\prime) n)) m\tag*{by (\ref{eqn-comp-gkl-2})}\\
&&
&\ccol{\ \xrightarrow{(LE_DL\varepsilon_{lc^\prime,n})_m}\ }&
&(LE_D((LE_Dlc^\prime) n)) m\\
&&
&\ccol{\ =\ }&
&((LRLE_DLE_Dlc^\prime) n) m\\
&&
&\ccol{\ =\ }&
&((LE_{RLD}LE_Dlc^\prime) n) m\\
&&
&\ccol{\ =\ }&
&((M_{LD}LE_Dlc^\prime) n) m
\tag*{by (\ref{eqn-comp-gkl-5})}\\
&&
&\ccol{\ =\ }&
&(LE_Dlc^\prime) (m\tensor n)\\
&&
&\ccol{\ \xrightarrow{(LE_Dlc^\prime)v}\ }&
&(LE_Dlc^\prime)m^\prime.\tag*{\qedhere}
\end{alignat*}
\end{defn}

\begin{prop}
\label{prop-comparison-kl-gm}
The 1-cell~$(K,k)$ of $\Emm$ satisfies the equations~$
(K,k)\circ(F_\cat{M},f_\gm{T})=(L,l)$ and
$(U_\cat{M},u_\gm{T})=(R,r)\circ(K,k)$.
Moreover, it is the unique such.
\end{prop}

\subsection{The 2-dimensional universality}

\subsubsection{Statement of the theorem}

We will show that there is a family of
isomorphisms of categories
\begin{multline*}
  \Emm\left(\left(\EM{\Cat}{[\cat{M},-]},
  \vcoalg{[\cat{M},\Klg]}{\klact}{\Klg}\right),\ 
  (\tcat{X},X)\right)\\
  \cong\quad 
  \Emm\big((\Cat,\cat{C}),\ (\tcat{X},X)\big)^{\Emm(([\cat{M},-],\gm{T}),\ (\tcat{X},X))}
\end{multline*}
2-natural in $(\tcat{X},X)\in \Emm$; cf.~Definition~\ref{defn-Kl-obj}.
We claim that the data 
    \[
    \begin{tikzpicture}
      \node (A)  at (0,0) {$\left(\EM{\Cat}{[\cat{M},-]},
        \vcoalg{[\cat{M},\Klg]}{\klact}{\Klg}\right)$};
      \node (B) at (-4,-2) {$\left(\EM{\Cat}{[\cat{M},-]},
        \vcoalg{[\cat{M},\Klg]}{\klact}{\Klg}\right)$};
      \node (C) at (-4,0) {$(\Cat,\cat{C})$};
      \node (D) at (-4,-4) {$(\Cat,\cat{C})$};
      
      \draw [->] (C) to [bend left=0]node [auto,labelsize] {$(F_\cat{M},f_\gm{T})$}(A);
      \draw [->] (B) to [bend left=0]node [auto,labelsize] {$(U_\cat{M},u_\gm{T})$}(C);
      \draw [->] (D) to [bend left=0]node [auto,labelsize] {$(F_\cat{M},f_\gm{T})$}(B);

      \draw [->] (B) to [bend right=20]node [auto,swap,labelsize] (Fpp) {$(\id{\EM{\Cat}{[\cat{M},-]}},\id{\klact})$}(A);
      \draw [2cell] (-2,-0.1) to node [auto,swap,labelsize]{$(E_\cat{M},\varepsilon_\gm{T})$}(-2,-1.6);

    \end{tikzpicture}
    \]
provides the universal right $\gm{T}$-module, in the sense 
that every right $\gm{T}$-module
    \[
    \begin{tikzpicture}
      \node (A)  at (0,0) {$(\tcat{X},X)$};
      \node (B) at (-3,-1.5) {$(\Cat,\cat{C})$};
      \node (C) at (-3,0) {$(\Cat,\cat{C})$};
      
      \draw [->] (C) to [bend left=0]node [auto,labelsize] {$(G,g)$}(A);
      \draw [->] (B) to [bend left=0]node [auto,labelsize] {$([\cat{M},-],\gm{T})$}(C);

      \draw [->] (B) to [bend right=20]node [auto,swap,labelsize] (Fpp) {$(G,g)$}(A);
      \draw [2cell] (-1.5,-0.1) to node [auto,swap,labelsize]{$(\Gamma,\gamma)$}(-1.5,-1.1);

    \end{tikzpicture}
    \]
i.e., an object of the category~$
\Emm((\Cat,\cat{C}),\ (\tcat{X},X))^{\Emm(([\cat{M},-],\gm{T}),\ (\tcat{X},X))}$,
factors uniquely as 
    \[
    \begin{tikzpicture}
      \node (E) at (4,0) {$(\tcat{X},X)$};
      \node (A)  at (0,0) {$\left(\EM{\Cat}{[\cat{M},-]},
              \vcoalg{[\cat{M},\Klg]}{\klact}{\Klg}\right)$};
      \node (B) at (-4,-2) {$\left(\EM{\Cat}{[\cat{M},-]},
              \vcoalg{[\cat{M},\Klg]}{\klact}{\Klg}\right)$};
      \node (C) at (-4,0) {$(\Cat,\cat{C})$};
      \node (D) at (-4,-4) {$(\Cat,\cat{C})$};

      \draw [->] (A) to [bend left=0]node [auto,labelsize] {$(\tilde{G},\tilde{g})$}(E);
      \draw [->] (C) to [bend left=0]node [auto,labelsize] {$(F_\cat{M},f_\gm{T})$}(A);
      \draw [->] (B) to [bend left=0]node [auto,labelsize] {$(U_\cat{M},u_\gm{T})$}(C);
      \draw [->] (D) to [bend left=0]node [auto,labelsize] {$(F_\cat{M},f_\gm{T})$}(B);

      \draw [->] (B) to [bend right=20]node [auto,swap,labelsize] (Fpp) {$(\id{\EM{\Cat}{[\cat{M},-]}},\id{\klact})$}(A);
      \draw [2cell] (-2,-0.1) to node [auto,swap,labelsize]{$(E_\cat{M},\varepsilon_\gm{T})$}(-2,-1.6);

    \end{tikzpicture}
    \]
and similarly every morphism of right $\gm{T}$-modules
    \[
    \begin{tikzpicture}
      \node (A)  at (0,0) {$(\tcat{X},X)$};
      \node (C) at (-4,0) {$(\Cat,\cat{C})$};
      
      \draw [->] (C) to [bend left=20]node (dom)[auto,labelsize] {$(G,g)$}(A);
      \draw [->] (C) to [bend right=20]node (cod)[auto,swap,labelsize] {$(G^\prime,g^\prime)$}(A);

      \draw [2cellnode] (dom) to node [auto,labelsize]{$(\Omega,\omega)$}(cod);

    \end{tikzpicture}
    \]
i.e., a morphism of the category~$
\Emm((\Cat,\cat{C}),\ (\tcat{X},X))^{\Emm(([\cat{M},-],\gm{T}),\ (\tcat{X},X))}$,
factors uniquely as 
    \[
    \begin{tikzpicture}
      \node (E) at (5,0) {$(\tcat{X},X)$};
      \node (A)  at (0,0) {$\left(\EM{\Cat}{[\cat{M},-]},
           \vcoalg{[\cat{M},\Klg]}{\klact}{\Klg}\right)$};
      \node (C) at (-4,0) {$(\Cat,\cat{C})$};

      \draw [->] ([yshift=3mm]A.east) to [bend left=20]node(dom) [auto,labelsize] {$(\tilde{G},\tilde{g})$} ([yshift=3mm]E.west);
      \draw [->] ([yshift=-3mm]A.east) to [bend right=20]node(cod) [auto,swap,labelsize] {$(\tilde{G^\prime},\tilde{g^\prime})$}([yshift=-3mm]E.west);
      \draw [->] (C) to [bend left=0]node [auto,labelsize] {$(F_\cat{M},f_\gm{T})$}(A);

      \draw [2cellnode] (dom) to node [auto,labelsize]{$(\tilde{\Omega},\tilde{\omega})$}(cod);

    \end{tikzpicture}
    \]
    
\subsubsection{The 1-dimensional aspect}
Suppose we have a right $\gm{T}$-module, i.e., a piece of data
    \[
    \begin{tikzpicture}
      \node (A)  at (0,0) {$(\tcat{X},X)$};
      \node (B) at (-3,-1.5) {$(\Cat,\cat{C})$};
      \node (C) at (-3,0) {$(\Cat,\cat{C})$};
      
      \draw [->] (C) to [bend left=0]node [auto,labelsize] {$(G,g)$}(A);
      \draw [->] (B) to [bend left=0]node [auto,labelsize] {$([\cat{M},-],\gm{T})$}(C);

      \draw [->] (B) to [bend right=20]node [auto,swap,labelsize] (Fpp) {$(G,g)$}(A);
      \draw [2cell] (-1.5,-0.1) to node [auto,swap,labelsize]{$(\Gamma,\gamma)$}(-1.5,-1.1);

    \end{tikzpicture}
    \]
in $\Emm$ satisfying
{\allowdisplaybreaks
\begin{align}
\label{eqn-univ-gkl-1}
\begin{tikzpicture}[baseline=-\the\dimexpr\fontdimen22\textfont2\relax ]
        \begin{scope} [shift={(0,0.75)}]
      \node (A)  at (0,0) {$\tcat{X}$};
      \node (B) at (-3,-1.5) {$\Cat$};
      \node (C) at (-3,0) {$\Cat$};
      \draw [->] (A) to [bend left=0]node [auto,swap,labelsize] {$G$}(C);
      \draw [->] (C) to [bend left=0]node (t)[auto,labelsize] {$[\cat{M},-]$}(B);
      \draw [->] (C) to [bend right=90]node(id) [auto,swap,labelsize] {$\id{\Cat}$}(B);
      %
      \draw [->] (A) to [bend left=20]node [auto,labelsize] (Fpp) {$G$}(B);
      \draw [2cellr] (-1.5,-0.1) to node [auto,labelsize]{$\Gamma$}(-1.5,-1.1);
      \draw [2cellnoder] (id) to node [auto,labelsize]{$H_\cat{M}$}(t);
        \end{scope}
\end{tikzpicture}
\quad&=\quad
\begin{tikzpicture}[baseline=-\the\dimexpr\fontdimen22\textfont2\relax ]
        \begin{scope} 
       \node (A)  at (0,0) {$\tcat{X}$};
       \node (C) at (-2,0) {$\Cat$};
       \draw [->] (A) to [bend left=0]node [auto,swap,labelsize] {$G$}(C);
        \end{scope}
\end{tikzpicture}\\
\label{eqn-univ-gkl-2}
\begin{tikzpicture}[baseline=-\the\dimexpr\fontdimen22\textfont2\relax ]
        \begin{scope} 
      \node (A)  at (0,0) {$[\cat{M},GX]$};
      \node (B) at (-3,0) {$[\cat{M},\cat{C}]$};
      \node (C) at (-5,0) {$\cat{C}$};
      \node (D)  at (-1.5,-1.5) {$GX$};
      \node (E) at (-3.5,1.5) {$\cat{C}$};
      \node (F) at (-0.5,1.5) {$GX$};
      \draw [->] (B) to [bend left=0]node [auto,labelsize] {$[\cat{M},g]$}(A);
      \draw [->] (C) to [bend left=0]node [auto,labelsize] {$\gm{T}$}(B);
      \draw [->] (D) to [bend right=20]node [auto,swap,labelsize] {$\Gamma_X$}(A);
      \draw [->] (C) to [bend right=20]node [auto,swap,labelsize] {$g$}(D);
      \draw [->] (B) to [bend right=20]node [auto,swap,labelsize] {$H_\cat{C}$}(E);
      \draw [->] (E) to [bend right=0]node [auto,labelsize] {$g$}(F);
      \draw [->] (C) to [bend left=20]node [auto,labelsize] {$\id{\cat{C}}$}(E);
      \draw [->] (A) to [bend right=20]node [auto,swap,labelsize] {$H_{GX}$}(F);
      \draw [2cell] (-1.5,-0.1) to node [auto,labelsize]{$\gamma$}(-1.5,-1.1);
      \draw [2cell] (-3.5,1.2) to node [auto,labelsize]{$\eta$}(-3.5,0.3);
        \end{scope}
\end{tikzpicture}
\quad&=\quad
\begin{tikzpicture}[baseline=-\the\dimexpr\fontdimen22\textfont2\relax ]
        \begin{scope} 
       \node (A)  at (0,0) {$G{X}$};
       \node (C) at (-2,0) {$\cat{C}$};
       \draw [->] (C) to [bend left=0]node [auto,labelsize] {$g$}(A);
        \end{scope}
\end{tikzpicture}\\
\label{eqn-univ-gkl-3}
\begin{tikzpicture}[baseline=-\the\dimexpr\fontdimen22\textfont2\relax ]
        \begin{scope} 
      \node (A)  at (0,0) {$\tcat{X}$};
      \node (B) at (-3,-1.5) {$\Cat$};
      \node (C) at (-4.5,0) {$\Cat$};
      \node (D) at (-3,1.5) {$\Cat$};
      \draw [->] (D) to [bend left=0]node [auto,labelsize] {$[\cat{M},-]$}(B);
      \draw [->] (D) to [bend left=-20]node (t)[auto,swap,labelsize] {$[\cat{M},-]$}(C);
      \draw [->] (C) to [bend left=-20]node (t)[auto,swap,labelsize] {$[\cat{M},-]$}(B);
      \draw [->] (A) to [bend left=-20]node [auto,swap,labelsize] (Fpp) {$G$}(D);
      \draw [->] (A) to [bend right=-20]node [auto,labelsize] (Fpp) {$G$}(B);
      \draw [2cellr] (-1.5,1) to node [auto,labelsize]{$\Gamma$}(-1.5,-1);
      \draw [2cellr] (-4.1,0) to node [auto,labelsize]{$M_\cat{M}$}(-3.1,0);
        \end{scope}
\end{tikzpicture}
\quad&=\quad 
\begin{tikzpicture}[baseline=-\the\dimexpr\fontdimen22\textfont2\relax ]
        \begin{scope} 
      \node (A)  at (0,0) {$\tcat{X}$};
      \node (B) at (-3,-1.5) {$\Cat$};
      \node (C) at (-4,0) {$\Cat$};
      \node (D) at (-3,1.5) {$\Cat$};
      \draw [->] (A) to [bend left=0]node [auto,swap,labelsize] {$G$}(C);
      \draw [->] (D) to [bend left=-20]node (t)[auto,swap,labelsize] {$[\cat{M},-]$}(C);
      \draw [->] (C) to [bend left=-20]node (t)[auto,swap,labelsize] {$[\cat{M},-]$}(B);
      \draw [->] (A) to [bend left=-20]node [auto,swap,labelsize] (Fpp) {$G$}(D);
      \draw [->] (A) to [bend right=-20]node [auto,labelsize] (Fpp) {$G$}(B);
      \draw [2cellr] (-3,1.2) to node [auto,labelsize]{$\Gamma$}(-3,0.1);
      \draw [2cellr] (-3,-0.1) to node [auto,labelsize]{$\Gamma$}(-3,-1.2);
        \end{scope}
\end{tikzpicture}
\end{align}
}%
\vspace{-20pt}
\begin{multline}
\label{eqn-univ-gkl-4}
\begin{tikzpicture}[baseline=-\the\dimexpr\fontdimen22\textfont2\relax ]
        \begin{scope} 
      \node (A)  at (0,0) {$[\cat{M},GX]$};
      \node (B) at (-4,0) {$[\cat{M},\cat{C}]$};
      \node (C) at (-7,0) {$\cat{C}$};
      \node (D)  at (-2,-1.5) {$GX$};
      \node (E) at (-4,1.5) {$[\cat{M},[\cat{M},\cat{C}]]$};
      \node (F) at (0,1.5) {$[\cat{M},[\cat{M},GX]]$};
      \node (G) at (-7,1.5) {$[\cat{M},\cat{C}]$};
      \draw [->] (B) to [bend left=0]node [auto,labelsize] {$[\cat{M},g]$}(A);
      \draw [->] (C) to [bend left=0]node [auto,swap,labelsize] {$\gm{T}$}(B);
      \draw [->] (D) to [bend right=20]node [auto,swap,labelsize] {$\Gamma_X$}(A);
      \draw [->] (C) to [bend right=15]node [auto,swap,labelsize] {$g$}(D);
      \draw [->] (B) to [bend right=0]node [auto,swap,labelsize] {$M_\cat{C}$}(E);
      \draw [->] (C) to [bend right=0]node [auto,labelsize] {$\gm{T}$}(G);
      \draw [->] (G) to [bend right=0]node [auto,labelsize] {$[\cat{M},\gm{T}]$}(E);
      \draw [->] (E) to [bend right=0]node [auto,labelsize] {$[\cat{M},[\cat{M},g]]$}(F);
      \draw [->] (A) to [bend right=0]node [auto,swap,labelsize] {$M_{GX}$}(F);
      \draw [2cell] (-2,-0.1) to node [auto,labelsize]{$\gamma$}(-2,-1.1);
      \draw [2cell] (-5.5,1.4) to node [auto,labelsize]{$\mu$}(-5.5,0.1);
        \end{scope}
\end{tikzpicture}\\
=\quad
\begin{tikzpicture}[baseline=-\the\dimexpr\fontdimen22\textfont2\relax ]
        \begin{scope} 
      \node (A)  at (0,0) {$[\cat{M},GX]$};
      \node (B) at (-4,0) {$[\cat{M},\cat{C}]$};
      \node (C) at (-7,0) {$\cat{C}$};
      \node (D)  at (-2,-1.5) {$GX$};
      \node (E) at (-4,1.5) {$[\cat{M},[\cat{M},\cat{C}]]$};
      \node (F) at (0,1.5) {$[\cat{M},[\cat{M}, GX]]$};
      \draw [->] (B) to [bend left=0]node [auto,labelsize] {$[\cat{M}, g]$}(A);
      \draw [->] (C) to [bend left=0]node [auto,labelsize] {$\gm{T}$}(B);
      \draw [->] (D) to [bend right=20]node [auto,swap,labelsize] {$\Gamma_X$}(A);
      \draw [->] (C) to [bend right=15]node [auto,swap,labelsize] {$g$}(D);
      \draw [->] (B) to [bend right=0]node [auto,labelsize] {$[\cat{M},\gm{T}]$}(E);
      \draw [->] (E) to [bend right=0]node [auto,labelsize] {$[\cat{M},[\cat{M},g]]$}(F);
      \draw [->] (A) to [bend right=0]node [auto,swap,labelsize] {$[\cat{M},\Gamma_{X}]$}(F);
      \draw [2cell] (-2,-0.1) to node [auto,labelsize]{$\gamma$}(-2,-1.1);
      \draw [2cell] (-2,1.4) to node [auto,labelsize]{$[\cat{M},\gamma]$}(-2,0.5);
        \end{scope}
\end{tikzpicture}
\end{multline}

Equations~(\ref{eqn-univ-gkl-1}) and (\ref{eqn-univ-gkl-3}) imply the 
unique factorization
    \[
    \begin{tikzpicture}[baseline=-\the\dimexpr\fontdimen22\textfont2\relax ]
            \begin{scope} [shift={(0,0.75)}]
      \node (A)  at (0,0) {$\tcat{X}$};
      \node (B) at (-3,-1.5) {$\Cat$};
      \node (C) at (-3,0) {$\Cat$};
          
      \draw [->] (A) to [bend left=0]node [auto,swap,labelsize] {$G$}(C);
      \draw [->] (C) to [bend left=0]node [auto,swap,labelsize] {$[\cat{M},-]$}(B);

      \draw [->] (A) to [bend right=-20]node [auto,labelsize] (Fpp) {$G$}(B);
      \draw [2cellr] (-1.5,-0.1) to node [auto,labelsize]{$\Gamma$}(-1.5,-1.1);
            \end{scope}
    \end{tikzpicture}
    \quad=\quad
    \begin{tikzpicture}[baseline=-\the\dimexpr\fontdimen22\textfont2\relax ]
            \begin{scope} [shift={(0,1.5)}]
      \node (E) at (3,0) {$\tcat{X}$};
      \node (A)  at (0,0) {$\EM{\Cat}{[\cat{M},-]}$};
      \node (B) at (-3,-1.5) {$\EM{\Cat}{[\cat{M},-]}$};
      \node (C) at (-3,0)  {$\Cat$};
      \node (D) at (-3,-3) {$\Cat$};

      \draw [->] (E) to [bend left=0]node [auto,swap,labelsize] {$\tilde{G}$}(A);
      \draw [->] (A) to [bend left=0]node [auto,swap,labelsize] {$F_\cat{M}$}(C);
      \draw [->] (C) to [bend left=0]node [auto,swap,labelsize] {$U_\cat{M}$}(B);
      \draw [->] (B) to [bend left=0]node [auto,swap,labelsize] {$F_\cat{M}$}(D);

      \draw [->] (A) to [bend right=-20]node [auto,labelsize] (Fpp) {$\id{\EM{\Cat}{[\cat{M},-]}}$}(B);
      \draw [2cellr] (-1.5,-0.1) to node [auto,labelsize]{$E_\cat{M}$}(-1.5,-1.1);
            \end{scope}
    \end{tikzpicture}
    \]
since the co-Eilenberg--Moore 2-category~$\EM{\Cat}{[\cat{M},-]}$
is the co-Eilenberg--Moore object in $\TCat$.
\begin{defn}
The 2-functor~$\tilde{G}\colon\tcat{X}\to\EM{\Cat}{[\cat{M},-]}$
is the mediating 2-functor
$X^\prime\mapsto\vcoalgp{[\cat{M},GX^\prime]}{\Gamma_{X^\prime}}{GX^\prime}$.
\end{defn}

Now it remains to construct a 1-cell~$\tilde{g}\colon
\vcoalgp{[\cat{M},\Klg]}{\klact}{\Klg}\to
\vcoalgp{[\cat{M},GX]}{\Gamma_{X}}{GX}$
of $\EM{\Cat}{[\cat{M},-]}$ which satisfies
\begin{align}
\label{eqn-univ-gkl-5}
\begin{tikzpicture}[baseline=-\the\dimexpr\fontdimen22\textfont2\relax ]
        \begin{scope} 
      \node (A)  at (0,0) {$GX$};
      \node (B) at (-2,0) {$\cat{C}$};
      \draw [->] (B) to [bend left=0]node [auto,labelsize] {${g}$}(A);
        \end{scope}
\end{tikzpicture}\quad
=\quad\begin{tikzpicture}[baseline=-\the\dimexpr\fontdimen22\textfont2\relax ]
        \begin{scope} 
      \node (A)  at (0,0) {$GX$};
      \node (B) at (-2,0) {$\Klg$};
      \node (C) at (-4,0) {$\cat{C}$};
      \draw [->] (B) to [bend left=0]node [auto,labelsize] {$\tilde{g}$}(A);
      \draw [->] (C) to [bend left=0]node [auto,labelsize] {$f_\gm{T}$}(B);
        \end{scope}
\end{tikzpicture}
\end{align}
\begin{multline}
\label{eqn-univ-gkl-6}
\begin{tikzpicture}[baseline=-\the\dimexpr\fontdimen22\textfont2\relax ]
        \begin{scope} [shift={(0,0.75)}]
      \node (A)  at (0,0) {$[\cat{M},GX]$};
      \node (B) at (-3,0) {$[\cat{M},\cat{C}]$};
      \node (C) at (-5,0) {$\cat{C}$};
      \node (D)  at (-1.5,-1.5) {$GX$};
      \draw [->] (B) to [bend left=0]node [auto,labelsize] {$[\cat{M},g]$}(A);
      \draw [->] (C) to [bend left=0]node [auto,labelsize] {$\gm{T}$}(B);
      \draw [->] (D) to [bend right=20]node [auto,swap,labelsize] {$\Gamma_X$}(A);
      \draw [->] (C) to [bend right=20]node [auto,swap,labelsize] {$g$}(D);
      \draw [2cell] (-1.5,-0.1) to node [auto,labelsize]{$\gamma$}(-1.5,-1.1);
        \end{scope}
\end{tikzpicture}\\
=\quad
\begin{tikzpicture}[baseline=-\the\dimexpr\fontdimen22\textfont2\relax ]
        \begin{scope} [shift={(0,0.75)}]  
      \node (E) at (3,0) {$[\cat{M},GX]$};
      \node (F) at (-6.5,0) {$\cat{C}$};
      \node (A)  at (0,0) {$[\cat{M},\Klg]$};
      \node (B) at (-3,0) {$[\cat{M},\cat{C}]$};
      \node (C) at (-5,0) {$\Klg$};
      \node (D)  at (-1.5,-1.5) {$\Klg$};
      \draw [->] (B) to [bend left=0]node [auto,labelsize] {$[\cat{M},f_\gm{T}]$}(A);
      \draw [->] (A) to [bend left=0]node [auto,labelsize] {$[\cat{M},\tilde{g}]$}(E);
      \draw [->] (C) to [bend left=0]node [auto,labelsize] {$u_\gm{T}$}(B);
      \draw [->] (F) to [bend left=0]node [auto,labelsize] {$f_\gm{T}$}(C);
      \draw [->] (D) to [bend right=20]node [auto,swap,labelsize] {$E_{\klact}$}(A);
      \draw [->] (C) to [bend right=20]node [auto,swap,labelsize] {$\id{\Klg}$}(D);
      \draw [2cell] (-1.5,-0.1) to node [auto,labelsize]{$\varepsilon_\gm{T}$}(-1.5,-1.1);
        \end{scope}
\end{tikzpicture}
\end{multline}
and show its uniqueness.

\begin{defn}
The 1-cell~$\tilde{g}\colon
\vcoalgp{[\cat{M},\Klg]}{\klact}{\Klg}\to
\vcoalgp{[\cat{M},GX]}{\Gamma_{X}}{GX}$
of $\EM{\Cat}{[\cat{M},-]}$ is the functor~$
\tilde{g}\colon\Klg\to GX$ defined as
\[
\tilde{g}(m,c)\quad\coloneqq\quad (\Gamma_Xgc) m
\]
on an object~$(m,c)$ and 
\[
\tilde{g}[n,\ m\tensor n\xrightarrow{v}m^\prime,\ c\xrightarrow{f}n\act c^\prime]
\quad\coloneqq\quad
((\Gamma_Xgc^\prime)v)
\circ
(\Gamma_X\gamma_{c^\prime,n})_m
\circ
(\Gamma_Xgf)_m
\]
on a morphism~$
[n,\,m\tensor n\xrightarrow{v}m^\prime,\,c\xrightarrow{f}n\act c^\prime]
\colon (m,c)\to (m^\prime,c^\prime)$.
To check the type of $\tilde{g}[n,v,f]$, observe
\begin{alignat*}{3}
&(\Gamma_Xgc) m&
&\ccol{\ \xrightarrow{(\Gamma_Xgf)_m}\ }&
&(\Gamma_Xg(n\act c^\prime)) m\\
&&
&\ccol{\ =\ }&
&(\Gamma_X(([\cat{M},g]\gm{T}c^\prime) n)) m\\
&&
&\ccol{\ \xrightarrow{(\Gamma_X\gamma_{c^\prime,n})_m}\ }&
&(\Gamma_X((\Gamma_Xgc^\prime) n)) m\\
&&
&\ccol{\ =\ }&
&((M_{GX}\Gamma_Xgc^\prime)n)m
\tag*{by (\ref{eqn-univ-gkl-3})}\\
&&
&\ccol{\ =\ }&
&(\Gamma_Xgc^\prime) (m\tensor n)\\
&&
&\ccol{\ \xrightarrow{(\Gamma_Xgc^\prime)v}\ }&
&(\Gamma_Xgc^\prime) m^\prime.\tag*{\qedhere}
\end{alignat*}
\end{defn}

In order to show the uniqueness of factorization,
we need the following calculational result.

\begin{lem}
\label{lem-decomp-klg}
Every morphism~$
[n,\,m\tensor n\xrightarrow{v}m^\prime,\,c\xrightarrow{f}n\act c^\prime]
\colon (m,c)\to (m^\prime,c^\prime)$
of $\Klg$ can be decomposed as
\begin{multline*}
[n,\,m\tensor n\xrightarrow{v}m^\prime,\,c\xrightarrow{f}n\act c^\prime]\\
=\quad
(m,c)
\xrightarrow{m\klact f_\gm{T} (f)}
(m,n\act c^\prime)
\xrightarrow{m\klact \varepsilon_{\gm{T},(\e,c^\prime),n}}
(m\tensor n,c^\prime)
\xrightarrow{v\klact f_\gm{T} (c^\prime)}
(m^\prime,c^\prime).
\end{multline*}
\end{lem}
\begin{pf}
As equivalence classes, the morphisms on the right hand side are
\begin{align*}
m\klact f_\gm{T} (f)
\ &=\ 
m\klact
[\e,\ 
\e\xrightarrow{\id{}}\e,\ 
c\xrightarrow{f}n\act c^\prime\xrightarrow{\eta_{n\act c^\prime}}\e\act n\act c^\prime]\\
\ &=\ 
[\e,\ 
m\xrightarrow{\id{}}m,\ 
c\xrightarrow{f}n\act c^\prime
\xrightarrow{\eta_{n\act c^\prime}}\e\act n\act c^\prime],\\
m\klact \varepsilon_{\gm{T},(\e,c^\prime),n}
\ &=\ 
m\klact
[n,\ 
n\xrightarrow{\id{}}n,\ 
n\act c^\prime\xrightarrow{\id{}}n\act c^\prime]\\
\ &=\ 
[n,\ 
m\tensor n\xrightarrow{\id{}}m\tensor n,\ 
n\act c^\prime\xrightarrow{\id{}}n\act c^\prime],\\
v\klact f_\gm{T} (c^\prime)
\ &=\ 
v\klact 
[\e,\ 
\e\xrightarrow{\id{}}\e,\ 
c^\prime\xrightarrow{\eta_{c^\prime}}\e\act c^\prime]\\
\ &=\ 
[\e,\ 
m\tensor n\xrightarrow{v}m^\prime,\ 
c^\prime\xrightarrow{\eta_{c^\prime}}\e\act c^\prime].
\end{align*}
The composite of the first two morphisms is
\begin{align*}
\big(m\klact \varepsilon_{\gm{T},(\e,c^\prime),n}\big)
&\circ
\big(m\klact f_\gm{T} (f)\big)\\
&=\ [n,\ 
m\tensor n\xrightarrow{\id{}}m\tensor n,\ 
c\xrightarrow{f}n\act c^\prime
\xrightarrow{\eta_{n\act c^\prime}}\e\act n\act c^\prime
\xrightarrow{\mu_{\e,n,c^\prime}}n\act c^\prime]\\
&=\ [n,\ 
m\tensor n\xrightarrow{\id{}}m\tensor n,\ 
c\xrightarrow{f}n\act c^\prime],
\end{align*}
so finally the composite of the three morphisms is
\begin{align*}
\big(v\klact f_\gm{T} (c^\prime)\big)
\circ
\big(m\klact &\varepsilon_{\gm{T},(\e,c^\prime),n}\big)
\circ
\big(m\klact f_\gm{T} (f)\big)\\
&=\ 
[n,\ 
m\tensor n\xrightarrow{v}m^\prime,\ 
c\xrightarrow{f}n\act c^\prime
\xrightarrow{n\act \eta_{c^\prime}}n\act \e\act c^\prime
\xrightarrow{\mu_{n,\e,c^\prime}}n\act c^\prime]\\
&=\ 
[n,\ 
m\tensor n\xrightarrow{v}m^\prime,\ 
c\xrightarrow{f}n\act c^\prime].\qedhere
\end{align*}
\end{pf}

\begin{prop}
The functor~$\tilde{g}$ defined above is indeed a 1-cell 
of $\EM{\Cat}{[\cat{M},-]}$ which satisfies
(\ref{eqn-univ-gkl-5}) and (\ref{eqn-univ-gkl-6}).
Moreover, it is the unique such.
\end{prop}
\begin{pf}
First observe that thanks to Lemma~\ref{lem-decomp-klg},
the uniqueness is obvious.
Indeed, equation~(\ref{eqn-univ-gkl-5}) determines 
the values of $\tilde{g}$ at $f_\gm{T}(f)$ and
$f_\gm{T}(c^\prime)$, and 
equation~(\ref{eqn-univ-gkl-6}) determines the value of
$\tilde{g}$ at
$\varepsilon_{\gm{T},(\e,c^\prime),n}
=\varepsilon_{\gm{T},f_\gm{T}(c^\prime),n}$.
Finally, the requirement that $\tilde{g}$ is a 1-cell of 
$\EM{\Cat}{[\cat{M},-]}$ enforces the equation~$
\tilde{g}(m\klact -)=(\Gamma_X\tilde{g}(-))m$.

That $\tilde{g}$ indeed satisfies the conditions is
straightforward to check.
$\tilde{g}$ is a 1-cell of $\EM{\Cat}{[\cat{M},-]}$
because, on objects
\begin{align*}
\tilde{g}(l\klact (m,c))
\ &=\ 
\tilde{g}(l\tensor m,c)\\
&=\ 
(\Gamma_Xgc)(l\tensor m)\\
&=\ 
((M_{GX}\Gamma_Xgc)m)l\\
&=\ 
(\Gamma_X((\Gamma_Xgc)m))l \tag*{by (\ref{eqn-univ-gkl-3})}\\
&=\ 
(\Gamma_X\tilde{g}(m,c))l,
\end{align*}
and on morphisms
\begin{align*}
\tilde{g}(u&\klact[n,\,m\tensor n\xrightarrow{v}m^\prime,\,
c\xrightarrow{f}n\act c^\prime])\\
&=\ 
\tilde{g}[n,\,l\tensor m\tensor n\xrightarrow{u\tensor v}l^\prime\tensor m^\prime,\,
c\xrightarrow{f}n\act c^\prime]\\
&=\ 
(\Gamma_Xgc^\prime)(u\tensor v)
\circ
(\Gamma_X\gamma_{c^\prime,n})_{l\tensor m}
\circ
(\Gamma_Xgf)_{l\tensor m}\\
&=\ 
(((M_{GX}\Gamma_Xgc^\prime)v)u)
\circ
(M_{GX}\Gamma_X\gamma_{c^\prime,n})_{m,l}
\circ
(M_{GX}\Gamma_Xgf)_{m,l}\\
&=\ 
\big(
((M_{GX}\Gamma_Xgc^\prime)v)
\circ
(M_{GX}\Gamma_X\gamma_{c^\prime,n})_m
\circ
(M_{GX}\Gamma_Xgf)_m
\big)u\\
&=\ 
\big(
\Gamma_X((\Gamma_Xgc^\prime)v)
\circ
\Gamma_X(\Gamma_X\gamma_{c^\prime,n})_m
\circ
\Gamma_X(\Gamma_Xgf)_m
\big)u
\tag*{by (\ref{eqn-univ-gkl-3})}\\
&=\ 
\big(\Gamma_X
\big(((\Gamma_Xgc^\prime)v)
\circ
(\Gamma_X\gamma_{c^\prime,n})_m
\circ
(\Gamma_Xgf)_m\big)\big)u\\
&=\ 
(\Gamma_X\tilde{g}[n,v,f])u
\end{align*}

$\tilde{g}$ satisfies (\ref{eqn-univ-gkl-5}) since
\begin{align*}
gc
\ &=\ 
H_{GX}\Gamma_Xgc\tag*{by (\ref{eqn-univ-gkl-1})}\\
&=\ 
(\Gamma_Xgc)\e\\
&=\ 
\tilde{g}(\e,c)\\
&=\ 
\tilde{g}f_\gm{T}(c)
\end{align*}
on objects and 
\begin{align*}
gf
\ &=\ 
\gamma_{c^\prime,\e}
\circ
g\eta_{c^\prime}
\circ gf\tag*{by (\ref{eqn-univ-gkl-2})}\\
&=\ 
(H_{GX}\Gamma_X\gamma_{c^\prime,\e})
\circ
(H_{GX}\Gamma_Xg(\eta_{c^\prime}\circ f))\tag*{by (\ref{eqn-univ-gkl-1})}\\
&=\ 
(\Gamma_X\gamma_{c^\prime,\e})_\e
\circ
(\Gamma_Xg(\eta_{c^\prime}\circ f))_{\e}\\
&=\ 
\tilde{g}[\e,\,\e\xrightarrow{\id{}}\e,\,
c\xrightarrow{f}c^\prime\xrightarrow{\eta_{c^\prime}}\e\act c^\prime]\\
&=\ 
\tilde{g}f_\gm{T}(f)
\end{align*}
on morphisms.

Finally, $\tilde{g}$ satisfies (\ref{eqn-univ-gkl-6}) since
\begin{align*}
\gamma_{c,n}
\ &=\ H_{GX}\Gamma_X\gamma_{c,n}\tag*{by (\ref{eqn-univ-gkl-1})}\\
&=\ (\Gamma_X\gamma_{c,n})_\e\\
&=\ \tilde{g}[n,\,n\xrightarrow{\id{}}n,\,
n\act c\xrightarrow{\id{}}n\act c]\\
&=\ \tilde{g}\varepsilon_{\gm{T},(\e,c),n}\\
&=\ ([\cat{M},\tilde{g}]\varepsilon_{\gm{T},f_\gm{T}(c)})_n
\end{align*}
holds.
\end{pf}

\subsubsection{The 2-dimensional aspect}
Suppose we have a morphism of right $\gm{T}$-modules, i.e., 
a 2-cell
    \[
    \begin{tikzpicture}
      \node (A)  at (0,0) {$(\tcat{X},X)$};
      \node (C) at (-4,0) {$(\Cat,\cat{C})$};
      
      \draw [->] (C) to [bend left=20]node (dom)[auto,labelsize] {$(G,g)$}(A);
      \draw [->] (C) to [bend right=20]node (cod)[auto,swap,labelsize] {$(G^\prime,g^\prime)$}(A);

      \draw [2cellnode] (dom) to node [auto,labelsize]{$(\Omega,\omega)$}(cod);

    \end{tikzpicture}
    \]
of $\Emm$ satisfying
\begin{align}
\label{eqn-univ-gkl-7}
    \begin{tikzpicture}[baseline=-\the\dimexpr\fontdimen22\textfont2\relax ]
            \begin{scope} [shift={(0,0.75)}]
      \node (A)  at (0,0) {$\tcat{X}$};
      \node (B) at (-3,-1.5) {$\Cat$};
      \node (C) at (-3,0) {$\Cat$};
      \draw [->] (A) to [bend left=-20]node [auto,swap,labelsize] {$G$}(C);
      \draw [->] (A) to [bend right=-20]node [auto,swap,very near end,labelsize] {$G^\prime$}(C);
      \draw [->] (C) to [bend left=0]node [auto,swap,labelsize] {$[\cat{M},-]$}(B);
      %
      \draw [->] (A) to [bend right=-20]node [auto,labelsize] (Fpp) {$G^\prime$}(B);
      \draw [2cellr] (-1.5,0.25) to node [auto,labelsize]{$\Omega$}(-1.5,-0.25);
      \draw [2cellr] (-1.5,-0.5) to node [auto,labelsize]{$\Gamma^\prime$}(-1.5,-1.1);
            \end{scope}
    \end{tikzpicture}
    \quad=\quad
    \begin{tikzpicture}[baseline=-\the\dimexpr\fontdimen22\textfont2\relax ]
            \begin{scope} [shift={(0,0.75)}]
      \node (A)  at (0,0) {$\tcat{X}$};
      \node (B) at (-3,-1.5) {$\Cat$};
      \node (C) at (-3,0) {$\Cat$};
      \draw [->] (A) to [bend left=-20]node [auto,swap,labelsize] {$G$}(C);
      \draw [->] (A) to [bend left=-20]node [auto,swap,near end,labelsize] {$G$}(B);
      \draw [->] (C) to [bend left=0]node [auto,swap,labelsize] {$[\cat{M},-]$}(B);
      %
      \draw [->] (A) to [bend right=-20]node [auto,labelsize] (Fpp) {$G^\prime$}(B);
      \draw [2cellr] (-1.5,0.25) to node [auto,labelsize]{$\Gamma$}(-1.5,-0.25);
      \draw [2cellr] (-1.5,-0.5) to node [auto,labelsize]{$\Omega$}(-1.5,-1.1);
            \end{scope}
    \end{tikzpicture}
\end{align}
\begin{multline}
\label{eqn-univ-gkl-8}
\begin{tikzpicture}[baseline=-\the\dimexpr\fontdimen22\textfont2\relax ]
        \begin{scope} 
      \node (A)  at (0,0) {$[\cat{M},G^\prime X]$};
      \node (B) at (-3,0) {$[\cat{M},\cat{C}]$};
      \node (C) at (-5,0) {$\cat{C}$};
      \node (D)  at (-1.5,-1.5) {$G^\prime X$};
      \node (F) at (-1,1.5) {$[\cat{M},GX]$};
      \draw [->] (B) to [bend left=0]node [auto,labelsize] {$[\cat{M},g^\prime]$}(A);
      \draw [->] (C) to [bend left=0]node [auto,labelsize] {$\gm{T}$}(B);
      \draw [->] (D) to [bend right=20]node [auto,swap,labelsize] {$\Gamma^\prime _X$}(A);
      \draw [->] (C) to [bend right=20]node [auto,swap,labelsize] {$g^\prime$}(D);
      \draw [->] (A) to [bend right=20]node [auto,swap,labelsize] {$[\cat{M},\Omega_X]$}(F);
      \draw [->] (B) to [bend left=20]node [auto,labelsize] {$[\cat{M},g]$}(F);
      \draw [2cell] (-1.5,-0.1) to node [auto,labelsize]{$\gamma^\prime$}(-1.5,-1.1);
      \draw [2cell] (-1,1.2) to node [auto,swap,labelsize]{$[\cat{M},\omega]$}(-1,0.4);
        \end{scope}
\end{tikzpicture}\\
=\quad
    \begin{tikzpicture}[baseline=-\the\dimexpr\fontdimen22\textfont2\relax ]
            \begin{scope} 
          \node (M)  at (0,0)  {$GX$};
          \node (N)  at (-3.8,1.05) {$[\cat{M},\cat{C}]$};
          \node (RR)  at (-5,0) {$\cat{C}$};
          \node (B)  at (-1.5,1.5) {$[\cat{M},GX]$};
          \node (E) at (-1.5,-1.5) {$G^\prime X$};
          \draw [->] (RR) to node  [auto,labelsize]      {$g$} (M);
          \draw [<-] (B) to [bend left=20]node  [auto,labelsize] {$\Gamma_X$}(M);
          \draw [->] (N) to [bend left=8]node  [auto,labelsize] {$[\cat{M}, g]$}(B);
          \draw [->] (RR) to [bend left=5]node  [auto,labelsize] {$\gm{T}$}(N);
          \draw [<-] (M) to [bend left=20]node  [auto,labelsize] {$\Omega_X$}(E);
          \draw [->] (RR) to [bend right=20]node  [auto,swap,labelsize] {$g^\prime$}(E);
          \draw [2cellnode] (B) to node[auto,labelsize] {$\gamma$}(-1.5,0);
          \draw [2cellnode] (-1.5,0) to node[auto,labelsize] {$\omega$}(E);
            \end{scope}
    \end{tikzpicture}
\end{multline}
Equation~(\ref{eqn-univ-gkl-7}) implies the unique factorization
    \[
    \begin{tikzpicture}[baseline=-\the\dimexpr\fontdimen22\textfont2\relax ]
      \node (E) at (3,0) {$\tcat{X}$};
      \node (A)  at (0,0) {$\Cat$};
      \draw [->] ([yshift=2mm]E.west) to [bend left=-20]node(dom) [auto,swap,labelsize] {${G}$} ([yshift=2mm]A.east);
      \draw [->] ([yshift=-2mm]E.west) to [bend right=-20]node(cod) [auto,labelsize] {${G^\prime}$}([yshift=-2mm]A.east);
      \draw [2cellnoder] (dom) to node [auto,labelsize]{${\Omega}$}(cod);

    \end{tikzpicture}
    \quad=\quad
    \begin{tikzpicture}[baseline=-\the\dimexpr\fontdimen22\textfont2\relax ]
      \node (E) at (3,0) {$\tcat{X}$};
      \node (A)  at (0,0) {$\EM{\Cat}{[\cat{M},-]}$};
      \node (C) at (-3,0) {$\Cat$};
      \draw [->] ([yshift=2mm]E.west) to [bend left=-20]node(dom) [auto,swap,labelsize] {$\tilde{G}$} ([yshift=2mm]A.east);
      \draw [->] ([yshift=-2mm]E.west) to [bend right=-20]node(cod) [auto,labelsize] {$\tilde{G^\prime}$}([yshift=-2mm]A.east);
      \draw [->] (A) to [bend left=0]node [auto,swap,labelsize] {$F_\cat{M}$}(C);
      \draw [2cellnoder] (dom) to node [auto,labelsize]{$\tilde{\Omega}$}(cod);

    \end{tikzpicture}
    \]
where the 2-natural transformation~$\tilde{\Omega}$ is defined as follows:

\begin{defn}
The 2-natural transformation~$\tilde{\Omega}\colon\tilde{G^\prime}\Rightarrow\tilde{G}$
has, as components, 1-cells~$
\tilde{\Omega}_{X^\prime}\colon
\vcoalgp{[\cat{M},G^\prime {X^\prime}]}{\Gamma^\prime_{X^\prime}}{G^\prime {X^\prime}}
\to\vcoalgp{[\cat{M},G{X^\prime}]}{\Gamma_{X^\prime}}{G{X^\prime}}$
of $\EM{\Cat}{[\cat{M},-]}$
given by $\tilde{\Omega}_{X^\prime}\coloneqq\Omega_{X^\prime}$
as functors.
\end{defn}

Next we proceed to construct a 2-cell
    \[
    \begin{tikzpicture}
      \node (a)  at (0,0) {$\vcoalgp{[\cat{M},G{X}]}{\Gamma_{X}}{G{X}}$};
      \node (ap) at (-6,0) {$\vcoalgp{[\cat{M},\Klg]}{\klact}{\Klg}$};
      \node (b)  at (-2,-3) {$\vcoalgp{[\cat{M},G^\prime {X}]}{\Gamma^\prime_{X}}{G^\prime {X}}$};
      
      \draw [->] (ap) to [bend right=0]node [auto,labelsize] {$\tilde{g}$}(a);
      \draw [->] (b)  to [bend right=20]node [auto,swap,labelsize] {$\tilde{\Omega}_X$}(a);
      \draw [->] (ap) to [bend right=20]node [auto,swap,labelsize] {$\tilde{g^\prime}$}(b);
      \draw [2cell] (-2,-0.1) to node [auto,labelsize]{$\tilde{\omega}$}(b);
    \end{tikzpicture}
    \]
of $\EM{\Cat}{[\cat{M},-]}$ satisfying
\begin{align}
\label{eqn-univ-gkl-9}
    \begin{tikzpicture}[baseline=-\the\dimexpr\fontdimen22\textfont2\relax ]
     \begin{scope} [shift={(0,0.75)}]
      \node (a)  at (0,0) {$GX$};
      \node (ap) at (-3,0) {$\cat{C}$};
      \node (b)  at (-1,-1.5) {$G^\prime X$};
      \draw [->] (ap) to [bend right=0]node [auto,labelsize] {$g$}(a);
      \draw [->] (b)  to [bend right=20]node [auto,swap,labelsize] {$\Omega_X$}(a);
      \draw [->] (ap) to [bend right=20]node [auto,swap,labelsize] {$g^\prime$}(b);
      \draw [2cell] (-1,-0.1) to node [auto,labelsize]{$\omega$}(b);
     \end{scope}
    \end{tikzpicture}
   \quad=\quad
    \begin{tikzpicture}[baseline=-\the\dimexpr\fontdimen22\textfont2\relax ]
     \begin{scope} [shift={(0,0.75)}]
      \node (a)  at (0,0) {$GX$};
      \node (ap) at (-3,0) {$\Klg$};
      \node (b)  at (-1,-1.5) {$G^\prime X$};
      \node (c)  at (-5,0) {$\cat{C}$};
      \draw [->] (ap) to [bend right=0]node [auto,labelsize] {$\tilde{g}$}(a);
      \draw [->] (b)  to [bend right=20]node [auto,swap,labelsize] {$\tilde{\Omega}_X$}(a);
      \draw [->] (ap) to [bend right=20]node [auto,swap,labelsize] {$\tilde{g^\prime}$}(b);
      \draw [->] (c) to [bend right=0]node [auto,labelsize] {$f_\gm{T}$}(ap);
      \draw [2cell] (-1,-0.1) to node [auto,labelsize]{$\tilde{\omega}$}(b);
     \end{scope}
    \end{tikzpicture}
\end{align}

\begin{defn}
The 2-cell~$\tilde{\omega}$ of $\EM{\Cat}{[\cat{M},-]}$ has
components given by
\[
\tilde{\omega}_{(m,c)}\quad\colon\quad
(\Gamma_Xgc)m
\xrightarrow{(\Gamma_X\omega_c)_{m}}
(\Gamma_X\Omega_Xg^\prime c)m
\overset{\text{(\ref{eqn-univ-gkl-7})}}=
\Omega_X((\Gamma^\prime_Xg^\prime c)m)
\]
at $(m,c)\in \Klg$.
\end{defn}

\begin{prop}
The natural transformation~$\tilde{\omega}$ defined above is indeed a
2-cell of $\EM{\Cat}{[\cat{M},-]}$ which satisfies
(\ref{eqn-univ-gkl-9}).
Moreover, it is the unique such.
\end{prop}
\begin{pf}
Let us first check that $\tilde{\omega}$ satisfies the conditions.
$\tilde{\omega}$ is a 2-cell of $\EM{\Cat}{[\cat{M},-]}$, i.e.,
$\tilde{\omega}_{l\klact (m,c)}=(\Gamma_X\tilde{\omega}_{(m,c)})_l$
holds, since 
\begin{align*}
\tilde{\omega}_{l\klact (m,c)}\ &=\ 
\tilde{\omega}_{(l\tensor m,c)}\\
&=\ 
(\Gamma_X\omega_c)_{l\tensor m}\\
&=\ 
(M_{GX}\Gamma_X\omega_c)_{m,l}\\
&=\ 
(\Gamma_X(\Gamma_X\omega_c)_m)_l\tag*{by (\ref{eqn-univ-gkl-3})}\\
&=\ 
(\Gamma_X\tilde{\omega}_{(m,c)})_l.
\end{align*}

$\tilde{\omega}$ satisfies (\ref{eqn-univ-gkl-9}) since
\begin{align*}
\omega_c\ 
&=\ H_{GX}\Gamma_X\omega_c\tag*{by (\ref{eqn-univ-gkl-1})}\\
&=\ (\Gamma_X\omega_c)_\e\\
&=\ \tilde{\omega}_{(\e,c)}\\
&=\ \tilde{\omega}_{f_\gm{T}(c)}.
\end{align*}

For uniqueness, observe that equation (\ref{eqn-univ-gkl-9})
determines the $(\e,c)$-component of $\tilde{\omega}$,
and the requirement that $\tilde{\omega}$ is a 2-cell 
of $\EM{\Cat}{[\cat{M},-]}$ implies
$\tilde{\omega}_{(m,c)}=(\Gamma_X\tilde{\omega}_{(\e,c)})_m$,
thus determining all the components.
\end{pf}

\begin{thm}
\label{thm-univ-kl-gm}
The object~$\left(\EM{\Cat}{[\cat{M},-]},
\vcoalg{[\cat{M},\Klg]}{\klact}{\Klg}\right)$ of 
$\Emm$ is the Kleisli object of the
graded monad~$\gm{T}$, considered as a monad
$([\cat{M},-],\gm{T})$ in $\Emm$ on $(\Cat,\cat{C})$.
\end{thm}


\section{The Eilenberg--Moore construction for indexed monads}
\label{sec-iem}
Let $\cat{B}$ and $\cat{C}$ be categories,
and $\im{T}$ a $\cat{B}$-indexed monad on $\cat{C}$.
The Eilenberg--Moore adjunction for $\im{T}$,
which lives in $\Epm$ (cf.~Section~\ref{subsec-im-mnd-epm}), lies above
the co-Eilenberg--Moore adjunction for the 2-comonad~$\cat{B}\times(-)$ on 
$\Cat$.
    \[
    \begin{tikzpicture}
      \node (A)  at (8,0) {$\left(\EM{\Cat}{\cat{B}\times(-)},
            \vcoalg{\cat{B}\times\EMi}{\pi}{\EMi}\right)$};
      \node (B) at (0,0) {$(\tcat{D},D)$};
      \node (C) at (4,2) {$(\Cat,\cat{C})$};
      
      \node at (2,1) [rotate=120,font=\large]{$\vdash$};
      \node at (6,1) [rotate=240,font=\large]{$\vdash$};
      
      \draw [->] (A) to [bend right=17]node [auto,labelsize,swap] {$(U^\cat{B},u^\im{T})$}(C);
      \draw [->] (C) to [bend right=15]node [auto,labelsize,swap] {$(F^\cat{B},f^\im{T})$}(A);
      
      \draw [->] (C) to [bend right=20]node [auto,swap,labelsize] {$(L,l)$}(B);
      \draw [->] (B) to [bend right=20]node [auto,swap,labelsize] {$(R,r)$}(C);
      
      \draw [->,dashed,rounded corners=8pt] (B)--(0,-1) to [bend right=0]node [auto,swap,labelsize]
       {$(K,k)$} (8,-1)--(A);
      
      \draw [rounded corners=15pt,myborder] (-2,-1.5) rectangle (10,2.5);
      \node at (9,2) {$\Epm$};

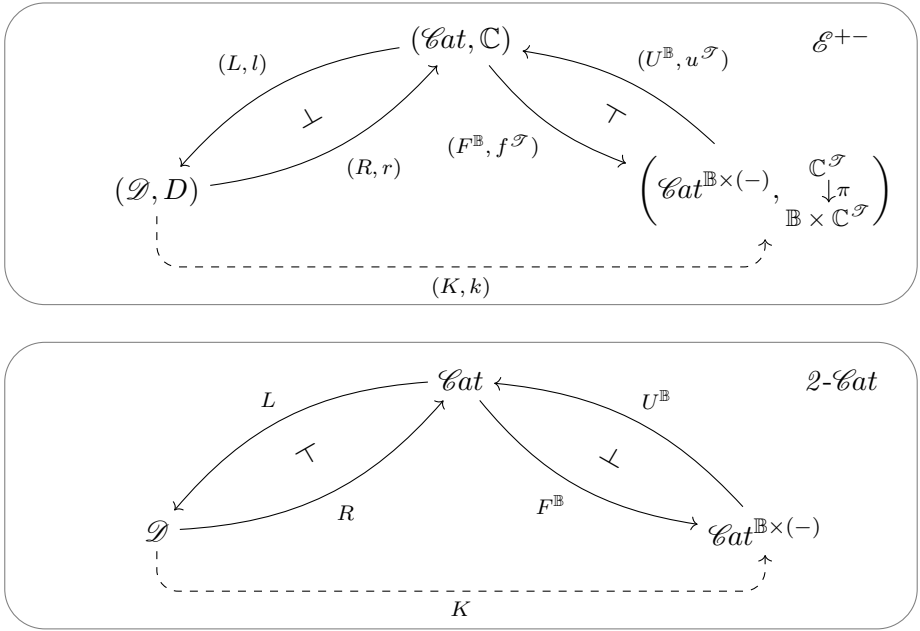
      \begin{scope}[shift={(0,-4.5)}]
       \node (A)  at (8,0) {$\EM{\Cat}{\cat{B}\times(-)}$};
       \node (B) at (0,0) {$\tcat{D}$};
       \node (C) at (4,2) {$\Cat$};
      
       \node at (2,1) [rotate=120,font=\large]{$\dashv$};
       \node at (6,1) [rotate=240,font=\large]{$\dashv$};
      
       \draw [->] (A) to [bend right=21]node [auto,labelsize,swap] {$U^\cat{B}$}(C);
       \draw [->] (C) to [bend right=21]node [auto,labelsize,swap] {$F^\cat{B}$}(A);
      
       \draw [->] (C) to [bend right=23]node [auto,swap,labelsize] (Fpp) {$L$}(B);
       \draw [->] (B) to [bend right=23]node [auto,swap,labelsize] (Fpp) {$R$}(C);
       
       \draw [->,dashed,rounded corners=8pt] (B)--(0,-0.8) to [bend right=0]node [auto,swap,labelsize]
              {$K$} (8,-0.8)--(A);
       \draw [rounded corners=15pt,myborder] (-2,-1.3) rectangle (10,2.5);
             \node at (9,2) {$\TCat$};
   
      \end{scope}
    \end{tikzpicture}
    \]

\subsection{The Eilenberg--Moore category}
Noting that the classical Eilenberg--Moore construction for 
ordinary monads in $\Cat$
defines a functor~$\cat{C}^{(-)}\colon \Mnd{\cat{C}}\to\Cat$,
it follows that from the indexed monad~$
\im{T}\colon \cat{B}^\opp\to\Mnd{\cat{C}}$,
we can obtain an \emph{indexed category}
\[
\cat{C}^{\imc{(-)}}\quad\colon\quad
\cat{B}^\opp\quad
\longrightarrow\quad\Mnd{\cat{C}}
\quad\longrightarrow\quad
\Cat.\]
The Eilenberg--Moore category~$
\EMi$ for $\im{T}$ is given as the total category of the Grothendieck
fibration over $\cat{B}$ corresponding to this indexed category.

\begin{defn}
Define the category $\EMi$ as follows:
\begin{itemize}
\item An object of $\EMi$ is a pair 
$\left(b,\valg{\imc{b}c}{\chi}{c}\right)$,
or more concisely $(b,\chi)$,
where $b$ is an object of $\cat{B}$ and 
$\valgp{\imc{b}c}{\chi}{c}$ is a $\im{T}_b$-algebra
(an object of $\EM{\cat{C}}{\imc{b}}$).
\item A morphism from $\left(b,\valg{\imc{b}c}{\chi}{c}\right)$ to 
$\left(b^\prime,\valg{\imc{b^\prime}c^\prime}{\chi^\prime}{c^\prime}\right)$
is a pair~$(u,h)$
where $u\colon b\to b^\prime$
is a morphism of $\cat{B}$ and 
$h\colon \valgp{\imc{b}c}{\chi}{c}\to
\valgp{\imc{b}c^\prime}{\chi^\prime\circ \imc{u,c^\prime}}{c^\prime}$
is a homomorphism of $\imc{b}$-algebras, i.e., a
morphism~$h\colon c\to c^\prime$ in $\cat{C}$ which makes the diagram
\[
\begin{tikzpicture}
      \node (TL)  at (0,1)  {$\imc{b}c$};
      \node (TR)  at (3,1)  {$\imc{b}c^\prime$};
      \node (BL)  at (0,-1) {$c$};
      \node (BR)  at (3,-1) {$c^\prime$};
      \node (MR)  at (3,0)  {$\imc{b^\prime}c^\prime$};
      
      \draw [->] (TL) to node (T) [auto,labelsize]      {$\imc{b}h$} (TR);
      \draw [->] (TR) to node (R) [auto,labelsize]      {$\imc{u,c^\prime}$} (MR);
      \draw [->] (MR) to node (R) [auto,labelsize]      {$\chi^\prime$} (BR);
      \draw [->] (TL) to node (L) [auto,swap,labelsize] {$\chi$}(BL);
      \draw [->] (BL) to node (B) [auto,swap,labelsize] {$h$}(BR);
\end{tikzpicture}
\]
commute.
\item The identity morphism on $\left(b,\valg{\imc{b}c}{\chi}{c}\right)$
is given by $(\id{b},\id{\chi})$.
\item For two composable morphisms
\begin{align*}
\left(b\xrightarrow{u}b^\prime,\ 
\valgp{\imc{b}c}{\chi}{c}\xrightarrow{h}
\valgp{\imc{b}c^\prime}{\chi^\prime\circ\imc{u,c^\prime}}{c^\prime}\right)
\quad&\colon\quad
(b,\chi)
\quad\longrightarrow\quad
(b^\prime,\chi^\prime),\\
\left(b^\prime\xrightarrow{u^\prime} b^{\prime\prime},\ 
\valgp{\imc{b^\prime} c^\prime}{\chi^\prime}{c^\prime}
\xrightarrow{h^\prime}
\valgp{\imc{b^\prime} c^{\prime\prime}}{\chi^{\prime\prime}\circ
\imc{u^\prime,c^{\prime\prime}} }{c^{\prime\prime}}\right)
\quad&\colon\quad
(b^\prime,\chi^\prime)
\quad\longrightarrow\quad
(b^{\prime\prime},\chi^{\prime\prime}),
\end{align*}
their composite is given by 
\begin{multline*}
\left(b\xrightarrow{u}b^\prime\xrightarrow{u^\prime}b^{\prime\prime},\ 
\valgp{\imc{b}c}{\chi}{c}\xrightarrow{h}
\valgp{\imc{b}c^\prime}{\chi^\prime\circ\imc{u,c^\prime}}{c^\prime}
\xrightarrow{h^\prime}
\valgp{\imc{b} c^{\prime\prime}}{\chi^{\prime\prime}\circ\imc{u^\prime 
\circ u,c^{\prime\prime}}}{c^{\prime\prime}}\right)\\
\colon\quad (b,\chi)\quad\longrightarrow\quad
(b^{\prime\prime},\chi^{\prime\prime}).\tag*{\qedhere}
\end{multline*}
\end{itemize}
\end{defn}

\begin{defn}
Define the functor~$
\pi\colon\EMi\to\cat{B}\times\EMi$
by $(b,\chi)\mapsto (b,(b,\chi))$ and $(u,h)\mapsto (u,(u,h))$.
\end{defn}

\subsection{The Eilenberg--Moore adjunction}

\subsubsection{The left adjoint}

We define the 1-cell
\[
(F^\cat{B},f^\im{T})\quad\colon\quad 
(\Cat,\cat{C})
\quad\longrightarrow\quad
\left(\EM{\Cat}{\cat{B}\times(-)},
\vcoalg{\cat{B}\times\EMi}{\pi}{\EMi}
\right)
\] 
of $\Epm$ as follows:

\begin{defn}
\label{defn-F-B-times}
The 2-functor~$F^\cat{B}\colon \Cat\to\EM{\Cat}{\cat{B\times(-)}}$
is the cofree 2-functor~$
\cat{X}\mapsto\vcoalgp{\cat{B}\times\cat{B}\times\cat{X}}{M^\cat{B}_\cat{X}}{\cat{B}\times\cat{X}}$
where $M^\cat{B}_\cat{X}=M_\cat{X}$
is the one defined in Definition~\ref{defn-B-times}.
\end{defn}

\begin{defn}
The 1-cell~$
f^\im{T}\colon 
\vcoalgp{\cat{B}\times\cat{B}\times\cat{C}}{M_\cat{C}}{\cat{B}\times\cat{C}}
\to
\vcoalgp{\cat{B}\times\EMi}{\pi}{\EMi}$
of $\EM{\Cat}{\cat{B}\times(-)}$
is the functor~$f^\im{T}\colon \cat{B}\times\cat{C}\to\EMi$ defined as
\[
f^\im{T}(b,c)\quad\coloneqq\quad\left(b,
\valg{\imc{b}\imc{b}c}{\mu_{b,c}}{\imc{b}c}\right)
\]
on an object~$(b,c)$ and 
\[
f^\im{T}(u,f)\quad\coloneqq\quad
\left(b\xrightarrow{u}b^\prime,\ 
\valgp{\imc{b}\imc{b}c}{\mu_{b,c}}{\imc{b}c}
\xrightarrow{\imc{u,c^\prime}\circ\imc{b}f}
\valgp{\imc{b}\imc{b^\prime}c^\prime}
{\mu_{b^\prime,c^\prime}\circ \imc{u,\imc{b^\prime}c^\prime}}
{\imc{b^\prime}c^\prime}\right)
\]
on a morphism~$(u,f)\colon (b,c)\to (b^\prime,c^\prime)$.
\end{defn}

\subsubsection{The right adjoint}

We define the 1-cell
\[
(U^\cat{B},u^\im{T})\quad\colon\quad 
\left(\EM{\Cat}{\cat{B}\times(-)},
\vcoalg{\cat{B}\times\EMi}{\pi}{\EMi}
\right)
\quad\longrightarrow\quad
(\Cat,\cat{C})
\] 
of $\Epm$ as follows:
\begin{defn}
\label{defn-U-B-times}
The 2-functor~$U^\cat{B}\colon \EM{\Cat}{\cat{B}\times(-)}\to\Cat$ is 
the forgetful 2-functor $
\vcoalgp{\cat{B}\times\cat{A}}{\alpha}{\cat{A}}\mapsto\cat{A}$.
\end{defn}

\begin{defn}
The functor~$u^\im{T}\colon \EMi\to\cat{C}$
is also given as the forgetful functor
$\left(b,\valgp{\imc{b}c}{\chi}{c}\right)\mapsto c$ and 
$(u,h)\mapsto h$.
\end{defn}

\subsubsection{The unit}

We define the 2-cell
    \[
    \begin{tikzpicture}
      \node (A)  at (0,0) {$(\Cat,\cat{C})$};
      \node (B) at (8,0) {$(\Cat,\cat{C})$};
      \node (C) at (4,-2) {$\left(\EM{\Cat}{\cat{B}\times(-)},
      \vcoalg{\cat{B}\times\EMi}{\pi}{\EMi}\right)$};
      
      \draw [->] (A) to [bend right=20]node [auto,swap,labelsize] {$(F^\cat{B},f^\im{T})$}(C);
      \draw [->] (C) to [bend right=20]node [auto,swap,labelsize] {$(U^\cat{B},u^\im{T})$}(B);

      \draw [->] (A) to [bend left=20]node [auto,labelsize] (Fpp) {$(\id{\Cat},\id{\cat{C}})$}(B);
      \draw [2cellnoder] (C) to node [auto,labelsize,swap]{$(H^\cat{B},\eta^\im{T})$}(Fpp);

    \end{tikzpicture}
    \]
of $\Epm$ as follows:

\begin{defn}
The 2-natural transformation~$
H^\cat{B}\colon U^\cat{B}\circ F^\cat{B}\Rightarrow\id{\Cat}$
has components~$
H^\cat{B}_\cat{X}=H_\cat{X}$
defined in Definition~\ref{defn-B-times}.
\end{defn}

\begin{defn}
The natural transformation
    \[
    \begin{tikzpicture}
      \node (L) at (0,0) {$\cat{B}\times\cat{C}$};
      \node (M) at (3,0) {$\EMi$};
      \node (R) at (6,0) {$\cat{C}$};
      \node (T) at (2,2) {$\cat{C}$};
      
      \draw [->] (L) to node [auto,swap,labelsize] {$f^\im{T}$}(M);
      \draw [->] (M) to node [auto,swap,labelsize] {$u^\im{T}$}(R);
      \draw [->] (L) to [bend left=20]node [auto,labelsize] {$H_\cat{C}$}(T);
      \draw [->] (T) to [bend left=20]node [auto,labelsize] {$\id{\cat{C}}$}(R);
            
      \draw [2cell] (T) to node [auto,labelsize]{$\eta^\im{T}$}(2,0.1);
    \end{tikzpicture}
    \]
has components~$\eta^{\im{T}}_{b,c}\colon c\to\imc{b}c$ 
given by the data of the indexed monad.
\end{defn}

\subsubsection{The counit}

We define the 2-cell
    \[
    \begin{tikzpicture}
      \node (A)  at (0,0) {$\left(\EM{\Cat}{\cat{B}\times(-)},
            \vcoalg{\cat{B}\times\EMi}{\pi}{\EMi}\right)$};
      \node (B) at (8,0) {$\left(\EM{\Cat}{\cat{B}\times(-)},
            \vcoalg{\cat{B}\times\EMi}{\pi}{\EMi}\right)$};
      \node (C) at (4,2) {$(\Cat,\cat{C})$};
      
      \draw [->] (A) to [bend left=20]node [auto,labelsize] {$(U^\cat{B},u^\im{T})$}(C);
      \draw [->] (C) to [bend left=20]node [auto,labelsize] {$(F^\cat{B},f^\im{T})$}(B);

      \draw [->] (A) to [bend right=20]node [auto,swap,labelsize] (Fpp) {$(\id{\EM{\Cat}{\cat{B}\times(-)}},\id{\pi})$}(B);
      \draw [2cellnoder] (Fpp) to node [auto,swap,labelsize]{$(E^\cat{B},\varepsilon^\im{T})$}(C);

    \end{tikzpicture}
    \]
of $\Epm$ as follows:

\begin{defn}
The 2-natural transformation~$E^\cat{B}\colon \id{\EM{\Cat}{\cat{B}\times(-)}}
\Rightarrow F^\cat{B}\circ U^\cat{B}$
has components~$
E^\cat{B}_\alpha\colon 
\vcoalgp{\cat{B}\times\cat{A}}{\alpha}{\cat{A}}
\to
\vcoalgp{\cat{B}\times\cat{B}\times\cat{A}}{M_\cat{A}}{\cat{B}\times\cat{A}}$
given by $E^\cat{B}_\alpha=E_\alpha\coloneqq\alpha\colon \cat{A}\to\cat{B}\times\cat{A}$.
\end{defn}

\begin{defn}
The 2-cell
    \[
    \begin{tikzpicture}
      \node (L) at (0,0) {$\vcoalgp{\cat{B}\times\EMi}{\pi}{\EMi}$};
      \node (R) at (-9,0) {$\vcoalgp{\cat{B}\times\EMi}{\pi}{\EMi}$};
      \node (T) at (-6,3) {$\vcoalgp{\cat{B}\times\cat{B}\times\EMi}{M_\EMi}{\cat{B}\times\EMi}$};
      
      \node (Mm) at (-2.5,2.2) {$\vcoalgp{\cat{B}\times\cat{B}\times\cat{C}}{M_\cat{C}}{\cat{B}\times\cat{C}}$};
      
      \draw [<-] (L) to [bend right=5]node [auto,swap,labelsize] {$f^\im{T}$}(Mm);
      \draw [<-] (Mm) to [bend right=5]node [auto,swap,labelsize,near end] {$\cat{B}\times u^\im{T}$}(T);

      \draw [<-] (L) to node [auto,labelsize] {$\id{\pi}$}(R);
      \draw [<-] (T) to [bend right=20]node [auto,swap,labelsize] {$E_\pi$}(R);
            
      \draw [2cell] (T) to node [auto,labelsize]{$\varepsilon^\im{T}$}(-6,0.1);
    \end{tikzpicture}
    \]
of $\EM{\Cat}{\cat{B}\times(-)}$ is the natural transformation
with the component at 
$\left(b,\valg{\imc{b}c}{\chi}{c}\right)\in\EMi$,
$\varepsilon^{\im{T}}_{(b,\chi)}\colon 
\left(b, \valg{\imc{b}\imc{b}c}{\mu_{b,c}}{\imc{b}c}\right)
\to
\left(b, \valg{\imc{b}c}{\chi}{c}\right)$,
given by
$\varepsilon^{\im{T}}_{(b,\chi)}\coloneqq (\id{b},\chi)$.
\end{defn}

\subsection{Comparison maps}

Suppose we have an adjunction~$
(L,l)\dashv(R,r)\colon(\tcat{D},D)\to(\Cat,\cat{C})$ in $\Epm$
with unit~$(H,\eta)\colon$ $(\id{\Cat},\id{\cat{C}})\Rightarrow(R,r)\circ(L,l)$ 
and counit~$(E,\varepsilon)\colon(L,l)\circ(R,r)\Rightarrow(\id{\tcat{D}},\id{D})$,
which gives a resolution of the monad $(\cat{B}\times(-),\im{T})$,
i.e., such that the following equations hold:
{\allowdisplaybreaks
\begin{align}
\label{eqn-comp-iem-1}
\Cat\xrightarrow{\cat{B}\times(-)}\Cat\quad
&=
\quad\Cat\xrightarrow{L}\tcat{D}\xrightarrow{R}\Cat\\
\label{eqn-comp-iem-2}
\cat{B}\times\cat{C}\xrightarrow{\im{T}}\cat{C}\quad&=\quad
RL\cat{C}\xrightarrow{Rl}RD\xrightarrow{r}\cat{C}\\
\label{eqn-comp-iem-3}
H^\cat{B}\quad&=\quad H\\
\label{eqn-comp-iem-4}
\eta^\im{T}\quad&=\quad\eta\\
\label{eqn-comp-iem-5}
\begin{tikzpicture}[baseline=-\the\dimexpr\fontdimen22\textfont2\relax ]
        \begin{scope} 
      \node (TL)  at (0,0)  {$\Cat$};
      \node (TT)  at (1.2,0.7) {$\Cat$};
      \node (RR)  at (2.4,0) {$\Cat$};
      \draw [->] (TL) to [bend left=20]node (L) [auto,labelsize] {$\cat{B}\times(-)$}(TT);
      \draw [->] (TT) to [bend left=20]node  [auto,labelsize] {$\cat{B}\times(-)$}(RR);
      \draw [->] (TL) to [bend right=20]node (B) [auto,swap,labelsize] {$\cat{B}\times(-)$} (RR);
      \draw [2cellnoder] (TT) to node[auto,labelsize] {$M^\cat{B}$}(B);
        \end{scope}
\end{tikzpicture}
\quad&=\quad
\begin{tikzpicture}[baseline=-\the\dimexpr\fontdimen22\textfont2\relax ]
        \begin{scope} 
      \node (TL)  at (0,0)  {$\Cat$};
      \node (TR)  at (1.2,0)  {$\tcat{D}$};
      \node (TT)  at (2.4,0.7) {$\Cat$};
      \node (BR)  at (3.6,0) {$\tcat{D}$};
      \node (RR)  at (4.8,0) {$\Cat$};
      \draw [->] (TL) to node (T) [auto,labelsize]      {$L$} (TR);
      \draw [->] (TR) to [bend left=20]node (R) [auto,labelsize]      {$R$} (TT);
      \draw [->] (TT) to [bend left=20]node (L) [auto,labelsize] {$L$}(BR);
      \draw [->] (BR) to node  [auto,labelsize] {$R$}(RR);
      \draw [->] (TR) to [bend right=20]node (B) [auto,swap,labelsize] {$\id{\tcat{D}}$} (BR);
      \draw [2cellnoder] (TT) to node[auto,labelsize] {$E$}(B);
        \end{scope}
\end{tikzpicture}
\end{align}
}%
\vspace{-20pt}
\begin{multline}
\label{eqn-comp-iem-6}
\begin{tikzpicture}[baseline=-\the\dimexpr\fontdimen22\textfont2\relax ]
        \begin{scope} 
      \node (L)  at (0,0)  {$\cat{B}\times\cat{C}$};
      \node (R)  at (6,0) {$\cat{C}$};
      \node (B)  at (2,2) {$\cat{B}\times\cat{B}\times\cat{C}$};
      \node (N) at (4.7,1.3) {$\cat{B}\times\cat{C}$};
      %
      \draw [->] (L) to node  [auto,swap,labelsize]      {$\im{T}$} (R);
      \draw [->] (L) to [bend left=20]node  [auto,labelsize] {$M_\cat{C}$}(B);
      \draw [->] (B) to [bend left=7]node  [auto,labelsize] {$\cat{B}\times\im{T}$}(N);
      \draw [->] (N) to [bend left=5]node  [auto,labelsize] {$\im{T}$}(R);
      \draw [2cellnode] (B) to node[auto,labelsize] {$\mu^\im{T}$}(2,0);
        \end{scope}
\end{tikzpicture}\\
=\quad
\begin{tikzpicture}[baseline=-\the\dimexpr\fontdimen22\textfont2\relax ]
        \begin{scope} [shift={(0,-0.75)}]
      \node (L)  at (0.75,0)  {$RL\cat{C}$};
      \node (M)  at (2.75,0)  {$RD$};
      \node (N)  at (6.5,1.2) {$RL\cat{C}$};
      \node (RR)  at (7.8,0) {$RD$};
      \node (RRR)  at (9.4,0) {$\cat{C}$};
      \node (B)  at (4.2,1.8) {$RLRD$};
      \draw [->] (L) to node  [auto,labelsize]      {$Rl$} (M);
      \draw [->] (M) to node  [auto,swap,labelsize]      {$\id{RD}$} (RR);
      \draw [->] (RR) to node  [auto,labelsize]      {$r$} (RRR);
      \draw [->] (M) to [bend left=20]node  [auto,labelsize] {$RE_D$}(B);
      \draw [->] (B) to [bend left=8]node  [auto,labelsize] {$RLr$}(N);
      \draw [->] (N) to [bend left=5]node  [auto,labelsize] {$Rl$}(RR);
      \draw [2cellnode] (B) to node[auto,labelsize] {$R\varepsilon$}(4.2,0);
        \end{scope}
\end{tikzpicture}
\end{multline}

The equations~(\ref{eqn-comp-iem-1}), (\ref{eqn-comp-iem-3}) and
(\ref{eqn-comp-iem-5}) imply the existence of the comparison
2-functor~$K$:

\begin{defn}
The 2-functor~$K\colon\tcat{D}\to\EM{\Cat}{\cat{B}\times(-)}$
is the comparison 2-functor
$D^\prime\mapsto \vcoalgp{\cat{B}\times RD^\prime}{RE_{D^\prime}}{RD^\prime}$.
\end{defn}

Before constructing the 1-cell~$k$, we introduce a notation 
to describe the structure maps of coalgebras for the 
2-comonad~$\cat{B}\times(-)$.
Given an object~$\vcoalgp{\cat{B}\times\cat{A}}{\alpha}{\cat{A}}$
of $\EM{\Cat}{\cat{B}\times(-)}$, 
we will write $\alpha =(\alpha_0 (-), \alpha_1 (-))$.
Note that by one of the axioms for $\cat{B}\times(-)$-coalgebras 
saying that 
\[
\begin{tikzpicture}[baseline=-\the\dimexpr\fontdimen22\textfont2\relax ]
      \node (TL)  at (0,0.75)  {$\cat{A}$};
      \node (TR)  at (2,0.75)  {$\cat{B}\times\cat{A}$};
      \node (BR)  at (2,-0.75) {$\cat{A}$};
      
      \draw [->] (TL) to node (T) [auto,labelsize]      {$\alpha$} (TR);
      \draw [->] (TR) to node (R) [auto,labelsize]      {$H_\cat{A}$} (BR);
      \draw [->] (TL) to node (L) [auto,swap,labelsize] {$\id{\cat{A}}$}(BR);
\end{tikzpicture}
\]
commutes, $\alpha_1=\id{\cat{A}}$ holds;
thus the only meaningful data of $\alpha$ is 
actually $\alpha_0\colon\cat{A}\to\cat{B}$, which is subject to
no axioms.
In fact, there is an isomorphism of 2-categories~$
\EM{\Cat}{\cat{B}\times(-)}\cong\Cat/\cat{B}$
given by $\vcoalgp{\cat{B}\times\cat{A}}{\alpha}{\cat{A}}\mapsto\vcoalgp{\cat{B}}{\alpha_0}{\cat{A}}$.

\begin{defn}
The 1-cell~$
k\colon\vcoalgp{\cat{B}\times RD}{RE_{D}}{RD}
\to
\vcoalgp{\cat{B}\times\EMi}{\pi}{\EMi}$
of $\EM{\Cat}{\cat{B}\times(-)}$
is defined as 
\[
k(d)\quad\coloneqq\quad\left(RE_{D,0}d,\valg{\imc{RE_{D,0}d}rd}{rR\varepsilon_d}{rd}\right)
\]
on an object~$d$. 
To check the type of the structure map, observe
\begin{alignat*}{3}
&\imc{RE_{D,0}d}rd&
&\ccol{\ =\ }&
&rRl(RE_{D,0}d,rd)\\
&&
&\ccol{\ =\ }&
&rRlRLr(RE_{D,0}d,d)\\
&&
&\ccol{\ =\ }&
&rRlRLrRE_Dd\\
&&
&\ccol{\ \xrightarrow{rR\varepsilon_d}\ }&
&rd.
\end{alignat*}
$k$ is defined as 
\begin{multline*}
k(w)\quad\coloneqq\quad
\left(RE_{D,0}d\xrightarrow{RE_{D,0}w}RE_{D,0}d^\prime,
{\vphantom{\valgp{R}{r}{r}}}\right.\\
\left. \valgp{\imc{RE_{D,0}d}rd}{rR\varepsilon_d}{rd}\xrightarrow{rw}
\valgp{\imc{RE_{D,0}d}rd^\prime}{rR\varepsilon_{d^\prime}\circ\imc{RE_{D,0}w,rd^\prime}}{rd^\prime}\right)
\end{multline*}
on a morphism~$w\colon d\to d^\prime$.
\end{defn}

\begin{prop}
The 1-cell~$(K,k)$ of $\Epm$ satisfies the equations~$
(K,k)\circ(L,l)=(F^\cat{B},f^\im{T})$ and
$(R,r)=(U^\cat{B},u^\im{T})\circ(K,k)$.
Moreover, it is the unique such.
\end{prop}

\subsection{The 2-dimensional universality}
\subsubsection{Statement of the theorem}
We will show that there is a family of
isomorphisms of categories
\begin{multline*}
  \Epm\left((\tcat{X},X),\ 
  \left(\EM{\Cat}{\cat{B}\times(-)},
              \vcoalg{\cat{B}\times\EMi}{\pi}{\EMi}\right)\right)\\
  \cong\quad 
  \Epm\big((\tcat{X},X),\ (\Cat,\cat{C})\big)^{\Epm((\tcat{X},X),\ (\cat{B}\times(-),\im{T}))}
\end{multline*}
2-natural in $(\tcat{X},X)\in \Epm$,
by showing that the data 
    \[
    \begin{tikzpicture}
      \node (A)  at (0,0) {$\left(\EM{\Cat}{\cat{B}\times(-)},
                    \vcoalg{\cat{B}\times\EMi}{\pi}{\EMi}\right)$};
      \node (B) at (4,-2) {$\left(\EM{\Cat}{\cat{B}\times(-)},
                    \vcoalg{\cat{B}\times\EMi}{\pi}{\EMi}\right)$};
      \node (C) at (4,0) {$(\Cat,\cat{C})$};
      \node (D) at (4,-4) {$(\Cat,\cat{C})$};
      
      \draw [->] (A) to [bend left=0]node [auto,labelsize] {$(U^\cat{B},u^\im{T})$}(C);
      \draw [->] (C) to [bend left=0]node [auto,labelsize] {$(F^\cat{B},f^\im{T})$}(B);
      \draw [->] (B) to [bend left=0]node [auto,labelsize] {$(U^\cat{B},u^\im{T})$}(D);

      \draw [->] (A) to [bend right=20]node [auto,swap,labelsize] (Fpp) {$(\id{\EM{\Cat}{\cat{B}\times(-)}},\id{\pi})$}(B);
      \draw [2cell] (2,-0.1) to node [auto,labelsize]{$(E^\cat{B},\varepsilon^\im{T})$}(2,-1.6);

    \end{tikzpicture}
    \]
provides the universal left $\im{T}$-module.

\subsubsection{The 1-dimensional aspect}
Suppose we have a left $\im{T}$-module, i.e., a diagram
    \[
    \begin{tikzpicture}
      \node (A)  at (0,0) {$(\tcat{X},X)$};
      \node (B) at (3,-1.5) {$(\Cat,\cat{C})$};
      \node (C) at (3,0) {$(\Cat,\cat{C})$};
      
      \draw [->] (A) to [bend left=0]node [auto,labelsize] {$(G,g)$}(C);
      \draw [->] (C) to [bend left=0]node [auto,labelsize] {$(\cat{B}\times(-),\im{T})$}(B);

      \draw [->] (A) to [bend right=20]node [auto,swap,labelsize] (Fpp) {$(G,g)$}(B);
      \draw [2cell] (1.5,-0.1) to node [auto,labelsize]{$(\Gamma,\gamma)$}(1.5,-1.1);

    \end{tikzpicture}
    \]
in $\Epm$ satisfying
{\allowdisplaybreaks
\begin{align}
\label{eqn-univ-iem-1}
\begin{tikzpicture}[baseline=-\the\dimexpr\fontdimen22\textfont2\relax ]
        \begin{scope} [shift={(0,0.75)}]
      \node (A)  at (0,0) {$\tcat{X}$};
      \node (B) at (3,-1.5) {$\Cat$};
      \node (C) at (3,0) {$\Cat$};
      \draw [->] (A) to [bend left=0]node [auto,labelsize] {$G$}(C);
      \draw [->] (C) to [bend left=0]node (t)[auto,swap,labelsize] {$\cat{B}\times(-)$}(B);
      \draw [->] (C) to [bend left=90]node(id) [auto,labelsize] {$\id{\Cat}$}(B);
      %
      \draw [->] (A) to [bend right=20]node [auto,swap,labelsize] (Fpp) {$G$}(B);
      \draw [2cellr] (1.5,-0.1) to node [auto,swap,labelsize]{$\Gamma$}(1.5,-1.1);
      \draw [2cellnoder] (id) to node [auto,swap,labelsize]{$H^\cat{B}$}(t);
        \end{scope}
\end{tikzpicture}
\quad&=\quad
\begin{tikzpicture}[baseline=-\the\dimexpr\fontdimen22\textfont2\relax ]
        \begin{scope} 
       \node (A)  at (0,0) {$\tcat{X}$};
       \node (C) at (2,0) {$\Cat$};
       \draw [->] (A) to [bend left=0]node [auto,labelsize] {$G$}(C);
        \end{scope}
\end{tikzpicture}\\
\label{eqn-univ-iem-2}
\begin{tikzpicture}[baseline=-\the\dimexpr\fontdimen22\textfont2\relax ]
        \begin{scope} 
      \node (A)  at (0,0) {$\cat{B}\times GX$};
      \node (B) at (3,0) {$\cat{B}\times\cat{C}$};
      \node (C) at (5,0) {$\cat{C}$};
      \node (D)  at (1.5,-1.5) {$GX$};
      \node (E) at (3.5,1.5) {$\cat{C}$};
      \node (F) at (0.5,1.5) {$GX$};
      \draw [->] (A) to [bend left=0]node [auto,labelsize] {$\cat{B}\times g$}(B);
      \draw [->] (B) to [bend left=0]node [auto,labelsize] {$\im{T}$}(C);
      \draw [<-] (A) to [bend right=20]node [auto,swap,labelsize] {$\Gamma_X$}(D);
      \draw [->] (D) to [bend right=20]node [auto,swap,labelsize] {$g$}(C);
      \draw [<-] (E) to [bend right=20]node [auto,swap,labelsize] {$H_\cat{C}$}(B);
      \draw [->] (F) to [bend right=0]node [auto,labelsize] {$g$}(E);
      \draw [->] (E) to [bend left=20]node [auto,labelsize] {$\id{\cat{C}}$}(C);
      \draw [<-] (F) to [bend right=20]node [auto,swap,labelsize] {$H_{GX}$}(A);
      \draw [2cell] (1.5,-0.1) to node [auto,labelsize]{$\gamma$}(1.5,-1.1);
      \draw [2cell] (3.5,1.2) to node [auto,labelsize]{$\eta$}(3.5,0.3);
        \end{scope}
\end{tikzpicture}
\quad&=\quad
\begin{tikzpicture}[baseline=-\the\dimexpr\fontdimen22\textfont2\relax ]
        \begin{scope} 
       \node (A)  at (0,0) {$G{X}$};
       \node (C) at (2,0) {$\cat{C}$};
       \draw [->] (A) to [bend left=0]node [auto,labelsize] {$g$}(C);
        \end{scope}
\end{tikzpicture}\\
\label{eqn-univ-iem-3}
\begin{tikzpicture}[baseline=-\the\dimexpr\fontdimen22\textfont2\relax ]
        \begin{scope} 
      \node (A)  at (0,0) {$\tcat{X}$};
      \node (B) at (3,-1.5) {$\Cat$};
      \node (C) at (4.5,0) {$\Cat$};
      \node (D) at (3,1.5) {$\Cat$};
      \draw [->] (D) to [bend left=0]node [auto,swap,labelsize] {$\cat{B}\times(-)$}(B);
      \draw [->] (D) to [bend left=20]node (t)[auto,labelsize] {$\cat{B}\times(-)$}(C);
      \draw [->] (C) to [bend left=20]node (t)[auto,labelsize] {$\cat{B}\times(-)$}(B);
      \draw [->] (A) to [bend left=20]node [auto,labelsize] (Fpp) {$G$}(D);
      \draw [->] (A) to [bend right=20]node [auto,swap,labelsize] (Fpp) {$G$}(B);
      \draw [2cellr] (1.5,1) to node [auto,swap,labelsize]{$\Gamma$}(1.5,-1);
      \draw [2cellr] (4.1,0) to node [auto,swap,labelsize]{$M^\cat{B}$}(3.1,0);
        \end{scope}
\end{tikzpicture}
\quad&=\quad 
\begin{tikzpicture}[baseline=-\the\dimexpr\fontdimen22\textfont2\relax ]
        \begin{scope} 
      \node (A)  at (0,0) {$\tcat{X}$};
      \node (B) at (3,-1.5) {$\Cat$};
      \node (C) at (4,0) {$\Cat$};
      \node (D) at (3,1.5) {$\Cat$};
      \draw [->] (A) to [bend left=0]node [auto,labelsize] {$G$}(C);
      \draw [->] (D) to [bend left=20]node (t)[auto,labelsize] {$\cat{B}\times(-)$}(C);
      \draw [->] (C) to [bend left=20]node (t)[auto,labelsize] {$\cat{B}\times(-)$}(B);
      \draw [->] (A) to [bend left=20]node [auto,labelsize] (Fpp) {$G$}(D);
      \draw [->] (A) to [bend right=20]node [auto,swap,labelsize] (Fpp) {$G$}(B);
      \draw [2cellr] (3,1.2) to node [auto,labelsize]{$\Gamma$}(3,0.1);
      \draw [2cellr] (3,-0.1) to node [auto,labelsize]{$\Gamma$}(3,-1.2);
        \end{scope}
\end{tikzpicture}
\end{align}
}%
\vspace{-20pt}
\begin{multline}
\label{eqn-univ-iem-4}
\begin{tikzpicture}[baseline=-\the\dimexpr\fontdimen22\textfont2\relax ]
        \begin{scope} 
      \node (A)  at (0,0) {$\cat{B}\times GX$};
      \node (B) at (4,0) {$\cat{B}\times\cat{C}$};
      \node (C) at (7,0) {$\cat{C}$};
      \node (D)  at (2,-1.5) {$GX$};
      \node (E) at (4,1.5) {$\cat{B}\times\cat{B}\times\cat{C}$};
      \node (F) at (0,1.5) {$\cat{B}\times\cat{B}\times GX$};
      \node (G) at (7,1.5) {$\cat{B}\times\cat{C}$};
      \draw [->] (A) to [bend left=0]node [auto,labelsize] {$\cat{B}\times g$}(B);
      \draw [->] (B) to [bend left=0]node [auto,swap,labelsize] {$\im{T}$}(C);
      \draw [<-] (A) to [bend right=20]node [auto,swap,labelsize] {$\Gamma_X$}(D);
      \draw [->] (D) to [bend right=15]node [auto,swap,labelsize] {$g$}(C);
      \draw [<-] (E) to [bend right=0]node [auto,swap,labelsize] {$M_\cat{C}$}(B);
      \draw [->] (G) to [bend right=0]node [auto,labelsize] {$\im{T}$}(C);
      \draw [->] (E) to [bend right=0]node [auto,labelsize] {$\cat{B}\times \im{T}$}(G);
      \draw [->] (F) to [bend right=0]node [auto,labelsize] {$\cat{B}\times\cat{B}\times g$}(E);
      \draw [<-] (F) to [bend right=0]node [auto,swap,labelsize] {$M_{GX}$}(A);
      \draw [2cell] (2,-0.1) to node [auto,labelsize]{$\gamma$}(2,-1.1);
      \draw [2cell] (5.5,1.4) to node [auto,labelsize]{$\mu$}(5.5,0.1);
        \end{scope}
\end{tikzpicture}\\
=\quad
\begin{tikzpicture}[baseline=-\the\dimexpr\fontdimen22\textfont2\relax ]
        \begin{scope} 
      \node (A)  at (0,0) {$\cat{B}\times GX$};
      \node (B) at (4,0) {$\cat{B}\times\cat{C}$};
      \node (C) at (7,0) {$\cat{C}$};
      \node (D)  at (2,-1.5) {$GX$};
      \node (E) at (4,1.5) {$\cat{B}\times\cat{B}\times\cat{C}$};
      \node (F) at (0,1.5) {$\cat{B}\times\cat{B}\times GX$};
      \draw [->] (A) to [bend left=0]node [auto,labelsize] {$\cat{B}\times g$}(B);
      \draw [->] (B) to [bend left=0]node [auto,labelsize] {$\im{T}$}(C);
      \draw [<-] (A) to [bend right=20]node [auto,swap,labelsize] {$\Gamma_X$}(D);
      \draw [->] (D) to [bend right=15]node [auto,swap,labelsize] {$g$}(C);
      \draw [->] (E) to [bend right=0]node [auto,labelsize] {$\cat{B}\times\im{T}$}(B);
      \draw [->] (F) to [bend right=0]node [auto,labelsize] {$\cat{B}\times\cat{B}\times g$}(E);
      \draw [<-] (F) to [bend right=0]node [auto,swap,labelsize] {$\cat{B}\times\Gamma_{X}$}(A);
      \draw [2cell] (2,-0.1) to node [auto,labelsize]{$\gamma$}(2,-1.1);
      \draw [2cell] (2,1.4) to node [auto,labelsize]{$\cat{B}\times\gamma$}(2,0.5);
        \end{scope}
\end{tikzpicture}
\end{multline}

Equations~(\ref{eqn-univ-iem-1}) and (\ref{eqn-univ-iem-3})
imply the unique factorization
    \[
    \begin{tikzpicture}[baseline=-\the\dimexpr\fontdimen22\textfont2\relax ]
            \begin{scope} [shift={(0,0.75)}]
      \node (A)  at (0,0) {$\tcat{X}$};
      \node (B) at (3,-1.5) {$\Cat$};
      \node (C) at (3,0) {$\Cat$};
          
      \draw [->] (A) to [bend left=0]node [auto,labelsize] {$G$}(C);
      \draw [->] (C) to [bend left=0]node [auto,labelsize] {$\cat{B}\times(-)$}(B);

      \draw [->] (A) to [bend right=20]node [auto,swap,labelsize] (Fpp) {$G$}(B);
      \draw [2cellr] (1.5,-0.1) to node [auto,labelsize]{$\Gamma$}(1.5,-1.1);
            \end{scope}
    \end{tikzpicture}
    \quad=\quad
    \begin{tikzpicture}[baseline=-\the\dimexpr\fontdimen22\textfont2\relax ]
            \begin{scope} [shift={(0,1.5)}]
      \node (E) at (-3,0) {$\tcat{X}$};
      \node (A)  at (0,0) {$\EM{\Cat}{\cat{B}\times(-)}$};
      \node (B) at (3,-1.5) {$\EM{\Cat}{\cat{B}\times(-)}$};
      \node (C) at (3,0)  {$\Cat$};
      \node (D) at (3,-3) {$\Cat$};

      \draw [->] (E) to [bend left=0]node [auto,labelsize] {$\tilde{G}$}(A);
      \draw [->] (A) to [bend left=0]node [auto,labelsize] {$U^\cat{B}$}(C);
      \draw [->] (C) to [bend left=0]node [auto,labelsize] {$F^\cat{B}$}(B);
      \draw [->] (B) to [bend left=0]node [auto,labelsize] {$U^\cat{B}$}(D);

      \draw [->] (A) to [bend right=20]node [auto,swap,labelsize] (Fpp) {$\id{\EM{\Cat}{\cat{B}\times(-)}}$}(B);
      \draw [2cellr] (1.5,-0.1) to node [auto,labelsize]{$E^\cat{B}$}(1.5,-1.1);
            \end{scope}
    \end{tikzpicture}
    \]
since the co-Eilenberg--Moore 2-category~$\EM{\Cat}{\cat{B}\times(-)}$
is the co-Eilenberg--Moore object in $\TCat$.

\begin{defn}
The 2-functor~$\tilde{G}\colon\tcat{X}\to\EM{\Cat}{\cat{B}\times(-)}$
is the mediating 2-functor
$X^\prime\mapsto\vcoalgp{\cat{B}\times GX^\prime}{\Gamma_{X^\prime}}{GX^\prime}$.
\end{defn}

It remains to construct a 1-cell~$\tilde{g}\colon
\vcoalgp{\cat{B}\times GX}{\Gamma_{X}}{GX}\to\vcoalgp{\cat{B}\times\EMi}{\pi}{\EMi}$
of $\EM{\Cat}{\cat{B}\times(-)}$ which satisfies
\begin{align}
\label{eqn-univ-iem-5}
\begin{tikzpicture}[baseline=-\the\dimexpr\fontdimen22\textfont2\relax ]
        \begin{scope} 
      \node (A)  at (0,0) {$GX$};
      \node (B) at (2,0) {$\cat{C}$};
      \draw [->] (A) to [bend left=0]node [auto,labelsize] {${g}$}(B);
        \end{scope}
\end{tikzpicture}\quad
=\quad\begin{tikzpicture}[baseline=-\the\dimexpr\fontdimen22\textfont2\relax ]
        \begin{scope} 
      \node (A)  at (0,0) {$GX$};
      \node (B) at (2,0) {$\EMi$};
      \node (C) at (4,0) {$\cat{C}$};
      \draw [->] (A) to [bend left=0]node [auto,labelsize] {$\tilde{g}$}(B);
      \draw [->] (B) to [bend left=0]node [auto,labelsize] {$u^\im{T}$}(C);
        \end{scope}
\end{tikzpicture}
\end{align}
\begin{multline}
\label{eqn-univ-iem-6}
\begin{tikzpicture}[baseline=-\the\dimexpr\fontdimen22\textfont2\relax ]
        \begin{scope} [shift={(0,0.75)}]
      \node (M)  at (0,0)  {$GX$};
      \node (N)  at (3.8,1.05) {$\cat{B}\times\cat{C}$};
      \node (RR)  at (5,0) {$\cat{C}$};
      \node (B)  at (1.5,1.5) {$\cat{B}\times GX$};
      \draw [->] (M) to node  [auto,swap,labelsize]      {$g$} (RR);
      \draw [->] (M) to [bend left=20]node  [auto,labelsize] {$\Gamma_X$}(B);
      \draw [->] (B) to [bend left=8]node  [auto,labelsize] {$\cat{B}\times g$}(N);
      \draw [->] (N) to [bend left=5]node  [auto,labelsize] {$\im{T}$}(RR);
      \draw [2cellnode] (B) to node[auto,labelsize] {$\gamma$}(1.5,0);
        \end{scope}
\end{tikzpicture}\\
=\quad
\begin{tikzpicture}[baseline=-\the\dimexpr\fontdimen22\textfont2\relax ]
        \begin{scope} [shift={(0,-0.75)}]
      \node (E) at (-2,0) {$GX$};
      \node (M)  at (0,0)  {$\EMi$};
      \node (N)  at (3.8,1.05) {$\cat{B}\times\cat{C}$};
      \node (RR)  at (5,0) {$\EMi$};
      \node (B)  at (1.5,1.5) {$\cat{B}\times \EMi$};
      \node (F) at (7,0) {$\cat{C}$};
      \draw [->] (E) to node  [auto,labelsize]      {$\tilde{g}$} (M);
      \draw [->] (M) to node  [auto,swap,labelsize]      {$\id{\EMi}$} (RR);
      \draw [->] (M) to [bend left=20]node  [auto,labelsize] {$E_\pi$}(B);
      \draw [->] (B) to [bend left=8]node  [auto,labelsize] {$\cat{B}\times u^\im{T}$}(N);
      \draw [->] (N) to [bend left=4]node  [auto,labelsize] {$f^\im{T}$}(RR);
      \draw [->] (RR) to node  [auto,labelsize]      {$u^\im{T}$} (F);
      \draw [2cellnode] (B) to node[auto,labelsize] {$\varepsilon^\im{T}$}(1.5,0);
        \end{scope}
\end{tikzpicture}
\end{multline}
and show its uniqueness.
To describe the definition of this functor concisely, let us write
$\Gamma_Xx=(\Gamma_{X,0} x,x)=(\Gamma x,x)\in\cat{B}\times GX$
in what follows.
\begin{defn}
The 1-cell~$\tilde{g}\colon
\vcoalgp{\cat{B}\times GX}{\Gamma_{X}}{GX}\to\vcoalgp{\cat{B}\times\EMi}{\pi}{\EMi}$
of $\EM{\Cat}{\cat{B}\times(-)}$ is the functor~$
\tilde{g}\colon GX\to\EMi$ defined as
\[
\tilde{g}(x)\quad\coloneqq\quad\left(\Gamma x,
\valg{\imc{\Gamma x}gx}{\gamma_{x}}{gx}\right)
\]
on an object~$x$ and 
\[
\tilde{g}(z)\quad\coloneqq\quad\left(\Gamma x\xrightarrow{\Gamma z}\Gamma x^\prime,
\valgp{\imc{\Gamma x}gx}{\gamma_{x}}{gx}
\xrightarrow{gz}
\valgp{\imc{\Gamma x}gx^\prime}
{\gamma_{x^\prime}\circ\imc{\Gamma z, gx^\prime}}
{gx^\prime}\right)
\]
on a morphism~$z\colon x\to x^\prime$.
\end{defn}

\begin{prop}
The functor~$\tilde{g}$ defined above is indeed a 1-cell of 
$\EM{\Cat}{\cat{B}\times(-)}$ which satisfies
(\ref{eqn-univ-iem-5}) and (\ref{eqn-univ-iem-6}).
Moreover, it is the unique such.
\end{prop}
\begin{pf}
First, $\tilde{g}$ is a 1-cell of
$\EM{\Cat}{\cat{B}\times(-)}$ because
\begin{align*}
(\id{\cat{B}}\times\tilde{g})\Gamma_Xx\ &=\ 
(\Gamma x,\tilde{g}x)\\
&=\ 
(\Gamma x, (\Gamma x,\gamma_x))\\
&=\ 
\pi \tilde{g}x
\end{align*} 
on objects and 
\begin{align*}
(\id{\cat{B}}\times\tilde{g})\Gamma_Xz\ &=\ 
(\Gamma z,\tilde{g}z)\\
&=\ 
(\Gamma z,(\Gamma z,gz))\\
&=\ 
\pi \tilde{g}z
\end{align*} 
on morphisms.

$\tilde{g}$ satisfies (\ref{eqn-univ-iem-5}) since
\begin{align*}
gx\ 
&=\ u^\im{T}\left(\Gamma x, \valg{\imc{\Gamma x}gx}{\gamma_{x}}{gx}\right)\\
&=\ u^\im{T}\tilde{g}x
\end{align*}
on objects and
\begin{align*}
gz\ &=\ u^\im{T}\left(\Gamma x\xrightarrow{\Gamma z}\Gamma x^\prime,
\valgp{\imc{\Gamma x}gx}{\gamma_{x}}{gx}
\xrightarrow{gz}
\valgp{\imc{\Gamma x}gx^\prime}
{\gamma_{x^\prime}\circ\imc{\Gamma z, gx^\prime}}
{gx^\prime}\right)\\
&=\ u^\im{T}\tilde{g}z
\end{align*}
on morphisms.

Finally, $\tilde{g}$ satisfies (\ref{eqn-univ-iem-6}) since
\begin{align*}
\gamma_x\ 
&=\ u^\im{T}(\id{\Gamma_X},\gamma_x)\\
&=\ u^\im{T}\varepsilon^\im{T}_{(\Gamma x, \gamma_x)}\\
&=\ u^\im{T}\varepsilon^\im{T}_{\tilde{g}x}.
\end{align*}

Let us now move on to the proof of uniqueness.
The requirement that $\tilde{g}$ is a 1-cell of 
$\EM{\Cat}{\cat{B}\times(-)}$ determines the
first components of $\tilde{g}x$.
(\ref{eqn-univ-iem-5}) determines the underlying object of
the algebra part of $\tilde{g}x$.
Finally, (\ref{eqn-univ-iem-6}) forces the structure map
of the algebra part of $\tilde{g}x$ to be $\gamma_x$,
thus completely specifies the definition of $\tilde{g}$.
\end{pf}

\subsubsection{The 2-dimensional aspect}
Suppose we have a morphism of left $\im{T}$-modules, i.e., 
a 2-cell
    \[
    \begin{tikzpicture}
      \node (A)  at (0,0) {$(\tcat{X},X)$};
      \node (C) at (4,0) {$(\Cat,\cat{C})$};
      
      \draw [->] (A) to [bend left=20]node (dom)[auto,labelsize] {$(G,g)$}(C);
      \draw [->] (A) to [bend right=20]node (cod)[auto,swap,labelsize] {$(G^\prime,g^\prime)$}(C);

      \draw [2cellnode] (dom) to node [auto,labelsize]{$(\Omega,\omega)$}(cod);

    \end{tikzpicture}
    \]
of $\Epm$ satisfying
\begin{align}
\label{eqn-univ-iem-7}
    \begin{tikzpicture}[baseline=-\the\dimexpr\fontdimen22\textfont2\relax ]
            \begin{scope} [shift={(0,0.75)}]
      \node (A)  at (0,0) {$\tcat{X}$};
      \node (B) at (3,-1.5) {$\Cat$};
      \node (C) at (3,0) {$\Cat$};
      \draw [->] (A) to [bend left=20]node [auto,labelsize] {$G$}(C);
      \draw [->] (A) to [bend right=20]node [auto,very near start,labelsize] {$G^\prime$}(C);
      \draw [->] (C) to [bend left=0]node [auto,labelsize] {$\cat{B}\times(-)$}(B);
      %
      \draw [->] (A) to [bend right=20]node [auto,swap,labelsize] (Fpp) {$G^\prime$}(B);
      \draw [2cellr] (1.5,0.25) to node [auto,labelsize]{$\Omega$}(1.5,-0.25);
      \draw [2cellr] (1.5,-0.5) to node [auto,labelsize]{$\Gamma^\prime$}(1.5,-1.1);
            \end{scope}
    \end{tikzpicture}
    \quad=\quad
    \begin{tikzpicture}[baseline=-\the\dimexpr\fontdimen22\textfont2\relax ]
            \begin{scope} [shift={(0,0.75)}]
      \node (A)  at (0,0) {$\tcat{X}$};
      \node (B) at (3,-1.5) {$\Cat$};
      \node (C) at (3,0) {$\Cat$};
      \draw [->] (A) to [bend left=20]node [auto,labelsize] {$G$}(C);
      \draw [->] (A) to [bend left=20]node [auto,near end,labelsize] {$G$}(B);
      \draw [->] (C) to [bend left=0]node [auto,labelsize] {$\cat{B}\times(-)$}(B);
      %
      \draw [->] (A) to [bend right=20]node [auto,swap,labelsize] (Fpp) {$G^\prime$}(B);
      \draw [2cellr] (1.5,0.25) to node [auto,labelsize]{$\Gamma$}(1.5,-0.25);
      \draw [2cellr] (1.5,-0.5) to node [auto,labelsize]{$\Omega$}(1.5,-1.1);
            \end{scope}
    \end{tikzpicture}
\end{align}
\begin{multline}
\label{eqn-univ-iem-8}
\begin{tikzpicture}[baseline=-\the\dimexpr\fontdimen22\textfont2\relax ]
        \begin{scope} 
      \node (A)  at (0,0) {$\cat{B}\times G^\prime X$};
      \node (B) at (3,0) {$\cat{B}\times\cat{C}$};
      \node (C) at (5,0) {$\cat{C}$};
      \node (D)  at (1.5,-1.5) {$G^\prime X$};
      \node (F) at (1,1.5) {$\cat{B}\times GX$};
      \draw [->] (A) to [bend left=0]node [auto,labelsize] {$\cat{B}\times g^\prime$}(B);
      \draw [->] (B) to [bend left=0]node [auto,labelsize] {$\im{T}$}(C);
      \draw [<-] (A) to [bend right=20]node [auto,swap,labelsize] {$\Gamma^\prime _X$}(D);
      \draw [->] (D) to [bend right=20]node [auto,swap,labelsize] {$g^\prime$}(C);
      \draw [<-] (F) to [bend right=20]node [auto,swap,labelsize] {$\cat{B}\times\Omega_X$}(A);
      \draw [->] (F) to [bend left=20]node [auto,labelsize] {$\cat{B}\times g$}(B);
      \draw [2cell] (1.5,-0.1) to node [auto,labelsize]{$\gamma^\prime$}(1.5,-1.1);
      \draw [2cell] (1,1.2) to node [auto,labelsize]{$\cat{B}\times\omega$}(1,0.4);
        \end{scope}
\end{tikzpicture}\\
=\quad
    \begin{tikzpicture}[baseline=-\the\dimexpr\fontdimen22\textfont2\relax ]
            \begin{scope} 
          \node (M)  at (0,0)  {$GX$};
          \node (N)  at (3.8,1.05) {$\cat{B}\times\cat{C}$};
          \node (RR)  at (5,0) {$\cat{C}$};
          \node (B)  at (1.5,1.5) {$\cat{B}\times GX$};
          \node (E) at (1.5,-1.5) {$G^\prime X$};
          \draw [->] (M) to node  [auto,labelsize]      {$g$} (RR);
          \draw [->] (M) to [bend left=20]node  [auto,labelsize] {$\Gamma_X$}(B);
          \draw [->] (B) to [bend left=8]node  [auto,labelsize] {$\cat{B}\times g$}(N);
          \draw [->] (N) to [bend left=5]node  [auto,labelsize] {$\im{T}$}(RR);
          \draw [->] (E) to [bend left=20]node  [auto,labelsize] {$\Omega_X$}(M);
          \draw [->] (E) to [bend right=20]node  [auto,swap,labelsize] {$g^\prime$}(RR);
          \draw [2cellnode] (B) to node[auto,labelsize] {$\gamma$}(1.5,0);
          \draw [2cellnode] (1.5,0) to node[auto,labelsize] {$\omega$}(E);
            \end{scope}
    \end{tikzpicture}
\end{multline}

Equation~(\ref{eqn-univ-iem-7}) implies the unique factorization
    \[
    \begin{tikzpicture}[baseline=-\the\dimexpr\fontdimen22\textfont2\relax ]
      \node (E) at (-3,0) {$\tcat{X}$};
      \node (A)  at (0,0) {$\Cat$};
      \draw [->] ([yshift=2mm]E.east) to [bend left=20]node(dom) [auto,labelsize] {${G}$} ([yshift=2mm]A.west);
      \draw [->] ([yshift=-2mm]E.east) to [bend right=20]node(cod) [auto,swap,labelsize] {${G^\prime}$}([yshift=-2mm]A.west);
      \draw [2cellnoder] (dom) to node [auto,labelsize]{${\Omega}$}(cod);

    \end{tikzpicture}
    \quad=\quad
    \begin{tikzpicture}[baseline=-\the\dimexpr\fontdimen22\textfont2\relax ]
      \node (E) at (-3,0) {$\tcat{X}$};
      \node (A)  at (0,0) {$\EM{\Cat}{\cat{B}\times(-)}$};
      \node (C) at (3,0) {$\Cat$};
      \draw [->] ([yshift=2mm]E.east) to [bend left=20]node(dom) [auto,labelsize] {$\tilde{G}$} ([yshift=2mm]A.west);
      \draw [->] ([yshift=-2mm]E.east) to [bend right=20]node(cod) [auto,swap,labelsize] {$\tilde{G^\prime}$}([yshift=-2mm]A.west);
      \draw [->] (A) to [bend left=0]node [auto,labelsize] {$U^\cat{B}$}(C);
      \draw [2cellnoder] (dom) to node [auto,labelsize]{$\tilde{\Omega}$}(cod);

    \end{tikzpicture}
    \]
where the 2-natural transformation~$\tilde{\Omega}$ is defined as follows:
\begin{defn}
The 2-natural transformation~$\tilde{\Omega}\colon\tilde{G^\prime}\Rightarrow\tilde{G}$
has, as components, 1-cells~$
\tilde{\Omega}_{X^\prime}\colon
\vcoalgp{\cat{B}\times G^\prime {X^\prime}}{\Gamma^\prime_{X^\prime}}{G^\prime {X^\prime}}\to
\vcoalgp{\cat{B}\times G{X^\prime}}{\Gamma_{X^\prime}}{G{X^\prime}}$
of $\EM{\Cat}{\cat{B}\times(-)}$
given by $\tilde{\Omega}_{X^\prime}\coloneqq\Omega_{X^\prime}$
as functors.
\end{defn}

Next we proceed to construct a 2-cell
    \[
    \begin{tikzpicture}
      \node (a)  at (0,0) {$\vcoalgp{\cat{B}\times G^\prime{X}}{\Gamma^\prime_{X}}{G^\prime{X}}$};
      \node (ap) at (6,0) {$\vcoalgp{\cat{B}\times\EMi}{\pi}{\EMi}$};
      \node (b)  at (2,3) {$\vcoalgp{\cat{B}\times G{X}}{\Gamma_{X}}{G{X}}$};
      
      \draw [->] (a) to [bend right=0]node [auto,swap,labelsize] {$\tilde{g^\prime}$}(ap);
      \draw [->] (a)  to [bend right=-20]node [auto,labelsize] {$\tilde{\Omega}_X$}(b);
      \draw [->] (b) to [bend right=-20]node [auto,labelsize] {$\tilde{g}$}(ap);
      \draw [2cellr] (2,0.1) to node [auto,swap,labelsize]{$\tilde{\omega}$}(b);
    \end{tikzpicture}
    \]
of $\EM{\Cat}{\cat{B}\times(-)}$ satisfying
\begin{align}
\label{eqn-univ-iem-9}
    \begin{tikzpicture}[baseline=-\the\dimexpr\fontdimen22\textfont2\relax ]
     \begin{scope} [shift={(0,-0.75)}]
      \node (a)  at (0,0) {$G^\prime X$};
      \node (ap) at (3,0) {$\cat{C}$};
      \node (b)  at (1,1.5) {$GX$};
      \draw [->] (a) to [bend right=0]node [auto,swap,labelsize] {$g^\prime$}(ap);
      \draw [->] (a)  to [bend right=-20]node [auto,labelsize] {$\Omega_X$}(b);
      \draw [->] (b) to [bend right=-20]node [auto,labelsize] {$g$}(ap);
      \draw [2cellr] (1,0.1) to node [auto,swap,labelsize]{$\omega$}(b);
     \end{scope}
    \end{tikzpicture}
   \quad=\quad
    \begin{tikzpicture}[baseline=-\the\dimexpr\fontdimen22\textfont2\relax ]
     \begin{scope} [shift={(0,-0.75)}]
      \node (a)  at (0,0) {$G^\prime X$};
      \node (ap) at (3,0) {$\EMi$};
      \node (b)  at (1,1.5) {$GX$};
      \node (c)  at (5,0) {$\cat{C}$};
      \draw [->] (a) to [bend right=0]node [auto,swap,labelsize] {$\tilde{g^\prime}$}(ap);
      \draw [->] (a)  to [bend right=-20]node [auto,labelsize] {$\tilde{\Omega}_X$}(b);
      \draw [->] (b) to [bend right=-20]node [auto,labelsize] {$\tilde{g}$}(ap);
      \draw [->] (ap) to [bend right=0]node [auto,labelsize] {$u^\im{T}$}(c);
      \draw [2cellr] (1,0.1) to node [auto,swap,labelsize]{$\tilde{\omega}$}(b);
     \end{scope}
    \end{tikzpicture}
\end{align}

\begin{defn}
To define the 2-cell~$\tilde{\omega}$ of $\EM{\Cat}{\cat{B}\times(-)}$,
with components
\[
\tilde{\omega}_x\quad\colon\quad
\left(\Gamma \Omega_Xx,
\valg{\imc{\Gamma\Omega_Xx}g\Omega_Xx}{\gamma_{g\Omega_Xx}}{g\Omega_Xx}\right)\quad
\longrightarrow\quad
\left(\Gamma^\prime x,
\valg{\imc{\Gamma^\prime x}g^\prime x}{\gamma^\prime_{g^\prime x}}{g^\prime x}\right)
\]
at $x\in G^\prime X$, first observe that by
equation~(\ref{eqn-univ-iem-7}),
$\Gamma\Omega_X x=\Gamma^\prime x$.
Now we define $\tilde{\omega}_x$ to be
\[
\tilde{\omega}_x\quad\coloneqq\quad
\left(
\Gamma\Omega_X x\xrightarrow{\id{}}\Gamma^\prime x,\,
\valgp{\imc{\Gamma\Omega_Xx}g\Omega_Xx}{\gamma_{g\Omega_Xx}}{g\Omega_Xx}
\xrightarrow{\omega_x}
\valgp{\imc{\Gamma^\prime x}g^\prime x}{\gamma^\prime_{g^\prime x}}{g^\prime x}
\right).\qedhere
\]
\end{defn}
\begin{prop}
The natural transformation~$\tilde{\omega}$ defined above is indeed a 
2-cell of $\EM{\Cat}{\cat{B}\times(-)}$ which satisfies
(\ref{eqn-univ-iem-9}). Moreover, it is the unique such.
\end{prop}
\begin{pf}
First, $\tilde{\omega}$ is a 2-cell of $\EM{\Cat}{\cat{B}\times(-)}$,
i.e., $(\id{\cat{B}}\times\tilde{\omega})_{\Gamma^\prime_X x}
=\pi\tilde{\omega}_x$ holds, since
\begin{align*}
(\id{\cat{B}}\times\tilde{\omega})_{\Gamma^\prime_X x}\ &=\ 
\left(\Gamma^\prime x\xrightarrow{\id{}}\Gamma^\prime x,\,
\tilde{\omega}_x\right)\\
&=\ \left(\Gamma\Omega_X x\xrightarrow{\id{}}\Gamma^\prime x,\,
\tilde{\omega}_x\right)\\
&=\ \pi\tilde{\omega}_x.
\end{align*}

That $\tilde{\omega}$ satisfies (\ref{eqn-univ-iem-9}) is also 
easy to see:
\begin{align*}
\omega_x\ 
&=\ u^\im{T}\left(
\Gamma\Omega_X x\xrightarrow{\id{}}\Gamma^\prime x,\,
\valgp{\imc{\Gamma\Omega_Xx}g\Omega_Xx}{\gamma_{g\Omega_Xx}}{g\Omega_Xx}
\xrightarrow{\omega_x}
\valgp{\imc{\Gamma^\prime x}g^\prime x}{\gamma^\prime_{g^\prime x}}{g^\prime x}
\right)\\
&=\ u^\im{T}\tilde{\omega}_x.
\end{align*}

For the uniqueness, observe that 
the requirement that $\tilde{\omega}$ is a 2-cell of 
$\EM{\Cat}{\cat{B}\times(-)}$ forces the first component
of $\tilde{\omega}_x$ to be the identity, and 
the requirement that $\tilde{\omega}$ satisfies 
(\ref{eqn-univ-iem-9}) determines the second component.
\end{pf}
\begin{thm}
\label{thm-univ-em-im}
The object~$\left(\EM{\Cat}{\cat{B}\times(-)},
\vcoalg{\cat{B}\times\EMi}{\pi}{\EMi}\right)$ of 
$\Epm$ is the Eilenberg--Moore object of the
indexed monad~$\im{T}$, considered as a monad
$(\cat{B}\times(-),\im{T})$ in $\Epm$ on $(\Cat,\cat{C})$.
\end{thm}

\section{Constructions for graded and indexed comonads}
The constructions introduced so far dualize to those for 
graded and indexed \emph{comonads} rather straightforwardly.
We briefly describe how the co-Eilenberg--Moore and co-Kleisli
categories look like.

\subsection{Co-Eilenberg--Moore categories of graded comonads}
Let us fix a strict monoidal category~$\cat{M}=(\cat{M},\tensor,\e)$,
a category~$\cat{C}$, and an $\cat{M}$-graded comonad~$\gm{S}$
on $\cat{C}$.
We write the functor part of $\gm{S}$ also as 
$\act\colon \cat{M}\times\cat{C}\to\cat{C}$ and 
use the infix notation.
As observed in Section~\ref{subsec-gicomnd-comnd},
graded comonads can be seen as comonads in the 2-category~$\Epm$;
the co-Eilenberg--Moore construction is performed inside $\Epm$.

\begin{defn}
Define the category~$\CoEMg$ as follows:
\begin{itemize}
\item An object of $\CoEMg$ is a 
\defemph{graded $\gm{S}$-coalgebra}, i.e.,
a pair~$(A,h)$ where $A\colon \cat{M}\to\cat{C}$ is a functor and 
$h$ is a natural transformation of type
\[
\begin{tikzpicture}
      \node (TL)  at (0,0.75)  {$\cat{M}\times\cat{M}$};
      \node (TR)  at (3,0.75)  {$\cat{M}\times\cat{C}$};
      \node (BL)  at (0,-0.75) {$\cat{M}$};
      \node (BR)  at (3,-0.75) {$\cat{C}$};
      
      \draw [->] (TL) to node (T) [auto,labelsize]      {$\cat{M}\times A$} (TR);
      \draw [->] (TR) to node (R) [auto,labelsize]      {$\gm{S}$} (BR);
      \draw [->] (TL) to node (L) [auto,swap,labelsize] {$\tensor$}(BL);
      \draw [->] (BL) to node (B) [auto,swap,labelsize] {$A$}(BR);
      
      \draw [2cellnoder] (R) to  [bend right=30] node [auto,swap,labelsize] {$h$}(B);
\end{tikzpicture}
\]
So the component of $h$ at $(m,n)\in\cat{M}\times\cat{M}$ is of type
\[
h_{m,n}\quad\colon \quad A_{m\tensor n}\quad \longrightarrow\quad m\act A_n.
\] 
These data are subject to the following axioms:
\begin{itemize}
\item $\begin{tikzpicture}[baseline=-\the\dimexpr\fontdimen22\textfont2\relax ]
      \node (TL)  at (0,0.75)  {$A_n$};
      \node (TR)  at (2,0.75)  {$\e\act A_n$};
      \node (BR)  at (2,-0.75) {$A_n$};
      
      \draw [->] (TL) to node (T) [auto,labelsize]      {$h_{\e, n}$} (TR);
      \draw [->] (TR) to node (R) [auto,labelsize]      {$\varepsilon_{A_n}$} (BR);
      \draw [->] (TL) to node (L) [auto,swap,labelsize] {$\id{A_n}$}(BR);
\end{tikzpicture}$
commutes for each object~$n$ of $\cat{M}$.
\item 
$\begin{tikzpicture}[baseline=-\the\dimexpr\fontdimen22\textfont2\relax ]
      \node (TL)  at (0,-0.75)  {$l\act m\act A_n$};
      \node (TR)  at (-4,-0.75)  {$(l\tensor m)\act A_n$};
      \node (BL)  at (0,0.75) {$l\act A_{m\tensor n}$};
      \node (BR)  at (-4,0.75) {$A_{l\tensor m\tensor n}$};
      
      \draw [<-] (TL) to node (T) [auto,labelsize]      {$\delta_{l,m,A_n}$} (TR);
      \draw [<-] (TR) to node (R) [auto,labelsize]      {$h_{l\tensor m,n}$} (BR);
      \draw [<-] (TL) to node (L) [auto,swap,labelsize] {$l\act h_{m,n}$}(BL);
      \draw [<-] (BL) to node (B) [auto,swap,labelsize] {$h_{l,m\tensor n}$}(BR);
\end{tikzpicture}$
commutes for each triple of objects~$l$, $m$, $n$ of $\cat{M}$.
\end{itemize}
\item
A morphism of $\CoEMg$ from $(A,h)$ to $(A^\prime,h^\prime)$ is a 
\defemph{homomorphism} of graded $\gm{S}$-coalgebras between them, i.e., a 
natural transformation
$\varphi\colon A\Rightarrow A^\prime$ making the diagram
\[
\begin{tikzpicture}
      \node (TL)  at (0,-0.75)  {$m\act A_n$};
      \node (TR)  at (3,-0.75)  {$m\act A^\prime_n$};
      \node (BL)  at (0,0.75) {$A_{m\tensor n}$};
      \node (BR)  at (3,0.75) {$A^\prime_{m\tensor n}$};
      
      \draw [->] (TL) to node (T) [auto,swap,labelsize]      {$m\act \varphi_{n}$} (TR);
      \draw [<-] (TR) to node (R) [auto,swap,labelsize]      {$h^\prime_{m,n}$} (BR);
      \draw [<-] (TL) to node (L) [auto,labelsize] {$h_{m,n}$}(BL);
      \draw [->] (BL) to node (B) [auto,labelsize] {$\varphi_{m\tensor n}$}(BR);
\end{tikzpicture}
\]
commute for each pair of objects~$m,n$ of $\cat{M}$.\qedhere
\end{itemize}
\end{defn}

\begin{defn}
Define the functor
\[
\emact\quad\colon\quad\cat{M}\times\CoEMg\quad\longrightarrow
\quad\CoEMg
\]
as follows:
\begin{itemize}
\item Given objects $p$ and $(A,h)$ of 
$\cat{M}$ and $\CoEMg$ respectively,
we define the graded $\gm{S}$-algebra~$p\emact (A,h)$ by
the precomposition of $(-)\tensor p\colon\cat{M}\to\cat{M}$:
\[
p\emact (A,h)\quad\coloneqq\quad((A_{n\tensor p})_{n\in\cat{M}},\ 
(h_{m,n\tensor p})_{m,n\in\cat{M}}).
\]
\item Given morphisms $u\colon p\to p^\prime$ and $\varphi\colon(A,h)\to(A^\prime,h^\prime)$
of $\cat{M}$ and $\CoEMg$ respectively, 
we define the homomorphism~$u\emact \varphi\colon
p\emact(A,h)\to p^\prime\emact(A^\prime,h^\prime)$
by setting the component~$(u\emact \varphi)_{n}\colon
A_{n\tensor p}\to A^\prime_{n\tensor p^\prime}$
to be either of the following two equivalent composites:
\[
\begin{tikzpicture}
      \node (TL)  at (0,0.75)  {$A_{n\tensor p}$};
      \node (TR)  at (3,0.75)  {$A^\prime_{n\tensor p}$};
      \node (BL)  at (0,-0.75) {$A_{n\tensor p^\prime}$};
      \node (BR)  at (3,-0.75) {$A^\prime_{n\tensor p^\prime}$};
      
      \draw [->] (TL) to node (T) [auto,labelsize]      {$\varphi_{n\tensor p}$} (TR);
      \draw [->] (TR) to node (R) [auto,labelsize]      {$A^\prime_{n\tensor u}$} (BR);
      \draw [->] (TL) to node (L) [auto,swap,labelsize] {$A_{n\tensor u}$}(BL);
      \draw [->] (BL) to node (B) [auto,swap,labelsize] {$\varphi_{n\tensor p^\prime}$}(BR);
\end{tikzpicture}
\qedhere
\]
\end{itemize}
\end{defn}

\begin{thm}
The object~$\left(\EM{\Cat}{\cat{M}\times(-)},
\valg{\cat{M}\times\CoEMg}{\emact}{\CoEMg}\right)$ of 
$\Epm$ is the co-Eilenberg--Moore object of the
graded comonad~$\gm{S}$, considered as a comonad
$(\cat{M}\times(-),\gm{S})$ in $\Epm$ on $(\Cat,\cat{C})$.
\end{thm}

Recall that a $(\cat{M},\tensor,\e)$-graded comonad on $\cat{C}$
is the same thing as a $(\cat{M}^\opp,\tensor,\e)$-graded monad
on $\cat{C}^\opp$.
Let $\bar{\gm{S}}$ denote the $(\cat{M}^\opp,\tensor,\e)$-graded monad
on $\cat{C}^\opp$ corresponding to $\gm{S}$. 
Now, Eilenberg--Moore categories for graded monads and 
co-Eilenberg--Moore categories for graded comonads are 
related to each other in the following way:
\[
\EM{(\cat{C}^\opp)}{\bar{\gm{S}}}\quad\cong \quad (\EM{\cat{C}}{\gm{S}})^\opp.
\]
Actually, one can say more: 
the canonical $\cat{M}^\opp\times(-)$-algebra structure~$
\bar{\emact}\colon\cat{M}^\opp\times\EM{(\cat{C}^\opp)}{\bar{\gm{S}}}
\to\EM{(\cat{C}^\opp)}{\bar{\gm{S}}}$ on $\EM{(\cat{C}^\opp)}{\bar{\gm{S}}}$
corresponds to the canonical $\cat{M}\times(-)$-algebra structure~$
\emact\colon\cat{M}\times\CoEMg
\to\CoEMg$ on $\CoEMg$ via this dualization: $\bar{\emact}\cong\emact^\opp$.

\subsection{Co-Kleisli categories of graded comonads}
Again we fix a strict monoidal category~$\cat{M}=(\cat{M},\tensor,\e)$,
a category~$\cat{C}$, and an $\cat{M}$-graded comonad~$\gm{S}$
on $\cat{C}$;
we continue to write the functor part of $\gm{S}$ also as 
$\act\colon \cat{M}\times\cat{C}\to\cat{C}$ and 
use the infix notation.
For the co-Kleisli construction for graded comonads, we use another
observation in Section~\ref{subsec-gicomnd-comnd} that
graded comonads can also be seen as comonads in the 2-category~$\Emp$.

\begin{defn}
Define the category~$\CoKlg$ as follows:
\begin{itemize}
\item An object of $\CoKlg$ is a pair~$(m,c)$ where $m$ and $c$
are objects of $\cat{M}$ and $\cat{C}$ respectively.
\item The set of morphisms from $(m,c)$ to $(m^\prime,c^\prime)$
is defined by the coend formula
\[
\CoKlg((m,c),(m^\prime,c^\prime))\quad\coloneqq\quad\int^{n\in \cat{M}}
\cat{M}(m,m^\prime\tensor n)\times\cat{C}(n\act c,c^\prime).
\]
Explicitly, a morphism~$(m,c)\to(m^\prime,c^\prime)$ is an 
equivalence class~$
[n,\,m\xrightarrow{v}m^\prime\tensor n,\,n\act c\xrightarrow{f}c^\prime]$
of tuples consisting of an object~$n\in\cat{M}$ and morphisms~$v,f$, 
where the equivalence relation is generated by 
\begin{multline*}
(n,\ m\xrightarrow{v}m^\prime\tensor n,\ 
n\act c \xrightarrow{w\act c}n^\prime\act c\xrightarrow{f}c^\prime)\\
\sim\quad
(n^\prime,\ m\xrightarrow{v}m^\prime\tensor n\xrightarrow{m^\prime\tensor w}
m^\prime \tensor n^\prime,\ 
n^\prime \act c\xrightarrow{f}c^\prime)
\end{multline*}
for each morphism~$w\colon n\to n^\prime$ of $\cat{M}$.
\item The identity morphism on $(m,c)$ is given by 
$[\e,\,m\xrightarrow{\id{m}}m\tensor \e,\,
\e\act c\xrightarrow{\varepsilon_c}c]$.
\item For two composable morphisms
\begin{align*}
[n,\ m\xrightarrow{v}m^\prime\tensor n,\ n\act c\xrightarrow{f}c^\prime]
\quad&\colon\quad
(m,c)\quad\longrightarrow\quad(m^\prime,c^\prime),\\
[n^\prime,\ m^\prime\xrightarrow{v^\prime}m^{\prime\prime}\tensor n^\prime,\ 
n^\prime\act c^\prime\xrightarrow{f^\prime}c^{\prime\prime}]
\quad&\colon\quad 
(m^\prime,c^\prime)\quad\longrightarrow\quad(m^{\prime\prime},c^{\prime\prime}),
\end{align*}
their composite is given by
\begin{multline*}
[n^\prime\tensor n,\ m\xrightarrow{v}m^\prime\tensor n
\xrightarrow{v^\prime\tensor n}m^{\prime\prime}\tensor n^\prime\tensor n,\\
(n^\prime\tensor n)\act c
\xrightarrow{\delta_{n^\prime,n,c}}
n^\prime\act n\act c
\xrightarrow{n^\prime\act f}
n^\prime\act c^{\prime}
\xrightarrow{f^\prime}
c^{\prime\prime}]\\
\quad\colon\quad
(m,c)\quad\longrightarrow\quad(m^{\prime\prime},c^{\prime\prime}).
\tag*{\qedhere}
\end{multline*}
\end{itemize}
\end{defn}

Note that the first component of the composite morphism 
defined above is $n^\prime\tensor n$, rather than $n\tensor n^\prime$.

\begin{defn}
Define the functor
\[
\klact\quad\colon\quad\CoKlg\quad\longrightarrow\quad[\cat{M},\CoKlg]
\]
as follows;
note that we will use the infix notation~$l\klact (m,c)$
to denote the value of the functor $\klact (m,c)\colon \cat{M}\to\CoKlg$
applied to $l\in\cat{M}$ and similarly for morphisms.
\begin{itemize}
\item Given objects~$l$ and $(m,c)$ of $\cat{M}$ and $\CoKlg$ respectively,
we define $l\klact(m,c)\coloneqq(l\tensor m,c)$.
\item Given morphisms~$u\colon l\to l^\prime$ and 
$[n,\,m\xrightarrow{v}m^\prime\tensor n,\,n\act c\xrightarrow{f} c^\prime]
\colon (m,c)\to (m^\prime,c^\prime)$
of $\cat{M}$ and $\CoKlg$ respectively,
we define
\begin{multline*}
u\klact 
[n,\,m\xrightarrow{v}m^\prime\tensor n,\,n\act c\xrightarrow{f} c^\prime]\\
\coloneqq\quad
[n,\ l\tensor m\xrightarrow{u\tensor v}l^\prime\tensor m^\prime\tensor n,\ 
n\act c\xrightarrow{f}c^\prime]\\
\colon\quad
(l\tensor m,c)\quad\longrightarrow\quad(l^\prime \tensor m^\prime,c^\prime).\tag*{\qedhere}
\end{multline*}
\end{itemize}
\end{defn}

\begin{thm}
The object~$\left(\EM{\Cat}{[\cat{M},-]},
\vcoalg{[\cat{M},\CoKlg]}{\klact}{\CoKlg}\right)$ of 
$\Emp$ is the co-Kleisli object of the
graded comonad~$\gm{S}$, considered as a comonad
$([\cat{M},-],\gm{S})$ in $\Emp$ on $(\Cat,\cat{C})$.
\end{thm}
The relation to the Kleisli construction for graded monads
is again given by 
$\Kl{(\cat{C}^\opp)}{\bar{\gm{S}}}\cong(\Kl{\cat{C}}{\gm{S}})^\opp$ and 
$\bar{\klact}\cong \klact^\opp$,
where $\bar{\klact}\colon
\Kl{(\cat{C}^\opp)}{\bar{\gm{S}}}\to[\cat{B}^\opp,\Kl{(\cat{C}^\opp)}{\bar{\gm{S}}}]$
and $\klact\colon\Kl{\cat{C}}{\gm{S}}\to
[\cat{B},\Kl{\cat{C}}{\gm{S}}]$.

\subsection{Co-Eilenberg--Moore categories of indexed comonads}
Let us fix categories~$\cat{B},\cat{C}$, and a $\cat{B}$-indexed
comonad~$\im{S}$ on $\cat{C}$.
Based on the observation in Section~\ref{subsec-gicomnd-comnd} that
indexed comonads can be considered as comonads in the 2-category~$\Epp$,
we construct the co-Eilenberg--Moore category of $\im{S}$ in $\Epp$.

\begin{defn}
Define the category $\CoEMi$ as follows:
\begin{itemize}
\item An object of $\CoEMi$ is a pair 
$\left(b,\vcoalg{\icc{b}c}{\chi}{c}\right)$,
or more concisely $(b,\chi)$,
where $b$ is an object of $\cat{B}$ and 
$\vcoalgp{\icc{b}c}{\chi}{c}$ is a $\icc{b}$-coalgebra
(an object of $\EM{\cat{C}}{\icc{b}}$).
\item A morphism from $\left(b,\vcoalg{\icc{b}c}{\chi}{c}\right)$ to 
$\left(b^\prime,\vcoalg{\icc{b^\prime}c^\prime}{\chi^\prime}{c^\prime}\right)$
is a pair~$(u,h)$
where $u\colon b\to b^\prime$
is a morphism of $\cat{B}$ and 
$h\colon \vcoalgp{\icc{b^\prime}c}{\icc{u,c}\circ\chi}{c}\to
\vcoalgp{\icc{b^\prime}c^\prime}{\chi^\prime}{c^\prime}$
is a homomorphism of $\icc{b^\prime}$-coalgebras, i.e., a
morphism~$h\colon c\to c^\prime$ in $\cat{C}$ which makes the diagram
\[
\begin{tikzpicture}
      \node (TL)  at (0,-1)  {$\icc{b^\prime}c^\prime$};
      \node (TR)  at (-3,-1)  {$\icc{b^\prime}c$};
      \node (BL)  at (0,1) {$c^\prime$};
      \node (BR)  at (-3,1) {$c$};
      \node (MR)  at (-3,0)  {$\icc{b}c$};
      
      \draw [<-] (TL) to node (T) [auto,labelsize]      {$\icc{b^\prime}h$} (TR);
      \draw [<-] (TR) to node (R) [auto,labelsize]      {$\icc{u,c}$} (MR);
      \draw [<-] (MR) to node (R) [auto,labelsize]      {$\chi$} (BR);
      \draw [<-] (TL) to node (L) [auto,swap,labelsize] {$\chi^\prime$}(BL);
      \draw [<-] (BL) to node (B) [auto,swap,labelsize] {$h$}(BR);
\end{tikzpicture}
\]
commute.
\item The identity morphism on $\left(b,\vcoalg{\icc{b}c}{\chi}{c}\right)$
is given by $(\id{b},\id{\chi})$.
\item For two composable morphisms
\begin{align*}
\left(b\xrightarrow{u}b^\prime,\ 
\vcoalgp{\icc{b^\prime}c}{\icc{u,c}\circ\chi}{c}\xrightarrow{h}
\vcoalgp{\icc{b^\prime}c^\prime}{\chi^\prime}{c^\prime}\right)
\quad&\colon\quad
(b,\chi)
\quad\longrightarrow\quad
(b^\prime,\chi^\prime),\\
\left(b^\prime\xrightarrow{u^\prime} b^{\prime\prime},\ 
\vcoalgp{\icc{b^{\prime\prime}} c^\prime}{\icc{u^\prime,c^\prime}\circ\chi^\prime}{c^\prime}
\xrightarrow{h^\prime}
\vcoalgp{\icc{b^{\prime\prime}} c^{\prime\prime}}{\chi^{\prime\prime}}{c^{\prime\prime}}\right)
\quad&\colon\quad
(b^\prime,\chi^\prime)
\quad\longrightarrow\quad
(b^{\prime\prime},\chi^{\prime\prime}),
\end{align*}
their composite is given by 
\begin{multline*}
\left(b\xrightarrow{u}b^\prime\xrightarrow{u^\prime}b^{\prime\prime},\ 
\vcoalgp{\icc{b}c}{\icc{u^\prime\circ u,c}\circ\chi}{c}\xrightarrow{h}
\vcoalgp{\icc{b}c^\prime}{\icc{u^\prime,c^\prime}\circ\chi^\prime}{c^\prime}
\xrightarrow{h^\prime}
\vcoalgp{\icc{b} c^{\prime\prime}}{\chi^{\prime\prime}}{c^{\prime\prime}}\right)\\
\colon\quad (b,\chi)\quad\longrightarrow\quad
(b^{\prime\prime},\chi^{\prime\prime}).\tag*{\qedhere}
\end{multline*}
\end{itemize}
\end{defn}

\begin{defn}
Define the functor~$
\pi\colon\CoEMi\to\cat{B}\times\CoEMi$
by $(b,\chi)\mapsto (b,(b,\chi))$ and $(u,h)\mapsto (u,(u,h))$.
\end{defn}
\begin{thm}
The object~$\left(\EM{\Cat}{\cat{B}\times(-)},
\vcoalg{\cat{B}\times\CoEMi}{\pi}{\CoEMi}\right)$ of 
$\Epp$ is the co-Eilenberg--Moore object of the
indexed comonad~$\im{S}$, considered as a comonad
$(\cat{B}\times(-),\im{S})$ in $\Epp$ on $(\Cat,\cat{C})$.
\end{thm}

A $\cat{B}$-indexed comonad on $\cat{C}$ corresponds to a 
$\cat{B}^\opp$-indexed monad on $\cat{C}^\opp$;
let $\bar{\im{S}}$ be the $\cat{B}^\opp$-indexed monad on $\cat{C}^\opp$
corresponding to $\im{S}$.
The relationship of the Eilenberg--Moore construction for
indexed monads and the co-Eilenberg--Moore construction for
indexed comonads is given by
$\EM{(\cat{C}^\opp)}{\bar{\im{S}}}\cong(\EM{\cat{C}}{\im{S}})^\opp$
and $\bar{\pi}\cong\pi^\opp$, where
$\bar{\pi}\colon\EM{(\cat{C}^\opp)}{\bar{\im{S}}}
\to \cat{B}^\opp\times\EM{(\cat{C}^\opp)}{\bar{\im{S}}}$
and 
$\pi\colon\EM{\cat{C}}{\im{S}}\to\cat{B}\times\EM{\cat{C}}{\im{S}}$.

Recall the isomorphism of 2-categories~$\EM{\Cat}{\cat{B}\times(-)}
\cong\Cat/\cat{B}$. 
In contrast to the phenomenon that Eilenberg--Moore categories for 
$\cat{B}$-indexed monads give rise to \emph{Grothendieck fibrations
over $\cat{B}$}
under this isomorphism, co-Eilenberg--Moore categories for
$\cat{B}$-indexed comonads give rise to \emph{Grothendieck opfibrations
over $\cat{B}$}.

\section*{Notes}
The definition of the Eilenberg--Moore category of a graded monad
has been suggested to me independently by Shin-ya Katsumata
and by Paul-Andr\'e Melli\`es.
I learned the definition of the Kleisli category of a graded 
monad, together with Lemma~\ref{lem-decomp-klg},
from Shin-ya Katsumata.
After learning these definition, I formulated and proved
the 2-dimensional universality of the Eilenberg--Moore and Kleisli 
categories of a graded monad (Theorems~\ref{thm-univ-em-gm}
and \ref{thm-univ-kl-gm}).
The contents of Sections~\ref{sec-em-gm} and \ref{sec-kl-gm} are 
included in the paper~\cite{fujii-katsumata-mellies}.

I defined the Eilenberg--Moore category of an indexed monad,
and proved its universality (Theorem~\ref{thm-univ-em-im}).
The construction turned out to be 
essentially the same as the one that appears in
\cite{maillard-mellies} (but not as the Eilenberg--Moore category).

\chapter{Applications of the constructions}
\label{chap-appl}

In this chapter, we present two applications of the 
constructions presented in the previous chapter.
The first one is discussed in Section~\ref{sec-decomp-lax-action}, and
is an application of the Eilenberg--Moore and  
Kleisli constructions for graded monads; 
we shall show that they can be understood as giving ways to
\emph{decompose} lax monoidal actions into strict monoidal actions and 
adjunctions.
Then, in Section~\ref{sec-constr-maillard-mellies}, 
we see how the other construction, the Eilenberg--Moore construction
for indexed monads, sheds new light to the previous work of 
Power~\cite{power:indexed,power:lce} and 
Maillard--Melli\`es~\cite{maillard-mellies}, by revealing 
the related notions which have been implicit in their work.

\section{Decomposition of lax monoidal actions}
\label{sec-decomp-lax-action}
Let us fix a strict monoidal category~$\cat{M}=(\cat{M},\tensor,\e)$
throughout this section.
For $\cat{C}$ a category, an $\cat{M}$-graded monad on $\cat{C}$ 
is equivalent to the notion known as 
\emph{lax action} of $\cat{M}$ on $\cat{C}$,
in the sense of lax algebras for 2-monads~\cite{bkp:2-monad}.
On the other hand, objects of the 2-category~$
\EM{\Cat}{\cat{M}\times(-)}\cong\EM{\Cat}{[\cat{M},-]}$
are naturally thought of as categories equipped with
\emph{strict actions} of $\cat{M}$.
In this section, we will explain that 
the Eilenberg--Moore and Kleisli constructions for graded monads
developed in Sections~\ref{sec-em-gm} and \ref{sec-kl-gm}
can also be understood as a 
result relating these different types of 
actions of a monoidal category;
in a certain sense, these constructions show that 
we can always \emph{reduce} the general notion of 
lax action to the more restrictive notion of strict action.

\subsection{Lax and strict actions}
Let us first fix the definitions of strict and lax actions.
\begin{defn}
Let $\cat{A}$ be a category.
A \defemph{strict $\cat{M}$-action on $\cat{A}$} is a functor
\[
\sact\quad\colon\quad\cat{M}\times\cat{A}\quad\longrightarrow\quad
\cat{A}
\]
satisfying the equalities
\[
a\quad=\quad\e\sact a,\qquad \quad
m\sact n\sact a\quad=\quad(m\tensor n)\sact a.
\]

A category~$\cat{A}$ equipped with a strict $\cat{M}$-action~$\sact$ 
on it is a \defemph{strict $\cat{M}$-category $
\valgp{\cat{M}\times\cat{A}}{\sact}{\cat{A}}$}.
\end{defn}

\begin{defn}
Let $\cat{C}$ be a category.
A \defemph{lax $\cat{M}$-action on $\cat{A}$} is a functor
\[
\act\quad\colon\quad\cat{M}\times\cat{C}\quad\longrightarrow
\quad\cat{C},
\]
together with a family of morphisms
\[
c\quad\overset{\eta_c}{\longrightarrow}\quad \e\act c,
\qquad\quad
m\act n\act c\quad\xrightarrow{\mu_{m,n,c}}\quad
(m\tensor n)\act c
\]
satisfying the suitable coherence axioms 
corresponding to those for graded monads.
\end{defn}

Therefore the notion of lax action can be obtained by 
\emph{relaxing} that of strict action, by systematically replacing 
equalities with morphisms, which in turn are subject to new 
coherence axioms; cf.~\cite{baez-dolan:categorification}.
Now, as a general phenomenon in the 2-monad theory, 
we have the following proposition.

\begin{prop}
\label{prop-transport}
Let $\valgp{\cat{M}\times\cat{A}}{\sact}{\cat{A}}$ be
a strict $\cat{M}$-category, $\cat{C}$ a category, and 
\[
    \begin{tikzpicture}
      \node (A)  at (0,0) {$\cat{C}$};
      \node (C) at (0,2) {$\cat{A}$};
      
      \draw [->] (A) to [bend left=50]node (dom)[auto,labelsize] {$l$}(C);
      \draw [->] (C) to [bend left=50]node (cod)[auto,labelsize] {$r$}(A);
      
      \node at (0,1) [rotate=180,font=\large]{$\vdash$};

    \end{tikzpicture}
\]
an adjunction between the category~$\cat{C}$ and 
the underlying category~$\cat{A}$ of 
$\valgp{\cat{M}\times\cat{A}}{\sact}{\cat{A}}$.
Then, the composite functor
\[
\cat{M}\times\cat{C}\quad\xrightarrow{\cat{M}\times l}\quad
\cat{M}\times\cat{A}\quad\xrightarrow{\sact}\quad
\cat{A}\quad\xrightarrow{r}\quad\cat{C}
\]
naturally has a structure of lax $\cat{M}$-action on $\cat{C}$.
\end{prop}

This leads us to the notion of \emph{resolution} of 
a lax monoidal action. 

\subsection{Resolutions}

\begin{defn} 
For a lax $\cat{M}$-action~$\act=(\act,\eta,\mu)$ 
on a category~$\cat{C}$,
define the category~$\Res{\act}$ as follows.
\begin{itemize}
\item An object of $\Res{\act}$ is a \defemph{resolution} of $\act$,
which is given by the following data:
\begin{itemize}
\item A strict $\cat{M}$-category~$\valgp{\cat{M}\times\cat{A}}{\sact}{\cat{A}}$.
\item An adjunction
\[
    \begin{tikzpicture}
      \node (A)  at (0,0) {$\cat{C}$};
      \node (C) at (0,2) {$\cat{A}$};
      
      \draw [->] (A) to [bend left=50]node (dom)[auto,labelsize] {$l$}(C);
      \draw [->] (C) to [bend left=50]node (cod)[auto,labelsize] {$r$}(A);
      
      \node at (0,1) [rotate=180,font=\large]{$\vdash$};

    \end{tikzpicture}
\]
\end{itemize}
These data must satisfy the condition that, 
by the procedure of Proposition~\ref{prop-transport}
they yield $(\act,\eta,\mu)$.
\item Suppose we have two resolutions of $\act$:
\begin{align*}
\rho\quad&=\quad
\left(\valgp{\cat{M}\times\cat{A}}{\sact}{\cat{A}},\ 
l\dashv r\colon\cat{A}\to\cat{C}\right),\\
\rho^\prime\quad&=\quad
\left(\valgp{\cat{M}\times\cat{A^\prime}}{\sact^\prime}{\cat{A^\prime}},\ 
l^\prime\dashv r^\prime\colon\cat{A}^\prime\to\cat{C}\right).
\end{align*}
A morphism of $\Res{\act}$ from $\rho$ to $\rho^\prime$ is a
morphism of strict $\cat{M}$-actions
\[
k\quad\colon \quad
\valgp{\cat{M}\times\cat{A}}{\sact}{\cat{A}}\quad\longrightarrow\quad
\valgp{\cat{M}\times\cat{A^\prime}}{\sact^\prime}{\cat{A^\prime}},
\]
i.e., a functor~$k\colon \cat{A}\to\cat{A^\prime}$ making the diagram
\[
\begin{tikzpicture}
      \node (TL)  at (0,0.75)  {$\cat{M}\times\cat{A}$};
      \node (TR)  at (3,0.75)  {$\cat{M}\times\cat{A^\prime}$};
      \node (BL)  at (0,-0.75) {$\cat{A}$};
      \node (BR)  at (3,-0.75) {$\cat{A^\prime}$};
      
      \draw [->] (TL) to node (T) [auto,labelsize]      {$\cat{M}\times k$} (TR);
      \draw [->] (TR) to node (R) [auto,labelsize]      {$\sact^\prime$} (BR);
      \draw [->] (TL) to node (L) [auto,swap,labelsize] {$\sact$}(BL);
      \draw [->] (BL) to node (B) [auto,swap,labelsize] {$k$}(BR);
\end{tikzpicture}
\]
commute, satisfying the following equations:
\begin{align*}
k\circ l\quad&=\quad l^\prime,\\
r\quad&=\quad r^\prime\circ k.\qedhere
\end{align*}
\end{itemize}
\end{defn}

A natural question to ask at this point is:
given an arbitrary lax $\cat{M}$-action~$\act$, does there exist a 
resolution of $\act$?
This problem generalizes the one of finding an adjunction
that generates an arbitrary monad (replace $\cat{M}$ by $1$),
of which there are two solutions obtained way back in 1960's,
one by Eilenberg--Moore~\cite{eilenberg-moore} and
one by Kleisli~\cite{kleisli}.
Actually, our constructions of Eilenberg--Moore and Kleisli
categories of graded monads, which generalize
constructions in \cite{eilenberg-moore} and
\cite{kleisli} respectively, provide answers to 
this generalized problem as well.

\subsection{The fibrational correspondence of adjunctions}
\label{subsec-fib-adj}
In order to solve the problem of finding resolutions of 
the lax $\cat{M}$-action~$\act$ via 
the Eilenberg--Moore and Kleisli constructions for graded monads,
we need the following 2-fibrational property~\cite{hermida:fib} 
of the 2-functor~$
\ppp\colon\Epp\to\TCat_2$ concerning the 
correspondence of adjunctions in the total 2-category and 
those in a fiber 2-category:
\begin{prop}
\label{prop-fibrational-adj-epp}
Let $\tcat{C}$ and $\tcat{A}$ be 2-categories, 
$L\dashv R\colon \tcat{A}\to\tcat{C}$ be a 2-adjunction
(adjunction in $\TCat$), 
and $C$ and $A$ be objects of $\tcat{C}$ and $\tcat{A}$
respectively.
Then there is a bijective correspondence between the following 
two notions:
\begin{itemize}
\item Adjunctions in $\Epp$ between $(\tcat{C},C)$ and $(\tcat{A},A)$
above $L\dashv R$:
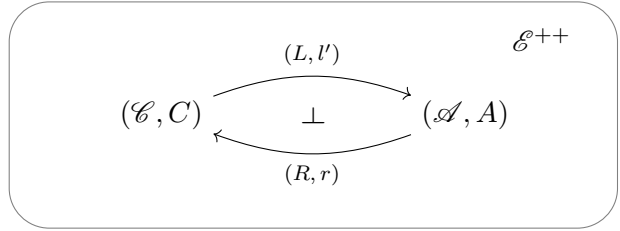
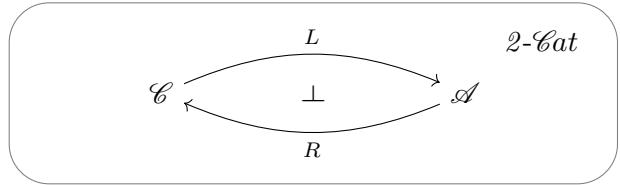
    \[
    \begin{tikzpicture}
      \node (B) at (0,0) {$(\tcat{C},C)$};
      \node (C) at (4,0) {$(\tcat{A},A)$};
      
      \node at (2,0) [rotate=90,font=\large]{$\vdash$};
      
      \draw [<-] (C) to [bend right=20]node [auto,swap,labelsize] {$(L,l^\prime)$}(B);
      \draw [<-] (B) to [bend right=20]node [auto,swap,labelsize] {$(R,r)$}(C);
      
      \draw [rounded corners=15pt,myborder] (-2,-1.5) rectangle (6,1.5);
      \node at (5,1) {$\Epp$};

     \begin{scope}[shift={(0,-3.2)}]
       \node (B) at (0,0) {$\tcat{C}$};
       \node (C) at (4,0) {$\tcat{A}$};
      
       \node at (2,0) [rotate=90,font=\large]{$\vdash$};
      
       \draw [<-] (C) to [bend right=23]node [auto,swap,labelsize] (Fpp) {$L$}(B);
       \draw [<-] (B) to [bend right=23]node [auto,swap,labelsize] (Fpp) {$R$}(C);
       \draw [rounded corners=15pt,myborder] (-2,-1.2) rectangle (6,1.2);
             \node at (5,0.7) {$\TCat$};
   
      \end{scope}
    \end{tikzpicture}
    \]
\item Adjunctions in $\tcat{C}$ between $C$ and $RA$:
    
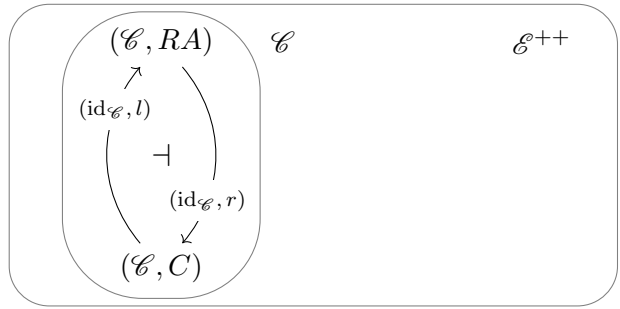
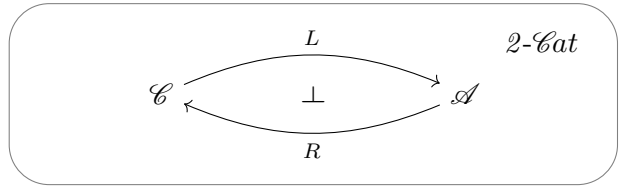
\[
    \begin{tikzpicture}
      \node (B) at (0,-1.5) {$(\tcat{C},C)$};
      \node (C) at (0,1.5) {$(\tcat{C},RA)$};
      
      \node at (0,0) [rotate=180,font=\large]{$\vdash$};
      
       \draw [rounded corners=30pt,myborder] (-1.3,-1.9) rectangle (1.3,1.9);
     \node at (1.6,1.5) {$\tcat{C}$};
      
      \draw [<-] (C) to [bend right=40]node [near start,labelsize,fill=white] {$(\id{\tcat{C}},l)$}(B);
      \draw [<-] (B) to [bend right=40]node [near start,labelsize,fill=white] {$(\id{\tcat{C}},r)$}(C);
      
      \draw [rounded corners=15pt,myborder] (-2,-2) rectangle (6,2);
      \node at (5,1.5) {$\Epp$};
   
      \begin{scope}[shift={(0,-3.7)}]
       \node (B) at (0,0) {$\tcat{C}$};
       \node (C) at (4,0) {$\tcat{A}$};
      
       \node at (2,0) [rotate=90,font=\large]{$\vdash$};
      
       \draw [<-] (C) to [bend right=23]node [auto,swap,labelsize] (Fpp) {$L$}(B);
       \draw [<-] (B) to [bend right=23]node [auto,swap,labelsize] (Fpp) {$R$}(C);
       \draw [rounded corners=15pt,myborder] (-2,-1.2) rectangle (6,1.2);
             \node at (5,0.7) {$\TCat$};
   
      \end{scope}
    \end{tikzpicture}
    \]
\end{itemize}
\end{prop}
\begin{pf}
Indeed, the former notion is given by the following data:
\begin{itemize}
\item A 1-cell~$l^\prime\colon LC\to A$ of $\tcat{A}$.
\item A 1-cell~$r\colon RA\to C$ of $\tcat{C}$.
\item A 2-cell~$\eta$ of $\tcat{C}$ of the following type:
    
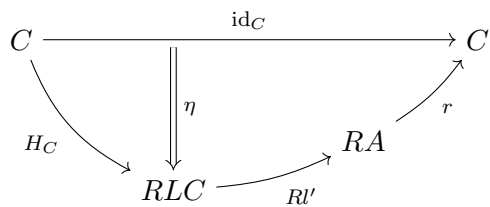
\[
    \begin{tikzpicture}
      \node (L) at (0,0) {$C$};
      \node (R) at (6,0) {$C$};
      \node (T) at (2,-2) {$RLC$};
      \node (N) at (4.5,-1.35) {$RA$};
      
      \draw [->] (L) to node [auto,labelsize] {$\id{{C}}$}(R);
      \draw [->] (L) to [bend right=20]node [auto,swap,labelsize] {$H_{C}$}(T);
      \draw [->] (T) to [bend right=10]node [auto,swap,labelsize] {$Rl^\prime$}(N);
      \draw [->] (N) to [bend right=10]node [auto,swap,labelsize] {$r$}(R);
            
      \draw [2cellr] (T) to node [auto,labelsize,swap]{$\eta$}(2,-0.1);
    \end{tikzpicture}
    \]
\item A 2-cell~$\varepsilon^\prime$ of $\tcat{A}$ of the following type:
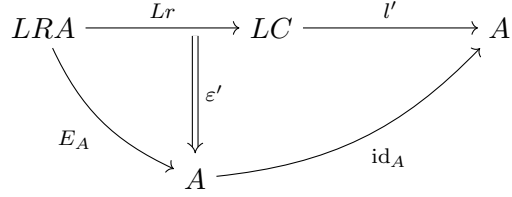
\[
\begin{tikzpicture}[baseline=-\the\dimexpr\fontdimen22\textfont2\relax ]
        \begin{scope} 
      \node (L)  at (0,0)  {$LRA$};
      \node (M)  at (3,0)  {$LC$};
      \node (R)  at (6,0) {$A$};
      \node (B)  at (2,-2) {$A$};
      \draw [->] (L) to node  [auto,labelsize]      {$Lr$} (M);
      \draw [->] (M) to node  [auto,labelsize]      {$l^\prime$} (R);
      \draw [->] (L) to [bend right=20]node  [auto,swap,labelsize] {$E_A$}(B);
      \draw [->] (B) to [bend right=20]node  [auto,swap,labelsize] {$\id{A}$}(R);
      \draw [2cellnode] (2,0) to node[auto,labelsize] {$\varepsilon^\prime$}(B);
        \end{scope}
\end{tikzpicture}
\]
\end{itemize}
Here, $H$ and $E$ are the unit and counit of the 2-adjunction~$
L\dashv R$ respectively. 

On the other hand, the latter notion is given by the following data:
\begin{itemize}
\item A 1-cell~$l\colon C\to RA$ of $\tcat{C}$.
\item A 1-cell~$r\colon RA\to C$ of $\tcat{C}$.
\item A 2-cell~$\eta$ of $\tcat{C}$ of the following type:
    
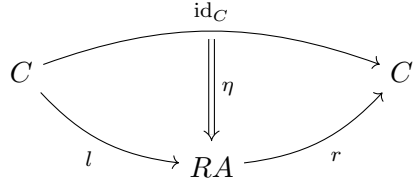
\[
    \begin{tikzpicture}
      \node (A)  at (0,0) {$C$};
      \node (B) at (5,0) {$C$};
      \node (C) at (2.5,-1.25) {$RA$};
      
      \draw [->] (A) to [bend right=20]node [auto,swap,labelsize] {$l$}(C);
      \draw [->] (C) to [bend right=20]node [auto,swap,labelsize] {$r$}(B);

      \draw [->] (A) to [bend left=20]node [auto,labelsize] (Fpp) {$\id{C}$}(B);
      \draw [2cellnoder] (C) to node [auto,labelsize,swap]{$\eta$}(Fpp);

    \end{tikzpicture}
    \]
\item A 2-cell~$\varepsilon$ of $\tcat{C}$ of the following type:
    
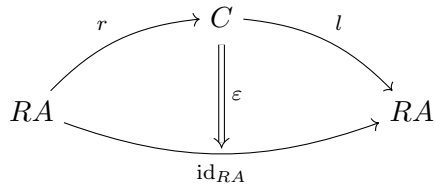
\[
    \begin{tikzpicture}
      \node (A)  at (0,0) {$RA$};
      \node (B) at (5,0) {$RA$};
      \node (C) at (2.5,1.25) {$C$};
      
      \draw [->] (A) to [bend right=-20]node [auto,labelsize] {$r$}(C);
      \draw [->] (C) to [bend right=-20]node [auto,labelsize] {$l$}(B);

      \draw [->] (A) to [bend left=-20]node [auto,swap,labelsize] (Fpp) {$\id{RA}$}(B);
      \draw [2cellnode] (C) to node [auto,labelsize]{$\varepsilon$}(Fpp);

    \end{tikzpicture}
    \]
\end{itemize}

Now the correspondence should be clear; $l^\prime$ and $\varepsilon^\prime$
correspond respectively to $l$ and $\varepsilon$ under the 
2-adjunction~$L\dashv R$.
That this correspondence preserves and reflects
the triangular identity is
straightforward to check.
\end{pf}
We also have an analogous result for 
$\pmm\colon\Emm\to\TCat_2^\op{1,2}$:
\begin{prop}
\label{prop-fibrational-adj-emm}
Let $\tcat{C}$ and $\tcat{A}$ be 2-categories, 
$L\dashv R\colon \tcat{C}\to\tcat{A}$ be a 2-adjunction
(adjunction in $\TCat$), 
and $C$ and $A$ be objects of $\tcat{C}$ and $\tcat{A}$
respectively.
Then there is a bijective correspondence between the following 
two notions:
\begin{itemize}
\item Adjunctions in $\Emm$ between $(\tcat{C},C)$ and $(\tcat{A},A)$
above $L\dashv R$:
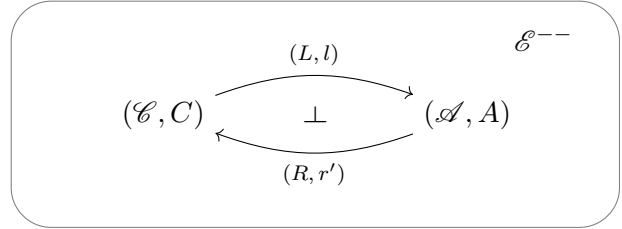
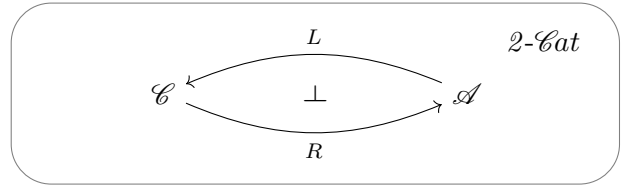
    \[
    \begin{tikzpicture}
      \node (B) at (0,0) {$(\tcat{C},C)$};
      \node (C) at (4,0) {$(\tcat{A},A)$};
      
      \node at (2,0) [rotate=90,font=\large]{$\vdash$};
      
      \draw [<-] (C) to [bend right=20]node [auto,swap,labelsize] {$(L,l)$}(B);
      \draw [<-] (B) to [bend right=20]node [auto,swap,labelsize] {$(R,r^\prime)$}(C);
      
      \draw [rounded corners=15pt,myborder] (-2,-1.5) rectangle (6,1.5);
      \node at (5,1) {$\Emm$};

     \begin{scope}[shift={(0,-3.2)}]
       \node (B) at (0,0) {$\tcat{C}$};
       \node (C) at (4,0) {$\tcat{A}$};
      
       \node at (2,0) [rotate=90,font=\large]{$\vdash$};
      
       \draw [->] (C) to [bend right=23]node [auto,swap,labelsize] (Fpp) {$L$}(B);
       \draw [->] (B) to [bend right=23]node [auto,swap,labelsize] (Fpp) {$R$}(C);
       \draw [rounded corners=15pt,myborder] (-2,-1.2) rectangle (6,1.2);
             \node at (5,0.7) {$\TCat$};
   
      \end{scope}
    \end{tikzpicture}
    \]
\item Adjunctions in $\tcat{C}$ between $C$ and $LA$:
    
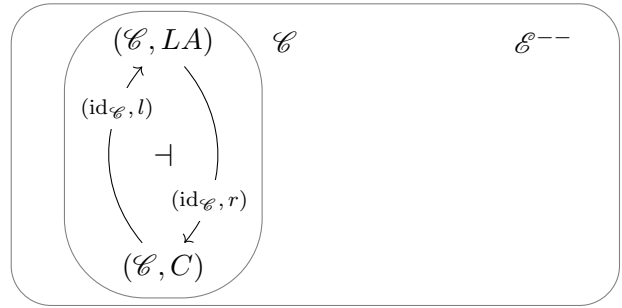
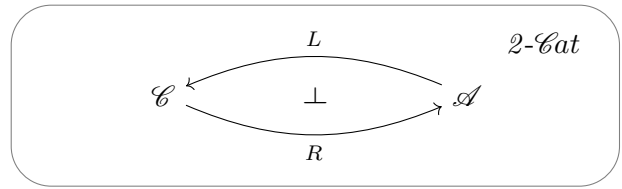
\[
    \begin{tikzpicture}
      \node (B) at (0,-1.5) {$(\tcat{C},C)$};
      \node (C) at (0,1.5) {$(\tcat{C},LA)$};
      
      \node at (0,0) [rotate=180,font=\large]{$\vdash$};
      
       \draw [rounded corners=30pt,myborder] (-1.3,-1.9) rectangle (1.3,1.9);
     \node at (1.6,1.5) {$\tcat{C}$};
      
      \draw [<-] (C) to [bend right=40]node [near start,labelsize,fill=white] {$(\id{\tcat{C}},l)$}(B);
      \draw [<-] (B) to [bend right=40]node [near start,labelsize,fill=white] {$(\id{\tcat{C}},r)$}(C);
      
      \draw [rounded corners=15pt,myborder] (-2,-2) rectangle (6,2);
      \node at (5,1.5) {$\Emm$};
   
      \begin{scope}[shift={(0,-3.7)}]
       \node (B) at (0,0) {$\tcat{C}$};
       \node (C) at (4,0) {$\tcat{A}$};
      
       \node at (2,0) [rotate=90,font=\large]{$\vdash$};
      
       \draw [->] (C) to [bend right=23]node [auto,swap,labelsize] (Fpp) {$L$}(B);
       \draw [->] (B) to [bend right=23]node [auto,swap,labelsize] (Fpp) {$R$}(C);
       \draw [rounded corners=15pt,myborder] (-2,-1.2) rectangle (6,1.2);
             \node at (5,0.7) {$\TCat$};
   
      \end{scope}
    \end{tikzpicture}
    \]
\end{itemize}
\end{prop}

\subsection{Existence of the terminal and initial resolutions}
Now we are ready to connect the Eilenberg--Moore and Kleisli
constructions for graded monads to the notion of 
resolutions of a lax action.
Let us fix an $\cat{M}$-graded monad~$\gm{T}$, or equivalently
a lax $\cat{M}$-action~$\act$, on a category $\cat{C}$.

We begin with the case of the 
Eilenberg--Moore construction.
Applying Proposition~\ref{prop-fibrational-adj-epp} to 
the free-forgetful 2-adjunction
\[
F^\cat{M}\ \dashv\ U^\cat{M}\quad\colon\quad
\EM{\Cat}{\cat{M}\times(-)}\quad\longrightarrow\quad\Cat,
\]
and objects~$\cat{C}$ and
$\valgp{\cat{M}\times\EMg}{\emact}{\EMg}$ of
$\Cat$ and $\EM{\Cat}{\cat{M}\times(-)}$ respectively,
we obtain from the Eilenberg--Moore adjunction for $\gm{T}$
\[
\begin{tikzpicture}
    \node (E) at (-6,0) {$(\Cat,\cat{C})$};
    \node (A)  at (0,0) {$\left(\EM{\Cat}{\cat{M}\times(-)},
        \valg{\cat{M}\times\EMg}{\emact}{\EMg}\right)$};

    \draw [<-] ([yshift=3mm]A.west) to [bend left=-20]node(dom) [auto,swap,labelsize] {$(F^\cat{M},f^\gm{T})$} ([yshift=3mm]E.east);
    \draw [->] ([yshift=-3mm]A.west) to [bend right=-20]node(cod) [auto,labelsize] {$(U^\cat{M},u^\gm{T})$}([yshift=-3mm]E.east);

    \draw [draw=none](dom) to node [rotate=90,font=\large]{$\vdash$} (cod);
\end{tikzpicture}
\]
the following resolution:

\begin{defn}
The \defemph{Eilenberg--Moore resolution} for the lax $\cat{M}$-action~$\act$ is 
given by the strict $\cat{M}$-category~$
\valgp{\cat{M}\times\EMg}{\emact}{\EMg}$
and the adjunction
\[
    \begin{tikzpicture}
      \node (A) at (0,0) {$\cat{C}$};
      \node (C) at (0,2) {$\EMg$};
      
      \draw [->] (A) to [bend left=50] (C);
      \draw [->] (C) to [bend left=50] (A);
      
      \node[labelsize] at (-1.45,1) {$f^\gm{T}\circ H^\cat{M}_\cat{C}$};
      \node[labelsize] at (1,1.05)  {$u^\gm{T}$};
      
      \node at (0,1) [rotate=180,font=\large]{$\vdash$};

    \end{tikzpicture}
\]
between categories~$\cat{C}$ and $\EMg$.
\end{defn}

Similarly, using Proposition~\ref{prop-fibrational-adj-emm}
to the forgetful-cofree 2-adjunction
\[
F_\cat{M}\ \dashv\ U_\cat{M}\quad\colon\quad
\Cat\quad\longrightarrow\quad\EM{\Cat}{[\cat{M},-]},
\]
and objects~$\cat{C}$ and
$\vcoalgp{[\cat{M},\Klg]}{\klact}{\Klg}$ of
$\Cat$ and $\EM{\Cat}{[\cat{M},-]}$ respectively,
it follows that the Kleisli adjunction for $\gm{T}$
\[
\begin{tikzpicture}
    \node (E) at (-6,0) {$(\Cat,\cat{C})$};
    \node (A)  at (0,0) {$\left(\EM{\Cat}{[\cat{M},-]},
        \vcoalg{[\cat{M},\Klg]}{\klact}{\Klg}\right)$};

    \draw [<-] ([yshift=3mm]A.west) to [bend left=-20]node(dom) [auto,swap,labelsize] {$(F_\cat{M},f_\gm{T})$} ([yshift=3mm]E.east);
    \draw [->] ([yshift=-3mm]A.west) to [bend right=-20]node(cod) [auto,labelsize] {$(U_\cat{M},u_\gm{T})$}([yshift=-3mm]E.east);

    \draw [draw=none](dom) to node [rotate=90,font=\large]{$\vdash$} (cod);
\end{tikzpicture}
\]
gives rise to:

\begin{defn}
The \defemph{Kleisli resolution} for the lax $\cat{M}$-action~$\act$ 
is given by the 
strict $\cat{M}$-category~$\vcoalgp{[\cat{M},\Klg]}{\klact}{\Klg}$
and the adjunction
\[
    \begin{tikzpicture}
      \node (A) at (0,0) {$\cat{C}$};
      \node (C) at (0,2) {$\Klg$};
      
      \draw [->] (A) to [bend left=50] (C);
      \draw [->] (C) to [bend left=50] (A);
      
      \node[labelsize] at (-1,1) {$f_\gm{T}$};
      \node[labelsize] at (1.5,1)  {$H_{\cat{M},\cat{C}}\circ u_\gm{T}$};
      
      \node at (0,1) [rotate=180,font=\large]{$\vdash$};

    \end{tikzpicture}
\]
between categories~$\cat{C}$ and $\Klg$.
\end{defn}

Moreover, as an easy corollary of the comparison theorems 
(Propositions~\ref{prop-comparison-em-gm} and 
\ref{prop-comparison-kl-gm}),
we may conclude:
\begin{thm}
The category~$\Res{\act}$ has both the terminal and initial objects,
given respectively by
the Eilenberg--Moore and Kleisli resolutions.
\end{thm}

\section{A construction of Maillard and Melli\`es}
\label{sec-constr-maillard-mellies}
Maillard and Melli\`es~\cite{maillard-mellies} introduced 
\emph{indexed monads}, which are actually somewhat more 
general than what we call indexed monads here.
They also introduced a construction which produces a  
2-fibration over a 2-category~$\tcat{B}$ for each 
$\tcat{B}$-indexed monad in their sense.
Interestingly, our construction of the Eilenberg--Moore categories of 
indexed monads turns out to constitute particular instances of 
their construction.
In this chapter, we see how the result of \cite{maillard-mellies}
connecting their construction to the notion of
\emph{model of an indexed Lawvere theory}
introduced by Power~\cite{power:lce,power:indexed}, can be 
understood in the light of notions related to the
Eilenberg--Moore construction for indexed monads.

\subsection{Indexed Lawvere theories and their models}

In his investigation of
the relationship between the \emph{global state monad} and the
\emph{local state monad}~\cite{plotkin-power}, 
Power \cite{power:indexed,power:lce} introduced the notion 
of \emph{indexed Lawvere theory}, whose definition is 
recalled below.

\begin{defn}
Let $\cat{B}$ be a category.
A \defemph{$\cat{B}$-indexed Lawvere theory} is a functor
\[
\il{L}\quad\colon\quad\cat{B}\quad\longrightarrow\quad\Law.\qedhere
\]
\end{defn}
Here, $\Law$ is the category consisting of Lawvere theories and 
maps of them; see \cite{hyland-power} for the detailed definition.
Thanks to the well-known correspondence of (ordinary) 
Lawvere theories and 
finitary monads on $\Set$~\cite{power:enriched}, i.e., the  
inclusion functor~$\iota\colon\Law\to\Mnd{\Set}^\opp$, 
it turns out that every 
$\cat{B}$-indexed Lawvere theory defines a $\cat{B}$-indexed
monad on $\Set$ (in our sense) by postcomposing $\iota$, 
as observed in \cite{maillard-mellies}.

Power also defined \emph{models} of an indexed Lawvere 
theory, generalizing the classical notion of model of 
a Lawvere theory.
\begin{defn}
\label{defn-ind-Law}
Let $\cat{B}$ be a category, $\il{L}$ a $\cat{B}$-indexed Lawvere theory,
and $\cat{C}$ a category with finite products. 
Define the category~$\mdl{\il{L}}{\cat{C}}$ as follows:
\begin{itemize}
\item An object of $\mdl{\il{L}}{\cat{C}}$ is 
a \defemph{model~$M$ of $\il{L}$ in $\cat{C}$}.
It consists of the following data:
\begin{itemize}
\item For each object~$b$ of $\cat{B}$, a model~$M_b$ of the Lawvere
theory~$\il{L}_b$ in $\cat{C}$, i.e., a finite product preserving 
functor~$M_b\colon \il{L}_b\to\cat{C}$.
\item For each morphism~$u\colon b\to b^\prime$ of $\cat{B}$, a 
natural transformation of type
\[
    \begin{tikzpicture}[baseline=-\the\dimexpr\fontdimen22\textfont2\relax ]
            \begin{scope} [shift={(0,0.75)}]
      \node (A)  at (0,0) {$\cat{C}$};
      \node (B) at (-3,-1.5) {$\il{L}_{b^\prime}$};
      \node (C) at (-3,0) {$\il{L}_b$};
          
      \draw [->] (C) to [bend left=0]node [auto,labelsize] {$M_b$}(A);
      \draw [->] (C) to [bend left=0]node [auto,swap,labelsize] {$\il{L}_u$}(B);

      \draw [->] (B) to [bend right=20]node [auto,swap,labelsize] (Fpp) {$M_{b^\prime}$}(A);
      \draw [2cell] (-1.5,-0.1) to node [auto,labelsize]{$M_u$}(-1.5,-1.1);
            \end{scope}
    \end{tikzpicture}
\]
\end{itemize}
These data are subject to the following axioms:
\begin{itemize}
\item $M_{\id{b}}=\id{M_b}$ for each object~$b$ of $\cat{B}$.
\item    $ \begin{tikzpicture}[baseline=-\the\dimexpr\fontdimen22\textfont2\relax ]
            \begin{scope} [shift={(0,0)}]
      \node (A)  at (0,0) {$\cat{C}$};
      \node (B) at (-3,-1.5) {$\il{L}_{b^{\prime\prime}}$};
      \node (C) at (-3,1.5) {$\il{L}_b$};
          
      \draw [->] (C) to [bend left=20]node [auto,labelsize] {$M_b$}(A);
      \draw [->] (C) to [bend left=0]node [auto,swap,labelsize] {$\il{L}_{u^\prime\circ u}$}(B);

      \draw [->] (B) to [bend right=20]node [auto,swap,labelsize] (Fpp) {$M_{b^{\prime\prime}}$}(A);
      \draw [2cell] (-1.5,1.1) to node [auto,labelsize]{$M_{u^\prime\circ u}$}(-1.5,-1.1);
            \end{scope}
    \end{tikzpicture}
    \ =\ 
        \begin{tikzpicture}[baseline=-\the\dimexpr\fontdimen22\textfont2\relax ]
                \begin{scope} [shift={(0,0)}]
          \node (A)  at (0,0) {$\cat{C}$};
          \node (B) at (-3,-1.5) {$\il{L}_{b^{\prime\prime}}$};
          \node (C) at (-3,0) {$\il{L}_{b^\prime}$};
          \node (D) at (-3,1.5) {$\il{L}_b$};
              
          \draw [->] (C) to [bend left=0]node [auto,near start,labelsize] {$M_{b^\prime}$}(A);
          \draw [->] (C) to [bend left=0]node [auto,swap,labelsize] {$\il{L}_{u^\prime}$}(B);
          \draw [->] (D) to [bend left=20]node [auto,labelsize] {$M_b$}(A);
          \draw [->] (D) to [bend left=0]node [auto,swap,labelsize] {$\il{L}_u$}(C);

          \draw [->] (B) to [bend right=20]node [auto,swap,labelsize] (Fpp) {$M_{b^{\prime\prime}}$}(A);
          \draw [2cell] (-1.5,-0.1) to node [auto,labelsize]{$M_{u^\prime}$}(-1.5,-1.1);
          \draw [2cell] (-1.5,1.1) to node [auto,labelsize]{$M_u$}(-1.5,0.1);
                \end{scope}
        \end{tikzpicture}$
        for each composable pair of morphisms~$u\colon
        b\to b^\prime,u^\prime\colon b^\prime\to b^{\prime\prime}$ 
        of $\cat{B}$.
\end{itemize}
\item 
A morphism~$\varphi\colon M\to M^\prime$ of $\mdl{\il{L}}{\cat{C}}$
is given by a family of 
natural transformations~$\varphi_b\colon M_b\Rightarrow M^\prime_b
\colon\il{L}_b\to\cat{C}$
for each $b\in\cat{B}$, such that
\[
    \begin{tikzpicture}[baseline=-\the\dimexpr\fontdimen22\textfont2\relax ]
            \begin{scope} [shift={(0,0.75)}]
      \node (A)  at (0,0) {$\cat{C}$};
      \node (B) at (-3,-1.5) {$\il{L}_{b^\prime}$};
      \node (C) at (-3,0) {$\il{L}_b$};
      \draw [<-] (A) to [bend left=-20]node [auto,swap,labelsize] {$M_b$}(C);
      \draw [<-] (A) to [bend right=-20]node [auto,near end,labelsize] {$M_b^\prime$}(C);
      \draw [->] (C) to [bend left=0]node [auto,swap,labelsize] {$\il{L}_u$}(B);
      %
      \draw [<-] (A) to [bend right=-20]node [auto,labelsize] (Fpp) {$M_{b^\prime}^\prime$}(B);
      \draw [2cell] (-1.5,0.25) to node [auto,labelsize]{$\varphi_b$}(-1.5,-0.3);
      \draw [2cell] (-1.5,-0.45) to node [auto,near start,labelsize]{$M_u^\prime$}(-1.5,-1.1);
            \end{scope}
    \end{tikzpicture}
    \quad=\quad
    \begin{tikzpicture}[baseline=-\the\dimexpr\fontdimen22\textfont2\relax ]
            \begin{scope} [shift={(0,0.75)}]
      \node (A)  at (0,0) {$\cat{C}$};
      \node (B) at (-3,-1.5) {$\il{L}_{b^\prime}$};
      \node (C) at (-3,0) {$\il{L}_b$};
      \draw [<-] (A) to [bend left=-20]node [auto,swap,labelsize] {$M_b$}(C);
      \draw [<-] (A) to [bend left=-15]node [auto,swap,near end,labelsize] {$M_{b^\prime}$}(B);
      \draw [->] (C) to [bend left=0]node [auto,swap,labelsize] {$\il{L}_u$}(B);
      %
      \draw [<-] (A) to [bend right=-20]node [auto,labelsize] (Fpp) {$M_{b^\prime}^\prime$}(B);
      \draw [2cell] (-1.5,0.25) to node [auto,labelsize]{$M_u$}(-1.5,-0.3);
      \draw [2cell] (-1.5,-0.5) to node [auto,labelsize]{$\varphi_{b^\prime}$}(-1.5,-1.1);
            \end{scope}
    \end{tikzpicture}
\]
holds for each morphism~$u\colon b\to b^\prime$ of $\cat{B}$.\qedhere
\end{itemize}
\end{defn}

One of the most striking results on the classical correspondence between 
Lawvere theories~${L}$ and finitary monads~$T_{L}=\iota(L)$ on 
$\Set$ says that 
there is an equivalence of categories between 
$\mdl{L}{\Set}$ and $\EM{\Set}{T_L}$.
Now we claim that this intimate relation between 
the \emph{category of models of $L$ in $\Set$}
and the
\emph{Eilenberg--Moore category of $T_L$}
generalizes to the \emph{indexed} setting,
in a somewhat nontrivial manner;
this is exactly what we intend to show in the current section.

\subsection{The base change adjoint triple}
\label{subsec-base-change}
Recall from Section~\ref{sec-iem} that given a $\cat{B}$-indexed monad~$
\im{T}$ on the category~$\cat{C}$, one may find its Eilenberg--Moore
object by considering $\im{T}$ as a monad in the 2-category~$\Epm$, and
is given as 
$\left(\EM{\Cat}{\cat{B}\times(-)},\vcoalg{\cat{B}\times\EMi}{\pi}{\EMi}\right)\in\Epm$.
As remarked in Section~\ref{sec-iem}, the 2-category~$
\EM{\Cat}{\cat{B}\times(-)}$ is isomorphic to the slice 2-category~$
\Cat/\cat{B}$,
by the 2-functor~$\vcoalgp{\cat{B}\times\cat{A}}{\alpha}{\cat{A}}\mapsto\vcoalgp{\cat{B}}{\alpha_0}{\cat{A}}$. 
Let us now remember the 
\emph{base change adjoint triple}
\[
    \begin{tikzpicture}[baseline=-\the\dimexpr\fontdimen22\textfont2\relax ]
            \begin{scope} [shift={(0,0.75)}]
      \node (A)  at (0,0) {$\Cat/\cat{B}$};
      \node (B) at (0,-2) {$\Cat$};
      
      \node [labelsize] at (-1.25,-1) {$\coprod_\cat{B}$};
      \node [labelsize] at (1.25,-1) {$\prod_\cat{B}$};
      
      \node at (-0.45,-1) [font=\large]{$\dashv$};
      \node at (0.45,-1) [font=\large]{$\dashv$};

      \draw [->] (A) to [bend left=-60](B);
      \draw [->] (B) to [bend left=0]node [fill=white,labelsize] {$\cat{B}^\ast$}(A);
      \draw [->] (A) to [bend left=60](B);
            \end{scope}
    \end{tikzpicture}
\]
connecting the 2-categories~$\Cat$ and $\Cat/\cat{B}$.
This notion allows us to nicely describe the 
relation between the category of models and the
Eilenberg--Moore category;
before recalling the precise definition of the adjoint triple,
we state our main theorem in the current section:
\begin{thm}
\label{thm-model-section}
Let $\cat{B}$ be a category and $\il{L}$ a $\cat{B}$-indexed 
Lawvere theory; note that $\il{L}$ determines a 
$\cat{B}$-indexed monad~$\iota \il{L}$ on $\Set$.
There is an equivalence of categories between
$\mdl{\il{L}}{\Set}$ and 
$\dprod\vcoalgp{\cat{B}}{\pi_0}{\EM{\Set}{\iota\il{L}}}$.
\end{thm}

\begin{defn}
The 2-functor~$\coprod_{\cat{B}}\colon\Cat/\cat{B}\to\Cat$ 
is, up to the isomorphism $\EM{\Cat}{\cat{B}\times(-)}
\cong \Cat/\cat{B}$, the forgetful 2-functor~$U^\cat{B}
\colon \EM{\Cat}{\cat{B}\times(-)}\to\Cat$ defined in
Definition~\ref{defn-U-B-times}:
$\vcoalgp{\cat{B}}{p}{\cat{E}}\mapsto \cat{E}$.
\end{defn}

\begin{defn}
The 2-functor~$\cat{B}^\ast\colon \Cat\to\Cat/\cat{B}$ is, up to
the isomorphism $\EM{\Cat}{\cat{B}\times(-)}
\cong \Cat/\cat{B}$, the cofree 2-functor~$F^\cat{B}
\colon \Cat\to\EM{\Cat}{\cat{B}\times(-)}$ defined in
Definition~\ref{defn-F-B-times}:
$\cat{X}\mapsto \vcoalgp{\cat{B}}{M^\cat{B}_{\cat{X},0}}{\cat{B}\times\cat{X}}$.
\end{defn}

So the 2-adjunction~$\coprod_\cat{B}\dashv \cat{B}^\ast\colon
\Cat\to\Cat/\cat{B}$ is, essentially, the 
co-Eilenberg--Moore 2-adjunction~$
U^\cat{B}\dashv F^\cat{B}\colon \Cat\to\EM{\Cat}{\cat{B}\times(-)}$
for the 2-comonad~$\cat{B}\times(-)$ on which we heavily relied
when performing the Eilenberg--Moore construction for indexed monads.
The following construction is new.

\begin{defn}
We define the 2-functor~$\prod_{\cat{B}}\colon\Cat/\cat{B}\to\Cat$ 
in the following way.
\begin{itemize}
\item Given an object~$\vcoalgp{\cat{B}}{p}{\cat{E}}$ of 
$\Cat/\cat{B}$, define the category~$\prod_{\cat{B}}\vcoalgp{\cat{B}}{p}{\cat{E}}$ as follows:
\begin{itemize}
\item Its object is a \defemph{section} of $p$, 
i.e., a functor~$s\colon \cat{B}\to\cat{E}$
satisfying $p\circ s=\id{\cat{B}}$.
\item Its morphism from $s$ to $s^\prime$ is a natural transformation~$
\psi\colon s\Rightarrow s^\prime\colon\cat{B}\to\cat{E}$
satisfying $p\hc \psi=\id{\cat{B}}$.
\end{itemize}
\item Given a morphism~$
f\colon\vcoalgp{\cat{B}}{p}{\cat{E}}\to\vcoalgp{\cat{B}}{p^\prime}{\cat{E}^\prime}$
of $\Cat/\cat{B}$, i.e., a functor~$f\colon\cat{E}\to\cat{E}^\prime$ 
making the diagram
\[
\begin{tikzpicture}
      \node (TL)  at (0,0.75)  {$\cat{E}$};
      \node (TR)  at (2,0.75)  {$\cat{E}^\prime$};
      \node (BR)  at (1,-0.75) {$\cat{B}$};
      
      \draw [->] (TL) to node (T) [auto,labelsize]      {$f$} (TR);
      \draw [->] (TR) to node (R) [auto,labelsize]      {$p^\prime$} (BR);
      \draw [->] (TL) to node (L) [auto,swap,labelsize] {$p$}(BR);
\end{tikzpicture}
\]
commute, define the functor~$\prod_{\cat{B}}f\colon
\prod_{\cat{B}}\vcoalgp{\cat{B}}{p}{\cat{E}}\to
\prod_{\cat{B}}\vcoalgp{\cat{B}}{p^\prime}{\cat{E}^\prime}$
as follows:
\begin{itemize}
\item It sends a section~$s$ of $p$ to the section~$f\circ s\colon
\cat{B}\to\cat{E}^\prime$ of $p^\prime$; observe that
$p^\prime\circ f\circ s=p\circ s=\id{\cat{B}}$.
\item It sends a morphism~$\varphi\colon s\to s^\prime$ of sections
of $p$ to $f\hc \varphi\colon f\circ s\to f\circ s^\prime$.
\end{itemize}
\item Given a 2-cell~$\alpha\colon f\Rightarrow f^\prime
\colon\vcoalgp{\cat{B}}{p}{\cat{E}}\to\vcoalgp{\cat{B}}{p^\prime}{\cat{E}^\prime}$ of $\Cat/\cat{B}$, i.e., a natural transformation~$\alpha
\colon f\Rightarrow f^\prime\colon\cat{E}\to\cat{E^\prime}$ 
satisfying the equation 
\[
\begin{tikzpicture}[baseline=-\the\dimexpr\fontdimen22\textfont2\relax ]
      \node (TL)  at (0,0.75)  {$\cat{E}$};
      \node (BR)  at (1,-0.75) {$\cat{B}$};
      
      \draw [->] (TL) to node (L) [auto,swap,labelsize] {$p$}(BR);
\end{tikzpicture}
\quad=\quad
\begin{tikzpicture}[baseline=-\the\dimexpr\fontdimen22\textfont2\relax ]
      \node (TL)  at (0,0.75)  {$\cat{E}$};
      \node (TR)  at (2,0.75)  {$\cat{E}^\prime$};
      \node (BR)  at (1,-0.75) {$\cat{B}$};
      
      \draw [->] (TL) to [bend left=25]node (dom) [auto,labelsize]      {$f$} (TR);
      \draw [->] (TL) to [bend right=25]node (cod) [auto,swap,labelsize]      {$f^\prime$} (TR);
      \draw [->] (TR) to node (R) [auto,labelsize]      {$p^\prime$} (BR);
      
      \draw [2cellnode] (dom) to node [auto,labelsize] {$\alpha$} (cod);
\end{tikzpicture}
\]
define the natural transformation~$\prod_{\cat{B}}\alpha
\colon \prod_{\cat{B}}f\Rightarrow\prod_{\cat{B}}f^\prime
\colon \prod_{\cat{B}}\vcoalgp{\cat{B}}{p}{\cat{E}}
\to\prod_{\cat{B}}\vcoalgp{\cat{B}}{p^\prime}{\cat{E}^\prime}$
as follows:
\begin{itemize}
\item Its $s$-component~$\left(\prod_{\cat{B}}\alpha\right)_s\colon
f\circ s\to f^\prime \circ s$ for a section $s$ of $p$ is 
given by $\left(\prod_{\cat{B}}\alpha\right)_s\coloneqq \alpha\hc s$.
\qedhere
\end{itemize}
\end{itemize}
\end{defn}

Or more concisely, 
\[
\dprod\quad\coloneqq
\quad\Cat/\cat{B}\left(\vcoalgp{\cat{B}}{\id{\cat{B}}}{\cat{B}},\ 
-\,\right)
\quad\colon\quad
\Cat/\cat{B}\quad\longrightarrow\quad\Cat.
\]

\subsection{Models as sections}
Let us finally prove Theorem~\ref{thm-model-section}, assuming
the classical equivalence between the category of models 
of a Lawvere theory in $\Set$ and the Eilenberg--Moore category 
of the corresponding finitary monad on $\Set$.
Below we give a more concrete description of the category~$ 
\dprod\vcoalgp{\cat{B}}{\pi_0}{\EM{\Set}{\iota\il{L}}}$.
To avoid too heavy notation we abbreviate the finitary 
monad~$T_{\iota\il{L}_b}$ on $\Set$ corresponding to the Lawvere 
theory~$\il{L}_b$ as $T_b$; similarly for the relevant monad morphisms.
\begin{itemize}
\item An object of 
$\dprod\vcoalgp{\cat{B}}{\pi_0}{\EM{\Set}{\iota\il{L}}}$ consists 
of the following data:
\begin{itemize}
\item For each object~$b$ of $\cat{B}$, a $T_b$-algebra~$
\valgp{T_b c_b}{\chi_b}{c_b}$.
\item For each morphism~$u\colon b\to b^\prime$ of $\cat{B}$, a
morphism of $T_b$-algebras
\[
h_u\quad\colon\quad\valgp{T_b c_b}{\chi_b}{c_b}
\quad\longrightarrow\quad
\valgp{T_b c_{b^\prime}}{\chi_{b^\prime}\circ 
T_{u,c_{b^\prime}}}{c_{b^\prime}}
\]
\end{itemize}
These data are subject to the following \emph{functoriality} axioms:
\begin{itemize}
\item $h_{\id{b}}=\id{\chi_b}$ for each object~$b$ of $\cat{B}$.
\item $h_{u^\prime\circ u}=h_{u^\prime}\circ h_u$ for each pair of
composable morphisms~$u,u^\prime$ of $\cat{B}$.
\end{itemize}
\item A morphism 
$\psi\colon \big((\chi_b),(h_u)\big)\to 
\big((\chi^\prime_b),(h^\prime_u)\big)$ of $\dprod\vcoalgp{\cat{B}}{\pi_0}{\EM{\Set}{\iota\il{L}}}$
is a family of morphisms of $T_b$-algebras
\[
\psi_b\quad\colon\quad
\valgp{T_b c_b}{\chi_b}{c_b}
\quad\longrightarrow\quad
\valgp{T_b c^\prime_b}{\chi^\prime_b}{c^\prime_b}
\]
for each~$b\in\cat{B}$, such that
the following \emph{naturality square}
\[
\begin{tikzpicture}
      \node (TL)  at (0,1.5)  {$\valgp{T_b c_b}{\chi_b}{c_b}$};
      \node (TR)  at (5,1.5)  {$\valgp{T_b c^\prime_b}{\chi^\prime_b}{c^\prime_b}$};
      \node (BL)  at (0,-1) {$\valgp{T_b c_b^\prime}{\chi_b^\prime\circ T_{u,c_{b^\prime}}}{c_b^\prime}$};
      \node (BR)  at (5,-1) {$\valgp{T_b c^\prime_{b^\prime}}{\chi^\prime_{b^\prime}\circ T_{u,c^\prime_{b^\prime}}}{c^\prime_{b^\prime}}$};
      
      \draw [->] (TL) to node (T) [auto,labelsize]      {$\psi_b$} (TR);
      \draw [->] (TR) to node (R) [auto,labelsize]      {$h^\prime_u$} (BR);
      \draw [->] (TL) to node (L) [auto,swap,labelsize] {$h_u$}(BL);
      \draw [->] (BL) to node (B) [auto,swap,labelsize] {$\psi_{b^\prime}$}(BR);
\end{tikzpicture}
\]
commutes for each morphism~$u\colon b\to b^\prime$ of $\cat{B}$.
\end{itemize}

Now it only remains to apply the classical equivalence of models of 
a Lawvere theory in $\Set$ and algebras of the corresponding
finitary monad on $\Set$, before we reach Definition~\ref{defn-ind-Law}.

\section*{Notes}
The idea of decomposing a lax action into a strict action together with
an adjunction is presented in the paper~\cite{fujii-katsumata-mellies}.
Among the authors of \cite{fujii-katsumata-mellies}, this material 
has been mainly developed by Shin-ya Katsumata, and indeed
it was in this context that he arrived at the definitions of 
Eilenberg--Moore and Kleisli categories of graded monads.
(These are then informed to me, as mentioned in Notes for
Chapter~\ref{chap-main-constr}.)
Then I noticed that thanks to the 2-fibrational
property of $\ppp\colon\Epp\to\TCat_2$ and 
$\pmm\colon\Emm\to\TCat_2^\op{1,2}$ 
concerning the correspondence of adjunctions 
(Section~\ref{subsec-fib-adj}), 
decompositions (called \emph{resolutions} here) can be seen as 
instances of the more familiar 2-categorical situation of 
adjunctions that generate monads.
The adjoint correspondence theorem 
(Propositions~\ref{prop-fibrational-adj-epp}
and \ref{prop-fibrational-adj-emm}) seem not to have been stated
explicitly in the literature as far as I am aware,
but the similar result (for a particular 2-fibration) has been shown 
as the main theorem in the paper~\cite{hermida:fib} by Hermida.
I conjecture that this property is possessed more generally by 
an arbitrary 2-fibration.

The main conceptual novelty in Section~\ref{sec-constr-maillard-mellies}
is the observation brought by our Eilenberg--Moore construction
for indexed monads that, 
the functor~$\vcoalgp{\cat{B}}{\pi_0}{\EMi}$ naturally lives in
the 2-category~$\Cat/\cat{B}$ (rather than other alternatives such 
as $\Fib{\cat{B}}$).
This enables us to identify the \emph{dependent product 2-functor}~$
\dprod\colon\Cat/\cat{B}\to\Cat$ as the right construction of 
the ``category of sections''.
Note that this actually corrects a subtle mistake in the paper~\cite{maillard-mellies}
concerning the definition of morphisms of the category of sections;
Maillard and Melli\`es defined their category of sections as a 
\emph{full} subcategory of the functor category,
but it is indeed our more restricted definition that
establishes an equivalence with the category of models of
an indexed Lawvere theory as defined by
Power~\cite{power:lce,power:indexed}.

\chapter{Conclusions and future work}
\label{chap-concl}
\section{Conclusions}
In this thesis we initiated a unified mathematical study on 
graded and indexed monads.
After providing in Chapter~\ref{chap-four-2-cats} the novel
2-categorical understanding of these notions 
(see also Appendix~\ref{apx-enriched-monad-tensor}),
we defined in Chapter~\ref{chap-main-constr} 
the following constructions and established their 
2-dimensional universality:
\begin{itemize}
\item The Eilenberg--Moore construction for graded monads.
\item The Kleisli construction for graded monads.
\item The Eilenberg--Moore construction for indexed monads.
\end{itemize}

These constructions are then applied in
Chapter~\ref{chap-appl} to two situations.

The first (Section~\ref{sec-decomp-lax-action}) 
deals with \emph{lax actions} of monoidal categories, and 
we showed that our Eilenberg--Moore and Kleisli constructions for
graded monads provide canonical \emph{resolutions} of
a lax action, canonical in the sense that they are
the \emph{terminal} and \emph{initial} ones respectively.
In doing so we encountered a theorem on
the correspondence of adjunctions (Section~\ref{subsec-fib-adj}),
whose nature seems to be {2-fibrational}.

As the second application (Section~\ref{sec-constr-maillard-mellies})
we reconstructed a proof of a beautiful theorem in 
\cite{maillard-mellies} providing the view that models 
of an indexed Lawvere theory~\cite{power:lce,power:indexed}
can be seen as \emph{sections}.
We employed the novel perspective brought to us by 
our Eilenberg--Moore construction for indexed monads, and 
gave a more conceptual construction of 
the \emph{category of sections} using the rightmost 2-functor~$\dprod$
constituting the famous \emph{base change adjoint triple} 
(Section~\ref{subsec-base-change}).
Not only our formulation now enables one to state more clearly 
in what sense the theorem generalizes the well-known close relationship 
between the category of models (in $\Set$) of a Lawvere theory and the 
Eilenberg--Moore category of the corresponding monad,
by providing clearer understanding of the situation
we were able to point out a missing condition in \cite{maillard-mellies} 
that should have been
posed on morphisms of sections
(see Notes for Chapter~\ref{chap-appl}).

\section{Directions for further research}
\subsection{The Kleisli construction for indexed monads}
One thing that is obviously missing from the current thesis is 
the \emph{Kleisli construction for indexed monads}.
Naturally we believe that the Kleisli construction should take place in
the 2-category~$\Emp$, using the observation in 
Section~\ref{subsec-im-mnd-emp}.
We conjecture the existence of a suitable 
construction completing the following picture:
    
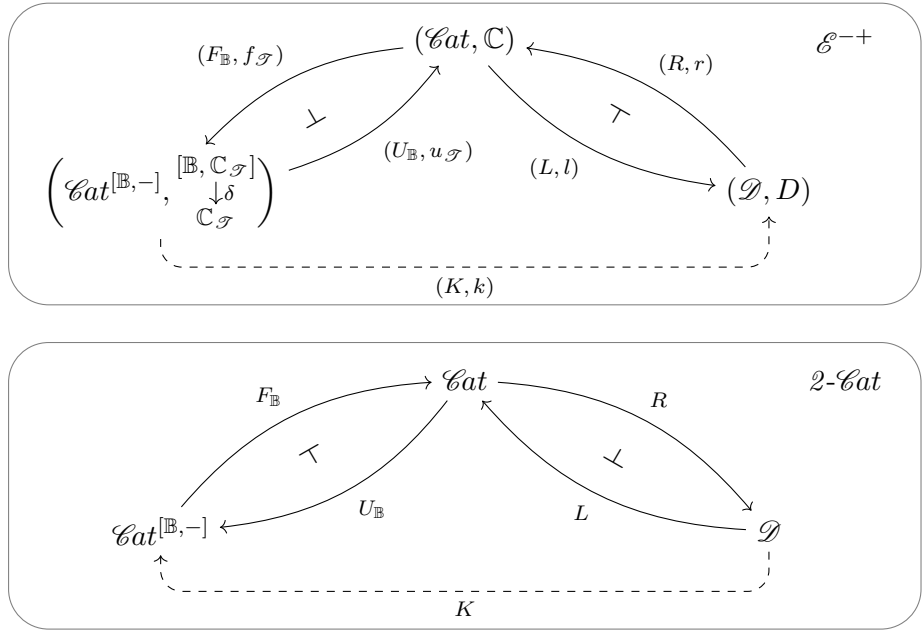
\[
    \begin{tikzpicture}
      \node (A)  at (0,0) {$\left(\EM{\Cat}{[\cat{B},-]},
            \valg{[\cat{B},\Kli]}{\delta}{\Kli}\right)$};
      \node (B) at (8,0) {$(\tcat{D},D)$};
      \node (C) at (4,2) {$(\Cat,\cat{C})$};
      
      \node at (2,1) [rotate=120,font=\large]{$\vdash$};
      \node at (6,1) [rotate=240,font=\large]{$\vdash$};
      
      \draw [->] (A) to [bend right=17]node [auto,labelsize,swap] {$(U_\cat{B},u_\im{T})$}(C);
      \draw [->] (C) to [bend right=20]node [auto,labelsize,swap] {$(F_\cat{B},f_\im{T})$}(A);
      
      \draw [->] (C) to [bend right=20]node [auto,swap,labelsize] {$(L,l)$}(B);
      \draw [->] (B) to [bend right=20]node [auto,swap,labelsize] {$(R,r)$}(C);
      
      \draw [->,dashed,rounded corners=8pt] (A)--(0,-1) to [bend right=0]node [auto,swap,labelsize]
       {$(K,k)$} (8,-1)--(B);
      
      \draw [rounded corners=15pt,myborder] (-2,-1.5) rectangle (10,2.5);
      \node at (9,2) {$\Emp$};

      \begin{scope}[shift={(0,-4.5)}]
       \node (A)  at (0,0) {$\EM{\Cat}{[\cat{B},-]}$};
       \node (B) at (8,0) {$\tcat{D}$};
       \node (C) at (4,2) {$\Cat$};
      
       \node at (2,1) [rotate=120,font=\large]{$\dashv$};
       \node at (6,1) [rotate=240,font=\large]{$\dashv$};
      
       \draw [<-] (A) to [bend right=23]node [auto,labelsize,swap] {$U_\cat{B}$}(C);
       \draw [<-] (C) to [bend right=23]node [auto,labelsize,swap] {$F_\cat{B}$}(A);
      
       \draw [<-] (C) to [bend right=23]node [auto,swap,labelsize] (Fpp) {$L$}(B);
       \draw [<-] (B) to [bend right=23]node [auto,swap,labelsize] (Fpp) {$R$}(C);
       
       \draw [<-,dashed,rounded corners=8pt] (A)--(0,-0.8) to [bend right=0]node [auto,swap,labelsize]
              {$K$} (8,-0.8)--(B);
       \draw [rounded corners=15pt,myborder] (-2,-1.3) rectangle (10,2.5);
             \node at (9,2) {$\TCat$};
   
      \end{scope}
    \end{tikzpicture}
    \]

All the attempts to define the object~$
\valgp{[\cat{B},\Kli]}{\delta}{\Kli}$ of $\EM{\Cat}{[\cat{B},-]}$
so far have failed.
One possible approach for this problem would be to understand
the other three constructions (Eilenberg--Moore and Kleisli
for graded monads, and Eilenberg--Moore for indexed monads) 
much more abstractly so that it is immediate how one can 
obtain the Kleisli construction for indexed monads.
The observation presented in Appendix~\ref{apx-enriched-monad-tensor}
might be valuable for this strategy.

\subsection{A 3-categorical study of graded and indexed monads}
The category~$\TCat$ of 2-categories, which we have employed when 
constructing the 2-categories~$\Epp$, $\Epm$, $\Emp$ and $\Emm$,
is inherently a \emph{3-category}. 
In fact, it may well be more natural to consider 
$\Epp$, $\Epm$, $\Emp$ and $\Emm$ as 3-categories as well;
indeed, there are fairly natural definitions of \emph{3-cells} of them.
They seem to arise as sub-3-categories of appropriate 
\emph{functor 3-categories}.

In this thesis we have confined ourselves to dealing only with 
\emph{strict} monoidal categories when considering 
graded monads.
Although this covers a large class of graded monads currently
employed, there do exist natural examples of graded monads 
graded by non-strict monoidal categories~\cite{mellies:tl-3}.
If one wishes to take an arbitrary monoidal category~$
\cat{M}=(\cat{M},\tensor,\e,\alpha,\lambda,\rho)$
as the parameter category of graded monads,
it seems inevitable to manipulate the 
\emph{pseudomonad}~$\cat{M}\times(-)$ and the
\emph{pseudocomonad}~$[\cat{M},-]$ on $\Cat$.
Just like monads live in 2-categories, pseudomonads live in
3-categories (or perhaps better: $\Gray$-categories),
and there is a work by Lack~\cite{lack}
which may naturally be thought of as the 3-dimensional version of 
Street's formal theory of monads~\cite{street:formal-theory-of-monads}.
In particular, Lack's work includes an abstract definition of the
\emph{object of pseudoalgebras} of a pseudomonad in a 
$\Gray$-category as a certain 3-dimensional limit, 
providing the 3-dimensional analogue of the Eilenberg--Moore object
of a monad.

Therefore there are evidences which support the claim that
our work would be done more properly
in the setting of 3-categories or 
$\Gray$-categories~\cite{gordon-power-street,gurski}.
As the world of 3-dimensional category theory still seems to remain 
rather unexplored, the upgraded 
\emph{3-categorical study of graded and indexed monads} could
bring important contribution to pure category theory, too.

\subsection{Categorical semantics of Bounded Linear Logic}
As mentioned at the end of Section~\ref{subsec-gm-defn},
the notion of graded comonad has been employed in the study of 
\emph{computational resources with parameters}.
Perhaps the most celebrated logical system dealing with 
parametrized computational resources is
\emph{Bounded Linear Logic}~\cite{girard-scedrov-scott},
which replaces the \emph{of course} modality~$\oc$ of 
Linear Logic~\cite{girard:linear} by a family of 
modalities~$\oc_x$ parametrized by 
\emph{resource polynomials}~$x$.
Bounded Linear Logic has recently been 
generalized~\cite{ghica-smith} so as to be able to
take an arbitrary semiring as the collection of parameters.

On the other hand, there has been a line of research seeking for
the appropriate categorical semantics for Linear Logic, such as
\cite{seely:linear,benton:lnl}, to name just two; 
see \cite{mellies:ll} for a nice survey.
The biggest challenge was the identification of 
a suitable categorical structure modeling the modality~$\oc$,
and the consensus reached is that $\oc$ should be modeled as a certain
\emph{comonad}.
Mathematical results on the co-Eilenberg--Moore and co-Kleisli
categories for comonads have been useful in clarifying the 
relationship between the various proposed semantics.

We expect that 
our co-Eilenberg--Moore and co-Kleisli
constructions for \emph{graded comonads} can be fruitfully 
employed in the study of 
\emph{categorical semantics of Bounded Linear Logic},
on which it looks that not much work has been done.

\subsection{Syntactical development}
From the viewpoint of the theory of computational effects,
the development presented in the current thesis remains
entirely \emph{semantical}.
The recent theoretical study of computational effects has been greatly 
benefited by the \emph{syntactical} approach,
in which the emphasis is placed on
\emph{Lawvere theories} rather than monads,
as is conspicuous for example in the seminal paper~\cite{plotkin-power}
by Plotkin and Power.

We believe that the development of suitable 
mathematical theories of \emph{graded Lawvere theories}
and \emph{indexed Lawvere theories} would be an important
contribution to the theory of 
\emph{computational effects with parameters}.
We do not know yet what a graded Lawvere theory means, and 
although Power~\cite{power:lce,power:indexed} has already defined
the notion of indexed Lawvere theory which possesses a nice 
relationship to our notion of \emph{indexed monad} 
(Section~\ref{sec-constr-maillard-mellies}), Power himself
makes it clear in his papers that his definition of 
indexed Lawvere theory is not a definitive one;
so this direction of research could bring us to 
a whole new world.

\subsection*{Notes}
I had a valuable discussion on possible definitions of the 
Kleisli category of an indexed monad with Kazuyuki Asada and 
Takeshi Tsukada.

The possibility of a \emph{3-categorical} approach has been in a sense 
evident as soon as I defined $\Epp$, 
but it was a series of enlightening discussions with
John Power which brought me much clearer insight into this.

The research theme on categorical semantics of Bounded Linear Logic
was suggested to me by Ichiro Hasuo.



\appendix

\chapter{Compositions in $\Epp$}
\label{apx-composition}
We describe in detail
how 1-cells and 2-cells in the 2-category $\Epp$ are
composed; compositions in the other three 2-categories~$
\Epm$, $\Emp$ and $\Emm$ are completely similar.

\section{Compositions of 1-cells}
Suppose we have the following diagram in $\Epp$: 
\[
(\tcat{A},A)\xrightarrow{(F,f)}(\tcat{B},B)\xrightarrow{(G,g)}(\tcat{C},C)
\]
Recall that $f\colon FA\to B$ can be thought of as a 1-cell in $\tcat{B}$ and 
$g\colon GB\to C$ as a 1-cell in $\tcat{C}$.
We define 
\[
(G,g)\circ (F,f)\quad\coloneqq\quad(GF,\ g\circ Gf)
\]
where the second component is the 1-cell 
\[
GFA\xrightarrow{Gf}GB\xrightarrow{\mathmakebox[1em]{g}}C
\]
in $\tcat{C}$. 

\section{Vertical compositions of 2-cells}
Suppose we have the following diagram in $\Epp$: 
    \[
    \begin{tikzpicture}
      \node (A)  at (0,0) {$(\tcat{A},A)$};
      \node (B) at (4,0) {$(\tcat{B},B)$};
      
      \draw [->] (A) to [bend left=70]node (F)  [auto,labelsize] {$(F,f)$}(B);
      \draw [->] (A) to [bend left=0]node [fill=white,labelsize] (Fp)  {$(F^\prime,f^\prime)$}(B);
      \draw [->] (A) to [bend right=70]node [auto,swap,labelsize] (Fpp) {$(F^{\prime\prime},f^{\prime\prime})$}(B);
      \draw [2cellnode] (F) to node [auto,labelsize]{$(\Theta,\alpha)$}(Fp);
      \draw [2cellnode] (Fp) to node [auto,labelsize]{$(\Theta^\prime,\alpha^\prime)$}(Fpp);
    \end{tikzpicture}
    \]
We regard $\alpha$ and $\alpha^\prime$ respectively as:
    \[
    \begin{tikzpicture}
      \node (a)  at (3.5,2) {$FA$};
      \node (ap) at (6.5,2) {$B$};
      \node (b)  at (4.5,0.5) {$F^\prime A$};
      \node at (7.5,1.25) {in $\tcat{B}$};
      
      \draw [->] (a) to [bend right=0]node [auto,labelsize] {$f$}(ap);
      \draw [->] (a)  to [bend right=20]node [auto,swap,labelsize] {$\Theta_A$}(b);
      \draw [->] (b) to [bend right=20]node [auto,swap,labelsize] {$f^\prime$}(ap);
      \draw [2cell] (4.5,1.9) to node [auto,labelsize]{$\alpha$}(b);
      
     \begin{scope}[shift={(7,0)}]
      \node (a)  at (3.5,2) {$F^\prime A$};
      \node (ap) at (6.5,2) {$B$};
      \node (b)  at (4.5,0.5) {$F^{\prime\prime} A$};
      \node at (7.5,1.25) {in $\tcat{B}$};
      
      \draw [->] (a) to [bend right=0]node [auto,labelsize] {$f^\prime$}(ap);
      \draw [->] (a)  to [bend right=20]node [auto,swap,labelsize] {$\Theta^\prime_A$}(b);
      \draw [->] (b) to [bend right=20]node [auto,swap,labelsize] {$f^{\prime\prime}$}(ap);
      \draw [2cell] (4.5,1.9) to node [auto,labelsize]{$\alpha^\prime$}(b);
      
     \end{scope}
    \end{tikzpicture}
    \]
Now define
\[
(\Theta^\prime,\alpha^\prime)\vc(\Theta,\alpha)\quad\coloneqq\quad
(\Theta^\prime\vc\Theta,\ (\alpha^\prime \hc \Theta_A)\vc \alpha)
\]
with the second component being the 2-cell
    
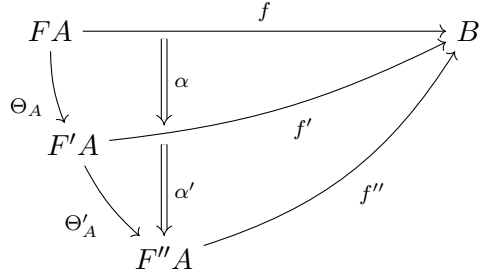
\[
    \begin{tikzpicture}
      \node (a)  at (0,2) {$FA$};
      \node (ap) at (5.5,2) {$B$};
      \node (b)  at (0.3,0.5) {$F^\prime A$};
      \node (ba)  at (1.5,-1) {$F^{\prime\prime} A$};
      
      \draw [->] (a) to [bend right=0]node [auto,labelsize] {$f$}(ap);
      \draw [->] (a)  to [bend right=10]node [auto,swap,labelsize] {$\Theta_A$}(b);
      \draw [->] (b)  to [bend right=10]node [auto,swap,labelsize] {$\Theta^\prime_A$}(ba);
      \draw [->] (b) to [bend right=10]node [auto,swap,labelsize] {$f^\prime$}(ap);
      \draw [->] (ba) to [bend right=20]node [auto,swap,labelsize] {$f^{\prime\prime}$}(ap);
      \draw [2cell] (1.5,1.9) to node [auto,labelsize]{$\alpha$}(1.5,0.75);
      \draw [2cell] (1.5,0.5) to node [auto,labelsize]{$\alpha^\prime$}(ba);
    \end{tikzpicture}
    \]
%
in $\tcat{B}$.
    
\section{Whiskerings}
Before describing the somewhat complicated horizontal compositions
of 2-cells in $\Epp$, we begin with the simpler situations of 
\emph{whiskerings}.

Suppose we have the following diagram in $\Epp$:
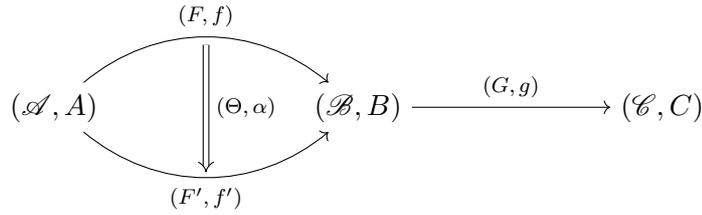
    \[
    \begin{tikzpicture}
      \node (A)  at (0,0) {$(\tcat{A},A)$};
      \node (B) at (4,0) {$(\tcat{B},B)$};
      \node (C) at (8,0) {$(\tcat{C},C)$};

      \draw [->] (A) to [bend left=40]node (F)  [auto,labelsize] {$(F,f)$}(B);
      \draw [->] (A) to [bend right=40]node [auto,swap,labelsize] (Fp) {$(F^{\prime},f^{\prime})$}(B);
      \draw [->] (B) to [bend left=0]node [auto,labelsize] (G) {$(G,g)$}(C);
      
      \draw [2cellnode] (F) to node [auto,labelsize]{$(\Theta,\alpha)$}(Fp);
    \end{tikzpicture}
    \]
We define
\[
(G,g)\hc(\Theta,\alpha)\quad\coloneqq\quad(G\hc \Theta,\ g\hc G\alpha)
\]
with the second component
    
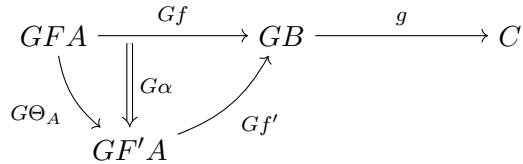
\[
    \begin{tikzpicture}
      \node (a)  at (3.5,2) {$GFA$};
      \node (ap) at (6.5,2) {$GB$};
      \node (b)  at (4.5,0.5) {$GF^\prime A$};
      \node (c)  at (9.5,2) {$C$};
      
      \draw [->] (a) to [bend right=0]node [auto,labelsize] {$Gf$}(ap);
      \draw [->] (a)  to [bend right=20]node [auto,swap,labelsize] {$G\Theta_A$}(b);
      \draw [->] (b) to [bend right=20]node [auto,swap,labelsize] {$Gf^\prime$}(ap);
      \draw [->] (ap) to [bend right=0]node [auto,labelsize] {$g$}(c);
      \draw [2cell] (4.5,1.9) to node [auto,labelsize]{$G\alpha$}(b);

    \end{tikzpicture}
    \]

Suppose we have the following diagram in $\Epp$:
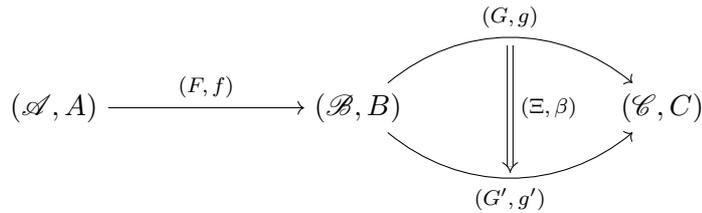
    \[
    \begin{tikzpicture}
      \node (A)  at (0,0) {$(\tcat{A},A)$};
      \node (B) at (4,0) {$(\tcat{B},B)$};
      \node (C) at (8,0) {$(\tcat{C},C)$};

      \draw [->] (A) to [bend left=0]node (F)  [auto,labelsize] {$(F,f)$}(B);
      \draw [->] (B) to [bend left=40]node [auto,labelsize] (G) {$(G,g)$}(C);
      \draw [->] (B) to [bend right=40]node [auto,swap,labelsize] (Gp) {$(G^{\prime},g^{\prime})$}(C);

      \draw [2cellnode] (G) to node [auto,labelsize]{$(\Xi,\beta)$}(Gp);
    \end{tikzpicture}
    \]
Define
\[
(\Xi,\beta)\hc(F,f)\quad\coloneqq\quad(\Xi\hc F,\ \beta\hc Gf)
\]
with the second component
    \[
    \begin{tikzpicture}
      \node (c)  at (0.5,2) {$GFA$};
      \node (cp)  at (1.5,0.5) {$G^\prime FA$};
      \node (a)  at (3.5,2) {$GB$};
      \node (ap) at (6.5,2) {$C$};
      \node (b)  at (4.5,0.5) {$G^\prime B$};

      \draw [->] (a) to [bend right=0]node [auto,labelsize] {$g$}(ap);
      \draw [->] (a)  to [bend right=20]node [auto,swap,labelsize] {$\Xi_B$}(b);
      \draw [->] (c)  to [bend right=20]node [auto,swap,labelsize] {$\Xi_{FA}$}(cp);
      \draw [->] (b) to [bend right=20]node [auto,swap,labelsize] {$g^\prime$}(ap);
      \draw [->] (c) to [bend right=0]node [auto,labelsize] {$Gf$}(a);
      \draw [->] (cp) to [bend right=0]node [auto,labelsize,swap] {$G^\prime f$}(b);
      \draw [2cell] (4.5,1.9) to node [auto,labelsize]{$\beta$}(b);

    \end{tikzpicture}
    \]
\section{Horizontal compositions of 2-cells}
Suppose we have the following diagram in $\Epp$:
    \[
    \begin{tikzpicture}
      \node (A)  at (0,0) {$(\tcat{A},A)$};
      \node (B) at (4,0) {$(\tcat{B},B)$};
      \node (C) at (8,0) {$(\tcat{C},C)$};

      \draw [->] (A) to [bend left=40]node (F)  [auto,labelsize] {$(F,f)$}(B);
      \draw [->] (A) to [bend right=40]node [auto,swap,labelsize] (Fp) {$(F^{\prime},f^{\prime})$}(B);
      \draw [->] (B) to [bend left=40]node [auto,labelsize] (G) {$(G,g)$}(C);
      \draw [->] (B) to [bend right=40]node [auto,swap,labelsize] (Gp) {$(G^{\prime},g^{\prime})$}(C);

      \draw [2cellnode] (F) to node [auto,labelsize]{$(\Theta,\alpha)$}(Fp);
      \draw [2cellnode] (G) to node [auto,labelsize]{$(\Xi,\beta)$}(Gp);
    \end{tikzpicture}
    \]
    
We should have, as in any 2-category, the following identity
\[
(\Xi,\beta)\hc(\Theta,\alpha)\quad=\quad\big((\Xi,\beta)\hc(F^\prime,f^\prime)\big)\vc\big((G,g)\hc(\Theta,\alpha)\big).
\]
Therefore we can take this as the definition:
\[
(\Xi,\beta)\hc(\Theta,\alpha)\quad\coloneqq\quad
\big(\Xi\hc\Theta
,\ 
\beta\hc G\alpha\big)
\]
with the second component depicted as follows:
    \[ 
    \begin{tikzpicture}
      \node (a)  at (0,2) {$GFA$};
      \node (ap) at (6,2) {$C$};
      \node (b)  at (1,0.5) {$GF^\prime A$};
      \node (ba)  at (2,-1) {$G^\prime F^{\prime} A$};
      \node (fb) at (3,2) {$GB$};
      \node (fpb) at (4,0.5) {$G^\prime B$};

      \draw [->] (a) to [bend right=0]node [auto,labelsize] {$Gf$}(fb);
      \draw [->] (fb) to [bend right=0]node [auto,labelsize] {$g$}(ap);
      \draw [->] (a)  to [bend right=20]node [auto,swap,labelsize] {$G\Theta_A$}(b);
      \draw [->] (b)  to [bend right=20]node [auto,swap,labelsize] {$\Xi_{F^\prime A}$}(ba);
      \draw [->] (b) to [bend right=20]node [auto,swap,near start,labelsize] {$Gf^\prime$}(fb);
      \draw [->] (ba) to [bend right=20]node [auto,swap,labelsize] {$G^\prime f^{\prime}$}(fpb);
      \draw [->] (fb) to [bend right=20]node [auto,swap,near end,labelsize] {$\Xi_B$}(fpb);
      \draw [->] (fpb) to [bend right=20]node [auto,swap,labelsize] {$g^{\prime}$}(ap);
      
      \draw [2cell] (1,1.9) to node [auto,labelsize]{$G\alpha$}(1,0.8);
      \draw [2cell] (4,1.9) to node [auto,labelsize]{$\beta$}(4,0.8);
    \end{tikzpicture}
    \]
    
Or alternatively, we can start from the following identity 
that also holds in any 2-category
\[
(\Xi,\beta)\hc(\Theta,\alpha)\quad=\quad\big((G^\prime,g^\prime)\hc(\Theta,\alpha)\big)\vc\big((\Xi,\beta)\hc(F,f)\big)
\]
and define
\[
(\Xi,\beta)\hc(\Theta,\alpha)\quad\coloneqq\quad
\big(\Xi\hc\Theta
,\ 
(g^\prime\hc G^\prime\alpha\hc\Xi_{FA})\vc(\beta\hc Gf)\big)
\]
with the second component
    \[
    \begin{tikzpicture}
      \node (c)  at (0.5,2) {$GFA$};
      \node (cp)  at (1.5,0.5) {$G^\prime FA$};
      \node (a)  at (3.5,2) {$GB$};
      \node (ap) at (6.5,2) {$C$};
      \node (b)  at (4.5,0.5) {$G^\prime B$};
      \node (d)  at (2.5,-1) {$G^\prime F^\prime A$};

      \draw [->] (a) to [bend right=0]node [auto,labelsize] {$g$}(ap);
      \draw [->] (a)  to [bend right=20]node [auto,swap,near start,labelsize] {$\Xi_B$}(b);
      \draw [->] (c)  to [bend right=20]node [auto,swap,labelsize] {$\Xi_{FA}$}(cp);
      \draw [->] (b) to [bend right=20]node [auto,swap,labelsize] {$g^\prime$}(ap);
      \draw [->] (c) to [bend right=0]node [auto,labelsize] {$Gf$}(a);
      \draw [->] (cp) to [bend right=0]node [auto,near start,labelsize] {$G^\prime f$}(b);
      \draw [->] (cp) to [bend right=20]node [auto,labelsize,swap] {$G^\prime\Theta_A$}(d);
      \draw [->] (d) to [bend right=20]node [auto,labelsize,swap] {$G^\prime f^\prime$}(b);
      \draw [2cell] (4.5,1.9) to node [auto,labelsize]{$\beta$}(b);
      \draw [2cell] (2.5,0.4) to node [auto,labelsize]{$G^\prime \alpha$}(d);
    \end{tikzpicture}
    \]
    
That these two definitions of vertical compositions
of 2-cells coincide is an immediate consequence of
the 2-naturality of $\Xi$.


\chapter{An abstract view of graded and indexed monads}
\label{apx-enriched-monad-tensor}

\section{Enriched (co)monads via (co)powers}
\label{apx-sec-enriched}
In developing the mathematical theory of graded and indexed monads,
we heavily employed the following 2-monads and 2-comonads on $\Cat$:
\begin{itemize}
\item The 2-monad~$\cat{M}\times(-)$ where $(\cat{M},\tensor,\e)$
is a strict monoidal category.
\item The 2-comonad~$[\cat{M},-]$ where $(\cat{M},\tensor,\e)$
is a strict monoidal category.
\item The 2-comonad~$\cat{B}\times(-)$ where $\cat{B}$ is a category.
\item The 2-monad~$[\cat{B},-]$ where $\cat{B}$ is a category.
\end{itemize}
In fact there is a general construction in enriched category theory~\cite{kelly:enriched}
of which these are instances;
the construction induces enriched monads and comonads 
on a suitably complete and cocomplete
enriched category from monoids and comonoids in the enriching category,
through \emph{powers} and \emph{copowers}.

Suppose that $\V=(\V,\bullet,[-,-],1)$ is a 
symmetric monoidal closed category
and $\ecat{A}$ a $\V$-category.
\begin{defn}
Given objects~$V$ and $A$ of $\V$ and $\ecat{A}$ respectively,
the \defemph{power of $A$ by $V$} is the object~$V\fork A$
of $\ecat{A}$ satisfying
\[
\ecat{A}(X,V\fork A)\quad\cong\quad [V,\ecat{A}(X,A)]
\]
$\V$-natural in $X\in\ecat{A}$.
A $\V$-category is said to be \defemph{powered} if it 
has all powers.
\end{defn}

\begin{defn}
Given objects~$V$ and $A$ of $\V$ and $\ecat{A}$ respectively,
the \defemph{copower of $A$ by $V$} is the object~$V\copow A$
of $\ecat{A}$ satisfying
\[
\ecat{A}(V\tensor A, X)\quad\cong\quad [V,\ecat{A}(A,X)]
\]
$\V$-natural in $X\in\ecat{A}$.
A $\V$-category is said to be \defemph{copowered} if it 
has all copowers.
\end{defn}

Suppose a $\V$-category~$\ecat{A}$ is both powered and copowered.
Now the general construction alluded above is the following:
\begin{itemize}
\item The functor~$M\tensor(-)\colon\ecat{A}\to\ecat{A}$ 
becomes a $\V$-monad on $\ecat{A}$ when
$(M,\ M\bullet M\xrightarrow{m}M,\ 1\xrightarrow{u}M)$
is a monoid in $\V$.
\item The functor~$M\fork(-)\colon\ecat{A}\to\ecat{A}$ 
becomes a $\V$-comonad on $\ecat{A}$ when
$(M,\ M\bullet M\xrightarrow{m}M,\ 1\xrightarrow{u}M)$
is a monoid in $\V$.
\item The functor~$B\tensor(-)\colon\ecat{A}\to\ecat{A}$ 
becomes a $\V$-comonad on $\ecat{A}$ when
$(B,\ B\xrightarrow{d}B\bullet B,\ B\xrightarrow{e}1)$
is a comonoid in $\V$.
\item The functor~$B\fork(-)\colon\ecat{A}\to\ecat{A}$ 
becomes a $\V$-monad on $\ecat{A}$ when
$(B,\ B\xrightarrow{d}B\bullet B,\ B\xrightarrow{e}1)$
is a comonoid in $\V$.
\end{itemize}

Note that 2-(co)monads are nothing but $\Cat$-(co)monads, 
and strict monoidal categories are nothing but monoids in
$(\Cat,\times,1)$.
Also, as in any Cartesian monoidal category,
a comonoid in $(\Cat,\times,1)$ is the same thing as an object of $\Cat$;
every category $\cat{B}$ admits a unique comonoid structure given by
$(\cat{B},\ \cat{B}\xrightarrow{\Delta_\cat{B}}\cat{B}\times\cat{B},\ 
\cat{B}\xrightarrow{t_\cat{B}} 1)$.
Combining these observations with the fact that
powers and copowers in $\Cat$ as a $\Cat$-category
are given by exponentials and Cartesian product respectively,
we now obtain a unified abstract explanation of 
the 2-(co)monads at the beginning of this section.

\section{Generalizing monads in $\Cat$ via lifting}

Now that we have understood abstractly the four kinds of 2-(co)monads 
used in studying graded and indexed monads, let us proceed to 
look again the relationship of graded and indexed monads and 
these 2-(co)monads.
Assuming the constructions of 
2-categories~$\Epp$, $\Epm$, $\Emp$ and $\Emm$,
a priori there are eight notions of \emph{generalized monad} on 
$\cat{C}\in\Cat$
obtained by lifting these 2-(co)monads:
\begin{enumerate}
\item Monads in $\Epp$ on $(\Cat,\cat{C})$ 
above the 2-monad $\cat{M}\times(-)$:
\begin{align*}
T\quad&\colon\quad\cat{M}\times\cat{C}\quad\longrightarrow\quad\cat{C}.
\intertext{\item Monads in $\Epp$ on $(\Cat,\cat{C})$ 
above the 2-monad $[\cat{B},-]$:}
T\quad&\colon\quad [\cat{B},\cat{C}]\quad\longrightarrow\quad\cat{C}.
\intertext{\item Monads in $\Epm$ on $(\Cat,\cat{C})$ 
above the 2-comonad $[\cat{M},-]$:}
T\quad&\colon\quad [\cat{M},\cat{C}]\quad\longrightarrow\quad\cat{C}.
\intertext{\item Monads in $\Epm$ on $(\Cat,\cat{C})$ 
above the 2-comonad $\cat{B}\times(-)$:}
T\quad&\colon\quad \cat{B}\times\cat{C}\quad\longrightarrow\quad\cat{C}.
\intertext{\item Monads in $\Emp$ on $(\Cat,\cat{C})$ 
above the 2-monad $\cat{M}\times(-)$:}
T\quad&\colon \quad\cat{C}\quad\longrightarrow\quad\cat{M}\times\cat{C}.
\intertext{\item Monads in $\Emp$ on $(\Cat,\cat{C})$ 
above the 2-monad $[\cat{B},-]$:}
T\quad&\colon\quad \cat{C}\quad\longrightarrow\quad[\cat{B},\cat{C}].
\intertext{\item Monads in $\Emm$ on $(\Cat,\cat{C})$ 
above the 2-comonad $[\cat{M},-]$:}
T\quad&\colon\quad \cat{C}\quad\longrightarrow\quad[\cat{M},\cat{C}].
\intertext{\item Monads in $\Emm$ on $(\Cat,\cat{C})$ 
above the 2-comonad $\cat{B}\times(-)$:}
T\quad&\colon\quad\cat{C}\quad\longrightarrow\quad\cat{B}\times\cat{C}.
\end{align*}
\end{enumerate}
Actually, it turns out that the notions 1 and 7 coincide and are that of 
\emph{graded monad},
and the notions 4 and 6 coincide and are that of \emph{indexed monad};
this coincidence is because of the general adjointness 
\begin{equation}
\label{eqn-apx-adj-power}
V\tensor (-)\quad\dashv\quad V\fork(-)
\end{equation}
between copowers and powers.

Therefore one can characterize the notions of 
graded and indexed monads exactly as those notions of
generalized monad on a category listed above, for which there are
two different ways of thinking about them thanks to the 
fundamental adjointness~(\ref{eqn-apx-adj-power}).
Note that indeed we used this dual view to construct (or at least 
try to construct)
\emph{both} Eilenberg--Moore and Kleisli 
categories.
\section*{Notes}
The materials contained in this chapter have occurred to  
my mind through discussions with John Power.
In particular I learned the construction of
enriched (co)monads from (co)monoids from him.

\end{document}